\documentclass[12pt,twoside]{article} 
\usepackage{amssymb,amsbsy,amsmath,amsfonts,amscd,times}
\newcommand{\C}{\mathbb{C}}
\newcommand{\N}{\mathbb{N}}
\newcommand{\Q}{\mathbb{Q}}

\date{17 march 2008} 

\begin{document} 

\centerline{\bf Int. J. Contemp. Math. Sciences, Vol. 3, 2008, no. 18, 861 - 933} 

\centerline{} 

\centerline{} 

\centerline{\Large{\bf Jets de Demailly-Semple d'Ordres 4 et 5}} 

\centerline{}

\centerline{\Large{\bf En Dimension 2}} 

\centerline{} 

\centerline{\bf {Jo\"el Merker}} 

\centerline{} 

\centerline{D\'epartement de Math\'ematiques et Applications} 

\centerline{\'Ecole Normale Sup\'erieure} 

\centerline{45 rue d'Ulm, F-75230 Paris Cedex 05, France} 

\centerline{merker@dma.ens.fr} 

\newtheorem{problem}{Probl\`eme}
\newtheorem{definition}{D\'efinition}
\newtheorem{theorem}{Th\'eor\`eme}
\newtheorem{proposition}{Proposition}
\newtheorem{lemma}{Lemme}
\newtheorem{assertion}{Assertion}
\newtheorem{corollary}{Corollaire}

\begin{abstract} 
In dimension 2 and for jets of order 4, there are 9 Demailly-Semple
invariant polynomials generated by bracketing the invariants for jets
of order $\leqslant 3$ (Demailly\,; Rousseau).  They share 9
fundamental syzygies and an Euler characteristic computation shows
that every entire holomorphic curve into a generic smooth complex
projective algebraic surface $X \subset P_3 ( \C)$ satisfies global
algebraic differential equations provided ${\rm deg}\, X \geqslant 9$.

For jets of order 5, the algebraic structure explodes in
complexity. Bracketing gives 36 invariants.  Removing redundancies
leaves 24 fundamental invariants, 11 of which, denoted $f_1'$,
$\Lambda^3$, $\Lambda_1^5$, $\Lambda_{ 1, 1}^7$, $\Lambda_{ 1, 1,
1}^9$, $M^8$, $M_1^{ 10}$, $K_{ 1, 1}^{ 12}$, $N^{ 12}$, $H_1^{ 14}$
and $F_{ 1, 1}^{ 16}$ and called {\sl bi-invariants}, are meaningful
for Euler characteristic estimates. Unexpectedly, 5 more appear: $X^{
18}$, $X^{ 19}$, $X^{ 21}$, $X^{ 23}$ and $X^{ 25}$.  The paper
contains an algorithm generating all such (complicated) bi-invariants.

\end{abstract} 



\bigskip

\hfill
\begin{minipage}[t]{10cm}
\baselineskip=0.32cm
{\tiny{\sf
I do not mean to suggest that all mathematical relations can be
perceived directly as obvious if they are visualised in the right
way\,\,---\,\,or merely that they can always be perceived in some
other way that is immediate to our intuitions. Far from it. Some
mathematical relations require long chains of reasoning before they
can be perceived with certainty. But the object of mathematical proof
is, in effect, to provide such chains of reasoning where each step is
indeed something that can be perceived as obvious. Consequently, the
endpoint of the reasoning is something that must be accepted as true,
even though it may not, in itself, be at all obvious.}
\hfill 
Sir Roger {\sc Penrose}, {\em Shadows of the Mind}, Oxford, 1994.}
\end{minipage}

\section*{\S1.~Introduction}

En dimension $\nu \leqslant 3$, la description de l'alg\`ebre
$\mathcal{ DS}_\nu^\kappa$ des polyn\^omes invariants de
Demailly-Semple n'est explicit\'ee dans la litt\'erature que pour les
jets d'ordre $\kappa \leqslant 3$ (\cite{ de1997, ro2007})\,; en
dimension $\nu = 2$ et \`a l'ordre $\kappa = 3$, on sait que
l'alg\`ebre $\mathcal{ DS}_2^3$ est engendr\'ee par 5 polyn\^omes
fondamentaux li\'es entre eux par une unique syzygie\,; en dimension 2
et \`a l'ordre 4, d'apr\`es un travail non publi\'e de Demailly ({\it
cf.} {\it e.g.} \cite{ de2007}), $\mathcal{ DS}_2^4$ est engendr\'ee
par 9 polyn\^omes invariants fondamentaux. Dans cet article, nous
\'etablissons ce r\'esultat, nous explicitons les syzygies
fondamentales entre ces 9 invariants\,\,---\,\,qui sont aussi au
nombre de 9\,\,---, nous effectuons un calcul de Riemann-Roch pour
estimer la caract\'eristique d'Euler du fibr\'e correspondant, et nous
en d\'eduisons que toute courbe holomorphe enti\`ere \`a valeurs dans
une surface projective alg\'ebrique complexe lisse $X^2 \subset P_3 (
\C)$ (tr\`es) g\'en\'erique
de degr\'e $d\geqslant 9$ satisfait des \'equations alg\'ebriques
globales non triviales d'ordre 4.

\smallskip

Pour $\kappa \geqslant 5$, la structure alg\'ebrique de $\mathcal{
DS}_2^\kappa$ explose en complexit\'e. Nous exposons un proc\'ed\'e
r\'ecursif\,: le {\sl crochet entre deux invariants}, qui semble
permettre ({\it cf.} \cite{ de2007}), en toute dimension et pour les
jets d'ordre quelconque, d'engendrer un syst\`eme fondamental de
polyn\^omes invariants, et aussi de trois autres proc\'ed\'es
r\'ecursifs\,: {\sl identit\'es de Jacobi}, {\sl identit\'es
pl\"ucke\-riennes d'ordre un}, et {\sl identit\'es pl\"ucke\-riennes
d'ordre deux}, qui d\'ecrivent exhaustivement le gigantesque id\'eal
des relations entre les invariants ainsi construits.

\smallskip

Ces quatre proc\'ed\'es (confirm\'es sur les cas connus) font
appara\^{\i}tre de mani\`ere saillante une explosion symbolique
incontr\^olable. Par exemple, pour le cas $\nu = 2$ et $\kappa = 5$
\'etudi\'e compl\`etement ici, on re\c coit 36 invariants bruts dont 12
exactement sont redondants, et on doit tenir compte de 210 syzygies
non redondantes de degr\'e $\leqslant 4$ entre ces invariants,
lesquelles se d\'eploient sur 13 pages manuscrites\footnote{\, Pour
les jets d'ordre $\kappa = 6$, nous ne nous sommes pas risqu\'e \`a
calculer les 325 invariants attendus, ni \`a entreprendre d'\'ecrire
les 14\,950 syzygies de degr\'e $\leqslant 5$ qui nous sont donn\'ees
automatiquement. }. En ne consid\'erant que les invariants stables par
l'action d'un certain sous-groupe unipotent de ${\sf GL}_2 ( \C)$, le
nombre de syzygies fondamentales se r\'eduit \`a 15, ce qui permet
d'effectuer un calcul de Riemann-Roch au niveau $\kappa = 5$.

\smallskip

Si l'on s'en tenait seulement aux estimations ainsi obtenues pour la
caract\'eristi\-que d'Euler du fibr\'e de Demailly-Semple $\mathcal{
DS}_{ 2, m}^\kappa T_X^*$, on pourrait penser qu'il faudrait
conna\^{\i}tre sa d\'ecomposition de Schur pour des niveaux $\kappa$ au
moins $\geqslant 20$, eu \'egard \`a la difficile conjecture
d'hyperbolicit\'e de Kobayashi concernant les surfaces complexes de
degr\'e $d\geqslant 5$ dans $P_3 ( \C)$ qui sont (tr\`es)
g\'en\'eriques (\cite{ de1997, ro2007, de2007}). Mais l'avenir dira si
cette approche ne devrait pas \^etre amend\'ee et r\'eorient\'ee d\`es
le niveau $\kappa = 5$, en tenant compte de la structure sp\'ecifique
de $\mathcal{ DS}_2^5$.

\smallskip

Nos r\'esultats principaux apparaissent dans les Sections~4, 5, 6, 7
et 8.

\section*{\S2.~Polyn\^omes invariants et diff\'erentiation compos\'ee}

\noindent{\bf Notations initiales.}
En dimension $\nu \geqslant 2$, le jet strict d'ordre $\kappa
\geqslant 1$ en un point fix\'e d'une application holomorphe locale $f =
( f_1, f_2, \dots, f_\nu)$ de $\C$ \`a valeurs dans $\C^\nu$ sera not\'e\,:
\[
j^\kappa f 
:=
\big(f_1', \dots, f_\nu', 
f_1'', \dots, f_\nu'', 
\dots\dots, f_1^{ (\kappa)},\dots,f_\nu^{(\kappa)}\big).
\]
 
\noindent{\bf Polyn\^omes invariants par reparam\'etrisation.}
Pour $\kappa \geqslant 1$ entier, on recherche les polyn\^omes ${\sf P}
= {\sf P } ( j^\kappa f)$ tels que\,:
\[
{\sf P}
\big(j^\kappa(f\circ \phi)\big)
=
(\phi')^m\,
{\sc P}\big(
(j^\kappa f)\circ\phi
\big),
\]
pour tout biholomorphisme local $\phi : U \to \phi ( U)$, o\`u $U
\subset \C$ est ouvert, et o\`u $m \geqslant 1$ est un entier que nous
appellerons {\sl poids}\, de ${\sf P}$. On note $\mathcal{ DS}_{ \nu,
m}^\kappa$ l'espace vectoriel constitu\'e par ces {\sl polyn\^omes
invariants par reparam\'etrisation}
(\cite{ de1997}). La r\'eunion $\mathcal{ DS}_\nu^\kappa :=
\bigoplus_{ m \geqslant 1} \, \mathcal{ DS}_{ \nu, m}^\kappa$ forme
une alg\`ebre gradu\'ee\,: $\mathcal{ DS}_{\nu, m_1}^\kappa \cdot \mathcal{
DS}_{\nu, m_2 }^\kappa \subset \mathcal{ DS}_{ \nu, m_1 +
m_2}^\kappa$. On travaillera toujours dans une fibre en un point $z
\in U$ qui n'appara\^{\i}tra pas dans les notations.

Lorsque $\kappa = 1$, les composantes $f_i'$ ($i = 1, \dots, \nu$) du
jet d'ordre un satisfont
\[
\big(
f_i\circ\phi
\big)'
=
\phi'\,f_i',
\]
et par cons\'equent, tout polyn\^ome ${\sf P} = {\sf P} \big( f_1', \dots,
f_\nu' \big)$ qui ne d\'epend que du jet d'ordre $1$ est invariant par
reparam\'etrisation.

\smallskip\noindent{\bf Diff\'erentiation compos\'ee jusqu'\`a l'ordre 5.} 
En posant $g_i := f_i \circ \phi$ pour $i = 1, \dots, \nu$, on calcule\,:
\[
\aligned
g_i'
&
=
\phi'f_i',
\\
g_i''
&
=
\phi''f_i'
+
{\phi'}^2f_i'',
\\
g_i'''
&
=
\phi'''f_i'
+
3\,\phi''\phi'f_i''
+
{\phi'}^3f_i''',
\\
g_i''''
&
=
\phi''''f_i'
+
4\,\phi'''\phi'f_i''
+
3\,{\phi''}^2f_i''
+
6\,\phi''{\phi'}^2f_i'''
+
{\phi'}^4f_i'''',
\\
g_i'''''
&
:=
\phi'''''f_i'
+
5\,\phi''''\phi'f_i''
+
10\,\phi'''\phi''f_i''
+
15\,{\phi''}^2\phi'f_i'''
+
\\
&
\ \ \ \ \ \ \ \ \ \ \ \ \ \ \ \ \ \ \ \ \ \ \ \
+
10\,\phi'''{\phi'}^2f_i'''
+
10\,\phi''{\phi'}^3f_i''''
+
{\phi'}^5f_i'''''.
\endaligned
\]
Moralement, on recherche tous les polyn\^omes possibles ${\sf P} = {\sf
P} \big( j^5 g \big)$ qui, lorsqu'on remplace $g_i'$, $g_i''$,
$g_i'''$, $g_i''''$ et $g_i'''''$ par ces valeurs, ont la vertu de
faire dispara\^{\i}tre toutes les d\'eriv\'ees intempestives $\phi''$,
$\phi'''$, $\phi''''$ et $\phi '''''$ d'ordre $\geqslant 2$ de $\phi$,
de telle sorte que ${\sf P} \big( j^5 g \big) = {\phi '}^m {\sf P}
\big( j^5 f)$ pour un $m \in \N$. Ce calcul d'\'elimination
heuristique fonctionne pour $\kappa = 2$ et $\kappa = 3$ en dimensions
$\nu = 2$ et $\nu = 3$, mais il se complexifie au-del\`a et nous ne
poursuivrons pas la recherche en prenant cette optique.

\smallskip

Donnons toutefois la formule g\'en\'erale de d\'erivation compos\'ee,
dite de {\sl Fa\`a di Bruno}, bien connue dans le cas classique d'une seule
variable $z \in \C$.

\def\thetheorem{\!}\begin{theorem} 
Pour tout entier $\kappa \geqslant 1$, la d\'eriv\'ee d'ordre $\kappa$
de chaque fonction compos\'ee $g_i (z) := f_i \circ \phi (z)$ $( 1
\leqslant i \leqslant \nu)$ par rapport \`a la variable $z \in \C$
s'exprime comme polyn\^ome \`a coefficients entiers en les d\'eriv\'ees
de $f_i$ et de $\phi${\rm \,:}
\[
\boxed{
\aligned
g_i^{(\kappa)}
& 
=
\sum_{d=1}^\kappa
\
\sum_{1\leqslant\lambda_1<\cdots<\lambda_d\leqslant\kappa}
\
\sum_{\mu_1\geqslant 1,\dots,\mu_d\geqslant 1}
\
\sum_{\mu_1\lambda_1+\cdots+\mu_d\lambda_d=\kappa}
\\
& \
\ \ \ \ \ 
\frac{\kappa !}{(\lambda_1!)^{\mu_1}\ \mu_1! 
\cdots
(\lambda_d!)^{\mu_d}\ \mu_d!}
\
\left(
\phi^{(\lambda_1)}
\right)^{\mu_1}
\cdots\cdots
\left(
\phi^{(\lambda_d)}
\right)^{\mu_d}
\
f_i^{(\mu_1+\cdots+\mu_d)}.
\endaligned
}
\]
\end{theorem}

Dans tout ce qui va suivre, par souci de simplicit\'e, nous nous
restreindrons dor\'enavant \`a la dimension $\nu = 2$\,; des
g\'en\'eralisations en dimension sup\'erieure appara\^{\i}tront
en temps voulu.

\section*{\S3.~D\'eterminants $2 \times 2$ et jets d'ordres 3 et 4}

\noindent{\bf Wronskien et g\'en\'eralisations.}
En examinant $g_1'$, $g_2'$, $g_1''$ et $g_2''$, on constate
l'invarian\-ce par reparam\'etrisation du {\sl wronskien}, d\'efini 
comme \'etant le d\'eterminant $2 \times
2$\,:
\[
\Delta^{1,2}
:=
\left\vert
\begin{array}{cc}
f_1'&f_2'
\\
f_1''&f_2''
\end{array}
\right\vert,
\]
et ce, gr\^ace au calcul \'el\'ementaire suivant\,:
\[
\aligned
\square^{1,2}
:=
\left\vert
\begin{array}{cc}
g_1'&g_2'
\\
g_1''&g_2''
\end{array}
\right\vert
&
=
\left\vert
\begin{array}{cc}
\phi'f_1'&\phi'f_2'
\\
\phi''f_1'+{\phi'}^2f_1''&\phi''f_2'+{\phi'}^2f_2''
\end{array}
\right\vert
\\
&
=
\left\vert
\begin{array}{cc}
\phi'f_1'&\phi'f_2'
\\
{\phi'}^2f_1''&{\phi'}^2f_2''
\end{array}
\right\vert
\\
&
=
\big(\phi'\big)^3
\Delta^{1,2}.
\endaligned
\]

\noindent
Son poids $m$ est \'egal \`a 3. Ensuite, en \'eliminant de mani\`ere analogue
$\phi'''$ et $\phi''$ parmi les six \'equations donnant $g_1'$, $g_2'$,
$g_1''$, $g_2''$, $g_1'''$ et $g_2'''$\,\,---\,\,ou bien en
proc\'edant d'une mani\`ere alternative\,\,---, on trouve les deux
invariants de poids $m = 5$\,:
\[
\small
\aligned
\big[
g_1'g_2'''
-
g_1'''g_2'
\big]
g_1'
-
3\,
\big[
g_1'g_2''
-
g_1''g_2'
\big]
g_1''
&
=
\big(
\phi'
\big)^5
\big[
f_1'f_2'''
-
f_1'''f_2'
\big]
f_1'
-
3\,
\big[
f_1'f_2''
-
f_1''f_2'
\big]
f_1'',
\\
\big[
g_1'g_2'''
-
g_1'''g_2'
\big]
g_2'
-
3\,
\big[
g_1'g_2''
-
g_1''g_2'
\big]
g_2''
&
=
\big(
\phi'
\big)^5
\big[
f_1'f_2'''
-
f_1'''f_2'
\big]
f_2'
-
3\,
\big[
f_1'f_2''
-
f_1''f_2'
\big]
f_2''.
\endaligned
\]

\noindent{\bf D\'eterminants $2 \times 2$ g\'en\'eralisant le wronskien.}
Il est commode de r\'e\'ecrire ces deux invariants de poids $m = 5$ sous
une forme contract\'ee en introduisant la notation\,:
\[
\Delta^{\alpha,\beta}
:=
\left\vert
\begin{array}{cc}
f_1^{(\alpha)}&f_2^{(\alpha)}
\\
f_1^{(\beta)}&f_2^{(\beta)}
\end{array}
\right\vert,
\]
pour tous entiers $\alpha, \beta \geqslant 1$, ce qui donne, pour $k =
1, 2$\,:
\[
\square^{1,3}\,g_k'
-
3\,\square^{1,2}\,g_k''
=
\big(\phi'\big)^5\,
\Big[
\Delta^{1,3}\,f_k'
-
3\,\Delta^{1,2}\,f_k''
\Big].
\]
Notons au passage la formule de d\'erivation bien connue qui sera utile
par la suite\,:
\[
\big[
\Delta^{\alpha,\beta}\big]'
=
\left\vert
\begin{array}{cc}
f_1^{(\alpha+1)}&f_2^{(\alpha)}
\\
f_1^{(\beta+1)}&f_2^{(\beta)}
\end{array}
\right\vert
+
\left\vert
\begin{array}{cc}
f_1^{(\alpha)}&f_2^{(\alpha+1)}
\\
f_1^{(\beta)}&f_2^{(\beta+1)}
\end{array}
\right\vert.
\]

\def\thelemma{\!}\begin{lemma} {\rm 
(\cite{ ro2007})} Le degr\'e de transcendance du corps engendr\'e par les
5 polyn\^omes invariants{\rm \,:}
\[
\boxed{
\aligned
f_1',
&
\ \ \ \ \
f_2',
\\
\Lambda^3
&
:=
\Delta^{1,2},
\\
\Lambda_1^5
&
:=
\Delta^{1,3}\,f_1'
-
3\,\Delta^{1,2}\,f_1''
\ \ \ \ \ \ 
\text{\it et}
\ \ \ \ \ \ 
\Lambda_2^5
:=
\Delta^{1,3}\,f_2'
-
3\,\Delta^{1,2}\,f_2''
\endaligned
}\,,
\]
au-dessus de $\C \big[ f_1', f_1'', f_1''', f_2', f_2'', f_2''' \big]$
est \'egal \`a 4, et pour pr\'eciser, les quatre polyn\^omes $f_1'$, $f_2'$,
$\Lambda_1^5$ et $\Lambda_2^5$ sont alg\'ebriquement ind\'ependants,
tandis que $\Lambda^3$ est quadratique sur $\C \big[ f_1', f_2',
\Lambda_1^5, \Lambda_2^5 \big]$ via la relation alg\'ebrique
imm\'ediatement v\'erifi\-able{\rm \,:}
\[
0
\equiv
f_2'\,\Lambda_1^5
-
f_1'\,\Lambda_2^5
-
3\,\Lambda^3\,\Lambda^3,
\]
et de plus, l'id\'eal des relations entre $f_1', f_2', \Lambda^3,
\Lambda_1^5, \Lambda_2^5$ est principal et se r\'eduit \`a cette unique
relation. 
\end{lemma}

Gr\^ace \`a ladite syzygie, on peut \'eliminer toutes les puissances de
$\Lambda^3$ sup\'erieures ou \'egales \`a 2 qui apparaissent dans un polyn\^ome
g\'en\'eral\,:
\[
\mathcal{P}
\big(
f_1',f_2',\Lambda^3,\Lambda_1^5,\Lambda_2^5
\big)
\]
exprim\'e en fonction de ces cinq polyn\^omes, et il ne reste alors que
des puissances de $\Lambda^3$ \'egales \`a $0$ ou \`a $1$. C'est un fait
remarquable que ces cinq polyn\^omes forment un syst\`eme g\'en\'erateurs,
comme l'\'enonce pr\'ecis\'ement le r\'esultat suivant.

\def\thetheorem{\!}\begin{theorem}
{\rm (\cite{ ro2007, de2007})} En dimension $\nu = 2$ et
au niveau $\kappa = 3$, tout polyn\^ome
${\sf P} \big( j^3 f \big)$ invariant par reparam\'etrisation
s'\'ecrit de mani\`ere unique{\rm \,:}
\[
{\sf P}
\big(
j^3f
\big)
=
\mathcal{P}
\big(
f_1',f_2',\Lambda_1^5,\Lambda_2^5
\big)
+
\Lambda^3\,
\mathcal{Q}
\big(
f_1',f_2',\Lambda_1^5,\Lambda_2^5
\big),
\]
avec des polyn\^omes quelconques 
$\mathcal{ P}$ et $\mathcal{ Q}$.
\end{theorem}

\noindent{\bf Travaux de calcul pour passer aux jets d'ordre
$\kappa = 4$ et $\kappa = 5$.} 

\begin{itemize}

\item[$\bullet$]
Trouver un syst\`eme de polyn\^omes invariants fondamentaux.

\item[$\bullet$]
Conna\^{\i}tre leur id\'eal des relations.

\item[$\bullet$]
Trouver une \'ecriture unique de tout polyn\^ome en les polyn\^omes
invariants fondamentaux.

\end{itemize}

\noindent{\bf Deux op\'erateurs de diff\'erentiation.}
Comment engendrer m\'ethodiquement une liste appropri\'ee de
polyn\^omes invariants fondamentaux pour les jets d'ordre $\kappa = 4$
ou $5$? Voici une premi\`ere id\'ee\,: si ${\sf P}$ est un polyn\^ome
invariant de {\sl poids}\, $m$, d\'efinissons la diff\'erentiation
``covariante''\,:
\[
{\sf P}_{;\,k}
:=
f_k'\,{\sf P}'
-
m\,
f_k''\,{\sf P},
\]
o\`u ${\sf P}' = {\sf P} \big( j^{\kappa + 1} f \big)$ s'obtient en
diff\'erentiant ${\sf P} = {\sf P} \big( j^\kappa f \big)$ par rapport
\`a la variable $z \in \C$.

\def\thelemma{\!}\begin{lemma} 
Ces deux op\'erateurs de diff\'erentiation $(\cdot)_{;\, 1}$ et
$(\cdot)_{;\, 2}$ satisfont la r\`egle de Leibniz{\rm \,:}
\[
\big(
{\sf P}\cdot{\sf Q}
\big)_{;\,k}
=
{\sf P}_{;\,k}\cdot{\sf Q}
+
{\sf P}\cdot{\sf Q}_{;\,k},
\]
et ils produisent\footnote{\, La d\'emonstration, laiss\'ee au lecteur
qui d\'esirerait anticiper, appara\^{\i}tra dans un instant comme
absorb\'ee par une observation plus g\'en\'erale. } des polyn\^omes
invariants par reparam\'etrisation ${\sf P}_{ \, ;1}$ et ${\sf P}_{ \,
;2}$ qui sont tous deux de poids $m + 2$. \hfill$\square$
\end{lemma}

\noindent{\bf Exemple.} On 
v\'erifie imm\'ediatement\,:
\[
\aligned
&
\big(f_2'\big)_{;\,1}
=
f_1'f_2''
-
f_1''f_2'
\equiv
\Lambda^3
\equiv
-
\big(f_1'\big)_{;\,2},
\\
&
{\Lambda^3}_{;\,i}
=
f_i'\,
\Delta^{1,3}
-
3\,f_i''\,\Delta^{1,2}
\equiv
\Lambda_i^5.
\endaligned
\]
Ensuite, \`a l'\'etage $\kappa = 4$, on est naturellement conduit \`a
introduire les quatre nouveaux invariants\,:
\[
\Lambda_{1,1}^7
:=
\big({\Lambda_1^5}\big)_{;\,1},
\ \ \ \ \ \ \ \
\Lambda_{1,2}^7
:=
\big({\Lambda_1^5}\big)_{;\,2},
\ \ \ \ \ \ \ \
\Lambda_{2,1}^7
:=
\big({\Lambda_2^5}\big)_{;\,1},
\ \ \ \ \ \ \ \
\Lambda_{2,2}^7
:=
\big({\Lambda_2^5}\big)_{;\,2},
\]
dont l'expression explicite sera fournie dans un instant.

\smallskip\noindent{\bf Produit crois\'e entre invariants.}
Comment donner corps \`a l'id\'ee qu'il 
doit exister des diff\'erentiations
covariantes, non seulement par rapport \`a 
$f_1'$ et $f_2'$, mais aussi par rapport
\`a n'importe quel invariant?

Supposons donc connus deux polyn\^omes homog\`enes 
invariants ${\sf P}$ de poids $m$ et
${\sf Q}$ de poids $n$\,:
\[
\aligned
{\sf P}
\big(
j^\kappa 
g
\big)
&
=
{\phi'}^m\,
{\sf P}
\big(
(j^\kappa f)
\circ
\phi
\big),
\\
{\sf Q}
\big(
j^\tau 
g
\big)
&
=
{\phi'}^n\,
{\sf Q}
\big(
(j^\tau f)
\circ
\phi
\big),
\endaligned
\]
o\`u l'on a pos\'e $g := f \circ \phi$. Diff\'erentier un polyn\^ome
par rapport \`a la variable $z \in \C$ revient \`a lui appliquer
l'op\'erateur de {\sl diff\'erentiation totale}\,:
\[
{\sf D}
:=
\sum_{\lambda\in\N}\,
\frac{\partial (\bullet)}{\partial f^{(\lambda)}}\,
\cdot\,
f^{(\lambda+1)},
\]
ce qui nous donne ici\,:
\[
\aligned
\big[
{\sf D}{\sf P}
\big]
\big(
j^{\kappa+1}g
\big)
&
=
m\,\phi''\,
{\phi'}^{m-1}\,
{\sf P}
\big(
(j^\kappa f)\circ\phi
\big)
+
{\phi'}^m\,\phi'\,
\big[
{\sf D}{\sf P}
\big]
\big(
(j^{\kappa+1}f)\circ\phi
\big)
\\
\big[
{\sf D}{\sf Q}
\big]
\big(
j^{\tau+1}g
\big)
&
=
n\,\phi''\,
{\phi'}^{n-1}\,
{\sf Q}
\big(
(j^\kappa f)\circ\phi
\big)
+
{\phi'}^n\,\phi'\,
\big[
{\sf D}{\sf Q}
\big]
\big(
(j^{\tau+1}f)\circ\phi
\big)
\endaligned
\]
et pour faire dispara\^{\i}tre la d\'eriv\'ee seconde $\phi''$, il suffit
d'effectuer un produit crois\'e, autrement dit de former le
d\'eterminant $2 \times 2$\,:
\[
\aligned
&
\left\vert
\begin{array}{cc}
\big[{\sf D}{\sf P}\big]\big(j^{\kappa+1}g\big)\
&\
m\,{\sf P}\big(j^\kappa g\big)
\\
\big[{\sf D}{\sf Q}\big]\big(j^{\tau+1}g\big)\
&\
n\,{\sf Q}\big(j^\tau g\big)
\end{array}
\right\vert
=
\\
&
=
\left\vert
\begin{array}{cc}
m\,\phi''\,
{\phi'}^{m-1}\,
{\sf P}
\big(
(j^\kappa f)\circ\phi
\big)
+
{\phi'}^{m+1}\,
\big[
{\sf D}{\sf P}
\big]
\big(
(j^{\kappa+1}f)\circ\phi
\big)\
&\
m\,{\phi'}^m\,{\sf P}\big((j^\kappa f)\circ\phi\big)
\\
n\,\phi''\,
{\phi'}^{n-1}\,
{\sf Q}
\big(
(j^\kappa f)\circ\phi
\big)
+
{\phi'}^{n+1}\,
\big[
{\sf D}{\sf Q}
\big]
\big(
(j^{\tau+1}f)\circ\phi
\big)
&\
n\,{\phi'}^n\,{\sf Q}\big((j^\kappa f)\circ\phi\big)
\end{array}
\right\vert
\\
&
\ \ \ \ \ \ \ \ \ \ \ \ \ \ \ \ \ \ \ \ \ \
=
\left\vert
\begin{array}{cc}
{\phi'}^{m+1}\,
\big[
{\sf D}{\sf P}
\big]
\big(
(j^{\kappa+1}f)\circ\phi
\big)\
&\
m\,{\phi'}^m\,{\sf P}\big((j^\kappa f)\circ\phi\big)
\\
{\phi'}^{n+1}\,
\big[
{\sf D}{\sf Q}
\big]
\big(
(j^{\tau+1}f)\circ\phi
\big)
&\
n\,{\phi'}^n\,{\sf Q}\big((j^\kappa f)\circ\phi\big)
\end{array}
\right\vert
\\
&
\ \ \ \ \ \ \ \ \ \ \ \ \ \ \ \ \ \ \ \ \ \
\ \ \ \ \ \ \ \ \ \ \ \ \ \ \ \ \ \ \ \ \ \
=
{\phi'}^{m+n+1}\,
\left\vert
\begin{array}{cc}
\big[{\sf D}{\sf P}\big]\big(j^{\kappa+1}f\big)\
&\
m\,{\sf P}\big(j^\kappa f\big)
\\
\big[{\sf D}{\sf Q}\big]\big(j^{\tau+1}f\big)\
&\
n\,{\sf Q}\big(j^\tau f\big)
\end{array}
\right\vert
\endaligned
\]
qui s'av\`ere ainsi constituer un nouvel invariant de poids $m + n +
1$.

\smallskip\noindent{\bf Notation crochet $\big[ \cdot, \, \cdot\big]$.}
Ainsi, toute paire d'invariants produit automatiquement un 
nouvel invariant\,:
\[
\boxed{
\big[
{\sf P},\,{\sf Q}
\big]
:=
n\,{\sf D}{\sf P}\cdot{\sf Q}
-
m\,{\sf P}\cdot{\sf D}{\sf Q}}\,,
\]
qui est \'evidemment antisym\'etrique 
par rapport au couple $\big( {\sf P}, {\sf Q} \big)$.

\smallskip\noindent{\bf Observation.}
Ces crochets pond\'er\'es g\'en\'eralisent les deux op\'erateurs de
diff\'erentia\-tion pr\'ec\'edents\,:
\[
{\sf P}_{;\,k}
\equiv
\big[
{\sf P},\,f_k'
\big].
\]
De plus, ils satisfont \`a la r\`egle de Leibniz\,:
\[
\big[
{\sf P},\,{\sf Q}{\sf R}
\big]
=
\big[
{\sf P},\,{\sf Q}
\big]\,{\sf R}
+
\big[{\sf P},\,{\sf R}\big]\,{\sf Q},
\]
de telle sorte que l'op\'erateur $\big[ \bullet, {\sf Q} \big]$, \`a
savoir\,: ${\sf P} \longmapsto \big[ {\sf P}, \, {\sf Q} \big]$, peut
\^etre consid\'er\'e comme un op\'erateur de d\'erivation.

\def\thelemma{\!}\begin{lemma}
Pour tout triplet $\big( {\sf P}, \, {\sf Q}, \, {\sf R} \big)$
d'invariants de poids $m, n, o$, l'identit\'e suivante de type Jacobi
est satisfaite{\rm \,:}
\def\theequation{$\mathcal{J}ac$}\begin{equation}
\boxed{
0
\equiv
\big[\big[{\sf P},\,{\sf Q}\big],\,{\sf R}\big]
+
\big[\big[{\sf R},\,{\sf P}\big],\,{\sf Q}\big]
+
\big[\big[{\sf Q},\,{\sf R}\big],\,{\sf P}\big]
}\,.
\end{equation}
\end{lemma}

\noindent{\em Preuve.}
D\'eveloppons le premier double crochet\,:
\[
\aligned
\big[\big[{\sf P},\,{\sf Q}\big],\,{\sf R}\big]
&
=
\big[
n\,{\sf D}{\sf P}\cdot{\sf Q}
-
m\,{\sf P}\cdot{\sf D}{\sf Q},\,
{\sf R}
\big]
\\
&
=
o\big(
n\,{\sf D}{\sf D}{\sf P}\cdot{\sf Q}
+
(n-m)\,{\sf D}{\sf P}\cdot{\sf D}{\sf Q}
-
m\,{\sf P}\cdot{\sf D}{\sf D}{\sf Q}
\big)
{\sf R}
-
\\
&
\ \ \ \
-
(m+n+1)
\big(
n\,{\sf D}{\sf P}\cdot{\sf Q}
-
m\,{\sf P}\cdot{\sf D}{\sf Q}
\big)
{\sf D}{\sf R}
\\
&
=
no\,{\sf D}{\sf D}{\sf P}\cdot{\sf Q}\cdot{\sf R}
+
(n-m)o\,{\sf D}{\sf P}\cdot{\sf D}{\sf Q}\cdot{\sf R}
-
mo\,{\sf P}\cdot{\sf D}{\sf D}{\sf Q}\cdot{\sf R}
-
\\
&
\ \ \ \
-
(m+n+1)n\,{\sf D}{\sf P}\cdot{\sf Q}\cdot{\sf D}{\sf R}
+
(m+n+1)m\,{\sf P}\cdot{\sf D}{\sf Q}\cdot{\sf D}{\sf R}.
\endaligned
\]
Il suffit alors de constater que l'annulation identique
de la somme suivante\,:
\[
\aligned
0
&
\equiv
no\,{\sf D}{\sf D}{\sf P}\cdot{\sf Q}\cdot{\sf R}
+
(n-m)o\,{\sf D}{\sf P}\cdot{\sf D}{\sf Q}\cdot{\sf R}
-
mo\,{\sf P}\cdot{\sf D}{\sf D}{\sf Q}\cdot{\sf R}
-
\\
&
\ \ \ \
-
(m+n+1)n\,{\sf D}{\sf P}\cdot{\sf Q}\cdot{\sf D}{\sf R}
+
(m+n+1)m\,{\sf P}\cdot{\sf D}{\sf Q}\cdot{\sf D}{\sf R}
+
\\
&
\ \ \ \
+
mn\,{\sf D}{\sf D}{\sf R}\cdot{\sf P}\cdot{\sf Q}
+
(m-o)n\,{\sf D}{\sf R}\cdot{\sf D}{\sf P}\cdot{\sf Q}
-
on\,{\sf R}\cdot{\sf D}{\sf D}{\sf P}\cdot{\sf Q}
-
\\
&
\ \ \ \
-
(o+m+1)m\,{\sf D}{\sf R}\cdot{\sf P}\cdot{\sf D}{\sf Q}
+
(o+m+1)o\,{\sf R}\cdot{\sf D}{\sf P}\cdot{\sf D}{\sf Q}
+\\
&
\ \ \ \
+
om\,{\sf D}{\sf D}{\sf Q}\cdot{\sf R}\cdot{\sf P}
+
(o-n)m\,{\sf D}{\sf Q}\cdot{\sf D}{\sf R}\cdot{\sf P}
-
nm\,{\sf Q}\cdot{\sf D}{\sf D}{\sf R}\cdot{\sf P}
-
\\
&
\ \ \ \
-
(n+o+1)o\,{\sf D}{\sf Q}\cdot{\sf R}\cdot{\sf D}{\sf P}
+
(n+o+1)n\,{\sf Q}\cdot{\sf D}{\sf R}\cdot{\sf D}{\sf P},
\endaligned
\]
est effectivement satisfaite.

\smallskip\noindent{\bf Exemple.} Avec ${\sf P} := f_1'$,
${\sf Q} := f_2'$ et ${\sf R} := \Lambda^3$, nous obtenons\,:
\[
\aligned
0
&
\equiv
\big[\big[f_1',\,f_2'\big],\,\Lambda^3\big]
+
\big[\big[\Lambda^3,\,f_1'\big],\,f_2'\big]
+
\big[\big[f_2',\,\Lambda^3\big],\,f_1'\big]
\\
&
\equiv
0
+
\big[\Lambda_1^5,\,f_2'\big]
-
\big[\Lambda_2^5,\,f_1'\big]
\\
&
\equiv
\Lambda_{1,2}^7
-
\Lambda_{2,1}^7.
\endaligned
\]
Cette relation sera confirm\'ee par les expressions explicites de
$\Lambda_{ 1, 2}^7$ et $\Lambda_{ 2, 1}^7$.

\smallskip\noindent{\bf Gen\`ese des invariants fondamentaux.}
Crucialement, il semblerait que l'on puisse engendrer tous les
polyn\^omes invariants en les jets d'un ordre $\kappa \geqslant 1$
quelconque, juste en calculant par r\'ecurrence tous les crochets
possibles, d'un \'etage de jets $\lambda$ \`a l'\'etage sup\'erieur
$\lambda + 1$. Cette id\'ee conjecturale (\cite{ de2007}), sur laquelle
nous donnerons plus de pr\'ecision ult\'erieurement, est renforc\'ee par le
fait que dans la th\'eorie classique des invariants pour une forme
binaire $\sum_{ i = 0 }^\kappa \, a_i\, x^i y^{ \kappa - i}$ de degr\'e
$\kappa$ par rapport \`a l'action lin\'eaire standard
de ${\sf SL}_2 (\C)$\,:
\[
\aligned
&
x\longmapsto 
\overline{x}
=
\alpha x+\beta y, 
\ \ \ \ \ \ \ \ \ \ \
y\longmapsto
\overline{y}
=
\gamma x+\delta y,
\ \ \ \ \ \ \ \ \ \ \ 
1
=
\alpha\delta-\beta\gamma,
\\
&
{\textstyle{\sum_{i=0}^\kappa}}\,a_i\,x^iy^{\kappa-i}\,
\longmapsto
{\textstyle{\sum_{i=0}^\kappa}}\,\overline{a}_i\,
\overline{x}^i\overline{y}^{\kappa-i},
\\
&
a_i
=
{\textstyle{\sum_{l=0}^\kappa}}\,\overline{a}_l\,
{\textstyle{\sum_{j=\max(0,i+l-\kappa)}^{\min(i,l)}}}\,
C_i^j\,
C_{\kappa-i}^{l-j}\,
\alpha^j\,
\beta^{l-j}\,
\gamma^{i-j}\,
\delta^{\kappa+j-i-l},
\ \ \ \ \ \ \ \ \
C_p^q
=
{\textstyle{
\frac{p!}{q!\,\,(p-q)!}}},
\endaligned
\]
on sait \'etablir que deux proc\'ed\'es alg\'ebriques
\'el\'ementaires, \`a savoir le ``processus $\Omega$'' et le
``processus $\sigma$'' ({\it cf.} \cite{ ol1999}) permettent
d'engendrer un syst\`eme fondamental de polyn\^omes ${\sf P} = {\sf P}
\big( a_0, a_1, \dots, a_\kappa \big)$ qui sont invariants\,:
\[
{\sf P}\big(
a_0,a_1,\dots,a_\kappa
\big)
=
{\sf P}
\big(
\overline{a}_0,\overline{a}_1,\dots,\overline{a}_\kappa
\big).
\]

\noindent{\bf Reconstitution par crochets des invariants 
connus.}
Pour passer des jets d'ordre 1 aux jets d'ordre 2, seul un crochet (au
signe pr\`es) peut \^etre form\'e\,:
\[
\big[f_1',\,f_2'\big]
=
-
\big[f_2',\,f_1'\big]
=
-
\Lambda^3.
\]
Pour passer des jets d'ordre 2 aux jets d'ordre 3, on peut former
trois crochets\,:
\[
\big[
\Lambda^3,\,f_1'
\big]
\ \ \ \ \ \ \ \ \ \ \ \ \
\big[
\Lambda^3,\,f_2'
\big]
\ \ \ \ \ \ \ \ \ \ \ \ \
\big[
\Lambda^3,\,\Lambda^3
\big],
\]
le dernier \'etant trivialement nul, et l'on v\'erifie imm\'ediatement
que les deux premiers fournissent les deux invariants propres \`a
l'\'etage $\kappa = 3$\,:
\[
\big[
\Lambda^3,\,f_i'
\big]
=
\Delta^{1,3}\,f_i'
-
3\,\Delta^{1,2}\,f_i''.
\]
Pour passer aux jets d'ordre 4, l'ensemble des crochets que l'on peut
former s'identi\-fie 
\`a la collection des d\'eterminants $2 \times 2$ de la
matrice matrice $2 \times 5$\,:
\[
\left\vert\!\left\vert
\begin{array}{ccccc}
f_1' \ \ & \ \ f_2' \ \ & \ \ 3\,\Lambda^3 \ \ & 
\ \ 5\,\Lambda_1^5 \ \ & \ \ 5\,\Lambda_2^5
\\
{\sf D}f_1' \ \ & \ \ {\sf D}f_2' \ \ & \ \
{\sf D}\Lambda^3 \ \ & \ \ {\sf D}\Lambda_1^5 \ \ & \ \
{\sf D}\Lambda_2^5
\end{array}
\right\vert\!\right\vert,
\]
ce qui fait au total de $C_{ 5}^2 = 10$ crochets, mais en tenant
compte du fait que nous connaissons d\'ej\`a les trois
mineurs\,\,---\,\,calcul\'es \`a l'\'etage $\kappa = 3$\,\,---\,\,de la
sous-matrice\,:
\[
\left\vert\!\left\vert
\begin{array}{ccc}
f_1' \ \ & \ \ f_2' \ \ & \ \ 3\,\Lambda^3
\\
{\sf D}f_1' \ \ & \ \ {\sf D}f_2' \ \ & \ \
{\sf D}\Lambda^3
\end{array}
\right\vert\!\right\vert,
\]
ce sont exactement sept nouveaux crochets qui apparaissent\,:
\[
\big[
\Lambda_i^5,\,f_j'
\big],
\ \ \ \ \ \ \ \ \
\big[
\Lambda_i^5,\,\Lambda^3
\big],
\ \ \ \ \ \ \ \ \
\big[
\Lambda_1^5,\,\Lambda_2^5
\big].
\]

\noindent{\bf Relations pl\"ucke\-riennes.}
Cependant, le calcul complet des crochets doit tenir compte des
relations de Pl\"ucker qui existent au niveau des variables initiales
des espaces de jets. En effet, l'id\'eal des relations pl\"ucke\-riennes
entre les $f_i^{ (\lambda)}$, $1\leqslant \lambda \leqslant 5$ et les
$\Delta^{ \alpha, \beta}$, $1 \leqslant \alpha < \beta \leqslant 5$
est engendr\'e par deux familles quadratiques de relations identiquement
satisfaites\,: dans la premi\`ere famille\,:
\[
0
\equiv
\Delta^{\beta,\gamma}\cdot f_i^{(\alpha)}
-
\Delta^{\alpha,\gamma}\cdot f_i^{(\beta)}
+
\Delta^{\alpha,\beta}\cdot f_i^{(\gamma)},
\]
$i$ est \'egal \`a $1$ ou \`a $2$, et les indices sup\'erieurs satisfont
$1\leqslant \alpha < \beta < \gamma \leqslant 5$, ce qui donne $10
\times 2 = 20$ relations\,; et dans la seconde famille\,:
\[
0
\equiv
\Delta^{\alpha,\delta}\cdot\Delta^{\beta,\gamma}
-
\Delta^{\alpha,\gamma}\cdot\Delta^{\beta,\delta}
+
\Delta^{\alpha,\beta}\cdot\Delta^{\gamma,\delta},
\]
les indices sup\'erieurs satisfont $1 \leqslant \alpha < \beta < \gamma <
\delta \leqslant 5$, ce qui donne $4$ relations.

En v\'erit\'e, seules les deux paires de relations suivantes, extraites de
la premi\`ere familles, seront utiles \`a l'\'etage des jets d'ordre 
$\kappa = 5$\,:
\[
\boxed{
\aligned
0
&
\equiv
\Delta^{2,3}\,f_i'
-
\Delta^{1,3}\,f_i''
+
\underline{
\Delta^{1,2}\,f_i'''}
\\
0
&
\equiv
\Delta^{2,4}\,f_i'
-
\Delta^{1,4}\,f_i''
+
\underline{
\Delta^{1,2}\,f_i''''}
\endaligned
}\,,
\]
et \`a l'\'etage $\kappa = 4$, seule la premi\`ere paire peut \^etre utilis\'ee,
tandis qu'aucune relation pl\"ucke\-rienne n'intervient aux \'etages $\kappa
\leqslant 3$. Il faut en outre attendre $\kappa = 6$ pour que la
premi\`ere relation de la seconde famille, \`a savoir\,: $0 \equiv \Delta^{
1,4}\, \Delta^{ 2,3} - \Delta^{ 1,3}\, \Delta^{ 2,4} + \Delta^{ 1,2}\,
\Delta^{ 3, 4}$ commence \`a interf\`erer, mais nous n'entreprendrons
pas l'\'etude de $\mathcal{ DS}_2^6$ dans cet article.

\smallskip\noindent{\bf Normalisations pr\'ealables des
diff\'erentielles totales.}
Ainsi, nous sommes conduits \`a normaliser ${\sf D} \Lambda_i^5$
avant de calculer $\Lambda_{ i,j}^7$, en d\'eveloppant tout d'abord\,:
\[
\aligned
{\sf D}\Lambda_i^5
&
=
\Delta^{1,4}\,f_i'
+
\Delta^{2,3}\,f_i'
+
\Delta^{1,3}\,f_i''
-
3\,\Delta^{1,3}\,f_i''
-
3,\Delta^{1,2}\,f_i'''
\\
&
=
\Delta^{1,4}\,f_i'
+
\Delta^{2,3}\,f_i'
-
2\,\Delta^{1,3}\,f_i''
-
3\,\Delta^{1,2}\,f_i''',
\endaligned
\]
expression dans laquelle nous pouvons remplacer $\underline{ \Delta^{
1,2}\, f_i''' }$ par $- \Delta^{ 2, 3}\, f_i' + \Delta^{ 1,3}\, f_i''$
afin d'\'eliminer toute pr\'esence de $f_i'''$, ce qui nous donne une
expression normalis\'ee et compacte ne contenant que trois termes\,:
\[
{\sf D}\Lambda_i^5
=
\Delta^{1,4}\,f_i'
+
4\,\Delta^{2,3}\,f_i'
-
5\,\Delta^{1,3}\,f_i''.
\]
Achevons donc le calcul de la premi\`ere famille de crochets $\big[
\Lambda_i^5,\, f_j' \big]$ en fournissant tous les d\'etails
interm\'ediaires\,:
\[
\aligned
\big[
\Lambda_i^5,\,f_j'
\big]
&
=
{\sf D}\Lambda_i^5\cdot f_j'
-
5\,\Lambda_i^5\cdot f_j''
\\
&
=
\Big(
\Delta^{1,4}\,f_i'
+
4\,\Delta^{2,3}\,f_i'f_j'
-
5\,\Delta^{1,3}
f_i''
\Big)\cdot f_j'
-
5\,\Big(
\Delta^{1,3}\,f_i'
-
3\,\Delta^{1,2}\,f_i''
\Big)\cdot f_j''
\\
&
=
\Delta^{1,4}\,f_i'f_j'
+
4\,\Delta^{2,3}\,f_i'f_j'
-
5\,\Delta^{1,3}\big(f_i''f_j'+f_i'f_j'')
+
15\,\Delta^{1,2}\,f_i''f_j''
\\
&
=:
\Lambda_{i,j}^7.
\endaligned
\]
Ici, la sym\'etrie indicielle $\Lambda_{ 1,2 }^7 = \Lambda_{ 2, 1}^7$
impos\'ee {\it a priori}\, par l'identit\'e de Jacobi montre que l'on
devrait se dispenser de $\Lambda_{ 2, 1 }^7$ (ou de $\Lambda_{ 1, 2
}^7$\big) dans une liste minimale de polyn\^omes invariants
fondamentaux.

\smallskip\noindent{\bf Invariant de poids 8.}
Ensuite, gr\^ace \`a notre normalisation pr\'ealable
de ${\sf D} \Lambda_i^5$, nous pouvons calculer
proprement chacun des deux crochets ($i = 1, 2$)\,:
\[
\small
\aligned
\big[
\Lambda_i^5,\,\Lambda^3
\big]
&
=
3\,{\sf D}\Lambda_i^5\cdot\Lambda^3
-
5\,\Lambda_i^5\cdot{\sf D}\Lambda^3
\\
&
=
\Big(
3\,\Delta^{1,4}\,f_i'
+
12\,\Delta^{2,3}\,f_i'
-
15\,\Delta^{1,3}\,f_i''
\Big)\cdot\Delta^{1,2}
-
\\
&
\ \ \ \ \ \ \ \ \ \ \ \
-
\Big(
5\,\Delta^{1,3}\,f_i'
-
15\,\Delta^{1,2}\,f_i''
\Big)\cdot\Delta^{1,3}
\\
&
=
3\,\Delta^{1,4}\,\Delta^{1,2}\,f_i'
+
12\,\Delta^{2,3}\,\Delta^{1,2}\,f_i'
-
5\,\Delta^{1,3}\,\Delta^{1,3}\,f_i'
\\
&
=
f_i'
\big(
3\,\Delta^{1,4}\,\Delta^{1,2}
+
12\,\Delta^{2,3}\,\Delta^{1,2}
-
5\,\Delta^{1,3}\,\Delta^{1,3}
\big)
\\
&
\equiv
f_i'\,M^8,
\endaligned
\]
o\`u le nouvel invariant 
\[
\boxed{
\aligned
M^8
:=
&\
\frac{1}{f_i'}\,
\big[
\Lambda_i^5,\,\Lambda^3
\big]
\\
=
&\
\
3\,\Delta^{1,4}\,\Delta^{1,2}
+
12\,\Delta^{2,3}\,\Delta^{1,2}
-
5\,\Delta^{1,3}\,\Delta^{1,3}
\endaligned
}
\] 
doit \^etre introduit, parce que le r\'esultat est divisible par $f_i'$.

\smallskip\noindent{\bf Question.} 
Pourquoi et comment doit-on \^etre conduit \`a diviser parfois les
crochets pour acc\'eder v\'eritablement \`a de nouveaux invariants
fondamentaux?

\smallskip\noindent{\bf Fin du passage \`a l'\'etage 
$\kappa = 4$} Enfin, calculons et examinons le dernier crochet
possible, \`a nouveau en fournissant scrupuleusement tous les d\'etails
interm\'ediaires\,:
\[
\small
\aligned
\big[
\Lambda_1^5,\,\Lambda_2^5
\big]
&
=
5\,{\sf D}\Lambda_1^5\cdot\Lambda_2^5
-
5\,\Lambda_1^5\cdot{\sf D}\Lambda_2^5
\\
&
=
5\Big(
\Delta^{1,4}\,f_1'
+
4\,\Delta^{2,3}\,f_1'
-
5\,\Delta^{1,3}\,f_1''
\Big)
\cdot
\Big(
\Delta^{1,3}\,f_2'
-
3\,\Delta^{1,2}\,f_2''
\Big)
-
\\
&
\ \ \ \ \ \ \ \ \ \
-
5\Big(
\Delta^{1,3}\,f_1'
-
3\,\Delta^{1,2}\,f_1''
\Big)\cdot
\Big(
\Delta^{1,4}\,f_2'
+
4\,\Delta^{2,3}\,f_2'
-
5\,\Delta^{1,3}\,f_2''
\Big)
\\
&
=
-
15\,\Delta^{1,4}\,\Delta^{1,2}\,f_1'f_2''
-
60\,\Delta^{2,3}\,\Delta^{1,2}\,f_1'f_2''
-
25\,\Delta^{1,3}\,\Delta^{1,3}\,f_1''f_2'
+
\\
&
\ \ \ \ 
+
15\,\Delta^{1,4}\,\Delta^{1,2}\,f_2'f_1''
+
60\,\Delta^{2,3}\,\Delta^{1,2}\,f_2'f_1''
+
25\,\Delta^{1,3}\,\Delta^{1,3}\,f_2''f_1'
\\
&
=
-
15\,\Delta^{1,4}\,\Delta^{1,2}\,\Delta^{1,2}
-
15\,\Delta^{2,3}\,\Delta^{1,2}\,\Delta^{1,2}
+
25\,\Delta^{1,3}\,\Delta^{1,3}\,\Delta^{1,2}
\\
&
=
-
5\,\Delta^{1,2}
\Big(
3\,\Delta^{1,4}\,\Delta^{1,2}
+
12\,\Delta^{2,3}\,\Delta^{1,2}
-
5\,\Delta^{1,3}\,\Delta^{1,3}
\Big)
\\
&
\equiv
-
5\,\Lambda^3\,M^8.
\endaligned
\]
Le r\'esultat \'etant multiple des deux invariants d\'ej\`a connus
$\Lambda^3$ et $M^8$, il n'apporte rien de nouveau. Toutefois,
conservons trace de la relation\,:
\[
\big[
\Lambda_1^5,\,
\Lambda_2^5
\big]
=
-
5\,\Lambda^3\,M^8.
\]

\def\theproposition{\!}\begin{proposition}
{\rm (\cite{ de2007})} 
En dimension $\nu = 2$, les neuf polyn\^omes{\rm \,:}
\[
\small
\aligned
&
f_1',\ \ \ \ \ \ \ \ \ \ \ \ \
f_2', \ \ \ \ \ \ \ \ \ \ \ \ \ \ \
\Lambda^3
:=
\Delta^{1,2},
\ \ \ \ \ \ \
\Lambda_1^5
&
:=
\Delta^{1,3}\,f_1'
-
3\,\Delta^{1,2}\,f_1'',
\\
\Lambda_2^5
&
:=
\Delta^{1,3}\,f_2'
-
3\,\Delta^{1,2}\,f_2'',
\endaligned
\]
\[
\small
\aligned
\Lambda_{1,1}^7
&
:=
\big(
\Delta^{1,4}
+
4\,
\Delta^{2,3}\big)\,f_1'f_1'
-
10\,\Delta^{1,3}\,f_1'f_1''
+
15\,\Delta^{1,2}\,f_1''f_1'',
\\
\Lambda_{1,2}^7
&
:=
\big(
\Delta^{1,4}
+
4\,
\Delta^{2,3}\big)\,f_1'f_2'
-
5\,\Delta^{1,3}\big(f_1''f_2'
+
f_2''f_1'\big)
+
15\,\Delta^{1,2}\,f_1''f_2'',
\\
\Lambda_{2,2}^7
&
:=
\big(\Delta^{1,4}
+
4\,\Delta^{2,3}\big)\,f_2'f_2'
-
10\,\Delta^{1,3}\,f_2'f_2''
+
15\,\Delta^{1,2}\,f_2''f_2'',
\endaligned
\]
\[
M^8
:=
3\,\Delta^{1,4}\,\Delta^{1,2}
+
12\,\Delta^{2,3}\,\Delta^{1,2}
-
5\,\Delta^{1,3}\,\Delta^{1,3}
\]
forment un syst\`eme g\'en\'erateur de polyn\^omes invariants par
reparam\'etrisation pour les jets d'ordre $\kappa = 4$.
\end{proposition}

Cette proposition sera englob\'ee dans un \'enonc\'e plus pr\'ecis dont la
preuve appara\^{\i}tra dans la Section~5.

\section*{ \S4.~Invariants fondamentaux pour les jets d'ordre 5} 

\noindent{\bf D\'enombrement des crochets.} Pour s'\'elever de
l'\'etage $\kappa = 4$ \`a l'\'etage $\kappa = 5$, l'ensemble des
crochets que l'on peut former s'identifie \`a la collection des
d\'etermi\-nants $2 \times 2$ de la matrice matrice $2 \times 9$\,:
\[
\small
\left\vert\!\left\vert
\begin{array}{ccccccccc}
f_1' \ \ & \ \ f_2' \ \ & \ \ 3\,\Lambda^3 \ \ & 
\ \ 5\,\Lambda_1^5 \ \ & \ \ 5\,\Lambda_2^5 \ \ & \ \ 
7\,\Lambda_{1,1}^7 \ \ & \ \ 7\,\Lambda_{1,2}^7 \ \ & \ \ 
7\,\Lambda_{2,2}^7 \ \ & \ \ 8\,M^8
\\
{\sf D}f_1' \ \ & \ \ {\sf D}f_2' \ \ & \ \
{\sf D}\Lambda^3 \ \ & \ \ {\sf D}\Lambda_1^5 \ \ & \ \
{\sf D}\Lambda_2^5 \ \ & \ \ {\sf D}\Lambda_{1,1}^7 \ \ & \ \ 
{\sf D}\Lambda_{1,2}^7 \ \ & \ \ {\sf D}\Lambda_{2,2}^7 \ \ & \ \
{\sf D}M^8
\end{array}
\right\vert\!\right\vert,
\]
ce qui fait au total de $C_{ 9}^2 = 36$ crochets, mais il n'y en a en
fait que $36 - 10 = 26$ \`a calculer, en tenant compte du fait que nous
connaissons d\'ej\`a les $C_5^2 = 10$ mineurs\,\,---\,\,calcul\'es dans la
section pr\'ec\'edente\,\,---\,\,de la sous-matrice\,:
\[
\left\vert\!\left\vert
\begin{array}{ccccc}
f_1' \ \ & \ \ f_2' \ \ & \ \ 3\,\Lambda^3 \ \ & 
\ \ 5\,\Lambda_1^5 \ \ & \ \ 5\,\Lambda_2^5
\\
{\sf D}f_1' \ \ & \ \ {\sf D}f_2' \ \ & \ \
{\sf D}\Lambda^3 \ \ & \ \ {\sf D}\Lambda_1^5 \ \ & \ \
{\sf D}\Lambda_2^5
\end{array}
\right\vert\!\right\vert.
\]

\smallskip\noindent{\bf Heuristique.}
Nous sommes par cons\'equent amen\'es \`a penser\footnote{\, Cette
id\'ee sera discut\'ee plus avant dans la Section~7. } que tout polyn\^ome
invariant ${\sf P} \big( j^5 f \big)$ en les jets d'ordre 5 est un
polyn\^ome en les neuf pr\'ec\'edents polyn\^omes fondamentaux\,:
\[
f_1',\ \ \ \ \
f_2',\ \ \ \ \
\Lambda^3,\ \ \ \ \ 
\Lambda_1^5,\ \ \ \ \
\Lambda_{1,1}^7,\ \ \ \ \
\Lambda_{1,2}^7,\ \ \ \ \
\Lambda_{2,2}^7,\ \ \ \ \
M^8,
\]
auxquels on ajoute tous ceux qui sont obtenus par crochets \`a
l'\'etage sup\'erieur, apr\`es simplification, normalisation
pl\"ucke\-rienne, division \'eventuelle, et suppression des invariants
redondants. On voit imm\'ediatement que les nouveaux crochets \`a \'etudier
se distribuent en huit familles\,:
\[
\begin{array}{cc}
\big[\Lambda_{i,j}^7,\,f_k'\big],
\ \ \ \ \ \ \ \ \ \ \ \ \ \ \ 
&
\ \ \ \ \ \ \ \ \ \ \ \ \ \ \ 
\big[M^8,\,f_i'\big],
\\
\big[\Lambda_{i,j}^7,\,\Lambda^3\big],
\ \ \ \ \ \ \ \ \ \ \ \ \ \ \ 
&
\ \ \ \ \ \ \ \ \ \ \ \ \ \ \ 
\big[M^8,\,\Lambda^3\big],
\\
\big[\Lambda_{i,j}^7,\,\Lambda_k^5\big],
\ \ \ \ \ \ \ \ \ \ \ \ \ \ \ 
&
\ \ \ \ \ \ \ \ \ \ \ \ \ \ \ 
\big[M^8,\,\Lambda_i^5\big],
\\
\big[\Lambda_{i,j}^7,\,\Lambda_{k,l}^7\big],
\ \ \ \ \ \ \ \ \ \ \ \ \ \ \ 
&
\ \ \ \ \ \ \ \ \ \ \ \ \ \ \ 
\big[M^8,\,\Lambda_{i,j}^7\big].
\end{array}
\]
Avant de calculer et d'examiner 
tous ces crochets\,\,---\,\,t\^ache substantielle s'il
en est\,\,---, reprenons en main la liste de tous les crochets pr\'ec\'edents
({\it i.e.} invariants \`a l'\'etage $\kappa = 4$),
en les \'ecrivant avec des indices\,:
\[
\small
\boxed{
\aligned
f_i'
&
\\
\Lambda^3
&
:=
\Delta^{1,2}
\\
\Lambda_i^5
&
:=
\Delta^{1,3}\,f_i'
-
3\,\Delta^{1,2}\,f_i''
\\
\Lambda_{i,j}^7
&
:=
\Delta^{1,4}\,f_i'f_j'
+
4\,\Delta^{2,3}\,f_i'f_j'
-
5\,\Delta^{1,3}\big(f_i''f_j'+f_i'f_j'')
+
15\,\Delta^{1,2}\,f_i''f_j''
\\
M^8
&
:=
3\,\Delta^{1,4}\,\Delta^{1,2}
+
12\,\Delta^{2,3}\,\Delta^{1,2}
-
5\,\Delta^{1,3}\,\Delta^{1,3}
\endaligned
}\,.
\]

\noindent{\bf Remarque sur le choix des notations.} 
Nous utiliserons syst\'ematiquement les grandes lettres, telles que
``$\Lambda$'' (particuli\`erement facile \`a \'ecrire \`a la main),
``$M$'', ``$H$'', {\it etc.}, parce que leur taille les rend
disponibles pour recevoir non seulement le nombre total de ``$'$''
en indice sup\'erieur (poids de l'invariant), mais aussi, en indices
inf\'erieurs, la suite {\it ordonn\'ee}\, de $1$ ou de $2$ dont d\'epend
chaque mon\^ome de l'invariant en question.

\smallskip\noindent{\bf 
Normalisation pr\'ealable des diff\'erentiations totales.}
Nous avons donc huit famil\-les de crochets \`a calculer, et pour cela,
nous travaillerons avec les repr\'esentations indici\'ees de nos
invariants connus \`a l'\'etage $\kappa = 4$. Auparavant, nous devons
calculer \`a l'avance les deux expressions d\'eriv\'ees ${\sf D} \Lambda_{
i,j}^7$ et ${\sf D} M^8$, et les normaliser en tenant compte des
identit\'es pl\"ucke\-riennes, comme nous l'avons expliqu\'e
ci-dessus. Calculons donc, en \'eliminant $\Delta^{ 1,2} \, f_i'''$ et
$\Delta^{ 1,2}\, f_j'''$ \`a la quatri\`eme ligne\,:
\[
\small
\aligned
{\sf D}\Lambda_{i,j}^7
&
=
\Delta^{1,5}\,f_i'f_j'
+
\Delta^{2,4}\,f_i'f_j'
+
\Delta^{1,4}\big(f_i''f_j'+f_i'f_j'')
+
\\
&
\ \ \ \ \
+
4\,\Delta^{2,4}\,f_i'f_j'
+
4\,\Delta^{2,3}\big(f_i''f_j'+f_i'f_j''\big)
-
5\,\Delta^{1,4}\big(f_i''f_j'+f_i'f_j''\big)
-
\\
&
\ \ \ \ \
-
5\,\Delta^{2,3}\big(f_i''f_j'+f_i'f_j''\big)
-
5\,\Delta^{1,3}\big(f_i'''f_j'+2\,f_i''f_j''+f_i'f_j'''\big)
+
15\,\Delta^{1,3}\,f_i''f_j''
+
\\
&
\ \ \ \ \ 
+
15\,\underline{\Delta^{1,2}\,f_i'''}f_j''
+
15\,\underline{\Delta^{1,2}}\,f_i''\underline{f_j'''}
\\
&
=
\Delta^{1,5}\,f_i'f_j'
+
5\,\Delta^{2,4}\,f_i'f_j'
-
4\,\Delta^{1,4}\big(f_i''f_j'+f_i'f_j''\big)
-
\Delta^{2,3}\big(f_i''f_j'+f_i'f_j''\big)
-
\\
&
\ \ \ \ \
-
5\,\Delta^{1,3}\big(f_i'''f_j'+f_i'f_j'''\big)
+
5\,\Delta^{1,3}\,f_i''f_j''
-
\\ 
&
\ \ \ \ \ 
-
15\,\Delta^{2,3}\,f_i'f_j''
+
15\,\Delta^{1,3}\,f_i''f_j''
-
15\,\Delta^{2,3}\,f_i''f_j''
+
15\,\Delta^{1,3}\,f_i''f_j''
\\
&
=
\Delta^{1,5}\,f_i'f_j'
+
5\,\Delta^{2,4}\,f_i'f_j'
-
4\,\Delta^{1,4}\big(f_i''f_j'+f_i'f_j''\big)
-
\\
&
\ \ \ \ \
-
16\,\Delta^{2,3}\big(f_i''f_j'+f_i'f_j''\big)
-
5\,\Delta^{1,3}\big(f_i'''f_j'+f_i'f_j'''\big)
+
35\,\Delta^{1,3}\,f_i''f_j''.
\endaligned
\]
Ensuite, le calcul de ${\sf D} M^8$ est imm\'ediat car il n'implique
aucune relation pl\"u\-cke\-rienne\,:
\[
\aligned
{\sf D}M^8
&
=
3\,\Delta^{1,5}\,\Delta^{1,2}
+
3\,\Delta^{2,4}\,\Delta^{1,2}
+
3\,\Delta^{1,4}\,\Delta^{1,3}
+
\\
&
\ \ \ \ \
+
12\,\Delta^{2,4}\,\Delta^{1,2}
+
12\,\Delta^{2,3}\,\Delta^{1,3}
-
10\,\Delta^{1,4}\,\Delta^{1,3}
-
10\,\Delta^{2,3}\,\Delta^{1,3}
\\
&
=
3\,\Delta^{1,5}\,\Delta^{1,2}
+
15\,\Delta^{2,4}\,\Delta^{1,2}
-
7\,\Delta^{1,4}\,\Delta^{1,3}
+
2\,\Delta^{2,3}\,\Delta^{1,3}.
\endaligned
\]

\noindent{\bf Tableau des diff\'erentiations totales.} 
En r\'esum\'e, 
nous obtenons les expressions normalis\'ees suivantes pour
nos invariants diff\'erenti\'es\,:
\[
\boxed{
\aligned
{\sf D}f_i'
&
=
f_i''
\\
{\sf D}\Lambda^3
&
=
\Delta^{1,3}
\\
{\sf D}\Lambda_i^5
&
=
\Delta^{1,4}\,f_i'
+
4\,\Delta^{2,3}\,f_i'
-
5\,\Delta^{1,3}\,f_i''
\\
{\sf D}\Lambda_{i,j}^7
&
=
\Delta^{1,5}\,f_i'f_j'
+
5\,\Delta^{2,4}\,f_i'f_j'
-
4\,\Delta^{1,4}\big(f_i''f_j'+f_i'f_j''\big)
-
\\
&
\ \ \ \ \
-
16\,\Delta^{2,3}\big(f_i''f_j'+f_i'f_j''\big)
-
5\,\Delta^{1,3}\big(f_i'''f_j'+f_i'f_j'''\big)
+
35\,\Delta^{1,3}\big(\,f_i''f_j''\big)
\\
{\sf D}M^8
&
=
3\,\Delta^{1,5}\,\Delta^{1,2}
+
15\,\Delta^{2,4}\,\Delta^{1,2}
-
7\,\Delta^{1,4}\,\Delta^{1,3}
+
2\,\Delta^{2,3}\,\Delta^{1,3}
\endaligned
}\,,
\]
et nous pouvons maintenant commencer \`a engendrer la table de
multiplication\,\,---\,\,pond\'er\'ee par le poids de nos
invariants\,\,---\,\,entre ces deux listes encadr\'ees, afin de
d\'ecouvrir de nouveaux polyn\^omes invariants fondamentaux \`a l'\'etage
$\kappa = 5$.

\smallskip\noindent{\bf Premi\`ere famille de crochets $\big[ \Lambda_{ i,j}^7, \,
f_k' \big]$.} Apr\`es un calcul direct que nous ne d\'etail\-lerons pas, mais
dans lequel les normalisations pl\"ucke\-riennes n'interviennent pas, nous
obtenons\,:
\[
\small
\aligned
\big[
\Lambda_{i,j}^7,\,f_k'
\big]
&
=
{\sf D}\Lambda_{i,j}^7\cdot f_k'
-
7\,\Lambda_{i,j}^7\cdot{\sf D}f_k'
\\
&
=
\Delta^{1,5}\,f_i'f_j'f_k'
+
5\,\Delta^{2,4}\,f_i'f_j'f_k'
-
\\
&
\ \ \ \ \
-
4\,\Delta^{1,4}\big(f_i''f_j'+f_i'f_j'')\,f_k'
-
7\,\Delta^{1,4}\,f_i'f_j'f_k''
-
\\
&
\ \ \ \ \
-
16\,\Delta^{2,3}\big(f_i''f_j'+f_i'f_j''\big)f_k'
-
28\,\Delta^{2,3}\,f_i'f_j'f_k''
-
\\
&
\ \ \ \ \
-
5\,\Delta^{1,3}\big(f_i'''f_j'+f_i'f_j''')\,f_k'
+
35\,\Delta^{1,3}\big(f_i''f_j''f_k'+f_i''f_j'f_k''+f_i'f_j''f_k''\big)
-
\\
&
\ \ \ \ \
-
105\,\Delta^{1,2}\,f_i''f_j''f_k''
\\
&
=:
\Lambda_{i,j,k}^9.
\endaligned
\]
Nous trouvons donc huit nouveaux invariants $\Lambda_{ 1,1,1}^9$,
$\Lambda_{ 1,1, 2}^9$, $\Lambda_{ 1,2, 1}^9$, $\Lambda_{ 1,2, 2}^9$,
$\Lambda_{ 2,1, 1}^9$, $\Lambda_{ 2, 1, 2}^9$, $\Lambda_{ 2,2, 1}^9$
et $\Lambda_{ 2,2, 2}^9$, qui ne s'expriment clairement pas en
fonction de ceux connus \`a l'\'etage $\kappa = 4$, \`a cause
par exemple de la pr\'esence
du d\'eterminant $\Delta^{ 1, 5}$ o\`u apparaissent $f_1 '''''$ et $f_2
'''''$.

\smallskip\noindent{\bf Observation.} 
Cependant, ces huit invariants ne sont pas ind\'ependants entre eux, ne
serait-ce que par h\'eritage de la sym\'etrie $\Lambda_{ 1,2}^7 =
\Lambda_{ 2, 1}^7$, qui implique les deux relations $\Lambda_{ 1, 2,
k}^9 = \Lambda_{ 2, 1, k}^9$, $k = 1, 2$. En fait, il y a quatre
relations ind\'ependantes, que l'on peut proposer
au lecteur de v\'erifier par un d\'eveloppement
direct\,:
\[
\aligned
\Lambda_{1,1,2}^9
&
=
\Lambda_{1,2,1}^9
-
f_1'\,M^8
\\
\Lambda_{1,2,1}^9
&
=
\Lambda_{2,1,1}^9
\\
\Lambda_{1,2,2}^9
&
=
\Lambda_{2,1,2}^9
\\
\Lambda_{2,2,1}^9
&
=
\Lambda_{2,1,2}^9
+
f_2'\,M^8.
\endaligned
\]
Toutefois, il est incontestablement pr\'ef\'erable d'obtenir
ces relations comme suit \`a partir de l'identit\'e de type Jacobi, 
en posant tout simplement ${\sf P} := f_i'$,
${\sf Q} := f_j'$ et ${\sf R} := \Lambda_k^5$\,:
\[
0
\equiv
\big[\big[f_i',\,f_j'\big],\,\Lambda_k^5\big]
+
\big[\big[\Lambda_k^5,\,f_i'\big],\,f_j'\big]
+
\big[\big[f_j',\,\Lambda_k^5\big],\,f_i'\big].
\]
Si l'on tient compte du fait que $\big[ f_i', \, f_j' \big]$ vaut $0$
ou $\pm \Lambda^3$ et si on utilise les relations $\big[ \Lambda_i^5,
\, \Lambda^3 \big] = f_i' \, M^8$, nos quatre relations en
d\'ecoulent. Ainsi, seuls quatre des huit crochets $\big[ \Lambda_{
i,j}^7, \, f_k' \big]$ peuvent \^etre fondamentaux, et on v\'erifie sans
peine que $\Lambda_{ 1, 1, 1}^9$, $\Lambda_{ 1, 2, 1}^9$, $\Lambda_{
2, 1, 2}^9$ et $\Lambda_{ 2, 2, 2}^9$ constituent bien de nouveaux
invariants qui sont ind\'ependants entre eux, parce que les trin\^omes
$f_1' f_1' f_1'$, $f_1' f_2' f_1'$, $f_2' f_1' f_2'$ et $f_2' f_2'
f_2'$ en facteur derri\`ere $\Delta^{ 1, 5} + 5\, \Delta^{ 2, 4}$ le sont.

\smallskip\noindent{\bf Deuxi\`eme famille de crochets $\big[ M^8, \, f_i' 
\big]$.} Le calcul, imm\'ediat, 
libre d'ambigu\"it\'e pl\"ucke\-rienne et ne n\'ecessitant aucune r\'eorganisation,
fournit le r\'esultat suivant\,:
\[
\aligned
\big[
M^8,\,f_i'
\big]
&
=
{\sf D}M^8\cdot f_i'
-
8\,M^8\cdot f_i''
\\
&
=
\big[
3\,\Delta^{1,5}\,\Delta^{1,2}
+
15\,\Delta^{2,4}\,\Delta^{1,2}
-
7\,\Delta^{1,4}\,\Delta^{1,3}
+
2\,\Delta^{2,3}\,\Delta^{1,3}
\big]\,f_i'
-
\\
&
\ \ \ \ \
-
\big[
24\,\Delta^{1,4}\,\Delta^{1,2}
+
96\,\Delta^{2,3}\,\Delta^{1,2}
-
40\,\Delta^{1,3}\,\Delta^{1,3}
\big]f_i''
\\
&
=:
M_i^{10}.
\endaligned
\]
Nous trouvons donc deux nouveaux invariants $M_1^{ 10}$ et $M_2^{ 10}$
de poids 10, qui ne s'expriment pas
en fonction de ceux d\'ej\`a connus, \`a cause de la
pr\'esence du produit de d\'eterminants 
$\Delta^{ 1, 5}\, \Delta^{ 1,2}$ o\`u 
appara\^{\i}t $f_1 ''''' f_2''$.

\smallskip\noindent{\bf Troisi\`eme famille de crochets $\big[ \Lambda_{ i,j}^7, \, 
\Lambda^3 
\big]$.}
Bien que le r\'esultat final n'apporte pas de nouvel invariant ({\it
cf. infra}), nous d\'etaillerons ce calcul\,:
\[
\small
\aligned
\big[
\Lambda_{i,j}^7,\,\Lambda^3
\big]
&
=
3\,{\sf D}\Lambda_{i,j}^7\cdot\Lambda^3
-
7\,\Lambda_{i,j}^7\cdot{\sf D}\Lambda^3
\\
&
=
3\,\Delta^{1,5}\,\Delta^{1,2}\,f_i'f_j'
+
15\,\Delta^{2,4}\,\Delta^{1,2}\,f_i'f_j'
-
12\,\Delta^{1,4}\,\Delta^{1,2}\big(f_i''f_j'+f_i'f_j''\big)
-
\\
&
\ \ \ \ \
-
48\,\Delta^{2,3}\,\Delta^{1,2}\big(f_i''f_j'+f_i'f_j''\big)
-
15\,\Delta^{1,3}\,\underline{\Delta^{1,2}\,f_i'''}f_j'
-
\\
&
\ \ \ \ \
-
15\,\Delta^{1,3}\,\underline{\Delta^{1,2}}\,f_i'\underline{f_j'''}
+
105\,\Delta^{1,3}\,\Delta^{1,2}\,f_i''f_j''
-
7\,\Delta^{1,4}\,\Delta^{1,3}\,f_i'f_j'
-
\\
&
\ \ \ \ 
-
28\,\Delta^{2,3}\,\Delta^{1,3}\,f_i'f_j'
+
35\,\Delta^{1,3}\,\Delta^{1,3}\big(f_i''f_j'+f_i'f_j''\big)
-
105\,\Delta^{1,3}\,\Delta^{1,2}\,f_i''f_j''.
\endaligned
\]
Nous utilisons la relation pl\"ucke\-rienne pour transformer les
deux termes soulign\'es, qui deviennent\,: 
\[
15\,\Delta^{2,3}\,\Delta^{1,3}\,f_i'f_j'
-
15\,\Delta^{1,3}\,\Delta^{1,3}\,f_i''f_j'
+
15\,\Delta^{2,3}\,\Delta^{1,3}\,f_i'f_j'
-
15\,\Delta^{1,3}\,\Delta^{1,3}\,f_i'f_j'',
\]
et ensuite, nous additionnons les mon\^omes \'egaux et nous
regroupons le tout dans un ordre naturel\,: 
\[
\small
\aligned
\big[
\Lambda_{i,j}^7,\,\Lambda^3
\big]
&
=
\Big(
3\,\Delta^{1,5}\,\Delta^{1,2}
+
15\,\Delta^{2,4}\,\Delta^{1,2}
-
7\,\Delta^{1,4}\,\Delta^{1,3}
+
2\,\Delta^{2,3}\,\Delta^{1,3}
\Big)
f_i'f_j'
+
\\
&
\ \ \ \ \
+
\Big(
-
12\,\Delta^{1,4}\,\Delta^{1,2}
-
48\,\Delta^{2,3}\,\Delta^{1,2}
+
20\,\Delta^{1,3}\,\Delta^{1,3}
\Big)
\big(f_i''f_j'+f_i'f_j''\big).
\endaligned
\]
Or, {\it cette troisi\`eme famille de crochets n'apporte aucun nouvel
invariant}. En effet, consid\'erons la famille d'identit\'es de Jacobi\,:
\[
\aligned
0
&
\equiv
\big[\big[\Lambda_i^5,\,f_j'\big],\,\Lambda^3\big]
+
\big[\big[\Lambda^3,\,\Lambda_i^5\big]\,f_j'\big]
+
\big[\big[f_j',\,\Lambda^3\big],\,\Lambda_i^5\big]
\\
&
\equiv
\big[\Lambda_{i,j}^7,\,\Lambda^3\big]
-
\big[f_i'\,M^8,\,f_j'\big]
-
\big[\Lambda_j^5,\,\Lambda_i^5\big].
\endaligned
\]
Si nous faisons tout d'abord $i = j$ ($=1$ ou $= 2$), en utilisant $\big[
f_i' \, M^8, \, f_i' \big] = f_i' \, \big[ M^8, \, f_i' \big] = f_i'
\, M_i^{ 10}$, nous obtenons les
deux relations\,:
\[
0
\equiv
\big[
\Lambda_{i,i}^7,\,\Lambda^3
\big]
-
f_i'\,M_i^{10}
\]
qui montrent que les deux crochets $\big[ \Lambda_{ i, i}^7, \,
\Lambda^3 \big]$ pour $i = 1, 2$ sont superflus. Si nous faisons
ensuite $i = 1$ et $j = 2$, nous obtenons\,:
\[
\aligned
0
&
\equiv
\big[
\Lambda_{1,2}^7,\,\Lambda^3
\big]
-
M^8\,\big[f_1',\,f_j'\big]
-
f_1'\,\big[M^8,\,f_2'\big]
+
\big[\Lambda_1^5,\,\Lambda_2^5\big]
\\
&
\equiv
\big[\Lambda_{1,2}^7,\,\Lambda^3\big]
+
M^8\,\Lambda^3
-
f_1'\,M_2^{10}
-
5\,\Lambda^3\,M^8
\\
&
\equiv
\big[\Lambda_{1,2}^7,\,\Lambda^3\big]
-
4\,\Lambda^3\,M^8
-
f_1'\,M_2^{10},
\endaligned
\]
ce qui montre que $\big[ \Lambda_{ 1,2}^7, \, \Lambda^3 \big]$ est
superflu. De la sym\'etrie indicielle $\Lambda_{ 1,2}^7 = \Lambda_{
2,1}^7$ on peut d\'eduire sans plus de calcul que le crochet $\big[
\Lambda_{ 2,1 }^7, \, \Lambda^3 \big]$ est lui aussi superflu, mais il
est instructif de faire quand m\^eme $i = 2$ et $j = 1$ dans l'identit\'e
g\'en\'erale\,:
\[
\aligned
0
&
\equiv
\big[\Lambda_{2,1}^7\,\Lambda^3\big]
-
M^8\,\big[f_2',\,f_1'\big]
-
f_2'\,\big[M^8,\,f_1'\big]
-
\big[\Lambda_1^5,\,\Lambda_2^5\big]
\\
&
\equiv
\big[\Lambda_{2,1}^7,\,\Lambda^3\big]
-
M^8\,\Lambda^3
-
f_2'\,M_1^{10}
+
5\,\Lambda^3\,M^8
\\
&
\equiv
\big[\Lambda_{2,1}^7,\,\Lambda^3\big]
+
4\,\Lambda^3\,M^8
-
f_2'\,M_1^{10}.
\endaligned
\]
Bien que $\Lambda_{ 1,2}^7 = \Lambda_{ 2,1}^7$, nous obtenons une
relation ind\'ependante de la pr\'ec\'edente, et par soustraction, nous
obtenons une nouvelle relation, que nous \'enon\c cons en passant\,:
\[
0
\equiv
f_2'\,M_1^{10}
-
f_1'\,M_2^{10}
-
8\,\Lambda^3\,M^8,
\]
ce qui anticipe un fait qui va se r\'ev\'eler crucial par la suite\,: les
invariants fondamentaux form\'es par crochets jouissent d'un tr\`es grand
nombre de relations alg\'ebriques, parfois appel\'ees {\sl syzygies},
qu'il est difficile d'englober dans une combinatoire uni\-fi\'ee.
Poursuivons toutefois pour l'instant notre pr\'eparation de tous les
invariants que l'on peut former par crochets.

\smallskip\noindent{\bf Quatri\`eme famille de crochets $\big[ M^8, \, 
\Lambda^3 \big]$.} Le calcul, ``{\it most elementary}'', donne\,:
\[
\aligned
\big[M^8,\,\Lambda^3\big]
&
=
3\,{\sf D}M^8\cdot\Lambda^3
-
8\,M^8\cdot{\sf D}\Lambda^3
\\
&
=
9\,\Delta^{1,5}\,\Delta^{1,2}\,\Delta^{1,2}
+
45\,\Delta^{2,4}\,\Delta^{1,2}\,\Delta^{1,2}
-
45\,\Delta^{1,4}\,\Delta^{1,3}\,\Delta^{1,2}
-
\\
&
\ \ \ \
-
90\,\Delta^{2,3}\,\Delta^{1,3}\,\Delta^{1,2}
+
40\,\Delta^{1,3}\,\Delta^{1,3}\,\Delta^{1,3}
\\
&
=:
N^{12}.
\endaligned
\]
C'est un nouvel invariant $N^{ 12}$ de poids 12 qui a la propri\'et\'e
remarquable de s'expri\-mer seulement en fonction des d\'eterminants
$\Delta^{ \alpha, \beta}$. En compagnie de $\Lambda^3$ et de $M^8$, il
jouera un r\^ole central dans l'\'elaboration d'une base de Gr\"obner pour
l'id\'eal des syzygies entre les invariants fondamentaux.

\smallskip\noindent{\bf Cinqui\`eme famille de crochets 
$\big[ \Lambda_{ i,j}^7, \, \Lambda_k^5 \big]$.} Cette fois-ci, nous ne
d\'etaillerons pas les calculs interm\'ediaires, puisque nous avons d\'ej\`a
\'evoqu\'e \`a pr\'esent tous les actes qui permettent de les accomplir. Nous
obtenons\,:
\[
\small
\aligned
\big[
\Lambda_{i,j}^7,\,\Lambda_k^5
\big]
&
=
5\,{\sf D}\Lambda_{i,j}^7\cdot\Lambda_k^5
-
7\,\Lambda_{i,j}^7\cdot{\sf D}\Lambda_k^5
\\
&
=
5\,\Delta^{1,5}\,\Delta^{1,3}\,f_i'f_j'f_k'
+
25\,\Delta^{2,4}\,\Delta^{1,3}\,f_i'f_j'f_k'
-
7\,\Delta^{1,4}\,\Delta^{1,4}\,f_i'f_j'f_k'
-
\\
&
\ \ \ \ \
-
56\,\Delta^{2,3}\,\Delta^{1,4}\,f_i'f_j'f_k'
-
112\,\Delta^{2,3}\,\Delta^{2,3}\,f_i'f_j'f_k'
-
15\,\Delta^{1,5}\,\Delta^{1,2}\,f_i'f_j'f_k''
-
\\
&
\ \ \ \ \
-
75\,\Delta^{2,4}\,\Delta^{1,2}\,f_i'f_j'f_k''
+
15\,\Delta^{1,4}\,\Delta^{1,3}\big(f_i''f_j'+f_i'f_j''\big)f_k'
+
\\
&
\ \ \ \ \
+
35\,\Delta^{1,4}\,\Delta^{1,3}\,f_i'f_j'f_k''
+
60\,\Delta^{2,3}\,\Delta^{1,3}\big(f_i''f_j'+f_i'f_j''\big)f_k'
-
\\
&
\ \ \ \ \
-
10\,\Delta^{2,3}\,\Delta^{1,3}\,f_i'f_j'f_k''
-
25\,\Delta^{1,3}\,\Delta^{1,3}\big(f_i'''f_j'+f_i'f_j'''\big)f_k'
+
\\
&
\ \ \ \ \
+
175\,\Delta^{1,3}\,\Delta^{1,3}\,f_i''f_j''f_k'
-
100\,\Delta^{1,3}\,\Delta^{1,3}\big(f_i''f_j'+f_i'f_j''\big)f_k''
+
\\
&
\ \ \ \ \
+
60\,\Delta^{1,4}\,\Delta^{1,2}\big(f_i''f_j'+f_i'f_j''\big)f_k''
-
105\,\Delta^{1,4}\,\Delta^{1,2}\,f_i''f_j''f_k'
+
\\
&
\ \ \ \ \
+
240\,\Delta^{2,3}\,\Delta^{1,2}\big(f_i''f_j'+f_i'f_j''\big)f_k''
-
420\,\Delta^{2,3}\,\Delta^{1,2}\,f_i''f_j''f_k'.
\endaligned
\]
Ces six invariants sont nouveaux, \`a ceci pr\`es qu'ils ne sont
pas ind\'ependants entre eux. En effet, $( \mathcal{ J} ac )$
donne\,:
\[
0
\equiv
\big[\big[f_i',\,\Lambda_j^5\big],\,\Lambda_k^5\big]
+
\big[\big[\Lambda_k^5,\,f_i'\big],\,\Lambda_j^5\big]
+
\big[\big[\Lambda_j^5,\,\Lambda_k^5\big],\,f_i'\big],
\]
relations qui se r\'eduisent \`a $0 \equiv 0$ lorsque 
$j = k$, mais qui fournissent deux relations non triviales lorsque
$j \neq k$, \`a savoir\,:
\[
\aligned
\big[\Lambda_{1,2}^7,\,\Lambda_1^5\big]
&
=
\big[\Lambda_{1,1}^7,\,\Lambda_2^5\big]
+
5\,M^8\,\Lambda_1^5
+
5\,\Lambda^3\,M_1^{10},
\\
\big[\Lambda_{1,2}^7,\,\Lambda_2^5\big]
&
=
\big[\Lambda_{2,2}^7,\,\Lambda_1^5\big]
-
5\,M^8\,\Lambda_2^5
-
5\,\Lambda^3\,M_2^{10}.
\endaligned
\]
Celles-ci nous permettent de n'avoir \`a consid\'erer que les quatre (au
lieu de six) nouveaux invariants\,:
\[
\aligned
K_{1,1,1}^{13}
&
:=
\big[\Lambda_{1,1}^7,\,\Lambda_1^5\big],
\ \ \ \ \ \ \ \ \ \ \ \ \ \ \ \ \ \ \ 
K_{1,1,2}^{13}
:=
\big[\Lambda_{1,1}^7,\,\Lambda_2^5\big],
\\
K_{2,2,1}^{13}
&
:=
\big[\Lambda_{2,2}^7,\,\Lambda_1^5\big],
\ \ \ \ \ \ \ \ \ \ \ \ \ \ \ \ \ \ \ 
K_{2,2,2}^{13}
:=
\big[\Lambda_{2,2}^7,\,\Lambda_2^5\big].
\endaligned
\]
Cependant, le travail n'est pas termin\'e. Puisque nous constatons que
$K_{ 1, 1, 1}^{ 13}$ est divisible par $f_1'$, nous devons introduire
l'invariant de poids 12\,:
\[
\small
\aligned
K_{1,1}^{12}
:=
&\
\frac{1}{f_1'}\,
\big[\Lambda_{1,1}^7,\,\Lambda_1^5\big]
\\
=
&\
f_1'f_1'
\Big(
5\,\Delta^{1,5}\,\Delta^{1,3}
+
25\,\Delta^{2,4}\,\Delta^{1,3}
-
7\,\Delta^{1,4}\,\Delta^{1,4}
-
56\,\Delta^{2,3}\,\Delta^{1,4}
-
\\
&
-112\,\Delta^{2,3}\,\Delta^{2,3}\Big)
+
f_1'f_1''\Big(
-15\,\Delta^{1,5}\,\Delta^{1,2}
-
75\,\Delta^{2,4}\,\Delta^{1,2}
+
\\
&
+
65\,\Delta^{1,4}\,\Delta^{1,3}
+
110\,\Delta^{2,3}\,\Delta^{1,3}
\Big)
+
f_1f_1'''\Big(
-50\,\Delta^{1,3}\,\Delta^{1,3}
\Big)
+
\\
&
+
f_1''f_1''\Big(
-25\,\Delta^{1,3}\,\Delta^{1,3}
+
15\,\Delta^{1,4}\,\Delta^{1,2}
+
60\,\Delta^{2,3}\,\Delta^{1,2}
\Big).
\endaligned
\]
Pareillement, $K_{ 2, 2, 2}^{ 13}$ \'etant divisible par $f_2'$, nous
devons introduire cet invariant d\'efini par
$K_{ 2, 2}^{ 12} := \frac{ 1}{ f_1'}\, \big[
\Lambda_{ 2, 2}^7, \, \Lambda_2^5 \big]$. Mais les deux invariants
restants, \`a savoir $K_{ 1, 1, 2}^{ 13}$ et $K_{ 2, 2, 1}^{ 13}$, ne
sont divisibles ni par $f_1'$ ni par $f_2'$. Or nous verrons dans la
suite qu'il est naturel d'introduire des ``polarisations'' de
certains invariant sp\'eciaux\,\,---\,\,tels que $K_{ 1, 1}^{
12}$\,\,---\,\,qui ne comportent que des ``$1$'' en indices inf\'erieurs
et que nous appelerons {\sl bi-invariants}, les ``polarisations''
de ces invariants sp\'eciaux consistant tout simplement \`a mettre des
``$1$'' et des ``$2$'' de toutes les mani\`eres possibles en indices
inf\'erieurs. Par exemple, les polarisations de $\Lambda_{ 1, 1}^7$
sont\,: $\Lambda_{ 1, 2}^7$, $\Lambda_{ 2, 1}^7$ et $\Lambda_{ 2,
2}^7$. Mais alors, comment donc les deux invariants $K_{ 1 , 1, 2}^{
13}$ et $K_{ 2, 2, 1}^{ 13}$ pourraient-ils \^etre obtenus par
polarisation, sachant qu'ils ont trois indices inf\'erieurs? Faut-il
revenir \`a $K_{ 1, 1, 1}^{ 13}$ et le polariser?

\def\thelemma{\!}\begin{lemma}
Si l'on introduit, pour tous $i, j$ appartenant \`a $\{ 1, 2\}$, les
quatre invariants de poids 12{\rm \,:}
\[
\small
\aligned
K_{i,j}^{12}
:=
&\
f_i'f_j'
\Big(
5\,\Delta^{1,5}\,\Delta^{1,3}
+
25\,\Delta^{2,4}\,\Delta^{1,3}
-
7\,\Delta^{1,4}\,\Delta^{1,4}
-
56\,\Delta^{2,3}\,\Delta^{1,4}
-
\\
&
-112\,\Delta^{2,3}\,\Delta^{2,3}\Big)
+
\frac{(f_i'f_j''+f_i''f_j')}{2}
\Big(
-15\,\Delta^{1,5}\,\Delta^{1,2}
-
75\,\Delta^{2,4}\,\Delta^{1,2}
+
\\
\endaligned
\]
\[
\small
\aligned
&
+
65\,\Delta^{1,4}\,\Delta^{1,3}
+
110\,\Delta^{2,3}\,\Delta^{1,3}
\Big)
+
\frac{(f_i'f_j'''+f_i'''f_j')}{2}
\Big(
-50\,\Delta^{1,3}\,\Delta^{1,3}
\Big)
+
\\
&
+
f_i''f_j''\Big(
-25\,\Delta^{1,3}\,\Delta^{1,3}
+
15\,\Delta^{1,4}\,\Delta^{1,2}
+
60\,\Delta^{2,3}\,\Delta^{1,2}
\Big),
\endaligned
\]
alors les quatre invariants $K_{ 1, 1, 1}^{ 13}$, $K_{ 1, 1, 2}^{
13}$, $K_{ 2, 2, 1}^{ 13}$ et $K_{ 2, 2, 2}^{ 13}$ pr\'ec\'edents sont
r\'eobtenus au moyen des quatre relations{\rm \,:}
\[
\aligned
K_{1,1,1}^{13}
&
=
f_1'\,K_{1,1}^{12},
\\
K_{1,1,2}^{13}
&
=
f_1'\,K_{1,2}^{12}
-
\frac{5}{2}\,\Lambda^3\,M_1^{10}
-
5\,\Lambda_1^5\,M^8,
\\
K_{2,2,1}^{13}
&
=
f_2'\,K_{2,1}^{12}
+
\frac{5}{2}\,\Lambda^3\,M_2^{10}
+
5\,\Lambda_2^5\,M^8,
\\
K_{2,2,2}^{13}
&
=
f_2'\,K_{2,2}^{12}.
\endaligned
\]
\end{lemma}

Gr\^ace \`a ce lemme bienvenu, nous pouvons donc introduire l'invariant
r\'eduit $K_{ 1, 1}^{ 12}$ accompagn\'e de ses trois polarisations $K_{ 1,
2}^{ 12}$, $K_{ 2, 1}^{ 12}$ et $K_{ 2, 2}^{ 12}$ (en fait $K_{ 1,2}^{
12} = K_{ 2, 1}^{ 12}$), et oublier purement et simplement
les quatre invariants de poids 13 qui nous \'etaient fournis par
crochets bruts.

\noindent{\em Preuve.}
En consid\'erant la permutation $1 \leftrightarrow 2$ des indices (noter que
les $\Delta^{ \alpha, \beta}$ changent de signe), il suffit d'\'etablir
la deuxi\`eme identit\'e. Nous r\'e\'ecrivons tout d'abord, en repartant de
l'expression obtenue pour $\big[ \Lambda_{ 1, 1}^7, \, \Lambda_2^5
\big]$\,:
\[
\small
\aligned
K_{1,1,2}^{13}
&
=
f_1'f_1'f_2'
\Big(
5\,\Delta^{1,5}\,\Delta^{1,3}
+
25\,\Delta^{2,4}\,\Delta^{1,3}
-
7\,\Delta^{1,4}\,\Delta^{1,4}
-
56\,\Delta^{2,3}\,\Delta^{1,4}
-
\\
&
-112\,\Delta^{2,3}\,\Delta^{2,3}\Big)
+
f_1'f_1'f_2''
\Big(
-15\,\Delta^{1,5}\,\Delta^{1,2}
-
75\,\Delta^{2,4}\,\Delta^{1,2}
+
\\
&
+
35\,\Delta^{1,4}\,\Delta^{1,3}
-
10\,\Delta^{2,3}\,\Delta^{1,3}
\Big)
+
f_1'f_1''f_2'
\Big(
30\,\Delta^{1,4}\,\Delta^{1,3}
+
120\,\Delta^{2,3}\,\Delta^{1,3}
\Big)
+
\\
&
+
f_1'f_1''f_2'
\Big(
-50\,\Delta^{1,3}\,\Delta^{1,3}
\Big)
+
f_1''f_1''f_2'
\Big(
175\,\Delta^{1,3}\,\Delta^{1,3}
-
105\,\Delta^{1,4}\,\Delta^{1,2}
-
\\
&
\ \ 
-
420\,\Delta^{2,3}\,\Delta^{1,2}
\Big)
+
f_1'f_1''f_2''
\Big(
-200\,\Delta^{1,3}\,\Delta^{1,3}
+
120\,\Delta^{1,4}\,\Delta^{1,2}
+
480\,\Delta^{2,3}\,\Delta^{1,2}
\Big).
\endaligned
\]
Ensuite, nous effectuons la soustraction\,:
\[
\footnotesize
\aligned
K_{1,1,1}^{13}
-
&\
f_1'\,K_{1,2}^{12}
=
\\
=
&\
f_1'f_1'f_2''
\Big(
-\frac{15}{2}\,\Delta^{1,5}\,\Delta^{1,2}
-
\frac{75}{2}\,\Delta^{2,4}\,\Delta^{1,2}
+
\frac{5}{2}\,\Delta^{1,4}\,\Delta^{1,3}
-
65\,\Delta^{2,3}\,\Delta^{1,3}
\Big)
+
\\
&
+
f_1'f_1''f_2'
\Big(
\frac{15}{2}\,\Delta^{1,5}\,\Delta^{1,2}
+
\frac{75}{2}\,\Delta^{2,4}\,\Delta^{1,2}
-
\frac{5}{2}\,\Delta^{1,4}\,\Delta^{1,3}
+
65\,\Delta^{2,3}\,\Delta^{1,3}
\Big)
+
\\
&
+
f_1'f_1'''f_2'
\Big(
-\frac{50}{2}\,\Delta^{1,3}\,\Delta^{1,3}
\Big)
+
\\
&
+
f_1'f_1'f_2'''
\Big(
\frac{50}{2}\,\Delta^{1,3}\,\Delta^{1,3}
\Big)
+
\\
&
+
f_1''f_1''f_2'
\Big(
175\,\Delta^{1,3}\,\Delta^{1,3}
-
105\,\Delta^{1,4}\,\Delta^{1,2}
-
420\,\Delta^{2,3}\,\Delta^{1,2}
\Big)
+
\\
&
+
f_1''f_1'f_2''
\Big(
-175\,\Delta^{1,3}\,\Delta^{1,3}
+
105\,\Delta^{1,4}\,\Delta^{1,2}
+
420\,\Delta^{2,3}\,\Delta^{1,2}
\Big).
\endaligned
\]
Remarquablement, les coefficients polynomiaux complexes \'etant oppos\'es
par paires de lignes qui se suivent, nous voyons des d\'eterminants $2
\times 2$ se reformer\,:
\[
\small
\aligned
K_{1,1,1}^{13}
-
&\
f_1'\,K_{1,2}^{12}
=
\\
=
&\
f_1'
\Big(
-\frac{15}{2}\,\Delta^{1,5}\,\Delta^{1,2}\,\Delta^{1,2}
-
\frac{75}{2}\,\Delta^{2,4}\,\Delta^{1,2}\,\Delta^{1,2}
+
\frac{5}{2}\,\Delta^{1,4}\,\Delta^{1,3}\,\Delta^{1,2}
-
\\
&
\ \ \ \ \
-
65\,\Delta^{2,3}\,\Delta^{1,2}\,\Delta^{1,2}
\Big)
+
f_1'\Big(
\frac{50}{2}\,\Delta^{1,3}\,\Delta^{1,3}\,\Delta^{1,2}
\Big)
+
\\
&
+
f_1''\Big(
105\,\Delta^{1,4}\,\Delta^{1,2}\,\Delta^{1,2}
+
420\,\Delta^{2,3}\,\Delta^{1,2}\,\Delta^{1,2}
-
175\,\Delta^{1,3}\,\Delta^{1,3}\,\Delta^{1,3}
\Big).
\endaligned
\]
Les deux premiers termes de la premi\`ere lignes ressemblant \`a ceux de
$\Lambda^3 \, M_1^{ 10}$, nous pouvons \'ecrire\,:
\[
\small
\aligned
&
=
-\frac{5}{2}\,\Lambda^3\,M_1^{10}
+
\\
&
+
f_1'\Big(
-\frac{30}{2}\,\Delta^{1,4}\,\Delta^{1,3}\Delta^{1,2}
-
60\,\Delta^{2,3}\,\Delta^{1,3}\,\Delta^{1,2}
+
\frac{50}{2}\,\Delta^{1,3}\,\Delta^{1,2}\,\Delta^{1,3}
\Big)
+
\\
&
+
f_1''\Big(
45\,\Delta^{1,4}\,\Delta^{1,2}\,\Delta^{1,2}
+
180\,\Delta^{2,3}\,\Delta^{1,2}\,\Delta^{1,2}
-
75\,\Delta^{1,3}\,\Delta^{1,3}\,\Delta^{1,2}
\Big).
\endaligned
\]
Et enfin, nous reconnaissons dans les termes restants l'expression
d\'evelopp\'ee de $- 5\, \Lambda_1^5 \, M^8$, ce qui nous donne bien la
relation annonc\'ee, que nous
r\'e\'ecrivons\,: $K_{ 1, 1, 2}^{ 13} - f_1' \, K_{ 1, 2}^{ 12} = -
\frac{ 5}{ 2}\, \Lambda^3 \, M_1^{ 10} - 5 \, \Lambda_1^5 \, M^8$.
\hfill$\square$

\smallskip\noindent{\bf Sixi\`eme famille de crochets 
$\big[ M^8, \, \Lambda_i^5 \big]$.} Par un calcul facile, court et sans
myst\`ere qui nous permer de reprendre haleine avant d'envisager la
septi\`eme et la plus complexe famille de crochets, nous obtenons\,:
\[
\small
\aligned
\big[
M^8,\,\Lambda_i^5
\big]
&
=
5\,{\sf D}M^8\cdot\Lambda_i^5
-
8\,M^8\cdot{\sf D}\Lambda_i^5
\\
&
=
\Big(
15\,\Delta^{1,5}\,\Delta^{1,3}\,\Delta^{1,2}
+
75\,\Delta^{2,4}\,\Delta^{1,3}\,\Delta^{1,2}
+
5\,\Delta^{1,4}\,\Delta^{1,3}\,\Delta^{1,3}
+
\\
&
\ \ \ \ \
+
170\,\Delta^{2,3}\,\Delta^{1,3}\,\Delta^{1,3}
-
24\,\Delta^{1,4}\,\Delta^{1,4}\,\Delta^{1,2}
-
192\,\Delta^{1,4}\,\Delta^{2,3}\,\Delta^{1,2}
-
\\
&
\ \ \ \ \
-
384\,\Delta^{2,3}\,\Delta^{2,3}\,\Delta^{1,2}
\Big)\,f_i'
+
\Big(
-
45\,\Delta^{1,5}\,\Delta^{1,2}\,\Delta^{1,2}
-
225\,\Delta^{2,4}\,\Delta^{1,2}\,\Delta^{1,2}
+
\\
&
\ \ \ \ \
+
225\,\Delta^{1,4}\,\Delta^{1,3}\,\Delta^{1,2}
+
450\,\Delta^{2,3}\,\Delta^{1,3}\,\Delta^{1,2}
-
200\,\Delta^{1,3}\,\Delta^{1,3}\,\Delta^{1,3}
\Big)\,f_i''
\\
&
=:
H_i^{14}.
\endaligned
\]
Le r\'esultat n'est divisible ni par $\Lambda_i^5$ (sinon ${\sf D}
\Lambda_i^5$ le serait), ni par $\Lambda^3$, ni par $f_i'$. Nous
trouvons donc deux nouveaux invariants $H_1^{ 14}$ et $H_2^{ 14}$ de
poids $14$.

\smallskip\noindent{\bf Septi\`eme famille de crochets $\big[
\Lambda_{ i,j}^7, \, \Lambda_{ k,l}^7 \big]$.} Le calcul complet, que
nous d\'etaillons dans la Section~9 parce qu'il est d\'elicat,
donne\,:
\[
\footnotesize
\aligned
\frac{
\big[
\Lambda_{i,j}^7,\,\Lambda_{k,l}^7
\big]}{7}
&
=
{\sf D}\Lambda_{i,j}^7\cdot\Lambda_{k,l}^7
-
\Lambda_{i,j}^7\cdot{\sf D}\Lambda_{k,l}^7
\\
&
=
\Big(
-
5\,\Delta^{1,5}\,\Delta^{1,3}
-
25\,\Delta^{2,4}\,\Delta^{1,3}
+
4\,\Delta^{1,4}\,\Delta^{1,4}
+
32\,\Delta^{1,4}\,\Delta^{2,3}
+
\\
&
\ \ \ \ \ \ \ \ \ \
+
64\,\Delta^{2,3}\,\Delta^{2,3}
\Big)
\big(
f_j'f_l'\,\Delta_{i,k}^{1,2}
+
f_i'f_k'\,\Delta_{j,l}^{1,2}
\big)
+
\\
\endaligned
\]
\[
\footnotesize
\aligned
&
\ \ \ \ \
+
\Big(
15\,\Delta^{1,5}\,\Delta^{1,2}
+
75\,\Delta^{2,4}\,\Delta^{1,2}
-
35\,\Delta^{1,4}\,\Delta^{1,3}
+
\\
&
\ \ \ \ \ \ \ \ \ \ \ \ \ \ \
+
10\,\Delta^{2,3}\,\Delta^{1,3}
\Big)
\big(
f_i'f_l''\,\Delta_{j,k}^{1,2}
+
f_k'f_j''\,\Delta_{i,l}^{1,2}
\big)
+
\\
&
\ \ \ \ \
+
\Big(
-5\,\Delta^{1,4}\,\Delta^{1,3}
-
20\,\Delta^{2,3}\,\Delta^{1,3}
\Big)
\big(
f_j'f_l'\,\Delta_{k,i}^{1,3}
+
f_i'f_k'\,\Delta_{l,j}^{1,3}
\big)
+
\\
&
\ \ \ \ \
+
\Big(
25\,\Delta^{1,3}\,\Delta^{1,3}
\Big)
\big(
f_j'f_l'\,\Delta_{k,i}^{2,3}
+
f_j'f_k'\,\Delta_{l,i}^{2,3}
+
f_i'f_l'\,\Delta_{k,j}^{2,3}
+
f_i'f_k'\,\Delta_{l,j}^{2,3}
\big)
+
\\
&
\ \ \ \ \
+
\Big(
-60\,\Delta^{1,4}\,\Delta^{1,2}
-
240\,\Delta^{2,3}\,\Delta^{1,2}
+
\\
&
\ \ \ \ \ \ \ \ \ \ \ \ \ \ \ 
+
100\,\Delta^{1,3}\,\Delta^{1,3}
\Big)
\big(
f_i''f_k''\,\Delta_{j,l}^{1,2}
+
f_j''f_l''\,\Delta_{i,k}^{1,2}
\big).
\endaligned
\]
Nous trouvons ainsi trois invariants de poids 15\,: 
\[
\small
\big[\Lambda_{1,1}^7,\,\Lambda_{1,2}^7\big],
\ \ \ \ \ \ \ \ \ \ \ \ \ \ \ \ \ \ \ \ \ \
\big[\Lambda_{1,1}^7,\,\Lambda_{2,2}^7\big],
\ \ \ \ \ \ \ \ \ \ \ \ \ \ \ \ \ \ \ \ \ \
\big[\Lambda_{1,2}^7,\,\Lambda_{2,2}^7\big],
\]
mais cependant, ces invariants s'expriment en fonction de ceux que
nous connaissons d\'ej\`a. En effet, sp\'ecialisons tout d'abord les
indices dans la formule g\'en\'erale et nettoyons les expressions
obtenues\,:
\[
\footnotesize
\aligned
\frac{\big[\Lambda_{1,1}^7,\,\Lambda_{1,2}^7\big]}{7}
&
=
f_1'f_1'
\Big(
-5\,\Delta^{1,5}\,\Delta^{1,3}\,\Delta^{1,2}
-
25\,\Delta^{2,4}\,\Delta^{1,3}\,\Delta^{1,2}
+
4\,\Delta^{1,4}\,\Delta^{1,4}\,\Delta^{1,2}
+
\\
&
\ \ \ \ \ \ \ \ \ \ \ \ \ \ \
+
32\,\Delta^{1,4}\,\Delta^{2,3}\,\Delta^{1,2}
+
64\,\Delta^{2,3}\,\Delta^{2,3}\,\Delta^{1,2}
+
\\
&
\ \ \ \ \ \ \ \ \ \ \ \ \ \ \
+
5\,\Delta^{1,4}\,\Delta^{1,3}\,\Delta^{1,3}
-
30\,\Delta^{2,3}\,\Delta^{1,3}\,\Delta^{1,3}\Big)
+
\\
&
\ \ \ \ \
+
f_1'f_1''\Big(
15\,\Delta^{1,5}\,\Delta^{1,2}\,\Delta^{1,2}
+
75\,\Delta^{2,4}\,\Delta^{1,2}\,\Delta^{1,2}
-
\\
&
\ \ \ \ \ \ \ \ \ \ \ \ \ \ \ \ \ \ \ \
-
35\,\Delta^{1,4}\,\Delta^{1,3}\,\Delta^{1,2}
+
10\,\Delta^{2,3}\,\Delta^{1,3}\,\Delta^{1,2}\Big)
+
\\
&
\ \ \ \ \ 
+
f_1''f_1''\Big(
-60\,\Delta^{1,4}\,\Delta^{1,2}\,\Delta^{1,2}
-
240\,\Delta^{2,3}\,\Delta^{1,2}\,\Delta^{1,2}
+
\\
&
\ \ \ \ \ \ \ \ \ \ \ \ \ \ \ \ \ \ \ \
+
100\,\Delta^{1,3}\,\Delta^{1,3}\,\Delta^{1,2}
\Big),
\endaligned
\]
\[
\footnotesize
\frac{\big[\Lambda_{1,2}^7,\,\Lambda_{2,2}^7\big]}{7}
=
\text{\rm Idem}
\big(
\text{\rm indice}\
{\bf 1}
\longleftrightarrow
\text{\rm indice}\
{\bf 2}
\big),
\]
\[
\footnotesize
\aligned
\frac{\big[\Lambda_{1,1}^7,\,\Lambda_{2,2}^7\big]}{7}
&
=
f_1'f_2'
\Big(
-10\,\Delta^{1,5}\,\Delta^{1,3}\,\Delta^{1,2}
-
50\,\Delta^{2,4}\,\Delta^{1,3}\,\Delta^{1,2}
+
8\,\Delta^{1,4}\,\Delta^{1,4}\,\Delta^{1,2}
+
\\
&
\ \ \ \ \ \ \ \ \ \ \ \ \ \ \
+
64\,\Delta^{1,4}\,\Delta^{2,3}\,\Delta^{1,2}
+
128\,\Delta^{2,3}\,\Delta^{2,3}\,\Delta^{1,2}
+
\\
&
\ \ \ \ \ \ \ \ \ \ \ \ \ \ \
+
10\,\Delta^{1,4}\,\Delta^{1,3}\,\Delta^{1,3}
-
60\,\Delta^{2,3}\,\Delta^{1,3}\,\Delta^{1,3}\Big)
+
\\
&
\ \ \ \ \
+
\Big(\frac{f_1'f_2''+f_1''f_2'}{2}\Big)\Big(
30\,\Delta^{1,5}\,\Delta^{1,2}\,\Delta^{1,2}
+
150\,\Delta^{2,4}\,\Delta^{1,2}\,\Delta^{1,2}
-
\\
&
\ \ \ \ \ \ \ \ \ \ \ \ \ \ \ \ \ \ 
\ \ \ \ \ \ \ \ \ \ \ \ \ \ \ \ \ \ \
-
70\,\Delta^{1,4}\,\Delta^{1,3}\,\Delta^{1,2}
+
20\,\Delta^{2,3}\,\Delta^{1,3}\,\Delta^{1,2}\Big)
+
\\
&
\ \ \ \ \ 
+
f_1''f_2''\Big(
-120\,\Delta^{1,4}\,\Delta^{1,2}\,\Delta^{1,2}
-
480\,\Delta^{2,3}\,\Delta^{1,2}\,\Delta^{1,2}
+
\\
&
\ \ \ \ \ \ \ \ \ \ \ \ \ \ \ \ \ \ \ \
+
200\,\Delta^{1,3}\,\Delta^{1,3}\,\Delta^{1,2}
\Big).
\endaligned
\]
Si nous examinons le polyn\^ome cubique en $\Delta$ 
qui est multiple de $f_1 ''
f_1''$ dans les deux derni\`eres lignes de l'expression de $ \frac{
1}{ 7}\, \big[ \Lambda_{ 1, 1}^7, \, \Lambda_{ 1,2}^7 \big]$, nous
reconnaissons par exemple $ - 20 \, M^8 \, \Delta^{ 1, 2}$, ce qui
constitue une co\"{\i}ncidence que nous devrions manifestement
exploiter, et ensuite, sans plus tenter de d\'ecrire l'asc\`ese visuelle
qui nous permet de deviner des relations alg\'ebriques entre de telles
expressions, nous trouvons les trois relations suivantes imm\'ediatement
v\'erifiables par d\'eveloppement\,:
\[
\aligned
0
&
\equiv
6\,\big[\Lambda_{1,1}^7,\,\Lambda_{1,2}^7\big]
+
35\,\Lambda_1^5\,M_1^{10}
+
f_1'\,H_1^{14},
\\
0
&
\equiv
6\,\big[\Lambda_{1,1}^7,\,\Lambda_{1,2}^7\big]
+
35\,\big(
\Lambda_1^5\,M_2^{10}
+
\Lambda_2^5\,M_1^{10} 
\big)
+
f_1'\,H_2^{14}
+
f_2'\,H_1^{14},
\\
0
&
\equiv
6\,\big[\Lambda_{1,2}^7,\,\Lambda_{2,2}^7\big]
+
35\,\Lambda_2^5\,M_2^{10}
+
f_2'\,H_2^{14},
\endaligned
\]
qui montrent que les trois invariants $\big[ \Lambda_{ 1, 1}^7, \,
\Lambda_{ 1, 2}^7 \big]$, $\big[ \Lambda_{ 1, 1}^7, \, \Lambda_{ 2, 2}^7
\big]$ et $\big[ \Lambda_{ 1, 2}^7, \, \Lambda_{ 2, 2}^7 \big]$ sont en fait
superflus. Bien qu'elle nous ait coût\'e de r\'eels efforts de calculs,
cette circonstance n'est pas sans nous d\'eplaire, puisque nous pouvons
ainsi r\'eduire de trois unit\'es le nombre d'invariants fondamentaux que
nous aurons \`a consid\'erer ult\'erieurement.

\smallskip\noindent{\bf Huiti\`eme famille de crochets $\big[ M^8, \, \Lambda_{ i,j}^7 
\big]$.} Le calcul, qui implique seulement quel\-ques normalisations
pl\"ucke\-riennes et bien sûr aussi de l'arithm\'etique formelle
\'el\'e\-men\-taire, fournit l'expression massive suivante, qui est en fait
compl\`etement simplifi\'ee\,:
\[
\footnotesize
\aligned
&
\big[
M^8,\,\Lambda_{i,j}^7
\big]
=
7\,{\sf D}M^8\cdot\Lambda_{i,j}^7
-
8\,M^8\cdot{\sf D}\Lambda_{i,j}^7
\\
&
=
\Big(
-3\,\Delta^{1,5}\,\Delta^{1,4}\,\Delta^{1,2}
-
15\,\Delta^{2,4}\,\Delta^{1,4}\,\Delta^{1,2}
-
12\,\Delta^{1,5}\,\Delta^{2,3}\,\Delta^{1,2}
+
\\
&
\ \ \ \ \ 
+
40\,\Delta^{1,5}\,\Delta^{1,3}\,\Delta^{1,3}
-
60\,\Delta^{2,4}\,\Delta^{2,3}\,\Delta^{1,2}
+
200\,\Delta^{2,4}\,\Delta^{1,3}\,\Delta^{1,3}
-
\\
&
\ \ \ \ \
-
49\,\Delta^{1,4}\,\Delta^{1,4}\,\Delta^{1,3}
-
422\,\Delta^{1,4}\,\Delta^{2,3}\,\Delta^{1,3}
-
904\,\Delta^{2,3}\,\Delta^{2,3}\,\Delta^{1,3}
\Big)f_i'f_j'
+
\\
&
\ \ \ \ \ 
+
\Big(
-105\,\Delta^{1,5}\,\Delta^{1,3}\,\Delta^{1,2}
-
525\,\Delta^{2,4}\,\Delta^{1,3}\,\Delta^{1,2}
+
205\,\Delta^{1,4}\,\Delta^{1,3}\,\Delta^{1,3}
-
\\
&
\ \ \ \ \
-
230\,\Delta^{2,3}\,\Delta^{1,3}\,\Delta^{1,3}
+
96\,\Delta^{1,4}\,\Delta^{1,4}\,\Delta^{1,2}
+
768\,\Delta^{1,4}\,\Delta^{2,3}\,\Delta^{1,2}
+
\\
&
\ \ \ \ \
+
1536\,\Delta^{2,3}\,\Delta^{2,3}\,\Delta^{1,2}
\Big)
\big(
f_i''f_j'+f_i'f_j''
\big)
+
\\
&
\ \ \ \ \
+
\Big(
-200\,\Delta^{1,3}\,\Delta^{1,3}\,\Delta^{1,3}
\Big)\big(
f_i'''f_j'+f_i'f_j'''
\big)
+
\\
&
\ \ \ \ \ 
+
\Big(
315\,\Delta^{1,5}\,\Delta^{1,2}\,\Delta^{1,2}
+
1575\,\Delta^{2,4}\,\Delta^{1,2}\,\Delta^{1,2}
-
1575\,\Delta^{1,4}\,\Delta^{1,3}\,\Delta^{1,2}
-
\\
&
\ \ \ \ \
-
3150\,\Delta^{2,3}\,\Delta^{1,3}\,\Delta^{1,2}
+
1400\,\Delta^{1,3}\,\Delta^{1,3}\,\Delta^{1,3}
\Big)f_i''f_j''.
\endaligned
\]
Nous obtenons ainsi trois nouveaux invariants de poids 16\,:
\[
\aligned
F_{1,1}^{16}
&
:=
\big[M^8,\,\Lambda_{1,1}^7\big],
\ \ \ \ \ \ \ \ \ \ \ \ \ \ \ \ \ \ \ \ \ \
F_{1,2}^{16}
:=
\big[M^8,\,\Lambda_{1,2}^7\big]
=
F_{2,1}^{16},
\\
F_{2,2}^{16}
&
:=
\big[M^8,\,\Lambda_{2,2}^7\big].
\endaligned
\]

\noindent{\bf Syst\`eme g\'en\'erateur pour les
jets d'ordre 5.} Il nous faut donc consid\'erer les vingt-cinq 
polyn\^o\-mes
invariants fondamentaux\,:
\[
\boxed{
\aligned
f_1',\ \
f_2',\ \
\Lambda^3,\ \
\Lambda_1^5,\ \
\Lambda_2^5,\ \
\Lambda_{1,1}^7,\ \
\Lambda_{1,2}^7,\ \
\Lambda_{2,2}^7,\ \
M^8,\
\\
\Lambda_{1,1,1}^9,\ \
\Lambda_{1,2,1}^9,\ \
\Lambda_{2,1,2}^9,\ \
\Lambda_{2,2,2}^9,\ \
M_1^{10},\ \
M_2^{10},\
\\
N^{12}, \ \
K_{1,1}^{12},\ \
K_{1,2}^{12},\ \
K_{2,1}^{12},\ \
K_{2,2}^{12},\
\\
H_1^{14},\ \
H_2^{14},\ \
F_{1,1}^{16},\ \
F_{1,2}^{16},\ \
F_{2,2}^{16}.\
\endaligned
}
\]

\noindent{\bf Probl\`eme.}
D\'ecrire explicitement l'id\'eal des relations entre ces vingt-cinq
polyn\^o\-mes et \'etablir que tout polyn\^ome invariant par reparam\'etrisation
${\sf P} \big( j^5 f \big)$ se repr\'e\-sente comme polyn\^ome $\mathcal{ P}
\big( f_1', \dots, M^8, \dots, H_1^{ 14}, \dots, F_{ 2, 2}^{ 16}
\big)$ en ces vingt-cinq invariants fondamentaux.

\section*{\S5.~Jets d'ordre 4 en dimension 2}

\noindent{\bf Neuf relations fondamentales}
Cette section et celle qui suit sont consacr\'ees \`a l'\'etude plus
accessible de $\mathcal{ DS}_2^4$. Pour les jets d'ordre 4, on
consid\`ere les neuf polyn\^omes invariants fondamentaux\,:
\[
\Big(\
f_1',\ \ \ \ \
f_2',\ \ \ \ \
\Lambda^3,\ \ \ \ \
\Lambda_1^5,\ \ \ \ \
\Lambda_2^5,\ \ \ \ \
\Lambda_{1,1}^7,\ \ \ \ \
\Lambda_{1,2}^7,\ \ \ \ \
\Lambda_{2,2}^7,\ \ \ \ \
M^8\
\Big),
\]
et avant de passer aux jets d'ordre cinq, nous allons 
d\'emontrer que tout polyn\^ome invariant par reparam\'etrisation 
${\sf P} \big( j^4 f \big)$ se repr\'esente comme
polyn\^ome de la forme
$\mathcal{ P} \big( f_1', \dots, \Lambda_2^5, \dots, M^8 \big)$.

Trois calculs distincts sur Maple\footnote{\, Les r\'esultats que l'on
re\c coit (apr\`es un quart d'heure de calcul environ) d\'ependent de
l'ordre monomial choisi\,; ils comprennent d'autres relations
superflues, {\it i.e.} d\'eduites des neuf fondamentales que nous
listons, parce que le logiciel doit effectuer algorithmiquement
l'op\'eration dite ``$S$-polyn\^ome'' sur tous les couples
d'identit\'es, suivie d'une division par les \'el\'ements pr\'esents, afin de
compl\'eter le calcul d'une base de Gr\"obner r\'eduite. L'auteur remercie
Erwan Rousseau de lui avoir communiqu\'e cette liste, ainsi que
l'invariant $M^8$.} conduisent aux neuf relations fondamentales
suivantes entre ces neuf invariants\,:
\[
\small
\left[
0
\overset{1}{\equiv}
f_2'\,\Lambda_1^5
-
f_1'\,\Lambda_2^5
-
3\,\Lambda^3\,\Lambda^3,
\right.
\]
\[
\small
\left[
\aligned
0
&
\overset{2}{\equiv}
f_2'\,\Lambda_{1,1}^7
-
f_1'\,\Lambda_{1,2}^7
-
5\,\Lambda^3\,\Lambda_1^5,
\\
0
&
\overset{3}{\equiv}
f_2'\,\Lambda_{1,2}^7
-
f_1'\,\Lambda_{2,2}^7
-
5\,\Lambda^3\,\Lambda_2^5,
\endaligned\right.
\]
\[
\small
\left[
\aligned
0
&
\overset{4}{\equiv}
f_1'\,f_1'\,M^8
-
3\,\Lambda^3\,\Lambda_{1,1}^7
+
5\,\Lambda_1^5\,\Lambda_1^5,
\\
0
&
\overset{5}{\equiv}
f_1'\,f_2'\,M^8
-
3\,\Lambda^3\,\Lambda_{1,2}^7
+
5\,\Lambda_1^5\,\Lambda_2^5,
\\
0
&
\overset{6}{\equiv}
f_2'\,f_2'\,M^8
-
3\,\Lambda^3\,\Lambda_{2,2}^7
+
5\,\Lambda_2^5\,\Lambda_2^5,
\endaligned\right.
\]
\[
\small
\left[
\aligned
0
&
\overset{7}{\equiv}
f_1'\,\Lambda^3\,M^8
-
\Lambda_1^5\,\Lambda_{1,2}^7
+
\Lambda_2^5\,\Lambda_{1,1}^7,
\\
0
&
\overset{8}{\equiv}
f_2'\,\Lambda^3\,M^8
-
\Lambda_1^5\,\Lambda_{2,2}^7
+
\Lambda_2^5\,\Lambda_{1,2}^7,
\endaligned\right.
\]
\[
\small
\left[
\aligned
0
&
\overset{9}{\equiv}
5\,\Lambda^3\,\Lambda^3\,M^8
-
\Lambda_{2,2}^7\,\Lambda_{1,1}^7
+
\Lambda_{1,2}^7\,\Lambda_{1,2}^7.
\endaligned\right.
\]
Dans le cas des jets d'ordre $\kappa = 3$, seule la premi\`ere relation
est pr\'esente\,; l'id\'eal des relations est principal et il constitue {\it
per se}\, la base de Gr\"obner ad\'equate.

\smallskip\noindent{\bf Question.} Comment retrouver 
et pr\'evoir l'existence de toutes ces relations~?

\smallskip\noindent{\bf Indication de r\'eponse.}
Au niveau $\kappa = 4$, une seule identit\'e de Jacobi peut \^etre form\'ee
et elle donne la relation d\'ej\`a connue $\Lambda_{ 1, 2}^7 = \Lambda_{
2, 1}^7$, qui est cependant triviale par rapport aux neufs identit\'es
list\'ees ci-dessus. Mais en d\'echiffrant ces neuf \'equations, ou bien en
observant que les deux familles d'identit\'es qui sont satisfaites dans
les alg\`ebres pl\"ucke\-riennes et que nous avons d\'ej\`a utilis\'ees
syst\'ematiquement pour normaliser l'expression d\'efinitive de nos
crochets, doivent n\'ecessairement et naturellement
``embo\^{\i}ter notre pas'' lorsque nous formons tous les d\'eterminants
$2 \times 2$ de la matrice
\[
\left\vert\!\left\vert
\begin{array}{ccccc}
f_1' \ \ & \ \ f_2' \ \ & \ \ 3\,\Lambda^3 \ \ & 
\ \ 5\,\Lambda_1^5 \ \ & \ \ 5\,\Lambda_2^5
\\
{\sf D}f_1' \ \ & \ \ {\sf D}f_2' \ \ & \ \
{\sf D}\Lambda^3 \ \ & \ \ {\sf D}\Lambda_1^5 \ \ & \ \
{\sf D}\Lambda_2^5
\end{array}
\right\vert\!\right\vert,
\] 
nous devinons\footnote{\, L'auteur a d'abord devin\'e et reconstitu\'e les
relations g\'en\'erales qui sous-tendent les neufs identit\'es pr\'ec\'edentes
avant d'arranger synoptiquement la formation de nouveaux invariants
par simple calcul de mineurs $2 \times 2$. }, ou nous constatons, que
le lemme g\'en\'eral suivant est vrai.

\def\thelemma{\!}\begin{lemma}
Pour tout quadruplet $\big( {\sf P}, \, {\sf Q},\, {\sf R}, \, {\sf S}
\big)$ d'invariants de poids $m, n, o, p$, les deux identit\'es suivantes
de type pl\"uck\'erien sont satisfaites{\rm \,:}
\def\theequation{$\mathcal{P}lck_1$}\begin{equation}
\boxed{
0
\equiv
m\,{\sf P}\,\big[{\sf Q},\,{\sf R}\big]
+
o\,{\sf R}\,\big[{\sf P},\,{\sf Q}\big]
+
n\,{\sf Q}\,\big[{\sf R},\,{\sf P}\big]
}\,,
\end{equation}
\def\theequation{$\mathcal{P}lck_2$}\begin{equation}
\boxed{
0
\equiv
\big[{\sf P},\,{\sf Q}\big]
\cdot
\big[{\sf R},\,{\sf S}\big]
+
\big[{\sf S},\,{\sf P}\big]
\cdot
\big[{\sf R},\,{\sf Q}\big]
+
\big[{\sf Q},\,{\sf S}\big]
\cdot
\big[{\sf R},\,{\sf P}\big]
}\,.
\end{equation}
\end{lemma}

\noindent{\em Preuve.}
Si en effet nous d\'eveloppons les deux derniers crochets\,:
\[
o\,{\sf R}
\big(
n\,{\sf D}{\sf P}\cdot{\sf Q}
-
m\,{\sf P}\cdot{\sf D}{\sf Q}
\big)
+
n\,{\sf Q}
\big(
m\,{\sf D}{\sf R}\cdot{\sf P}
-
o\,{\sf R}\cdot{\sf D}{\sf P}
\big),
\]
les deux termes extr\^emes s'annihilent, tandis que
les deux termes centraux\,:
\[
-
om\,{\sf R}\cdot{\sf P}\cdot{\sf D}{\sf Q}
+
nm\,{\sf Q}\cdot{\sf D}{\sf R}\cdot{\sf P}
\equiv
-
m\,{\sf P}\,\big[{\sf Q},\,{\sf R}\big]
\]
reconstituent l'oppos\'e du premier terme de la premi\`ere identit\'e
\thetag{ $ \mathcal{ P}lck_1$}. La seconde \thetag{ $\mathcal{ P}lck_2
$} n'est qu'une reformulation de la relation pl\"ucke\-rienne fondamentale
qui est satisfaite par les six mineurs $2 \times 2$ d'une matrice de
taille $2 \times 4$. On la v\'erifie en d\'eveloppant les trois produits
de d\'eterminants $2 \times 2$, ce qui produit 12 termes constitu\'es de 6
couples s'annihilant. \hfill$\square$

\smallskip\noindent{\bf Reconstitution des 9 syzygies.}
Il est tr\`es remarquable que les identit\'es \thetag{ $ \mathcal{ P}
lck_1$} permettent de reconstituer les huit premi\`eres (parmi neuf) des
identit\'es list\'ees (voir ci-dessous), et que la neuvi\`eme identit\'e
``$\overset{ 9}{ \equiv}$'' puisse \^etre obtenue comme l'une des
identit\'es \thetag{ $ \mathcal{ P} lck_2$}.

En effet, au niveau pr\'ec\'edent $\kappa = 3$, nous avions cinq polyn\^omes
invariants fondamentaux\,:
\[
\Big(\
f_1'\ \ \ \ \
f_2'\ \ \ \ \
\Lambda^3\ \ \ \ \
\Lambda_1^5\ \ \ \ \
\Lambda_2^5\
\Big).
\]
Par cons\'equent, le nombre d'identit\'es fondamentales 
\thetag{ $\mathcal{ P} lck_1$}
possibles est \'egal \`a $C_5^3 = 10$, 
et nous les \'ecrivons dans l'ordre suivant\,:
\[
\aligned
0
&
\overset{a}{\equiv}
f_1'\,\big[f_2',\,\Lambda^3\big]
+
3\,\Lambda^3\,\big[f_1',\,f_2'\big]
+
f_2'\,\big[\Lambda^3,\,f_1'\big],
\\
0
&
\overset{b}{\equiv}
f_1'\,\big[f_2',\,\Lambda_1^5\big]
+
5\,\Lambda_1^5\,\big[f_1',\,f_2'\big]
+
f_2'\,\big[\Lambda_1^5,\,f_1'\big],
\\
0
&
\overset{c}{\equiv}
f_1'\,\big[f_2',\,\Lambda_2^5\big]
+
5\,\Lambda_2^5\,\big[f_1',\,f_2'\big]
+
f_2'\,\big[\Lambda_2^5,\,f_1'\big],
\\
0
&
\overset{d}{\equiv}
f_1'\big[\Lambda^3,\,\Lambda_1^5\big]
+
5\,\Lambda_1^5\,\big[f_1',\,\Lambda^3\big]
+
3\,\Lambda^3\,\big[\Lambda_1^5,\,f_1'\big],
\\
0
&
\overset{e}{\equiv}
f_1'\,\big[\Lambda^3,\,\Lambda_2^5\big]
+
5\,\Lambda_2^5\,\big[f_1',\,\Lambda^3\big]
+
3\,\Lambda^3\,\big[\Lambda_2^5,\,f_1'\big],
\\
0
&
\overset{f}{\equiv}
f_1'\,\big[\Lambda_1^5,\,\Lambda_2^5\big]
+
5\,\Lambda_2^5\,\big[f_1',\,\Lambda_1^5\big]
+
5\,\Lambda_1^5\,\big[\Lambda_2^5,\,f_1'\big],
\\
0
&
\overset{g}{\equiv}
f_2'\,\big[\Lambda^3,\,\Lambda_1^5\big]
+
5\,\Lambda_1^5\,\big[f_2',\,\Lambda^3\big]
+
3\,\Lambda^3\,\big[\Lambda_1^5,\,f_2'\big],
\\
0
&
\overset{h}{\equiv}
f_2'\,\big[\Lambda^3,\,\Lambda_2^5\big]
+
5\,\Lambda_2^5\,\big[f_2',\,\Lambda^3\big]
+
3\,\Lambda^3\,\big[\Lambda_2^5,\,f_2'\big],
\\
0
&
\overset{i}{\equiv}
f_2'\,\big[\Lambda_1^5,\,\Lambda_2^5\big]
+
5\,\Lambda_2^5\,\big[f_2',\,\Lambda_1^5\big]
+
5\,\Lambda_1^5\,\big[\Lambda_2^5,\,f_2'\big],
\\
0
&
\overset{j}{\equiv}
3\,\Lambda^3\,\big[\Lambda_1^5,\,\Lambda_2^5\big]
+
5\,\Lambda_2^5\,\big[\Lambda^3,\,\Lambda_1^5\big]
+
5\,\Lambda_1^5\,\big[\Lambda_2^5,\,\Lambda^3\big].
\endaligned
\]
De m\^eme, le nombre d'identit\'es pl\"ucke\-riennes \thetag{
$\mathcal{ P}lck_2$} possibles est \'egal au 
nombre $C_5^4 = 5$\,:
\[
\aligned
0
&
\overset{k}{\equiv}
\big[f_1',\,f_2'\big]\cdot\big[\Lambda^3,\,\Lambda_1^5\big]
+
\big[\Lambda_1^5,\,f_1'\big]\cdot\big[\Lambda^3,\,f_2'\big]
+
\big[f_2',\,\Lambda_1^5\big]\cdot\big[\Lambda^3,\,f_1'\big],
\\
0
&
\overset{l}{\equiv}
\big[f_1',\,f_2'\big]\cdot\big[\Lambda^3,\,\Lambda_1^5\big]
+
\big[\Lambda_2^5,\,f_1'\big]\cdot\big[\Lambda^3,\,f_2'\big]
+
\big[f_2',\,\Lambda_2^5\big]\cdot\big[\Lambda^3,\,f_1'\big],
\\
0
&
\overset{m}{\equiv}
\big[f_1',\,f_2'\big]\cdot\big[\Lambda_1^5,\,\Lambda_2^5\big]
+
\big[\Lambda_2^5,\,f_1'\big]\cdot\big[\Lambda_1^5,\,f_2'\big]
+
\big[f_2',\,\Lambda_2^5\big]\cdot\big[\Lambda_1^5,\,f_1'\big],
\\
0
&
\overset{n}{\equiv}
\big[f_1',\,\Lambda^3\big]\cdot\big[\Lambda_1^5,\,\Lambda_2^5\big]
+
\big[\Lambda_2^5,\,f_1'\big]\cdot\big[\Lambda_1^5,\,\Lambda^3\big]
+
\big[\Lambda^3,\,\Lambda_2^5\big]\cdot\big[\Lambda_1^5,\,f_1'\big],
\\
0
&
\overset{o}{\equiv}
\big[f_2',\,\Lambda^3\big]\cdot\big[\Lambda_1^5,\,\Lambda_2^5\big]
+
\big[\Lambda_2^5,\,f_1'\big]\cdot\big[\Lambda_1^5,\,\Lambda^3\big]
+
\big[\Lambda^3,\,\Lambda_2^5\big]\cdot\big[\Lambda_1^5,\,f_2'\big].
\endaligned
\]
\`A pr\'esent, 
nous pouvons r\'e\'ecrire
tous ces crochets bruts en utilisant les notations que
nous avons introduites pour 
d\'esigner nos neuf invariants, tout d'abord dans les
dix identit\'es \thetag{ $\mathcal{ P}lck_1$}\,:
\[
\aligned
0
&
\overset{a}{\equiv}
-f_1'\,\Lambda_2^5
-
3\,\Lambda^3\,\Lambda^3
+
f_2'\,\Lambda_1^5,
\\
0
&
\overset{b}{\equiv}
-f_1'\,\Lambda_{1,2}^7
-
5\,\Lambda^3\,\Lambda_1^5
+
f_2'\,\Lambda_{1,1}^7,
\\
0
&
\overset{c}{\equiv}
-f_1'\,\Lambda_{2,2}^7
-
5\,\Lambda^3\,\Lambda_2^5
+
f_2'\,\Lambda_{2,1}^7,
\\
0
&
\overset{d}{\equiv}
-f_1'f_1'\,M^8
-
5\,\Lambda_1^5\,\Lambda_1^5
+
3\,\Lambda^3\,\Lambda_{1,1}^7,
\\
0
&
\overset{e}{\equiv}
-f_1'f_2'\,M^8
-
5\,\Lambda_1^5\,\Lambda_2^5
+
3\,\Lambda^3\,\Lambda_{2,1}^7,
\\
0
&
\overset{f}{\equiv}
-5\,f_1'\,\Lambda^3\,M^8
-
5\,\Lambda_2^5\,\Lambda_{1,1}^7
+
5\,\Lambda_1^5\,\Lambda_{2,1}^7,
\\
0
&
\overset{g}{\equiv}
-f_2'f_1'\,M^8
-
5\,\Lambda_2^5\,\Lambda_1^5
+
3\,\Lambda^3\,\Lambda_{1,2}^7,
\\
0
&
\overset{h}{\equiv}
-f_2'f_2'\,M^8
-
5\,\Lambda_2^5\,\Lambda_2^5
+
3\,\Lambda^3\,\Lambda_{2,2}^7,
\\
0
&
\overset{i}{\equiv}
-5\,f_2'\,\Lambda^3\,M^8
-
5\,\Lambda_2^5\,\Lambda_{1,2}^7
+
5\,\Lambda_1^5\,\Lambda_{2,2}^7,
\\
0
&
\overset{j}{\equiv}
-3\,\Lambda^3\,\Lambda^3\,M^8
-
\Lambda_2^5\,f_1'\,M^8
+
\Lambda_1^5\,f_2'\,M^8\,;
\endaligned
\]
et ensuite dans les cinq identit\'es \thetag{ $\mathcal{ P} lck_2$}\,:
\[
\aligned
0
&
\overset{k}{\equiv}
\Lambda^3\,f_1'\,M^8
+
\Lambda_{1,1}^7\,\Lambda_2^5
-
\Lambda_{1,2}^7\,\Lambda_1^5,
\\
0
&
\overset{l}{\equiv}
\Lambda^3\,f_2'\,M^8
+
\Lambda_{2,1}^7\,\Lambda_2^5
-
\Lambda_{2,2}^7\,\Lambda_1^5,
\\
0
&
\overset{m}{\equiv}
5\,\Lambda^3\,\Lambda^3\,M^8
+
\Lambda_{2,1}^7\,\Lambda_{1,2}^7
-
\Lambda_{2,2}^7\,\Lambda_{1,1}^7,
\\
0
&
\overset{n}{\equiv}
5\,\Lambda_1^5\,\Lambda^3\,M^8
+
\Lambda_{2,1}^7\,f_1'\,M^8
-
\Lambda_{1,1}^7\,f_2'\,M^8,
\\
0
&
\overset{o}{\equiv}
5\,\Lambda_2^5\,\Lambda^3\,M^8
+
\Lambda_{2,2}^7\,f_1'\,M^8
-
f_2'\,\Lambda_{1,2}^7\,M^8.
\endaligned
\]
Ici, en admettant bien sûr
que $\Lambda_{ 1,2} = \Lambda_{ 2,1}$, on
constate que\,:

\begin{itemize}

\smallskip\item[$\bullet$]
``$\overset{ a}{ \equiv}$'' fournit 
``$\overset{ 1}{ \equiv}$''\,;

\smallskip\item[$\bullet$]
``$\overset{ b}{ \equiv}$'' fournit
``$\overset{ 2}{ \equiv}$''\,;

\smallskip\item[$\bullet$]
``$\overset{ c}{ \equiv}$'' fournit
``$\overset{ 3}{ \equiv}$''\,;

\smallskip\item[$\bullet$]
``$\overset{ d}{ \equiv}$'' fournit
``$\overset{ 4}{ \equiv}$''\,;

\smallskip\item[$\bullet$]
``$\overset{ e}{ \equiv}$'' fournit
``$\overset{ 5}{ \equiv}$''\,;

\smallskip\item[$\bullet$]
``$\overset{ f}{ \equiv}$'' fournit
``$\overset{ 7}{ \equiv}$''\,;

\smallskip\item[$\bullet$]
``$\overset{ g}{ \equiv}$'' est redondant avec 
``$\overset{ e}{ \equiv}$''\,;

\smallskip\item[$\bullet$]
``$\overset{ h}{ \equiv}$'' fournit 
``$\overset{ 6}{ \equiv}$''\,;

\smallskip\item[$\bullet$]
``$\overset{ i}{ \equiv}$'' fournit
``$\overset{ 8}{ \equiv}$''\,;

\smallskip\item[$\bullet$]
``$\overset{ j}{ \equiv}$'' redouble 
``$\overset{ 1}{ \equiv}$'' en la multipliant par $M^8$\,;

\smallskip\item[$\bullet$]
``$\overset{ k}{ \equiv}$'' redouble 
``$\overset{ f}{ \equiv}$''\,;

\smallskip\item[$\bullet$]
``$\overset{ l}{ \equiv}$'' redouble
``$\overset{ i}{ \equiv}$''\,;

\smallskip\item[$\bullet$]
``$\overset{ m}{ \equiv}$'' fournit la derni\`ere identit\'e manquante 
``$\overset{ 9}{ \equiv}$''\,;

\smallskip\item[$\bullet$]
``$\overset{ n}{ \equiv}$'' redouble 
``$\overset{ b}{ \equiv}$'' en la multipliant par $M^8$\,;

\smallskip\item[$\bullet$]
``$\overset{ o}{ \equiv}$'' redouble 
``$\overset{ c}{ \equiv}$'' en la multipliant par $M^8$.

\end{itemize}\smallskip

\noindent{\bf Conclusion.} Toutes les identit\'es alg\'ebriques
entre les invariants fondamentaux que nous
avons trouv\'ees \`a l'aide de Maple
pour les niveaux
$\kappa = 3$ et $\kappa = 4$ peuvent en 
fait \^etre obtenues m\'ecaniquement gr\^ace aux trois familles
fondamentales de syzygies\,:
\[
\boxed{
\aligned
(\mathcal{J}ac)\,:
\ \ \ \ \ \ \ \ \ \ \ \
0
&
\equiv
\big[\big[{\sf P},\,{\sf Q}\big],\,{\sf R}\big]
+
\big[\big[{\sf R},\,{\sf P}\big],\,{\sf Q}\big]
+
\big[\big[{\sf Q},\,{\sf R}\big],\,{\sf P}\big],
\\
(\mathcal{R}\mathcal{F})\,:
\ \ \ \ \ \ \ \ \ \ \ \
0
&
\equiv
m\,{\sf P}\,\big[{\sf Q},\,{\sf R}\big]
+
o\,{\sf R}\,\big[{\sf P},\,{\sf Q}\big]
+
n\,{\sf Q}\,\big[{\sf R},\,{\sf P}\big],
\\
(\mathcal{P}lck)\,:
\ \ \ \ \ \ \ \ \ \ \ \
0
&
\equiv
\big[{\sf P},\,{\sf Q}\big]
\cdot
\big[{\sf R},\,{\sf S}\big]
+
\big[{\sf S},\,{\sf P}\big]
\cdot
\big[{\sf R},\,{\sf Q}\big]
+
\big[{\sf Q},\,{\sf S}\big]
\cdot
\big[{\sf R},\,{\sf P}\big]
\endaligned
}\,.
\]

\noindent{\bf Gen\`ese des syzygies.}
Crucialement, il semblerait que l'on puisse engendrer toutes les
relations entre tous les polyn\^omes invariants que l'on construit
r\'ecursivement par crochets, juste en d\'eveloppant par r\'ecurrence
toutes les identit\'es \thetag{ $\mathcal{ J}ac$}, \thetag{ $\mathcal{
P}lck_1$} et \thetag{ $\mathcal{ P}lck_2$} possibles lorsqu'on passe
d'un \'etage de jets $\lambda$ \`a l'\'etage sup\'erieur $\lambda +
1$. Cette id\'ee conjecturale est renforc\'ee par le fait que dans la
th\'eorie classique des invariants, il existe aussi trois proc\'ed\'es
automatiques qui engendrent l'id\'eal des relations entre les
invariants, et on d\'emontre rigoureusement (\cite{ ol1999}) que tel est
bien le cas, sans toutefois poursuivre l'\'etude plus avant, afin de
trouver des bases de Gr\"obner signifiantes d'un point combinatoire,
ou afin de d\'evoiler des harmonies formelles encore inconnues qui
montreraient explicitement en quoi l'alg\`ebre des invariants est de
Cohen-Macaulay, ce qui est toujours le cas pour
un groupe r\'eductif (\cite{ ol1999}).

\section*{ \S6.~D\'ecomposition en repr\'esentations de Schur} 

\noindent{\bf Motivation.} 
La cohomologie des fibr\'es de Schur sur une vari\'et\'e projective lisse
\'etant connue ({\it cf.} {\it e.g.} \cite{ ro2007} et {\it voir}\,
la Section~8 ci-dessous), nous cherchons maintenant \`a d\'ecomposer en
repr\'esentations irr\'eductibles de Schur les gradu\'es de poids $m$ de nos
deux alg\`ebres d'invariants $\mathcal{ DS}_2^4$ et $\mathcal{ DS}_2^5$.

\smallskip\noindent{\bf Action lin\'eaire diagonale sur les jets.} 
\`A cette fin, sur l'espa\-ce des jets d'ordre $\kappa$ en dimension
deux muni des coordonn\'ees $\big( f_1', f_2', \dots, f_1^{ (\lambda)},
f_2^{ ( \lambda)}, \dots, f_1^{ (\kappa)}, f_2^{ ( \kappa)} \big)$,
consid\'erons ({\it cf.} \cite{ ro2007}) l'action du groupe lin\'eaire \`a
deux dimensions ${\sf GL}_2 ( \C)$\,\,---\,\,constitu\'e des matrices $2
\times 2$ de la forme
\[
\text{\sc w}
:=
\left(
\begin{array}{cc}
t & v
\\
u & w
\end{array}
\right),
\]
o\`u $t, u, v, w \in \C$ satisfont $tw - uv \neq 0$\,\,---\,\,qui est
d\'efinie diagonalement par la m\^eme transformation 
\'evidente sur chaque \'etage de
jets\,:
\[
\aligned
{\sf w}\cdot
f_1^{(\lambda)}
&
:=
t\,f_1^{(\lambda)}
+
v\,f_2^{(\lambda)},
\\
{\sf w}\cdot
f_2^{(\lambda)}
&
:=
u\,f_1^{(\lambda)}
+
w\,f_2^{(\lambda)},
\endaligned
\]
pour tout $\lambda$ tel que $1 \leqslant \lambda \leqslant \kappa$.

\smallskip\noindent{\bf D\'ecompositions de Schur.}
La th\'eorie classique des repr\'esentations du groupe lin\'eaire permet
alors de d\'ecomposer toute repr\'esentation 
de ${\sf GL}_2 ( \C)$ comme
somme directe de repr\'esentations d'un certain type, dites {\sl de
Schur}, que l'on rep\`ere facilement en recherchant tous les vecteurs
qui sont invariants par un certain sous-groupe de ${\sf GL}_2 ( \C)$.
\'Enon\c cons ce que la th\'eorie g\'en\'erale donne dans le cas qui nous
int\'eresse.

\smallskip\noindent{\bf D\'efinition.} 
Un polyn\^ome invariant par reparam\'etrisation ${\sf P} \big( j^\kappa f
\big)$ est appel\'e {\sl bi-invariant}\, s'il est un {\sl vecteur de
plus haut poids}\, pour cette repr\'esentation, c'est-\`a-dire s'il est
invariant par l'action du sous-groupe unipotent ${\sf U}_2 ( \C)$
constitu\'e des matrices de la forme\,:
\[
\text{\sc u}
:=
\left(
\begin{array}{cc}
1&0
\\
u&1
\end{array}
\right).
\]
Autrement dit, un invariant simple ${\sf P}$ satisfait ${\sf P} \big(
j^\kappa ( f \circ \phi) \big) = ( \phi')^m \, {\sf P} \big( (
j^\kappa f) \circ \phi \big)$ pour un certain $m \geqslant 1$ et c'est
un {\sl bi-invariant} si l'on a de plus\,:
\[
\boxed{
{\sf P}^{2\times{\rm inv}}\big(\text{\sc u}\cdot j^\kappa f\big)
=
{\sf P}^{2\times{\rm inv}}\big(j^\kappa f\big)
}\,,
\]
pour toute matrice unipotente $\text{\sc u} \in {\sf U}_2 ( \C)$. 

\smallskip\noindent{\bf Exemples.} Puisque l'on a trivialement
$\text{\sc u} \cdot f_1' = f_1'$ et $\text{\sc u} \cdot f_1'' =
f_1''$, et aussi\,:
\[
\text{\sc u}\cdot\Delta^{\alpha,\beta}
=
\left\vert
\begin{array}{cc}
f_1^{(\alpha)}\ \ &f_2^{(\alpha)}+u\,f_1^{(\alpha)}
\\
f_1^{(\beta)}\ \ &f_2^{(\beta)}+u\,f_1^{(\beta)}
\end{array}
\right\vert
=
\Delta^{\alpha,\beta},
\]
nous voyons imm\'ediatement que $f_1'$, $\Lambda^3$, $\Lambda_1^5$,
$\Lambda_{ 1, 1}^7$ et $M^8$ sont des bi-invariants (au nombre de
cinq), tandis que les quatre invariants restants, \`a savoir $f_2'$,
$\Lambda_2^5$, $\Lambda_{ 1, 2}^7$ et $\Lambda_{ 2, 2}^7$ ne sont pas
bi-invariants.

\smallskip\noindent{\bf Rep\'erage des repr\'esentations de Schur.} 
D'apr\`es la th\'eorie g\'en\'erale, \`a tout vecteur ${\sf P}^{ 2\times {\rm
inv}}$ de plus haut poids correspond alors une et une seule
repr\'esentation de Schur $\Gamma^{ (l_1, l_2 )}$, o\`u les deux
entiers $l_1$ et $l_2$ satisfaisant $l_1 \geqslant
l_2$ sont ais\'ement rep\'er\'es comme \'etant les exposants des deux
\'el\'ements diagonaux qui apparaissent dans la valeur propre
\[
{\sf t}\cdot{\sf P}^{2\times{\rm inv}}
=
\left(
\begin{array}{cc}
t & 0 
\\
0 & w
\end{array}
\right)
\cdot
{\sf P}^{2\times{\rm inv}}
=
{\sf t}^{l_1}\,
{\sf w}^{l_2}\,
{\sf P}^{2\times{\rm inv}}
\]
dont jouit le bi-invariant\,\,---\,\,qui est n\'ecessairement vecteur
propre\,\,---\,\,par rapport au sous-groupe des matrices $2 \times 2$
diagonales. Tr\`es concr\`etement, l'entier $l_1$ compte le nombre
total d'indices inf\'erieurs ``${(\cdot)}_1$'' qui interviennent dans
chaque mon\^ome du bi-invariant en question, et de m\^eme, l'entier
$l_2$ compte le nombre d'indices ``${(\cdot)}_2$'', et puisque
chaque $\Delta^{ \alpha, \beta}$ contribue pour exactement un indice
``${(\cdot)}_1$'' et un indice ``${(\cdot)}_2$'', il est imm\'ediatement
clair que nous avons la correspondance suivante entre bi-invariants et
repr\'esentations de Schur\,:
\[
\aligned
&
f_1'
\longleftrightarrow
\Gamma^{(1,0)},
\ \ \ \ \ \ \ \ \ \
\Lambda^3
\longleftrightarrow
\Gamma^{(1,1)},
\\
\Lambda_1^5
\longleftrightarrow
\Gamma^{(2,1)},&
\ \ \ \ \ \ \ \ \ \
\Lambda_{1,1}^7
\longleftrightarrow
\Gamma^{(3,1)},
\ \ \ \ \ \ \ \ \ \
M^8
\longleftrightarrow
\Gamma^{(2,2)}.
\endaligned
\]

\noindent{\bf Fait d'exp\'erience.} 
La d\'etermination directe des bi-invariants en dimension $\nu = 2$ pour
les jets d'ordre $\kappa = 4$ ou $\kappa = 5$ est {\it beaucoup moins
coûteuse en calcul}\, que la d\'etermination de la totalit\'e des
invariants par reparam\'etrisation. Voici en effet
le premier de nos deux r\'esultats principaux.

\def\thetheorem{\!\!}\begin{theorem}
Pour les jets d'ordre 4 en dimension 2, tout bi-invariant de poids $m$, 
${\sf P}^{ 2
\times {\rm inv}} \big( j^4 f_1, \, j^4 f_2 \big)$,
s'\'ecrit sous forme unique{\rm \,:}
\[
{\sf P}^{2\times{\rm inv}}\big(j^4f\big)
=
\mathcal{Q}^{2\times{\rm inv}}
\big(f_1',\Lambda^3,\Lambda_{1,1}^7,M^8\big)
+
\Lambda_1^5\,\mathcal{R}^{2\times{\rm inv}}
\big(f_1',\Lambda^3,\Lambda_{1,1}^7,M^8\big),
\]
o\`u $\mathcal{ Q}^{2 \times {\rm inv }}$ et $\mathcal{ R}^{ 2 \times
{\rm inv}}$ sont deux polyn\^omes absolument arbitraires en leurs
arguments qui sont de poids $m$ et de poids $m - 5$,
respectivement. De plus, l'id\'eal des relations entre les cinq
bi-invariants fondamentaux{\rm \,:}
\[
\Big(\
f_1'\ \ \ \ \ 
\Lambda^3\ \ \ \ \ 
\Lambda_1^5\ \ \ \ \ 
\Lambda_{1,1}^7 \ \ \ \ \ 
M^8\
\Big),
\] 
est principal, et pour pr\'eciser, il est engendr\'e par
l'unique\footnote{\, Ce fait a \'et\'e confirm\'e par un calcul de l'id\'eal
des relations sur Maple. } relation{\rm \,:}
\[
0
\equiv
f_1'f_1'M^8
-
3\,\Lambda^3\Lambda_{1,1}^7
+
5\,\Lambda_1^5\Lambda_1^5.
\]
Par cons\'equent, une base de l'espace vectoriel des polyn\^omes de poids
$m$ invariant par reparam\'etrisation et par rapport \`a l'action de ${\sf
U}_2 ( \C)$ est constitu\'ee de l'ensemble des mon\^omes{\rm \,:}
\[
\aligned
&
(f_1')^\alpha\,\big(\Lambda^3\big)^\beta\,
\big(\Lambda_{1,1}^7\big)^\gamma\,\big(M^8\big)^\delta,
\ \ \ \ \ \ \ \ \
\text{\rm avec}\ \
\alpha+3\beta+7\gamma+8\delta=m,
\ \
\text{\rm et\,:}
\\
&
\Lambda_1^5\,(f_1')^\alpha\,\big(\Lambda^3\big)^\beta\,
\big(\Lambda_{1,1}^7\big)^\gamma\,\big(M^8\big)^\delta,
\ \ \ \ \ \ \ \ \
\text{\rm avec}\ \
\alpha+3\beta+7\gamma+8\delta=m-5,
\endaligned
\]
et chacun de ces deux mon\^omes correspond respectivement aux deux
repr\'esenta\-tions de Schur{\rm \,:}
\[
\Gamma^{\alpha+\beta+3\gamma+2\delta,\,\,\beta+\gamma+2\delta}
\ \ \ \ \ \ \ \ \ \ \ \ \ \ \ \ \ \ 
\text{\it et}
\ \ \ \ \ \ \ \ \ \ \ \ \ \ \ \ \ \ 
\Gamma^{2+\alpha+\beta+3\gamma+2\delta,\,\,1+\beta+\gamma+2\delta}.
\]
\end{theorem}

\noindent{\bf Cons\'equences.}
Avant d'entreprendre la d\'emonstration de ce premier th\'eor\`eme, notons
que l'alg\`ebre compl\`ete $\mathcal{ DS}_2^4$ des invariants par
reparam\'etrisation s'obtient maintenant facilement en regardant
l'orbite, par l'action du groupe complet ${\sf GL}_2 ( \C)$, de chacun
de nos cinq bi-invariants\,; on constate d'ailleurs qu'il suffit de
consid\'erer l'action des matrices de la forme
\[
\text{\sc v}
:=
\left(
\begin{array}{cc}
1 & v
\\
0 & 1
\end{array}
\right),
\]
qui nous fournissent imm\'ediatement\,:
\[
\aligned
\text{\sc v}\cdot f_1'
=
&\
f_1'
+
v\,f_2',
\ \ \ \ \ \ \ \ \ \ \ \ \ \ \ \
\text{\sc v}\cdot\Lambda_1^5
=
\Lambda_1^5
+
v\,\Lambda_2^5,
\\
&
\text{\sc v}\cdot\Lambda_{1,1}^7
=
\Lambda_{1,1}^7
+
2v\,\Lambda_{1,2}^7
+
v^2\,\Lambda_{2,2}^7,
\endaligned
\]
et de cette mani\`ere, non seulement nous engendrons facilement les
quatre invariants fondamentaux non bi-invariants que nous connaissions
d\'ej\`a, mais encore\,\,---\,\,et c'est l\`a qu'appara\^{\i}t {\it une strat\'egie
crucialement simplifi\'ee que nous r\'e-exploiterons ult\'erieurement pour
l'\'etude de $\mathcal{ DS}_2^5$}\,\,---, nous d\'eduisons que l'orbite des
polyn\^omes arbitraires en $\big( f_1', \, \Lambda^3, \, \Lambda_1^5, \,
\Lambda_{ 1, 1}^7, \, M^8 \big)$ est juste constitu\'ee des polyn\^omes en
les neuf invariants par reparam\'etrisation que nous avions engendr\'es en
calculant m\'ethodi\-quement des crochets.

\def\thecorollary{\!}\begin{corollary}
Pour les jets d'ordre quatre en dimension deux, l'alg\`ebre
$\mathcal{DS}_2^4$ des polyn\^omes invariants par reparam\'etrisation est
polynomialement engendr\'ee par les neuf invariants fondamentaux $\big(
f_1', \, f_2', \, \Lambda^3, \, \Lambda_1^5, \, \Lambda_2^5, \,
\Lambda_{ 1, 1}^7, \, \Lambda_{ 1, 2}^7, \, \Lambda_{ 2, 2}^7, \, M^8
\big)$.
\end{corollary}

\noindent{\bf Restrictions.}
Toutefois, cette mani\`ere \'economique de proc\'eder\,\,---\,\,\'etude
exclusive et exhaustive des bi-invariants suivie de la d\'eduction
raccourcie d'une description partielle de l'alg\`ebre compl\`ete des
invariants\,\,---\,\,ne fournit pas de description pr\'ecise de
$\mathcal{ DS}_2^4$, c'est-\`a-dire notamment qu'elle ne fournit pas une
\'ecriture unique, en tenant compte des 9 syzygies fondamentales, de
tout polyn\^ome g\'en\'eral de la forme\,:
\[
\mathcal{P}
\big(
f_1',\,f_2',\,\Lambda^3,\,\Lambda_1^5,\,\Lambda_2^5,\,
\Lambda_{1,1}^7,\,\Lambda_{1,2}^7,\,\Lambda_{2,2}^7,\,M^8
\big),
\]
et qui plus est, il serait impossible d'obtenir un tel r\'esultat
complet, et ce pour une raison profonde, \`a savoir que l'orbite
par ${\sf GL}_2 ( \C)$ de l'unique syzygie existant entre les
bi-invariants ne couvre pas l'ensemble des neuf syzygies fondamentales
qui existent entre les invariants complets.

Heureusement, puisque seule la d\'ecomposition en repr\'esentations
irr\'eductibles de Schur pr\'esente un v\'eritable sens alg\'ebrique, et
aussi, puisque nous aurons seulement besoin de cette d\'ecomposition
pour conduire nos calculs de caract\'eristique d'Euler dans la
Section~8, il est en v\'erit\'e essentiellement inutile de poursuivre plus
avant l'\'etude de $\mathcal{ DS}_2^4$. Nous confierons quand m\^eme au
lecteur d\'esireux de s'exercer \`a ma\^{\i}triser les bases de Gr\"obner le
soin d'\'etablir l'\'enonc\'e suivant, ou d'autres \'enonc\'es analogues qu'il
pourrait formuler en choisissant \`a sa guise des ordres monomiaux
diff\'erents.

\def\theproposition{\!}\begin{proposition}
Tout polyn\^ome en les neuf invariants fondamentaux 
s'\'ecrit de mani\`ere unique comme suit{\rm \,:}
\[
\aligned
&
\mathcal{P}
\big(f_1',\,f_2',\,\Lambda^3,\,\Lambda_{1,1}^7,\,\Lambda_{2,2}^7\big)
+
\Lambda^3\,
\mathcal{Q}
\big(f_1',\,f_2',\,\Lambda^3,\,\Lambda_{1,1}^7,\,\Lambda_{2,2}^7\big)
+
\\
&
+
\Lambda_1^5\,
\mathcal{R}
\big(f_1',\,f_2',\,\Lambda^3,\,\Lambda_{1,1}^7,\,\Lambda_{2,2}^7\big)
+
\Lambda_2^5\,
\mathcal{S}
\big(f_1',\,f_2',\,\Lambda^3,\,\Lambda_{1,1}^7,\,\Lambda_{2,2}^7\big)
+
\\
&
+
\Lambda_{1,2}^7\,
\mathcal{T}
\big(f_1',\,f_2',\,\Lambda^3,\,\Lambda_{1,1}^7,\,\Lambda_{2,2}^7\big)
+
\Lambda^3\Lambda_{1,2}^7\,
\mathcal{U}
\big(f_1',\,f_2',\,\Lambda^3,\,\Lambda_{1,1}^7,\,\Lambda_{2,2}^7\big),
\endaligned
\]
o\`u $\mathcal{ P}$, $\mathcal{ Q}$, $\mathcal{ R}$, 
$\mathcal{ S}$, $\mathcal{ T}$ et $\mathcal{ U}$ sont
des polyn\^omes arbitraires en leurs arguments.
\end{proposition}

\noindent{\bf D\'emonstration du premier th\'eor\`eme.}
Par d\'efinition de l'invariance par reparam\'e\-trisation d'un polyn\^ome
${\sf P} = {\sf P} \big( j^4 f)$ de poids $m$, on a\,:
\[
{\sf P}\big(
j^4(f\circ\phi)\big) 
= 
{\phi'}^m\,{\sf P}\big((j^4f)
\circ\phi\big),
\]
pour tout biholomorphisme local $\phi$ de $\C$. En suivant une astuce
de \cite{ ro2007}, nous allons appliquer cette formule \`a $\phi :=
f_1^{ - 1}$ en supposant l'inversibilit\'e, d'o\`u $\phi '= \frac{ 1}{
f_1'} \circ f_1^{ - 1}$. On a tout d'abord trivialement $\big( f_1
\circ f_1^{ - 1} \big) ' = {\rm Id}$, d'o\`u $\big( f_1 \circ f_1^{ -
1} \big)^{ ( \lambda)} = 0$ pour tout $\lambda \geqslant 2$ puis, par
des calculs directs dont la teneur est d\'ej\`a \'elucid\'ee par notre
connaissance pr\'ealable des invariants $\Lambda^3$, $\Lambda_1^5$ et
$\Lambda_{ 1,1}^7$\,:
\[
\aligned
\big(f_2\circ f_1^{-1}\big)'
&
=
\frac{f_2'}{f_1'}\circ f_1^{-1},
\\
\big(f_2\circ f_1^{-1}\big)''
&
=
\frac{\Lambda^3}{(f_1')^3}\circ f_1^{-1},
\\
\big(f_2\circ f_1^{-1}\big)'''
&
=
\frac{\Lambda_1^5}{(f_1')^5}\circ f_1^{-1},
\\
\big(f_2\circ f_1^{-1}\big)''''
&
=
\frac{\Lambda_{1,1}^7}{(f_1')^7}\circ f_1^{-1}.
\endaligned
\]
Par cons\'equent, tout polyn\^ome ${\sf P} ( j^4 f)$
invariant par reparam\'etrisation 
satisfait\,:
\[
{\sf P}
\bigg(
1,\frac{f_2'}{f_1'},
0,\frac{\Lambda^3}{(f_1')^3},
0,\frac{\Lambda_1^5}{(f_1')^5},
0,\frac{\Lambda_{1,1}^7}{(f_1')^7}
\bigg)
\circ f_1^{-1}
=
\Big(
\frac{1}{f_1'}\circ f_1^{-1}
\Big)^m\,
{\sf P}
\big((j^4f\circ f_1^{-1})\big).
\]
Recomposons imm\'ediatement avec $f_1$ pour faire dispara\^{\i}tre
$f_1^{ -1}$. Si ensuite nous \'ecrivons le polyn\^ome de d\'epart
${\sf P} = {\sf P} \big( j^4 f \big)$ sous la forme g\'en\'erale
suivante\,:
\[
\footnotesize
\aligned
\sum_{a_1+a_2+2b_1+2b_2+\cdots+4d_2=m}\,
{\sf p}_{a_1a_2\cdots d_2}\cdot\,
(f_1')^{a_1}(f_2')^{a_2}
(f_1'')^{b_1}(f_2'')^{b_2}
(f_1''')^{c_1}(f_2''')^{c_2}
(f_1'''')^{d_1}(f_2'''')^{d_2},
\endaligned
\]
avec des coefficients ${\sf p}_{ a_1 a_2 \cdots d_2} \in \C$,
l'identit\'e obtenue \`a l'instant nous permet alors d'obtenir une
repr\'esentation g\'en\'erale de ${\sf P}$\,:
\[
\small
\aligned
{\sf P}\big(j^4f\big)
&
=
(f_1')^m\,
{\sf P}
\bigg(
1,\frac{f_2'}{f_1'},
0,\frac{\Lambda^3}{(f_1')^3},
0,\frac{\Lambda_1^5}{(f_1')^5},
0,\frac{\Lambda_{1,1}^7}{(f_1')^7}
\bigg)
\\
&
=
(f_1')^m\,
\sum_{a_1+a_2+2b_2+3c_2+4d_2=m}\,
{\sf p}_{a_1a_20b_20c_20d_2}\,
\frac{1^{a_1}(f_2')^{a_2}\big(\Lambda^3\big)^{b_2}
\big(\Lambda_1^5\big)^{c_2}
\big(\Lambda_{1,1}^7\big)^{d_2}}{
(f_1')^{a_2+3b_2+5c_2+7d_2}}
\\
&
\ \ \ \ \ \ \ \ \ \ \ \ \ \ \ \ \ \ \ \ \ \ \ \ \ \ \ \
\ \ \ \ \ \ \ \ \ \ \ \ \ \ \ \ \ \ \ \ \ \ \ \ \ \ \ \
\in
\C\big[f_1',f_2',\Lambda^3,\Lambda_1^5,\Lambda_{1,1}^7\big]
\Big[
\frac{1}{f_1'}
\Big],
\endaligned
\]
qui est presque polynomiale, \`a ceci pr\`es qu'on s'autorise \`a diviser
par $f_1'$. Calculons alors l'ordre maximal en $\frac{ 1}{ f_1'}$
de cette expression rationnelle\,:
\[
\footnotesize
\aligned
\max_{a_1+a_2+2b_2+3c_2+4d_2=m}\,
\Big(a_2+3b_2+5c_2+7d_2-m\Big)
&
=
\max_{a_2+2b_2+3c_2+4d_2=m}
\Big(b_2+2c_2+3d_2\Big)
\\
&
=
\frac{1}{2}\cdot\,
\max_{2b_2+3c_2+4d_2=m}
\Big(2b_2+4c_2+6d_2\Big)
\\
&
=
\frac{m}{2}
+
\frac{1}{2}\cdot\,
\max_{3c_2+4d_2=m}
\Big(c_2+2d_2\Big)
\\
&
=
\frac{3}{4}\,m.
\endaligned
\]
Ainsi, tout polyn\^ome ${\sf P} \big( j^4 f \big) \in \mathcal{ DS}_{ 2,
m}^4$ est de la forme\,:
\[
\sum_{-\frac{3}{4}m\leqslant a\leqslant m}\,
(f_1')^a\,\mathcal{P}_a
\big(
f_2',\Lambda^3,\Lambda_1^5,\Lambda_{1,1}^7
\big).
\]
Cependant, toutes les expressions rationnelles de cette forme ne
conviennent pas\,: $\frac{ \Lambda^3 \, \Lambda^3}{ f_1' f_1'}$ avec
$m=4$ ne se simplifie pas pour produire un vrai polyn\^ome appartenant
$\mathcal{ DS}_{ 2, m}^4$, bien que cette expression rationnelle soit
invariante par reparam\'etrisation. Toutefois, l'\'enonc\'e suivant, 
que nous transf\'erons directement aux jets d'ordre quelconque,
est clair.

\def\thelemma{\!}\begin{lemma}
Tout polyn\^ome ${\sf P} \big( j^\kappa f)$ en le jet strict
\[
j^\kappa f
:=
\Big(
f_1',\,f_2',\
f_1'',\,f_2'',\dots\dots,\,
f_1^{(\kappa)},\,f_2^{(\kappa)}
\Big)
\]
d'ordre $\kappa \geqslant 1$ d'une application holomorphe locale $f =
(f_1, f_2) : \C \to \C^2$ qui est invariant par reparam\'etrisation,
{\rm i.e.} qui appartient \`a $\mathcal{ DS }_{ 2, m}^\kappa$, peut \^etre
repr\'esent\'e sous la forme{\rm \,:}
\[
{\sf P}\big(j^\kappa f)
=
\sum_{-\frac{\kappa-1}{\kappa}m\leqslant a\leqslant m}\,
(f_1')^a\,\mathcal{P}_a
\big(
f_2',\Lambda^3,\Lambda_1^5,\Lambda_{1,1}^7,\Lambda_{1,1,1}^9,
\dots,
\Lambda_{1,\dots,1}^{2\kappa-1}
\big),
\]
avec certains polyn\^omes $\mathcal{P}_a$ de poids $m - a$. R\'eciproquement,
toute expression rationnelle de cette forme qui s'av\`ere \^etre
polynomiale en $j^\kappa f$ quand on simplifie num\'erateur et
d\'enominateur appartient \`a $\mathcal{ DS }_{ 2, m}^\kappa$.
\hfill$\square$
\end{lemma}

Ici, on consid\`ere comme pr\'ec\'edemment $\Lambda_{ 1, 1, 1}^9 := \big[
\Lambda_{ 1, 1}^7, \, f_1' \big]$, et on introduit g\'en\'eralement par
r\'ecurrence $\Lambda_{ 1, \dots, 1, 1}^{ 2 \lambda - 1} := \big[
\Lambda_{ 1, \dots, 1}^{ 2 \lambda - 3}, \, f_1' \big]$ pour $3
\leqslant \lambda \leqslant \kappa$.

\smallskip

D\'ecrivons maintenant les polyn\^omes ${\sf P} = {\sf P} \big( j^\kappa f
\big)$ de $\mathcal{ DS}_{ 2, m}^\kappa$ \'ecrits sous une telle forme
rationnelle qui sont invariants par l'action de ${\sf U}_2 ( \C)$. Par
d\'efinition, $\text{\sc u}\cdot {\sf P} = {\sf P}$, {\it i.e.}
explicitement\,:
\[
{\sf P}
\Big(
f_1',\,f_2'+u\,f_1',
f_1'',\,f_2''+u\,f_1'',
\dots,
f_1^{(\kappa)},\,f_2^{(\kappa)}+u\,f_1^{(\kappa)}
\Big)
=
{\sf P}
\big(
j^\kappa f
\big),
\]
pour tout $u \in \C$. De mani\`ere \'equivalente,
\[
\frac{d}{du}\,
{\sf P}
\Big(
f_1',\,f_2'+u\,f_1',
f_1'',\,f_2''+u\,f_1'',
\dots,
f_1^{(\kappa)},\,f_2^{(\kappa)}+u\,f_1^{(\kappa)}
\Big)
\equiv
0,
\]
ce qui revient \`a dire que ${\sf P}$ est annul\'e identiquement par le
champ de vecteurs
\[
\underline{\mathcal{U}}
:=
f_1'\,\frac{\partial}{\partial f_2'}
+
f_1''\,\frac{\partial}{\partial f_2''}
+
\cdots
+
f_1^{(\kappa)}\,\frac{\partial}{\partial f_2^{(\kappa)}},
\]
{\it i.e.} que l'on a\,: $0 \equiv \underline{ \mathcal{ U}} \, {\sf P}$.
On constate ensuite imm\'ediatement que\,:
\[
\text{\sc u}\cdot
\Delta^{\alpha,\beta}
=
f_1^{(\alpha)}\big(
f_2^{(\beta)}+u\,f_1^{(\beta)}\big)
-
f_1^{(\beta)}\big(
f_2^{(\alpha)}+u\,f_1^{(\alpha)}
\big)
=
\Delta^{\alpha,\beta},
\]
{\it i.e.}\,: $0 \equiv \underline{ \mathcal{ U}} \, \Delta^{ \alpha,
\beta}$, d'o\`u nous d\'eduisons\,:
\[
\text{\sc u}\cdot\Lambda^3
=
\Lambda^3,
\ \ \ \ \ \ \
\text{\sc u}\cdot\Lambda_1^5
=
\Lambda_1^5,
\ \ \ \ \ \ \
\text{\sc u}\cdot\Lambda_{1,1}^7
=
\Lambda_{1,1}^7,
\ \ \ \ \ \ \
\text{\sc u}\cdot\Lambda_{1,1,1}^9
=
\Lambda_{1,1}^7,
\ \ \ \ \ \ \
\text{\it etc.}
\]
En appliquant donc cette d\'erivation $\underline{ \mathcal{ U}}$ \`a la
repr\'esentation rationnelle d'un polyn\^ome quelconque ${\sf P} \big(
j^\kappa f \big) \in \mathcal{ DS}_{ 2, m}^\kappa$ obtenue \`a
l'instant, nous voyons que l'\'equation $0 \equiv \underline{ \mathcal{
U}} \, {\sf P}$ est satisfaite si et seulement si chaque $\mathcal{ P
}_a$ est ind\'ependant de $f_2'$. Nous pouvons donc r\'esumer comme suit
le r\'esultat obtenu.

\def\thelemma{\!}\begin{lemma}
Tout polyn\^ome ${\sf P}^{2\times {\rm inv}} \big( j^\kappa f)$ qui est
invariant par reparam\'etrisation et qui est invariant par rapport \`a
l'action unipotente
de ${\sf U}_2 ( \C)$ peut \^etre repr\'esent\'e sous la forme{\rm \,:}
\[
{\sf P}^{2\times {\rm inv}}\big(j^\kappa f)
=
\sum_{-\frac{\kappa-1}{\kappa}m\leqslant a\leqslant m}\,
(f_1')^a\,\mathcal{P}_a^{2\times {\rm inv}}
\big(
\Lambda^3,\Lambda_1^5,\Lambda_{1,1}^7,\Lambda_{1,1,1}^9,
\dots,
\Lambda_{1,\dots,1}^{2\kappa-1}
\big),
\]
avec certains polyn\^omes $\mathcal{ P}_a^{2\times {\rm inv}}$ de poids $m -
a$. R\'eciproquement, toute expression rationnelle de cette forme, si
elle s'av\`ere \^etre polynomiale en $j^\kappa f$ quand on simplifie
num\'erateur et d\'enominateur, appartient n\'ecessairement \`a $\mathcal{ DS
}_{2, m}^\kappa$ et constitue un bi-invariant v\'eritable.
\hfill$\square$
\end{lemma}

\smallskip

En revenant \`a pr\'esent aux jets d'ordre 4, utilisons la
relation 
\[
0
\equiv
f_1'f_1'\,M^8
-
3\,\Lambda^3\,\Lambda_{1,1}^7
+
5\,\Lambda_1^5\,\Lambda_1^5
\]
pour \'eliminer toutes les puissances de $\Lambda_1^5$ qui sont
sup\'erieures ou \'egales \`a 2 dans chaque polyn\^ome $\mathcal{ P
}_a^{2 \times {\rm inv }}$ et
r\'eorganisons le tout en puissances de $f_1'$. Nous obtenons ainsi
une nouvelle repr\'esentation\,:
\[
{\sf P}^{2\times {\rm inv}}\big(j^4f\big)
=
\sum_{-\frac{3}{4}m\leqslant a\leqslant m}\,
(f_1')^a
\Big[
\mathcal{ Q}_a^{2\times{\rm inv}}
\big(\Lambda^3,\Lambda_{1,1}^7,M^8)
+
\Lambda_1^5\,\mathcal{ R}_a^{2\times{\rm inv}}
\big(\Lambda^3,\Lambda_{1,1}^7,M^8)
\Big],
\]
avec certains polyn\^omes $\mathcal{ Q}_a^{2\times{\rm inv}}$ 
de poids $m-a$ et $\mathcal{ R}_a^{2\times{\rm inv}}$ de
poids $m - 5 - a$.

\smallskip

Maintenant, c'est un fait remarquable qu'une telle repr\'esentation
(dans laquelle on a tenu compte de l'id\'eal des relations entre $f_1'$,
$\Lambda^3$, $\Lambda_1^5$, $\Lambda_{ 1, 1}^7$ et $M^8$\big) doit
n\'ecessairement ne faire appara\^{\i}tre que des puissances positives de
$f_1'$, et donc \^etre automatiquement polynomiale\,: nous l'affirmons.

\smallskip

En effet, si tel n'\'etait pas le cas, en prenant pour $a$ l'exposant le
plus n\'egatif tel que $\mathcal{ P}_a^{2\times{\rm inv}} 
+ \Lambda_1^5 \, \mathcal{ R}_a^{2\times{\rm inv}} \not \equiv
0$ et en chassant le d\'enominateur $\frac{ 1}{ (f_1')^{ - a}}$, nous
obtiendrions une \'equation de la forme
\[
\mathcal{ Q}_a^{2\times{\rm inv}}
\big(\Lambda^3,\Lambda_{1,1}^7,M^8\big)
+
\Lambda_1^5\,
\mathcal{ R}_a^{2\times{\rm inv}}
\big(\Lambda^3,\Lambda_{1,1}^7,M^8\big)
=
{\rm O}(f_1'),
\]
qui s'annulerait lorsque $f_1'$ est \'egal\'e \`a z\'ero,
circonstance qui est exclue par le lemme suivant.

\def\thelemma{\!}\begin{lemma}
\'Etant donn\'e deux polyn\^omes quelconques $\mathcal{ Q}$ et $
\mathcal{ R}$ de trois variables complexes, l'identit\'e{\rm \,:}
\[
0
\equiv
\mathcal{Q}\big(\Lambda^3,\Lambda_{1,1}^7,M^8\big)
+
\Lambda_1^5\,\mathcal{R}\big(\Lambda^3,\Lambda_{1,1}^7,M^8\big)
\Big\vert_{f_1'=0}
\]
est identiquement satisfaite dans $\C \big[ f_2', f_1'', f_2'',
f_1''', f_1''', f_1'''', f_2''''\big]$, lorsque $f_1 '$ est
\'ega\-l\'e \`a z\'ero, si et seulement si $\mathcal{Q }$ et
$\mathcal{ R}$ sont identiquement nuls.
\end{lemma}

\noindent{\em Preuve.}
D\'eveloppons en effet tout d'abord 
cette identit\'e suivant les puissances de $M^8$\,:
\[
0
\equiv
\sum_{k\geqslant 0}\,
\big(M^8\big)^k
\Big[
\mathcal{ Q}_k\big(\Lambda^3,\Lambda_{1,1}^7\big)
+
\Lambda_1^5\,\mathcal{ R}_k\big(\Lambda^3,\Lambda_{1,1}^7\big)
\Big]
\Big\vert_{f_1'=0},
\]
la somme \'etant bien entendu finie.
Lorsque $f_1' = 0$, les expressions\,:
\[
\aligned
\Lambda^3\big\vert_{f_1'=0}
&
=
-f_1''f_2',
\\
\Lambda_1^5\big\vert_{f_1'=0}
&
=
3(f_1''f_2')f_1'',
\\
\Lambda_{1,1}^7\big\vert_{f_1'=0}
&
=
-15\,(f_1''f_2')f_1''f_1'',
\\
M^8\big\vert_{f_1'=0}
&
=
3(f_1''''f_2')(f_1''f_2')
-
12(f_1'''f_2''-f_1''f_2''')(f_1''f_2')
-
5(f_1'''f_2')(f_1'''f_2').
\endaligned
\]
montrent que $M^8 \big\vert_{ f_1' = 0}$ est alg\'ebriquement
ind\'ependant des trois polyn\^omes $\Lambda^3 \big \vert_{ f_1' = 0}$,
$\Lambda_1^5 \big \vert_{ f_1' = 0}$, $\Lambda_{ 1,1}^7 \big \vert_{
f_1' = 0}$. Nous en d\'eduisons que
\[
0 
\equiv
\mathcal{ Q}_k\big(\Lambda^3,\,\Lambda_{1,1}^7\big)
+
\Lambda_1^5\,\mathcal{ R}_k\big(\Lambda^3,\,\Lambda_{1,1}^7\big)
\Big\vert_{f_1'=0}
\]
pour tout $k$. L'\'enonc\'e suivant permet alors de conclure.
\hfill$\square$

\def\thelemma{\!}\begin{lemma}
L'identit\'e polynomiale{\rm \,:}
\[
0 
\equiv 
\mathcal{S}\big(\Lambda^3,\,\Lambda_{1,1}^7\big)
+
\Lambda_1^5\,\mathcal{T}\big(\Lambda^3,\,\Lambda_{1,1}^7\big)
\Big\vert_{f_1'=0}
\]
est satisfaite si et seulement si $\mathcal{ S} = \mathcal{ T} = 0$.
\end{lemma}

\noindent{\em Preuve.}
Pour simplifier, introduisons les deux variables alg\'ebriquement
ind\'epen\-dantes $x := -f_1'' f_2'$ et $y := f_1''$, de telle sorte que
\[
\Lambda^3\big\vert_{f_1'=0}
=
x,
\ \ \ \ \ \ \ \ \
\Lambda_1^5\big\vert_{f_1'=0}
=
3\,yx,
\ \ \ \ \ \ \ \ \
\Lambda_{1,1}^7\big\vert_{f_1'=0}
=
15\,y^2x.
\]
En d\'eveloppant $\mathcal{ S}$ et $\mathcal{ T}$ en s\'erie de mon\^omes et en
regroupant les termes selon les puissances de $y$, on obtient une
identit\'e\,:
\[
0
\equiv
{\textstyle{\sum_l}}\,y^{2l}\,\Big(
{\textstyle{\sum_k}}\,15^l\,{\sf s}_{kl}\,x^{k+l}
\Big)
+
3\,
{\textstyle{\sum_l}}\,y^{2l+1}\,\Big(
{\textstyle{\sum_k}}\,15^l\,{\sf t}_{kl}\,x^{k+l+1}
\Big)
\]
qui se d\'eploie n\'ecessairement en deux collections d'identit\'es\,: 
\[
0
\equiv
{\textstyle{\sum_k}}\,15^l\,{\sf s}_{kl}\,x^{k+l}
\ \ \ \ \ \ \ \ \ \
\text{\rm et}
\ \ \ \ \ \ \ \ \ \
0
\equiv
{\textstyle{\sum_k}}\,15^l\,{\sf t}_{kl}\,x^{k+l+1}
\]
index\'ees par $l$, lesquelles impliquent enfin manifestement
l'annulation de tous les coefficients
${\sf s}_{ kl}$ et ${\sf t}_{
kl}$.
\hfill$\square$

Ainsi le lemme implique que $\mathcal{ Q}_a^{2 \times{\rm inv }} +
\Lambda_1^5 \, \mathcal{ R}_a^{2 \times{\rm inv }} \equiv 0$,
contradiction. Donc en conclusion, l'expression obtenue\,:
\[
{\sf P}^{2\times{\rm inv}}
\big(j^4f\big)
=
\sum_{0\leqslant a\leqslant m}\,
(f_1')^m\,
\Big[
\mathcal{Q}_a^{2\times{\rm inv}}
\big(\Lambda^3,\Lambda_{1,1}^7,M^8\big)
+
\Lambda_1^5\,
\mathcal{R}_a^{2\times{\rm inv}}
\big(\Lambda^3,\Lambda_{1,1}^7,M^8\big)
\Big]
\]
ne fait intervenir que des puissances positives de $f_1'$\,: c'est donc
un vrai polyn\^ome, et tout polyn\^ome de cette sorte est manifestement
invariant par reparam\'etrisation et par rapport \`a l'action de ${\sf U}_2
(\C)$. Le th\'eor\`eme est d\'emontr\'e.
\hfill$\square$

\smallskip\noindent{\bf Remarque sur le degr\'e de transcendance.}
Observons au passage que les quatre polyn\^omes fondamentaux $f_1'$,
$\Lambda^3$, $\Lambda_{ 1,1}^7$ et $M^8$, dont les puissances
apparaissent de mani\`ere quelconque dans ${\sf P}^{2 \times{ \rm inv}}
\big( j^4f \big)$, sont en fait alg\'ebriquement ind\'epen\-dants
(heureusement!). En effet, en partant des expressions compl\`etes\,:
\[
\aligned
f_1'
&
=
f_1'
\\
\Lambda^3
&
=
\Delta^{1,2}
\\
\Lambda_{1,1}^7
&
=
\Delta^{1,4}\,f_1'f_1'
+
4\,\Delta^{2,3}\,f_1'f_1'
-
10\,\Delta^{1,3}\,f_1'f_1''
+
15\,\Delta^{1,2}\,f_1''f_1''
\\
M^8
&
=
3\,\Delta^{1,4}\,\Delta^{1,2}
+
12\,\Delta^{2,3}\,\Delta^{1,2}
-
5\,\Delta^{1,3}\,\Delta^{1,3},
\endaligned
\]
si nous
introduisons la combinaison alg\'ebrique\,:
\[
\aligned
\widetilde{M}^8
:=
&\
M^8
-
3\,
\frac{\Lambda_{1,1}^7\,\Lambda^3}{f_1'f_1'}
\\
=
&\
-5\,\Delta^{1,3}\,\Delta^{1,3}
+
30\,\Delta^{1,3}\,\Delta^{1,2}\,
\frac{f_1''}{f_1'}
-
45\,\Delta^{1,2}\,\Delta^{1,2}\,
\frac{f_1''}{f_1'}\,\frac{f_1''}{f_1'},
\endaligned
\]
nous voyons imm\'ediatement que $f_1'$, $\Lambda^3$, $\Lambda_{ 1,1}^7$
et $\widetilde{ M}^8$ sont alg\'ebriquement ind\'epen\-dants, puisque leur
expression, de type triangulaire dans les variables de jets, fait
successivement appara\^{\i}tre $f_1'$, $f_1''$, $f_1'''$ et $f_1''''$.

\def\thelemma{\!}\begin{lemma}
Au dessus de $\C \big[ f_1', f_2', f_1'', f_2'', f_1''', f_2''',
f_1'''', f_2'''' \big]$, le degr\'e de transcendance du corps engendr\'e
par les cinq bi-invariants $f_1'$, $\Lambda^3$, $\Lambda_1^5$,
$\Lambda_{ 1, 1}^7$ et $M^8$ est \'egal \`a 4, tandis que celui du corps
engendr\'e par les neuf invariants $f_1'$, $f_2'$, $\Lambda^3$,
$\Lambda_1^5$, $\Lambda_2^5$, $\Lambda_{ 1, 1}^7$, $\Lambda_{ 1,
2}^7$, $\Lambda_{ 2, 2}^7$ et $M^8$ est \'egal \`a 5.
\end{lemma}

\section*{\S7.~Jets d'ordre 5 en dimension 2}

\noindent{\bf Id\'eal des relations.} 
Nous pouvons donc maintenant poursuivre notre \'etude des polyn\^omes
invariants par reparam\'etrisation au niveau des jets d'ordre $\kappa =
5$. \`A cet \'etage, parmi les ving-cinq invariants que nous avons
calcul\'es et normalis\'es dans la Section~4, onze d'entre eux sont
bi-invariants de mani\`ere \'evidente, \`a savoir ceux qui ne comportent que
des ``$(\cdot)_1$'' en indice inf\'erieur\,:
\[
f_1',\ \ \ \ \
\Lambda^3,\ \ \ \ \
\Lambda_1^5,\ \ \ \ \
\Lambda_{1,1}^7,\ \ \ \ \
M^8,\ \ \ \ \
\Lambda_{1,1,1}^9,\ \ \ \ \
M_1^{10},\ \ \ \ \
N^{12},\ \ \ \ \ 
K_{1,1}^{12},\ \ \ \ \
H_1^{14},\ \ \ \ \ 
F_{1,1}^{16}.
\]
Sachant que nous avons d\'ej\`a syst\'ematiquement tenu compte de l'identit\'e
de Jacobi toutes les fois qu'elle nous permettait de r\'eduire le nombre
d'invariants ind\'ependants qui doivent \^etre envisag\'es, l'id\'eal des
relations qui existe entre nos bi-invariants est alors maintenant
construit en \'ecrivant m\'ethodiquement les $C_5^3 = 10$ relations
\thetag{ $\mathcal{ P}lck_1$} que l'on peut former en s\'electionnant
trois colonnes abitraires de la matrice\,:
\[
\left\vert\!\left\vert
\begin{array}{ccccc}
f_1' \ \ & \ \ 3\,\Lambda^3 \ \ & 
\ \ 5\,\Lambda_1^5 \ \ & \ \ 7\,\Lambda_{1,1}^5 \ \ &
\ \ 8\,M^8
\\
{\sf D}f_1' \ \ & \ \
{\sf D}\Lambda^3 \ \ & \ \ {\sf D}\Lambda_1^5 \ \ & \ \
{\sf D}\Lambda_{1,1}^7 \ \ & \ \ {\sf D}M^8
\end{array}
\right\vert\!\right\vert,
\]
et aussi les $C_5^4 = 5$ relations \thetag{ $\mathcal{ P}lck_2$}
du deuxi\`eme type associ\'ees \`a chaque choix de
quatre colonnes de cette m\^eme matrice, ce qui nous donne\,:
\[
\aligned
0
&
\overset{8}{\equiv}
f_1'\,\big[\Lambda^3,\,\Lambda_1^5\big]
+
5\,\Lambda_1^5\,\big[f_1',\,\Lambda^3\big]
+
3\,\Lambda^3\,\big[\Lambda_1^5,\,f_1'\big],
\\
0
&
\overset{10}{\equiv}
f_1'\,\big[\Lambda^3,\,\Lambda_{1,1}^7\big]
+
7\,\Lambda_{1,1}^7\,\big[f_1',\,\Lambda^3\big]
+
3\,\Lambda^3\,\big[\Lambda_{1,1}^7,\,f_1'\big],
\\
0
&
\overset{13}{\equiv}
f_1'\,\big[\Lambda^3,\,M^8\big]
+
8\,M^8\,\big[f_1',\,\Lambda^3\big]
+
3\,\Lambda^3\,\big[M^8,\,f_1'\big],
\\
0
&
\overset{15}{\equiv}
f_1'\,\big[\Lambda_1^5,\,\Lambda_{1,1}^7\big]
+
7\,\Lambda_{1,1}^7\,\big[f_1',\,\Lambda_1^5\big]
+
5\,\Lambda_1^5\,\big[\Lambda_{1,1}^7,\,f_1'\big],
\\
0
&
\overset{18}{\equiv}
f_1'\,\big[\Lambda_1^5,\,M^8\big]
+
8\,M^8\,\big[f_1',\,\Lambda_1^5\big]
+
5\,\Lambda_1^5\,\big[M^8,\,f_1'\big],
\\
0
&
\overset{25}{\equiv}
f_1'\,\big[\Lambda_{1,1}^7,\,M^8\big]
+
8\,M^8\,\big[f_1',\,\Lambda_{1,1}^7\big]
+
7\,\Lambda_{1,1}^7\,\big[M^8,\,f_1'\big],
\endaligned
\]
\[
\aligned
0
&
\overset{51}{\equiv}
3\,\Lambda^3\,\big[\Lambda_1^5,\,\Lambda_{1,1}^7\big]
+
7\,\Lambda_{1,1}^7\,\big[\Lambda^3,\,\Lambda_1^5\big]
+
5\,\Lambda_1^5\,\big[\Lambda_{1,1}^7,\,\Lambda^3\big],
\\
0
&
\overset{54}{\equiv}
3\,\Lambda^3\,\big[\Lambda_1^5,\,M^8\big]
+
8\,M^8\,\big[\Lambda^3,\,\Lambda_1^5\big]
+
5\,\Lambda_1^5\,\big[M^8,\,\Lambda^3\big],
\\
0
&
\overset{61}{\equiv}
3\,\Lambda^3\,\big[\Lambda_{1,1}^7,\,M^8\big]
+
8\,M^8\,\big[\Lambda^3,\,\Lambda_{1,1}^7\big]
+
7\,\Lambda_{1,1}^7\,\big[M^8,\,\Lambda^3\big],
\\
0
&
\overset{71}{\equiv}
5\,\Lambda_1^5\,\big[\Lambda_{1,1}^7,\,M^8\big]
+
8\,M^8\,\big[\Lambda_1^5,\,\Lambda_{1,1}^7\big]
+
7\,\Lambda_{1,1}^7\,\big[M^8,\,\Lambda_1^5\big],
\endaligned
\]
\[
\aligned
0
&
\overset{23'}{\equiv}
\big[f_1',\,\Lambda^3\big]\cdot
\big[\Lambda_1^5,\,\Lambda_{1,1}^7\big]
+
\big[\Lambda_{1,1}^7,\,f_1'\big]\cdot
\big[\Lambda_1^5,\,\Lambda^3\big]
+
\big[\Lambda^3,\,\Lambda_{1,1}^7\big]\cdot
\big[\Lambda_1^5,\,f_1'\big],
\\
0
&
\overset{26'}{\equiv}
\big[f_1',\,\Lambda^3\big]\cdot
\big[\Lambda_1^5,\,M^8\big]
+
\big[M^8,\,f_1'\big]\cdot
\big[\Lambda_1^5,\,\Lambda^3\big]
+
\big[\Lambda^3,\,M^8\big]\cdot
\big[\Lambda_1^5,\,f_1'\big],
\\
0
&
\overset{33'}{\equiv}
\big[f_1',\,\Lambda^3\big]\cdot
\big[\Lambda_{1,1}^7,\,M^8\big]
+
\big[M^8,\,f_1'\big]\cdot
\big[\Lambda_{1,1}^7,\,\Lambda^3\big]
+
\big[\Lambda^3,\,M^8\big]\cdot
\big[\Lambda_{1,1}^7,\,f_1'\big],
\\
0
&
\overset{43'}{\equiv}
\big[f_1',\,\Lambda_1^5\big]\cdot
\big[\Lambda_{1,1}^7,\,M^8\big]
+
\big[M^8,\,f_1'\big]\cdot
\big[\Lambda_{1,1}^7,\,\Lambda_1^5\big]
+
\big[\Lambda_1^5,\,M^8\big]\cdot
\big[\Lambda_{1,1}^7,\,f_1'\big],
\\
0
&
\overset{98'}{\equiv}
\big[\Lambda^3,\,\Lambda_1^5\big]\cdot
\big[\Lambda_{1,1}^7,\,M^8\big]
+
\big[M^8,\,\Lambda^3\big]\cdot
\big[\Lambda_{1,1}^7,\,\Lambda_1^5\big]
+
\big[\Lambda_1^5,\,M^8\big]\cdot
\big[\Lambda_{1,1}^7,\,\Lambda^3\big].
\endaligned
\]

\noindent{\bf D\'enombrement des syzygies
compl\`etes.} L'id\'eal complet des relations \thetag{ $\mathcal{ P}lck_1$}
et \thetag{ $\mathcal{ P}lck_2$} existant entre les ving-cinq
invariants comporte\,:
\[
C_9^3
+
C_9^4
=
{\bf 84}
+
{\bf 126}
=
{\bf 210},
\]
relations que nous avons patiemment d\'evelopp\'ees sur treize pages
manuscrites, mais que nous renon\c cons \`a recopier dans ce fichier
\LaTeX, pour la simple raison qu'il suffit, comme nous l'avons
argument\'e, d'\'etudier seulement les bi-invariants. Sur les 15 signes
``$\equiv$'' donnant les syzygies qui existent entre les
bi-invariants, nous conservons, pour m\'emoire, la num\'erotation de nos
$84 + 126$ \'equations manuscrites. 

\smallskip\noindent{\bf \'Enonc\'e.} Voici maintenant l'\'enonc\'e que
nous devrions attendre comme constituant notre deuxi\`eme r\'esultat principal.

\def\thetheorem{\!}\begin{theorem}
Pour les jets d'ordre 5 en dimension 2, tout bi-invariant de poids $m$, 
${\sf
P}^{2\times {\rm inv}} \big( j^5 f_1, \, j^5 f_2 \big)$,
s'exprime polynomialement en fonction de onze polyn\^omes
fondamentaux{\rm \,:}
\[
\boxed{
\aligned
f_1'\ \ \ \ \
\Lambda^3\ \ \ \ \
\Lambda_1^5\ \ \ \ \
\Lambda_{1,1}^7\ \ \ \ \
M^8
\\
\Lambda_{1,1,1}^9\ \ \ 
M_1^{10}\ \ \ 
N^{12}\ \ \ 
K_{1,1}^{12}
\\
H_1^{14}\ \ \ \
F_{1,1}^{16}
\endaligned
}\,
\]
{\it qui sont donn\'es explicitement en fonction de $j^5 f$ par les
formules calcul\'ees \`a la Section~4, et dont l'id\'eal des relations est
constitu\'e des quinze \'equations de degr\'e $\leqslant 3$ suivantes{\rm
\,:}}
\[
\aligned
0
&
\overset{8}{\equiv}
-f_1'f_1'\,M^8
-
5\,\Lambda_1^5\,\Lambda_1^5
+
3\,\Lambda^3\,\Lambda_{1,1}^7,
\\
0
&
\overset{10}{\equiv}
-f_1'f_1'\,M_1^{10}
-
7\,\Lambda_1^5\,\Lambda_{1,1}^7
+
3\,\Lambda^3\,\Lambda_{1,1,1}^9,
\\
0
&
\overset{13}{\equiv}
-f_1'\,N^{12}
-
8\,\Lambda_1^5\,M^8
+
3\,\Lambda^3\,M_1^{10},
\\
0
&
\overset{15}{\equiv}
-f_1'f_1'\,K_{1,1}^{12}
-
7\,\Lambda_{1,1}^7\,\Lambda_{1,1}^7
+
5\,\Lambda_1^5\,\Lambda_{1,1,1}^9,
\\
0
&
\overset{18}{\equiv}
-f_1'\,H_1^{14}
-
8\,\Lambda_{1,1}^7\,M^8
+
5\,\Lambda_1^5\,M_1^{10},
\\
0
&
\overset{25}{\equiv}
-f_1'\,F_{1,1}^{16}
-
8\,M^8\,\Lambda_{1,1,1}^9
+
7\,\Lambda_{1,1}^7\,M_1^{10},
\endaligned
\]
\[
\aligned
0
&
\overset{51}{\equiv}
-3\,\Lambda^3\,K_{1,1}^{12}
-
7\,\Lambda_{1,1}^7\,M^8
+
5\,\Lambda_1^5\,M_1^{10},
\\
0
&
\overset{54}{\equiv}
-3\,\Lambda^3\,H_1^{14}
-
8\,f_1'\,M^8\,M^8
+
5\,\Lambda_1^5\,N^{12},
\\
0
&
\overset{61}{\equiv}
-3\,\Lambda^3\,F_{1,1}^{16}
-
8\,f_1'\,M^8\,M_1^{10}
+
7\,\Lambda_{1,1}^7\,N^{12},
\\
0
&
\overset{71}{\equiv}
-5\,\Lambda_1^5\,F_{1,1}^{16}
-
8\,f_1'\,M^8\,K_{1,1}^{12}
+
7\,\Lambda_{1,1}^7\,H_1^{14},
\endaligned
\]
\[
\aligned
0
&
\overset{23'}{\equiv}
\Lambda_1^5\,K_{1,1}^{12}
+
M^8\,\Lambda_{1,1,1}^9
-
\Lambda_{1,1}^7\,M_1^{10},
\\
0
&
\overset{26'}{\equiv}
\Lambda_1^5\,H_1^{14}
+
f_1'\,M^8\,M_1^{10}
-
\Lambda_{1,1}^7\,N^{12},
\\
0
&
\overset{33'}{\equiv}
\Lambda_1^5\,F_{1,1}^{16}
+
f_1'\,M_1^{10}\,M_1^{10}
-
\Lambda_{1,1,1}^9\,N^{12},
\\
0
&
\overset{43'}{\equiv}
\Lambda_{1,1}^7\,F_{1,1}^{16}
+
f_1'\,M_1^{10}\,K_{1,1}^{12}
-
\Lambda_{1,1,1}^9\,H_1^{14},
\\
0
&
\overset{98'}{\equiv}
M^8\,F_{1,1}^{16}
+
N^{12}\,K_{1,1}^{12}
-
M_1^{10}\,H_1^{14}.
\endaligned
\]
\end{theorem}

Par souci de ne pas alourdir exag\'er\'ement l'\'enonc\'e de ce
th\'eor\`eme nous repoussons \`a la Section~8 l'\'enonc\'e pr\'ecis
qui donne les sommes directes de repr\'esentations irr\'eductibles de
Schur permettant d'entreprendre un calcul de caract\'eristique
d'Euler.

\smallskip

L'action des matrices de la forme\,:
\[
\text{\sc v}
:=
\left(
\begin{array}{cc}
1 & v
\\
0 & 1
\end{array}
\right)
\]
faisant ``rena\^{\i}tre par polarisation'' les 14 invariants qui ne
sont pas bi-invariants, nous pouvons en d\'eduire une description
partielle, mais suffisante pour notre objectif, de $\mathcal{
DS}_2^5$.

\def\thecorollary{\!}\begin{corollary}
Tout polyn\^ome ${\sf P} \big( j^5 f_1, j^5 f_2 \big)$ invariant
par reparam\'etrisation s'exprime polynomialement en fonction de
vingt-cinq invariants fondamentaux{\rm \,:}
\[
\boxed{
\aligned
f_1'\ \ \ \ \
f_2'\ \ \ \ \
\Lambda^3\ \ \ \ \
\Lambda_1^5\ \ \ \ \
\Lambda_2^5\ \ \ \ \
\Lambda_{1,1}^7\ \ \ \ \
\Lambda_{1,2}^7\ \ \ \ \
\Lambda_{2,2}^7\ \ \ \ \
M^8\
\\
\Lambda_{1,1,1}^9\ \ \ \ \
\Lambda_{1,2,1}^9\ \ \ \ \
\Lambda_{2,1,2}^9\ \ \ \ \
\Lambda_{2,2,2}^9\ \ \ \ \
M_1^{10}\ \ \ \ \
M_2^{10}\
\\
N^{12}\ \ \ \ \
K_{1,1}^{12}\ \ \ \ \
K_{1,2}^{12}\ \ \ \ \
K_{2,1}^{12}\ \ \ \ \
K_{2,2}^{12}\
\\
H_1^{14}\ \ \ \ \
H_2^{14}\ \ \ \ \
F_{1,1}^{16}\ \ \ \ \
F_{1,2}^{16}\ \ \ \ \
F_{2,2}^{16}\
\endaligned
}\,,
\]
qui sont donn\'es par les formules explicites normalis\'ees{\rm \,:}
\[
\small
\aligned
f_i'
&
\\
\Lambda^3
&
:=
\Delta^{1,2}
\\
\Lambda_i^5
&
:=
\Delta^{1,3}\,f_i'
-
3\,\Delta^{1,2}\,f_i''
\\
\Lambda_{i,j}^7
&
:=
\Delta^{1,4}\,f_i'f_j'
+
4\,\Delta^{2,3}\,f_i'f_j'
-
5\,\Delta^{1,3}\big(f_i''f_j'+f_i'f_j'')
+
15\,\Delta^{1,2}\,f_i''f_j''
\\
M^8
&
:=
3\,\Delta^{1,4}\,\Delta^{1,2}
+
12\,\Delta^{2,3}\,\Delta^{1,2}
-
5\,\Delta^{1,3}\,\Delta^{1,3}
\endaligned
\]
\[
\small
\aligned
\Lambda_{i,j,k}^9
&
:=
\Delta^{1,5}\,f_i'f_j'f_k'
+
5\,\Delta^{2,4}\,f_i'f_j'f_k'
-
\\
&
\ \ \ \ \
-
4\,\Delta^{1,4}\big(f_i''f_j'+f_i'f_j'')\,f_k'
-
7\,\Delta^{1,4}\,f_i'f_j'f_k''
-
\\
&
\ \ \ \ \
-
16\,\Delta^{2,3}\big(f_i''f_j'+f_i'f_j''\big)f_k'
-
28\,\Delta^{2,3}\,f_i'f_j'f_k''
-
\\
&
\ \ \ \ \
-
5\,\Delta^{1,3}\big(f_i'''f_j'+f_i'f_j''')\,f_k'
+
35\,\Delta^{1,3}\big(f_i''f_j''f_k'+f_i''f_j'f_k''+f_i'f_j''f_k''\big)
-
\\
&
\ \ \ \ \
-
105\,\Delta^{1,2}\,f_i''f_j''f_k'',
\endaligned
\]
\[
\small
\aligned
M_i^{10}
&
:=
\big[
3\,\Delta^{1,5}\,\Delta^{1,2}
+
15\,\Delta^{2,4}\,\Delta^{1,2}
-
7\,\Delta^{1,4}\,\Delta^{1,3}
+
2\,\Delta^{2,3}\,\Delta^{1,3}
\big]\,f_i'
-
\\
&
\ \ \ \ \
-
\big[
24\,\Delta^{1,4}\,\Delta^{1,2}
+
96\,\Delta^{2,3}\,\Delta^{1,2}
-
40\,\Delta^{1,3}\,\Delta^{1,3}
\big]f_i'',
\endaligned
\]
\[
\small
\aligned
N^{12}
&
:=
9\,\Delta^{1,5}\,\Delta^{1,2}\,\Delta^{1,2}
+
45\,\Delta^{2,4}\,\Delta^{1,2}\,\Delta^{1,2}
-
45\,\Delta^{1,4}\,\Delta^{1,3}\,\Delta^{1,2}
-
\\
&
\ \ \ \
-
90\,\Delta^{2,3}\,\Delta^{1,3}\,\Delta^{1,2}
+
40\,\Delta^{1,3}\,\Delta^{1,3}\,\Delta^{1,3},
\endaligned
\]
\[
\footnotesize
\aligned
K_{i,j}^{12}
&
:=
f_i'f_j'
\Big(
5\,\Delta^{1,5}\,\Delta^{1,3}
+
25\,\Delta^{2,4}\,\Delta^{1,3}
-
7\,\Delta^{1,4}\,\Delta^{1,4}
-
56\,\Delta^{2,3}\,\Delta^{1,4}
-
112\,\Delta^{2,3}\,\Delta^{2,3}
\Big)
+
\\
&
\ \ \ \ \
+
\frac{(f_i'f_j''+f_i''f_j')}{2}
\Big(
-15\,\Delta^{1,5}\,\Delta^{1,2}
-
75\,\Delta^{2,4}\,\Delta^{1,2}
+
65\,\Delta^{1,4}\,\Delta^{1,3}
+
110\,\Delta^{2,3}\,\Delta^{1,3}
\Big)
+
\\
&
\ \ \ \ \
+
\frac{(f_i'f_j'''+f_i'''f_j')}{2}
\Big(
-50\,\Delta^{1,3}\,\Delta^{1,3}
\Big)
+
\\
&
\ \ \ \ \
+
f_i''f_j''
\Big(
-25\,\Delta^{1,3}\,\Delta^{1,3}
+
15\,\Delta^{1,4}\,\Delta^{1,2}
+
60\,\Delta^{2,3}\,\Delta^{1,2}
\Big),
\endaligned
\]
\[
\footnotesize
\aligned
H_i^{14}
&
:=
\Big(
15\,\Delta^{1,5}\,\Delta^{1,3}\,\Delta^{1,2}
+
75\,\Delta^{2,4}\,\Delta^{1,3}\,\Delta^{1,2}
+
5\,\Delta^{1,4}\,\Delta^{1,3}\,\Delta^{1,3}
+
\\
&
\ \ \ \ \
+
170\,\Delta^{2,3}\,\Delta^{1,3}\,\Delta^{1,3}
-
24\,\Delta^{1,4}\,\Delta^{1,4}\,\Delta^{1,2}
-
192\,\Delta^{1,4}\,\Delta^{2,3}\,\Delta^{1,2}
-
\\
&
\ \ \ \ \
-
384\,\Delta^{2,3}\,\Delta^{2,3}\,\Delta^{1,2}
\Big)\,f_i'
+
\Big(
-
45\,\Delta^{1,5}\,\Delta^{1,2}\,\Delta^{1,2}
-
225\,\Delta^{2,4}\,\Delta^{1,2}\,\Delta^{1,2}
+
\\
&
\ \ \ \ \
+
225\,\Delta^{1,4}\,\Delta^{1,3}\,\Delta^{1,2}
+
450\,\Delta^{2,3}\,\Delta^{1,3}\,\Delta^{1,2}
-
200\,\Delta^{1,3}\,\Delta^{1,3}\,\Delta^{1,3}
\Big)\,f_i'',
\endaligned
\]
\[
\footnotesize
\aligned
F_{i,j}^{16}
&
:=
\Big(
-3\,\Delta^{1,5}\,\Delta^{1,4}\,\Delta^{1,2}
-
15\,\Delta^{2,4}\,\Delta^{1,4}\,\Delta^{1,2}
-
12\,\Delta^{1,5}\,\Delta^{2,3}\,\Delta^{1,2}
+
\\
&
\ \ \ \ \ 
+
40\,\Delta^{1,5}\,\Delta^{1,3}\,\Delta^{1,3}
-
60\,\Delta^{2,4}\,\Delta^{2,3}\,\Delta^{1,2}
+
200\,\Delta^{2,4}\,\Delta^{1,3}\,\Delta^{1,3}
-
\\
&
\ \ \ \ \
-
49\,\Delta^{1,4}\,\Delta^{1,4}\,\Delta^{1,3}
-
422\,\Delta^{1,4}\,\Delta^{2,3}\,\Delta^{1,3}
-
904\,\Delta^{2,3}\,\Delta^{2,3}\,\Delta^{1,3}
\Big)f_i'f_j'
+
\\
&
\ \ \ \ \ 
+
\Big(
-105\,\Delta^{1,5}\,\Delta^{1,3}\,\Delta^{1,2}
-
525\,\Delta^{2,4}\,\Delta^{1,3}\,\Delta^{1,2}
+
205\,\Delta^{1,4}\,\Delta^{1,3}\,\Delta^{1,3}
-
\\
&
\ \ \ \ \
-
230\,\Delta^{2,3}\,\Delta^{1,3}\,\Delta^{1,3}
+
96\,\Delta^{1,4}\,\Delta^{1,4}\,\Delta^{1,2}
+
768\,\Delta^{1,4}\,\Delta^{2,3}\,\Delta^{1,2}
+
\\
&
\ \ \ \ \
+
1536\,\Delta^{2,3}\,\Delta^{2,3}\,\Delta^{1,2}
\Big)
\big(
f_i''f_j'+f_i'f_j''
\big)
+
\\
&
\ \ \ \ \
+
\Big(
-200\,\Delta^{1,3}\,\Delta^{1,3}\,\Delta^{1,3}
\Big)\big(
f_i'''f_j'+f_i'f_j'''
\big)
+
\\
&
\ \ \ \ \ 
+
\Big(
315\,\Delta^{1,5}\,\Delta^{1,2}\,\Delta^{1,2}
+
1575\,\Delta^{2,4}\,\Delta^{1,2}\,\Delta^{1,2}
-
1575\,\Delta^{1,4}\,\Delta^{1,3}\,\Delta^{1,2}
-
\\
&
\ \ \ \ \
-
3150\,\Delta^{2,3}\,\Delta^{1,3}\,\Delta^{1,2}
+
1400\,\Delta^{1,3}\,\Delta^{1,3}\,\Delta^{1,3}
\Big)f_i''f_j'',
\endaligned
\]
o\`u les indices $i$, $j$ et $k$ appartiennent \`a $\{ 1, 2 \}$. 
\end{corollary}

\smallskip\noindent{\bf Remarque.} 
Ce deuxi\`eme th\'eor\`eme ainsi que son corollaire doivent \^etre restreints
\`a la sous-alg\`ebre engendr\'ee par les crochets, laquelle s'organise de
mani\`ere coh\'erente (\cite{ de2007}) pour former un sous-fibr\'e du fibr\'e
des jets de Demailly-Semple au-dessus d'une surface projective
alg\'ebrique complexe $X^2 \subset P_3 ( \C)$. Nous montrerons en effet
\`a la fin de cette Section~7 que d\`es les jets d'ordre $\kappa = 5$, il
existe des invariants fondamentaux suppl\'ementaires qui ne sont pas
obtenus par crochets, et qui s'ajoutent aux nombreux invariants qui
apparaissent d\'ej\`a dans les deux \'enonc\'es pr\'ec\'edents\,; de plus, il
pourrait exister une infinit\'e d'invariants par reparam\'etrisation ainsi
que de bi-invariants qui sont fondamentaux, ce qui contredirait {\it a
fortiori}\, la pr\'esomption informelle d'apr\`es laquelle les invariants
form\'es par crochets engendrent $\mathcal{ DS}_2^5$. Ce ph\'enom\`ene est
d'autant plus troublant qu'au niveau $\kappa = 4$, comme nous l'avons
d\'emontr\'e, tous les invariants sont engendr\'es par crochets. Peut-\^etre
existe-il des liens math\'ematiques profonds entre ce ph\'enom\`ene
inattendu et le fait que $d = 5$ soit aussi le seuil critique optimal
attendu pour la Kobayashi-hyperbolicit\'e des surfaces g\'en\'eriques
$X^2 \subset P_3 ( \C)$. L'avenir le dira.

\smallskip\noindent{\bf Strat\'egie.} 
Insistons sur le fait que la strat\'egie de d\'emonstration que nous
allons entreprendre afin de tenter d'\'etablir, comme au niveau $\kappa
= 4$, que seuls les invariants form\'es par crochets existent, {\it
aurait n\'ecessairement dû aboutir si l'alg\`ebre des (bi)invariants
engendr\'es par crochet avait \underline{\rm a priori} co\"{\i}ncid\'e
avec l'alg\`ebre compl\`ete des invariants de Demailly-Semple}\,: ce fait
sera argument\'e apr\`es la fin de nos raisonnements. Ce n'est donc pas
cette strat\'egie qui est en cause, mais la r\'ealit\'e math\'ematique, et
m\^eme si cette derni\`ere contredit parfois nos attentes, il nous faut
bien admettre et reconna\^{\i}tre que c'est elle, et seulement elle qui
agit en ma\^{\i}tre, partout et toujours. Et puisque cette r\'ealit\'e pr\'ec\`ede
d'une certaine mani\`ere dans ses grandes lignes l'exploration et la
recherche, notamment lorsqu'il s'agit de structures alg\'ebriques, nous
n'aurions jamais pu aboutir \`a une telle conclusion n\'egative sans
entreprendre de consid\'erables efforts de calcul. C'est pourquoi nous
convions maintenant notre lecteur \`a d\'ecouvrir comment nous comptons
g\'en\'eraliser au niveau $\kappa = 5$ notre d\'emonstration qui \'etait
valable pour les jets d'ordre 4, avant de d\'evoiler les interstices
dans lesquelles s'insinuent de nombreux invariants fondamentaux
suppl\'ementaires (peut-\^etre une infinit\'e) qui ne sont pas engendr\'es par
crochets.

\smallskip\noindent{\bf D\'emonstration du second th\'eor\`eme.}
Partons de l'expression rationnelle
que nous avons obtenue pour tout polyn\^ome bi-invariant
\[
\sum_{-\frac{4}{5}m\leqslant a\leqslant m}\,
(f_1')^a\,\mathcal{P}_a
\big(
\Lambda^3,\Lambda_1^5,\Lambda_{1,1}^7,\Lambda_{1,1,1}^9
\big),
\]
expression dans laquelle entrent des puissances n\'egatives de $f_1'$.
Le raisonnement que nous avons \'elabor\'e pour les jets d'ordre 4 va se
g\'en\'eraliser ici, au prix d'une complication suppl\'ementaire mais
in\'evitable, parce que le recours aux bases de Gr\"obner est en g\'en\'eral
incontournable pour les id\'eaux de polyn\^omes \`a plusieurs variables qui
ne sont pas principaux.

\smallskip\noindent{\bf Ghost rationality.}
Observons que les six premi\`eres syzygies ``$\overset{ 8 }{ \equiv
}$'', ``$\overset{ 10 }{ \equiv }$'', ``$\overset{ 13 }{ \equiv }$'',
``$\overset{ 15 }{ \equiv }$'', ``$\overset{ 18 }{ \equiv }$'' et
``$\overset{ 25 }{ \equiv }$'' entre nos onze bi-invariants
fondamentaux font appara\^{\i}tre en premi\`ere place les six bi-invariants
$M^8$, $\Lambda_{ 1, 1, 1}^9$, $M_1^{ 10}$, $N^{ 12}$, $K_{ 1, 1}^{
12}$, $H_1^{ 14}$ et $F_{ 1, 1}^{ 16}$ que nous connaissons d\'ej\`a, mais
qui sont invisibles dans le d\'eveloppement en puissances positives et
n\'egatives de $f_1'$ que nous venons de rappeler \`a l'instant, ce
dernier ne constituant que la toute premi\`ere \'etape de la
d\'emonstration. Interpr\'etons donc ces six bi-invariants en les
qualifiant intuitivement de ``termes fant\^omes'' ``cach\'es''
derri\`ere $f_1'$ ou derri\`ere $f_1' f_1'$. Heuristiquement parlant, ce
pourrait tout \`a fait \underline{\^etre parce que}\footnote{\, (et m\^eme,
``ce {\it \underline{devrait} tr\`es vraisemblablement \underline{\^etre parce
que}}\ldots'', si nous avions d\'ej\`a achev\'e la d\'emonstration de notre
deuxi\`eme th\'eor\`eme!)}, dans toute expression purement polynomiale en
nos onze bi-invariants vers laquelle se dirige notre d\'emonstration\,:
\[
\mathcal{P}
\big(
f_1',\Lambda^3,\Lambda_1^5,\Lambda_{1,1}^7,M^8,
\Lambda_{1,1,1}^9,M_1^{10},N^{12},K_{1,1}^{12},H_1^{14},F_{1,1}^{16}
\big),
\] 
l'on peut remplacer ces six bi-invariants sp\'eciaux par leur expression
rationnelle en fonction seulement de $\Lambda^3$, de $\Lambda_1^5$, de
$\Lambda_{ 1, 1}^7$ et de $\Lambda_{ 1, 1, 1}^9$\,:
\[
\aligned
M^8
&
=
\frac{3\,\Lambda^3\,\Lambda_{1,1}^7
-
5\,\Lambda_1^5\,\Lambda_1^5}
{f_1'f_1'},
\\
M_1^{10}
&
=
\frac{3\,\Lambda^3\,\Lambda_{1,1,1}^9
-
7\,\Lambda_1^5\Lambda_{1,1}^7}
{f_1'f_1'},
\\
N^{12}
&
=
\frac{-45\,\Lambda^3\,\Lambda_1^5\,\Lambda_{1,1}^7
+
40\,\Lambda_1^5\,\Lambda_1^5\,\Lambda_1^5}
{f_1'f_1'f_1'},
\\
K_{1,1}^{12}
&
=
\frac{5\,\Lambda_1^5\,\Lambda_{1,1,1}^9
-
7\,\Lambda_{1,1}^7\,\Lambda_{1,1}^7}
{f_1'f_1'},
\\
H_1^{14}
&
=
\frac{-24\,\Lambda^3\,\Lambda_{1,1}^7\,\Lambda_{1,1}^7
+
5\,\Lambda_1^5\,\Lambda_1^5\,\Lambda_{1,1}^7
+
15\,\Lambda^3\,\Lambda_1^5\,\Lambda_{1,1,1}^9}
{f_1'f_1'f_1'},
\\
F_{1,1}^{16}
&
=
\frac{-3\,\Lambda^3\,\Lambda_{1,1}^7\,\Lambda_{1,1,1}^9
+
40\,\Lambda_1^5\,\Lambda_1^5\,\Lambda_{1,1,1}^9
-
49\,\Lambda_1^5\,\Lambda_{1,1}^7\,\Lambda_{1,1}^7}
{f_1'f_1'f_1'},
\endaligned
\]
ce qui impose manifestement des divisions par $f_1' f_1'$ ou par $f_1'
f_1' f_1'$, ce pourrait donc bien \^etre, disions-nous, pour cette seule
et simple raison que notre repr\'esentation initiale, (trop) facile \`a
obtenir, d'un bi-invariant sous la forme 
\[
\sum_{ - \frac{ 4}{ 5} \, m
\leqslant a \leqslant m}\, ( f_1' )^a \, \mathcal{ P}_a \big(
\Lambda^3, \Lambda_1^5, \Lambda_{ 1, 1}^7, \Lambda_{ 1, 1, 1}^9 \big)
\]
faisait in\'evitablement appara\^{\i}tre des puissances n\'egatives de
$f_1'$. Et pour \'eliminer ces d\'enominateurs, rien d'autre ne
s'offrirait \`a nous que d'{\it injecter les six bi-inva\-riants fant\^omes
dans l'expression rationnelle initiale}. Voil\`a\,: nous avons d\'evoil\'e une
nouvelle id\'ee essentielle qui s'av\`erera pertinente et efficiente pour
l'\'etude des invariants de Demailly-Semple \`a un ordre quelconque.

\smallskip\noindent{\bf \'Elimination des puissances n\'egatives de 
$f_1'$.}
En effet, rappelons-nous tout d'abord que dans le cas des jets d'ordre
4, apr\`es avoir inject\'e le seul bi-invariant ``fant\^ome'' existant, \`a
savoir $M^8$, nous sommes parvenus \`a \'eliminer les puissances n\'egatives
de $f_1'$ gr\^ace \`a une normalisation pr\'ealable de tout polyn\^ome
$\mathcal{ P} = \mathcal{ P} \big( \Lambda^3, \Lambda_1^5, \Lambda_{
1, 1}^7, M^8 \big)$ sous la forme $\mathcal{ Q} \big( \Lambda^3,
\Lambda_{ 1, 1}^7, M^8 \big) + \Lambda_1^5 \, \mathcal{ R} \big(
\Lambda^3, \Lambda_{ 1, 1}^7, M^8 \big)$, ce qui \'etait fort
\'el\'ementaire, sachant que l'id\'eal des relations est principal (la
th\'eorie des bases de Gr\"obner est vide dans ce cas), la fin de
l'argument reposant seulement sur le fait que lorsqu'on pose $f_1' =
0$, aucune relation polynomiale non triviale du type
\[
0
\equiv
\mathcal{Q}\big(\Lambda^3,\Lambda_1^5,M^8)
+
\Lambda_1^5\,\mathcal{R}\big(\Lambda^3,\Lambda_1^5,M^8\big)
\Big\vert_{f_1'=0}
\]
ne peut \^etre satisfaite. 

\smallskip\noindent{\bf Observation g\'en\'erale cruciale et poursuite de la d\'emonstration.}
{\it Aussi est-ce seulement l'id\'eal des relations entre les
bi-invariants restreints \`a l'hypersurface $\{ f_1 ' = 0 \}$ qui semble
compter}. L'enjeu, ici, apr\`es avoir inject\'e les six bi-invariants
fant\^omes $M^8$, $M_1^{ 10}$, $N^{ 12}$, $K_{ 1, 1}^{ 12}$, $H_1^{ 14}$
et $F_{ 1, 1}^{ 16}$ qui \'etaient cach\'es derri\`ere des puissances
positives de $f_1'$, ce qui nous donne ais\'ement une expression
g\'en\'erale du type\,:
\[
\sum_{-\frac{4}{5}m\leqslant a\leqslant m}\,
(f_1')^a\,\mathcal{P}_a
\big(
\Lambda^3,\Lambda_1^5,\Lambda_{1,1}^7,M^8,\Lambda_{1,1,1}^9,
M_1^{10},N^{12},K_{1,1}^{12},H_1^{14},F_{1,1}^{16}
\big),
\]
dans laquelle nous supposerons, puisque nous perdons tout contr\^ole
apr\`es injection des six bi-invariants suppl\'ementaires, que les nouveaux
polyn\^omes $\mathcal{ P}_a$ sont arbitraires de poids $m - a$, l'enjeu
alors semble \^etre de {\it parvenir \`a produire une {\sf \'ecriture
normalis\'ee en fonction des syzygies} pour repr\'esenter de mani\`ere
unique tout polyn\^ome $\mathcal{ P}$ de cette esp\`ece, de fa\c con \`a ce
que toute identit\'e du type}\,:
\[
0
\equiv
\text{\sf \'Ecriture unique}\
\Big\{
\mathcal{P}\big(
\Lambda^3,\Lambda_1^5,\dots,H_1^{14},F_{1,1}^{16}
\big)
\Big\}
\Big\vert_{f_1'=0},
\]
implique que le polyn\^ome $\mathcal{ P}$ est en fait identiquement nul.
Alors l'argument d'\'eli\-mination de la puissance maximalement n\'egative
de $f_1'$, le tout suivi de la restriction \`a $\{ f_1' = 0 \}$, cet
argument que nous avions utilis\'e avec succ\`es pour les jets d'ordre 4
fonctionnera \`a nouveau ici sans modification, et une r\'ecurrence
imm\'ediate montrera, comme pour les jets d'ordre 4, que les puissances
n\'egatives de $f_1'$ n'existent pas, ce que nous d\'esirions obtenir pour
achever la d\'emonstration du th\'eor\`eme.

\smallskip

Deux remarques pour mettre un terme \`a ces consid\'erations heuristiques
destin\'ees seulement \`a d\'evoiler nettement nos id\'ees en usant du
langage sp\'eculatif qui nous a servi de guide pour les \'elaborer.
Premi\`erement, il est clair que ce plan de d\'emonstration doit
fonctionner en toute g\'en\'eralit\'e pour des jets d'ordre quelconque
$\kappa \geqslant 4$, et nous montrerons en temps voulu qu'il
fonctionne aussi en dimension $\nu \geqslant 3$. Ensuite,
notons\,\,---\,\,puisqu'il est de l'essence des math\'ematiques d'\^etre
``truff\'ees d'obstacles''\,\,---\,\,que le saut en difficult\'e, lorsqu'on
passe des jets d'ordre 4 aux jets d'ordre 5, est presque trop
consid\'erable pour une intuition de g\'en\'eralit\'e habitu\'ee aux r\'ecurrences
r\'eguli\`eres et aux combinatoires qui d\'evoilent progressivement leurs
structures\,: on passe en effet brutalement de une syzygie \`a quinze, et
m\^eme de neuf \`a deux cents dix, pour ce qui concerne les invariants
complets\,; comment alors ne pas \'eprouver le sentiment que la complexit\'e
alg\'ebrique de ce probl\`eme explose {\it trop}\, rapidement?

\smallskip\noindent{\bf Restriction des syzygies.}
Poser $f_1' = 0$, comme nous devons maintenant le faire,
nous donne les 15 \'equations r\'eduites\,:
\[
\small
\aligned
0
&
\equiv
-5\,\Lambda_1^5\,\Lambda_1^5
+
3\,\Lambda^3\,\Lambda_{1,1}^7
\Big\vert_{f_1'=0},
\\
0
&
\equiv
-7\,\Lambda_1^5\,\Lambda_{1,1}^7
+
3\,\Lambda^3\,\Lambda_{1,1,1}^9
\Big\vert_{f_1'=0},
\\
0
&
\equiv
-8\,\Lambda_1^5\,M^8
+
3\,\Lambda^3\,M_1^{10}
\Big\vert_{f_1'=0},
\\
0
&
\equiv
-7\,\Lambda_{1,1}^7\,\Lambda_{1,1}^7
+
5\,\Lambda_1^5\,\Lambda_{1,1,1}^9
\Big\vert_{f_1'=0},
\\
0
&
\equiv
-8\,\Lambda_{1,1}^7\,M^8
+
5\,\Lambda_1^5\,M_1^{10}
\Big\vert_{f_1'=0},
\\
0
&
\equiv
-8\,M^8\,\Lambda_{1,1,1}^9
+
7\,\Lambda_{1,1}^7\,M_1^{10}\Big\vert_{f_1'=0},
\endaligned
\]
\[
\aligned
0
&
\equiv
-3\,\Lambda^3\,K_{1,1}^{12}
-
7\,\Lambda_{1,1}^7\,M^8
+
5\,\Lambda_1^5\,M_1^{10}
\Big\vert_{f_1'=0},
\\
0
&
\equiv
-3\,\Lambda^3\,H_1^{14}
+
5\,\Lambda_1^5\,N^{12}
\Big\vert_{f_1'=0},
\\
0
&
\equiv
-3\,\Lambda^3\,F_{1,1}^{16}
+
7\,\Lambda_{1,1}^7\,N^{12}
\Big\vert_{f_1'=0},
\\
0
&
\equiv
-5\,\Lambda_1^5\,F_{1,1}^{16}
+
7\,\Lambda_{1,1}^7\,H_1^{14}
\Big\vert_{f_1'=0},
\endaligned
\]
\[
\aligned
0
&
\equiv
\Lambda_1^5\,K_{1,1}^{12}
+
M^8\,\Lambda_{1,1,1}^9
-
\Lambda_{1,1}^7\,M_1^{10}
\Big\vert_{f_1'=0},
\\
0
&
\equiv
\Lambda_1^5\,H_1^{14}
-
\Lambda_{1,1}^7\,N^{12}
\Big\vert_{f_1'=0},
\\
0
&
\equiv
\Lambda_1^5\,F_{1,1}^{16}
-
\Lambda_{1,1,1}^9\,N^{12}
\Big\vert_{f_1'=0},
\\
0
&
\equiv
\Lambda_{1,1}^7\,F_{1,1}^{16}
-
\Lambda_{1,1,1}^9\,H_1^{14}
\Big\vert_{f_1'=0},
\\
0
&
\equiv
M^8\,F_{1,1}^{16}
+
N^{12}\,K^{12}
-
M_1^{10}\,H_1^{14}
\Big\vert_{f_1'=0}.
\endaligned
\]

\noindent{\bf Base de Gr\"obner.}
En choisissant l'ordre purement lexicographique (\cite{ clo2007}) sur
les mon\^omes de $\C\big[ f_1', \Lambda^3, \dots, H_1^{ 14}, F_{ 1, 1}^{
16} \big]$ qui est d\'eduit de l'ordre suivant sur les mon\^o\-mes
\'el\'ementaires restreints\,:
\[
\Lambda^3
>
\Lambda_1^5
>
\Lambda_{1,1}^7
>
M^8
>
\Lambda_{1,1,1}^9
>
M_1^{10}
>
N^{12}
>
K_{1,1}^{12}
>
H_1^{14}
>
F_{1,1}^{16},
\]
(nous sous-entendons ici la mention ``$(\cdot )\vert_{ f_1 ' = 0}$''),
Maple nous donne la base de Gr\"ob\-ner r\'eduite suivante pour l'id\'eal
complet des syzygies entre nos dix invariants restreints \`a $\{ f_1 ' =
0 \}$, laquelle est constitu\'ee de 21 \'equations\,:
\[
\small
\aligned
0
&
\overset{1}{\equiv}
-7\,H_1^{14}\,H_1^{14}
+
5\,\underline{N^{12}\,F_{1,1}^{16}}
\Big\vert_{f_1'=0},
\ \ \ \ \ \ \ \ \ \ \ \ \ \ \ \
0
\overset{11}{\equiv}
-
\Lambda_{1,1,1}^9\,N^{12}
+
\underline{\Lambda_1^5\,F_{1,1}^{16}}
\Big\vert_{f_1'=0},
\\
0
&
\overset{2}{\equiv}
-56\,K_{1,1}^{12}\,H_1^{14}
+
5\,\underline{M_1^{10}\,F_{1,1}^{16}}
\Big\vert_{f_1'=0},
\ \ \ \ \ \ \ \ \ \ \ \ \
0
\overset{12}{\equiv}
-\Lambda_{1,1}^7\,N^{12}
+
\underline{\Lambda_1^5\,H_1^{14}}
\Big\vert_{f_1'=0},
\\
0
&
\overset{3}{\equiv}
-8\,N^{12}\,K_{1,1}^{12}
+
\underline{M_1^{10}\,H_1^{14}}
\Big\vert_{f_1'=0},
\ \ \ \ \ \ \ \ \ \ \ \ \ \ \ \ \
0
\overset{13}{\equiv}
-M^8\,\Lambda_{1,1,1}^9
+
7\,\underline{\Lambda_1^5\,K_{1,1}^{12}}
\Big\vert_{f_1'=0},
\\
0
&
\overset{4}{\equiv}
-7\,N^{12}\,K_{1,1}^{12}
+
\underline{M^8\,F_{1,1}^{16}}
\Big\vert_{f_1'=0},
\ \ \ \ \ \ \ \ \ \ \ \ \ \ \ \ \ \ \
0
\overset{14}{\equiv}
-8\,\Lambda_{1,1}^7\,M^8
+
5\,\underline{\Lambda_1^5\,M_1^{10}}
\Big\vert_{f_1'=0},
\\
0
&
\overset{5}{\equiv}
-5\,M_1^{10}\,N^{12}
+
8\,\underline{M^8\,H_1^{14}}
\Big\vert_{f_1'=0},
\ \ \ \ \ \ \ \ \ \ \ \ \ \ \ \
0
\overset{15}{\equiv}
-7\,\Lambda_{1,1}^7\,\Lambda_{1,1}^7
+
5\,\underline{\Lambda_1^5\,\Lambda_{1,1,1}^9}
\Big\vert_{f_1'=0},
\\
0
&
\overset{6}{\equiv}
-5\,M_1^{10}\,M_1^{10}
+
64\,\underline{M^8\,K_{1,1}^{12}}
\Big\vert_{f_1'=0},
\ \ \ \ \ \ \ \ \ \ \ \ \
0
\overset{16}{\equiv}
-7\,\Lambda_{1,1}^7\,N^{12}
+
3\,\underline{\Lambda^3\,F_{1,1}^{16}}
\Big\vert_{f_1'=0},
\\
0
&
\overset{7}{\equiv}
-\Lambda_{1,1,1}^9\,H_1^{14}
+
\underline{\Lambda_{1,1}^7\,F_{1,1}^{16}}
\Big\vert_{f_1'=0},
\ \ \ \ \ \ \ \ \ \ \ \ \ \ \ \ \ \ \
0
\overset{17}{\equiv}
-5\,\Lambda_1^5\,N^{12}
+
3\,\underline{\Lambda^3\,H_1^{14}}
\Big\vert_{f_1'=0},
\\
0
&
\overset{8}{\equiv}
-5\,\Lambda_{1,1,1}^9\,N^{12}
+
7\,\underline{\Lambda_{1,1}^7\,H_1^{14}}
\Big\vert_{f_1'=0},
\ \ \ \ \ \ \ \ \ \ \ \ \ \ \ \
0
\overset{18}{\equiv}
-\Lambda_{1,1}^7\,M^8
+
3\,\underline{\Lambda^3\,K_{1,1}^{12}}
\Big\vert_{f_1'=0},
\\
0
&
\overset{9}{\equiv}
-5\,\Lambda_{1,1,1}^9\,M_1^{10}
+
56\,\underline{\Lambda_{1,1}^7\,K_{1,1}^{12}}
\Big\vert_{f_1'=0},
\ \ \ \ \ \ \ \ \ \
0
\overset{19}{\equiv}
-8\,\Lambda_1^5\,M^8
+
3\,\underline{\Lambda^3\,M_1^{10}}
\Big\vert_{f_1'=0},
\\
0
&
\overset{10}{\equiv}
-8\,M^8\,\Lambda_{1,1,1}^9
+
7\,\underline{\Lambda_{1,1}^7\,M_1^{10}}
\Big\vert_{f_1'=0},
\ \ \ \ \ \ \ \ \ \ \ \ \
0
\overset{20}{\equiv}
-7\,\Lambda_1^5\,\Lambda_{1,1}^7
+
3\,\underline{\Lambda^3\,\Lambda_{1,1,1}^9}
\Big\vert_{f_1'=0},
\\
&
\ \ \ \ \ \ \ \ \ \ \ \ \ \ \ \ \ \ \ \ \ \ \ \ \ \ \ \ \ \ \
\ \ \ \ \ \ \ \ \ \ \ \ \ \ \ \ \ \ \ \ \ \ \ \ \ \ \ \ \ \ \
\ \ \ \ \ \ \ \ \
0
\overset{21}{\equiv}
-5\,\Lambda_1^5\,\Lambda_1^5
+
3\,\underline{\Lambda^3\,\Lambda_{1,1}^7}
\Big\vert_{f_1'=0},
\endaligned
\]
toutes d\'eduites de nos 15 syzygies restreintes \`a $\{ f_1' = 0 \}$, et
dont l'ensemble recelle une combinatoire d'une simplicit\'e inattendue
qui va se d\'evoiler \`a nous dans un instant. Nous avons soulign\'e les
mon\^omes de t\^ete pour en extraire l'id\'eal monomial associ\'e ({\it
voir}\, {\it infra}). 

Bien que cette base de Gr\"obner nous ait \'et\'e procur\'ee par Maple, il n'est
pas n\'ecessaire que nous nous en remettions au calcul formel
\'electronique pour assurer la rigueur du r\'esultat, puisqu'il est ici
tr\`es ais\'e de v\'erifier\,:

\begin{itemize}

\smallskip\item[$\square$]
que ces 21 \'equations sont effectivement cons\'equence de nos
quinze syzygies r\'eduites\,; 

\smallskip\item[$\square$]
que ces 21 \'equations forment effectivement une
base de Gr\"obner.

\end{itemize}\smallskip

Le premier point se v\'erifie sans difficult\'e\,; le paragraphe ci-dessous
o\`u nous donnons l'expression des bi-invariants restreints \`a $\{ f_1' =
0 \}$ permet d'ailleurs de proc\'eder tr\`es rapidement. Pour ce qui
est du deuxi\`eme point, il nous suffit d'appli\-quer l'un 
des nombreux crit\`eres
caract\'erisant les bases de Gr\"obner (\cite{ clo2007}), 
d'apr\`es lequel chaque S-polyn\^ome
entre deux \'equations quelconques 
doit appartenir \`a l'i\-d\'eal engendr\'e par les
21 polyn\^omes, et \`a cette fin, la t\^ache de calcul manuel est
miraculeusement facilit\'ee par le fait que chacune de ces 21 \'equations
ne comporte que deux termes, avec \`a chaque fois un signe ``$-$'' et un
signe ``$+$'', ces deux termes \'etant chacun monomiaux et qui plus est,
de degr\'e deux, ce qui fait que {\it chaque S-polyn\^ome entre deux
\'equations ne poss\`ede encore que deux termes monomiaux}. Par exemple,
si on \'elimine les mon\^omes de t\^ete entre ``$\overset{ 20}{ \equiv}$''
et ``$\overset{ 21}{ \equiv}$'' en multipliant par des mon\^omes
appropri\'es et en soustrayant\,:
\[
\aligned
0
&
\equiv
\big(
-7\,\Lambda_1^5\,\Lambda_{1,1}^7
+
3\,\underline{\Lambda^3\,\Lambda_{1,1,1}^9}
\big)\cdot
\underline{\Lambda_{1,1}^7}
-
\big(
-5\,\Lambda_1^5\,\Lambda_1^5
+
3\,\underline{\Lambda^3\,\Lambda_{1,1}^7}
\big)\cdot
\underline{\Lambda_{1,1,1}^9}
\Big\vert_{f_1'=0},
\\
&
\equiv
-7\,\Lambda_1^5\,\Lambda_{1,1}^7\,\Lambda_{1,1}^7
+
5\,\Lambda_1^5\,\Lambda_1^5\,\Lambda_{1,1,1}^9,
\endaligned
\]
on constate que le S-polyn\^ome obtenu appartient bien \`a notre id\'eal,
puisqu'il co\"{\i}nci\-de avec l'\'equation ``$\overset{ 15}{ \equiv}$''
multipli\'ee par $\Lambda_1^5$. Les 209 autres S-polyn\^omes restants
se traitent de la m\^eme mani\`ere, \`a la main, en moins
de deux heures, apr\`es \'epuration pr\'ealable des notations.

\smallskip\noindent{\bf Expression des bi-invariants restreints \`a $\{ f_1 ' = 0 \}$.}
Mais avant de poursuivre, il est instructif d'\'ecrire, d'examiner et de
commenter la liste de nos onze bi-invariants restreints\,:
\[
\aligned
f_1'
\big\vert_0
&
=
0,
\\
\Lambda^3
\big\vert_0
&
=
-f_1''f_2'
=:
\Delta_0^{1,2},
\\
\Lambda_1^5
\big\vert_0
&
=
-3\,\Delta_0^{1,2}\,f_1'',
\\
\Lambda_{1,1}^7
\big\vert_0
&
=
15\,\Delta_0^{1,2}\,f_1''f_1'',
\\
M^8
\big\vert_0
&
=
3\,\Delta_0^{1,4}\,\Delta_0^{1,2}
+
12\,\Delta_0^{2,3}\,\Delta_0^{1,2}
-
5\,\Delta_0^{1,3}\,\Delta^{1,3},
\\
\Lambda_{1,1,1}^9
\big\vert_0
&
=
-105\,\Delta_0^{1,2}\,f_1''f_1''f_1'',
\\
M_1^{10}
\big\vert_0
&
=
-\Big[
24\,\Delta_0^{1,4}\,\Delta_0^{1,2}
+
96\,\Delta_0^{2,3}\,\Delta_0^{1,2}
-
40\,\Delta_0^{1,3}\,\Delta_0^{1,3}
\Big]\,f_1'',
\\
N^{12}
\big\vert_0
&
=
9\,\Delta_0^{1,5}\,\Delta_0^{1,2}\,\Delta_0^{1,2}
+
45\,\Delta_0^{2,4}\,\Delta_0^{1,2}\,\Delta_0^{1,2}
-
45\,\Delta_0^{1,4}\,\Delta_0^{1,3}\,\Delta_0^{1,2}
-
\\
&
\ \ \ \ \ \ \ \ \ \ \ \ \ \ \ \ \ \ \ \ \ \ \ \ \ \ \ \ \ \ \ \
-
90\,\Delta_0^{2,3}\,\Delta_0^{1,3}\,\Delta_0^{1,2}
+
40\,\Delta_0^{1,3}\,\Delta_0^{1,3}\,\Delta_0^{1,3},
\\
K_{1,1}^{12}
\big\vert_0
&
=
\Big[
15\,\Delta_0^{1,4}\,\Delta_0^{1,2}
+
60\,\Delta_0^{2,3}\,\Delta_0^{1,2}
-
25\,\Delta_0^{1,3}\,\Delta_0^{1,3}
\Big]\,
f_1''f_1'',
\\
H_1^{14}
\big\vert_0
&
=
\Big[
-45\,\Delta_0^{1,5}\,\Delta_0^{1,2}\,\Delta_0^{1,2}
-
225\,\Delta_0^{2,4}\,\Delta_0^{1,2}\,\Delta_0^{1,2}
+
225\,\Delta_0^{1,4}\,\Delta_0^{1,3}\,\Delta_0^{1,2}
+
\\
&
\ \ \ \ \ \ \ \ \ \ \ \ \ \ \ \ \ \ \ \ \ \ \ \ \ \
\ \ \ \ \ \ \ \
+
450\,\Delta_0^{2,3}\,\Delta_0^{1,3}\,\Delta_0^{1,2}
-
200\,\Delta_0^{1,3}\,\Delta_0^{1,3}\,\Delta_0^{1,3}
\Big]\,f_1'',
\\
F_{1,1}^{16}
\big\vert_0
&
=
\Big[
315\,\Delta_0^{1,5}\,\Delta_0^{1,2}\,\Delta_0^{1,2}
+
1575\,\Delta_0^{2,4}\,\Delta_0^{1,2}\,\Delta_0^{1,2}
-
1575\,\Delta_0^{1,4}\,\Delta_0^{1,3}\,\Delta_0^{1,2}
-
\\
&
\ \ \ \ \ \ \ \ \ \ \ \ \ \ \ \ \ \ \ \ \ \ \ \ \ \ \ 
-
3150\,\Delta_0^{2,3}\,\Delta_0^{1,3}\,\Delta_0^{1,2}
+
1400\,\Delta_0^{1,3}\,\Delta_0^{1,3}\,\Delta_0^{1,3}
\Big]\,
f_1''f_1''.
\endaligned
\]

\noindent{\bf Divisions wronskiennes.}
En comparant la deuxi\`eme et la troisi\`eme \'equation, nous obtenons par
exemple $f_1'' = - \frac{ 1}{ 3}\, \frac{ \Lambda_1^5 \vert_0 }{
\Lambda^3 \vert_0}$, et aussi $f_1 '' f_1'' = \frac{ 1}{ 15}\, \frac{
\Lambda_{ 1, 1}^7 \vert_0 }{ \Lambda^3 \vert_0 }$ si l'on compare la
deuxi\`eme et la quatri\`eme ligne, et en poursuivant ces observations,
nous pouvons \'ecrire\,:
\[
\aligned
&
\underline{\Lambda^3}\big\vert_0,
\\
&
\underline{\Lambda_1^5}\big\vert_0,
\\
&
\Lambda_{1,1}^7\big\vert_0
=
{\textstyle{\frac{5}{3}\,\frac{\Lambda_1^5\vert_0\,\Lambda_1^5\vert_0}{
\Lambda^3\vert_0}}},
\\
&
\underline{M^8}\big\vert_0,
\\
&
\Lambda_{1,1,1}^9\big\vert_0
=
{\textstyle{
\frac{7}{3}\,\frac{\Lambda_1^5\vert_0\,\Lambda_{1,1}^7\vert_0}
{\Lambda^3\vert_0}
=
\frac{35}{9}\,
\frac{\Lambda_1^5\vert_0\,\Lambda_1^5\vert_0\,\Lambda_1^5\vert_0}
{\Lambda^3\vert_0\,\Lambda^3\vert_0}
}},
\\
&
M_1^{10}\big\vert_0
=
{\textstyle{
\frac{8}{3}\,\frac{M^8\vert_0\,\Lambda_1^5\vert_0}
{\Lambda^3\vert_0}
}},
\\
&
\underline{N^{12}}\big\vert_0,
\\
&
K_{1,1}^{12}\big\vert_0
=
{\textstyle{
\frac{1}{3}\,\frac{M^8\vert_0\,\Lambda_{1,1}^7\vert_0}
{\Lambda^3\vert_0}
=
\frac{5}{9}\,\frac{M^8\vert_0\,\Lambda_1^5\vert_0\,\Lambda_1^5\vert_0}
{\Lambda^3\vert_0\,\Lambda^3\vert_0}
}},
\\
&
H_1^{14}\big\vert_0
=
{\textstyle{
\frac{5}{3}\,\frac{N^{12}\vert_0\,\Lambda_1^5\vert_0}
{\Lambda^3\vert_0}
}},
\\
&
F_{1,1}^{16}\big\vert_0
=
{\textstyle{
\frac{7}{3}\,\frac{\Lambda_{1,1}^7\vert_0\,N^{12}\vert_0}
{\Lambda^3\vert_0}
=
\frac{35}{9}\,\frac{\Lambda_1^5\vert_0\,\Lambda_1^5\vert_0\,N^{12}\vert_0}
{\Lambda^3\vert_0\,\Lambda^3\vert_0}
}},
\endaligned
\]
o\`u nous soulignons quatre bi-invariants restreints qui apparaissent
fondamentaux, \`a savoir $\underline{ \Lambda^3} \big\vert_0$,
$\underline{ \Lambda_1^5} \big\vert_0$, $\underline{ M^8 }\big
\vert_0$ et $\underline{ N^{ 12}} \big\vert_0$, puisque les six
autres s'expriment en fonction d'eux, apr\`es restriction \`a $\{ f_1' = 0
\}$, lorsqu'on autorise \`a diviser par le wronskien. Un examen
imm\'ediat de l'expression compl\`ete de ces quatre bi-invariants
restreints fondamentaux $\underline{ \Lambda^3} \big\vert_0$,
$\underline{ \Lambda_1^5 } \big\vert_0$, $\underline{ M^8 }\big
\vert_0$ et $\underline{ N^{ 12 }} \big\vert_0$ montre qu'ils sont
alg\'ebriquement ind\'ependants, puisqu'ils incorporent successivement
$f_1''$, $f_1'''$, $f_1''''$ et $f_1'''''$. Cette ind\'ependance
mutuelle resservira ult\'erieurement.

\smallskip\noindent{\bf Triangle harmonieux des mon\^omes de t\^ete.}
Comme nous le constatons en examinant notre base de Gr\"obner, les 21
bin\^omes de t\^ete s'organisent, lorsqu'on 
les range par ordre (lexicographique) croissant, en un triangle
remarquable\,:
\[
\begin{array}{cccccc}
\Lambda^3\,\Lambda_{1,1}^7 &
>\Lambda^3\,\Lambda_{1,1,1}^9 &
>\Lambda^3\,M_1^{10} &
>\Lambda^3\,K_{1,1}^{12} &
>\Lambda^3\,H_1^{14} &
>\Lambda^3\,F_{1,1}^{16}
\\
&
>\Lambda_1^5\,\Lambda_{1,1,1}^9 &
>\Lambda_1^5\,M_1^{10} &
>\Lambda_1^5\,K_{1,1}^{12} &
>\Lambda_1^5\,H_1^{14} &
>\Lambda_1^5\,F_{1,1}^{16}
\\
& &
>\Lambda_{1,1}^7\,M_1^{10} &
>\Lambda_{1,1}^7\,K_{1,1}^{12} &
>\Lambda_{1,1}^7\,H_1^{14} &
>\Lambda_{1,1}^7\,F_{1,1}^{16} 
\\
& & & 
>M^8\,K_{1,1}^{12} &
>M^8\,H_1^{14} &
>M^8\,F_{1,1}^{16} 
\\
& & & &
>M_1^{10}\,H_1^{14} &
>M_1^{10}\,F_{1,1}^{16}
\\
& & & & &
>N^{12}\,F_{1,1}^{16}
\end{array}
\]
(on sous-entend la mention ``$\vert_0$'' dans ce diagramme) dans
lequel nous reconnaissons, \`a la place des colonnes, les six
bi-invariants restreints $\Lambda_{ 1, 1}^7 \big \vert_0$, $\Lambda_{
1, 1, 1}^9 \big \vert_0$, $M_1^{ 10} \big \vert_0$, $K_{ 1, 1}^{ 12}
\big \vert_0$, $H_1^{ 14} \big \vert_0$ et $F_{ 1, 1}^{ 16} \big
\vert_0$, qui s'expriment rationnellement en fonction des quatre
bi-invariants restreints fondamentaux $\underline{ \Lambda^3}
\big\vert_0$, $\underline{ \Lambda_1^5 } \big\vert_0$, $\underline{
M^8 }\big \vert_0$ et $\underline{ N^{ 12 }} \big\vert_0$.

\smallskip\noindent{\bf Normalisation modulo les syzygies restreintes.}
Nous pouvons maintenant \'enoncer et d\'emontrer le lemme sur lequel
repose la fin de la d\'emonstration de notre second th\'eor\`eme.

\def\thelemma{\!\!}\begin{lemma}
Tout polyn\^ome arbitraire en les dix bi-invariants restreints{\rm \,:}
\[
\mathcal{P}
\big(
\Lambda^3\big\vert_0,\,
\Lambda_1^5\big\vert_0,\,
\Lambda_{1,1}^7\big\vert_0,\,
M^8\big\vert_0,\,
\Lambda_{1,1,1}^9\big\vert_0,\,
M_1^{10}\big\vert_0,\,
N^{12}\big\vert_0,\,
K_{1,1}^{12}\big\vert_0,\,
H_1^{14}\big\vert_0,\,
F_{1,1}^{16}\big\vert_0
\big)
\]
s'\'ecrit de mani\`ere unique, en tenant compte des 21 syzygies
gr\"obn\'eris\'ees ci-dessus, sous la forme unique{\rm \,:}
\[
\footnotesize
\aligned
&
\mathcal{P}_0\big(
\Lambda^3\big\vert_0,\,
\Lambda_1^5\big\vert_0,\,
M^8\big\vert_0,\,
N^{12}\big\vert_0
\big)
+
\Lambda_{1,1}^7\big\vert_0\,
\mathcal{Q}_0\big(
\Lambda_1^5\big\vert_0,\,
\Lambda_{1,1}^7\big\vert_0,\,
M^8\big\vert_0,\,
N^{12}\big\vert_0
\big)
+
\\
&
+
\Lambda_{1,1,1}^9\big\vert_0\,
\mathcal{R}_0\big(
\Lambda_{1,1}^7\big\vert_0,\,
M^8\big\vert_0,\,
\Lambda_{1,1,1}^9\big\vert_0,\,
N^{12}\big\vert_0
\big)
+
M_1^{10}\big\vert_0\,
\mathcal{S}_0\big(
M^8\big\vert_0,\,
\Lambda_{1,1,1}^9\big\vert_0,\,
M_1^{10}\big\vert_0,\,
N^{12}\big\vert_0
\big)
+
\\
&
+
K_{1,1}^{12}\big\vert_0\,
\mathcal{T}_0\big(
\Lambda_{1,1,1}^9\big\vert_0,\,
M_1^{10}\big\vert_0,\,
N^{12}\big\vert_0,\,
K_{1,1}^{12}\big\vert_0
\big)
+
H_1^{14}\big\vert_0\,
\mathcal{U}_0\big(
\Lambda_{1,1,1}^9\big\vert_0,\,
N^{12}\big\vert_0,\,
K_{1,1}^{12}\big\vert_0,\,
H_1^{14}\big\vert_0
\big)
+
\\
&
+
F_{1,1}^{16}\big\vert_0\,
\mathcal{V}_0\big(
\Lambda_{1,1,1}^9\big\vert_0,\,
K_{1,1}^{12}\big\vert_0,\,
H_1^{14}\big\vert_0,\,
F_{1,1}^{16}\big\vert_0
\big),
\endaligned
\]
o\`u $\mathcal{ P}_0$, $\mathcal{ Q}_0$, $\mathcal{ R}_0$, $\mathcal{
S}_0$, $\mathcal{ T}_0$, $\mathcal{ U}_0$ et $\mathcal{ V}_0$ sont des
polyn\^omes absolument arbitraires en leurs quatre arguments, et pour
pr\'eciser, toute relation du type{\rm \,:}
\[
0
\equiv
\mathcal{P}_0
+
\Lambda_{1,1}^7\big\vert_0\,\mathcal{Q}_0
+
\Lambda_{1,1,1}^9\big\vert_0\,\mathcal{R}_0
+
M_1^{10}\big\vert_0\,\mathcal{S}_0
+
K_{1,1}^{12}\big\vert_0\,\mathcal{T}_0
+
H_1^{14}\big\vert_0\,\mathcal{U}_0
+
F_{1,1}^{16}\big\vert_0\,\mathcal{V}_0
\]
qui est identiquement satisfaite dans $\C \big[ f_2', f_1'', f_2'',
f_1''', f_2''', f_1'''', f_2'''', f_1''''', f_2''''' \big]$ lors\-qu'on
remplace les dix bi-invariants restreints par leur expression en
fonction de $j^5 f \big\vert_0$, implique n\'ecessairement que les sept
polyn\^omes $\mathcal{ P}_0$, $\mathcal{ Q}_0$, $\mathcal{ R}_0$,
$\mathcal{ S}_0$, $\mathcal{ T}_0$, $\mathcal{ U}_0$ et $\mathcal{
V}_0$ s'annulent tous identiquement.
\end{lemma}

\smallskip\noindent{\bf D\'emonstration du lemme principal.}
D'apr\`es la th\'eorie \'el\'ementaire des bases de Gr\"obner, une base de
l'espace vectoriel quotient
\[
\C\Big[
\Lambda^3\big\vert_0,\,
\Lambda_1^5\big\vert_0,\,
\dots,\,
H_1^{14}\big\vert_0,\,
F_{1,1}^{16}\big\vert_0
\Big]
\Big/
\big(
\text{\rm 21 syzygies restreintes}
\big)
\]
est constitut\'ee de tous les mon\^omes
\[
\footnotesize
\aligned
\Big(
\Lambda^3\big\vert_0
\Big)^a
\Big(
\Lambda_1^5\big\vert_0
\Big)^b
\Big(
\Lambda_{1,1}^7\big\vert_0
\Big)^c
\Big(
M^8\big\vert_0
\Big)^d
&
\Big(
\Lambda_{1,1,1}^9\big\vert_0
\Big)^e
\Big(
M_1^{10}\big\vert_0
\Big)^f
\Big(
N^{12}\big\vert_0
\Big)^g
\\
&
\ \ 
\Big(
K_{1,1}^{12}\big\vert_0
\Big)^h
\Big(
H_1^{14}\big\vert_0
\Big)^i
\Big(
F_{1,1}^{16}\big\vert_0
\Big)^j
\endaligned
\]
qui {\it n'appartiennent pas}\, \`a l'id\'eal monomial engendr\'e
par les 21 mon\^omes de t\^ete que nous avons dispos\'es en triangle, o\`u les
exposants $a, b, c, d, e, f, g, h, i, j \in \N$ sont des entiers
positifs ou nuls. Or un tel mon\^ome appartient \`a cet id\'eal monomial si
et seulement si il est divisible par l'un des 21 mon\^omes de t\^ete, ce
qui revient \`a dire que le d\'eca-indice
\[
\big(
a,\,
b,\,
c,\,
d,\,
e,\,
f,\,
g,\,
h,\,
i,\,
j
\big)
\in\N^{10}
\]
appartient \`a la r\'eunion des 21 sous-ensembles suivants de $\N^{ 10}$\,:

\begin{tiny}
\[
\aligned
\{a\geqslant 1\}\cap \{c\geqslant 1\}
\bigcup
\{a\geqslant 1\}\cap \{e\geqslant 1\}
\bigcup
\{a\geqslant 1\}\cap \{f\geqslant 1\}
\bigcup
\{a\geqslant 1\}\cap \{h\geqslant 1\}
\bigcup
\{a\geqslant 1\}\cap \{i\geqslant 1\}
\bigcup
\{a\geqslant 1\}\cap \{j\geqslant 1\}
\\
\bigcup
\{b\geqslant 1\}\cap \{e\geqslant 1\}
\bigcup
\{b\geqslant 1\}\cap \{f\geqslant 1\}
\bigcup
\{b\geqslant 1\}\cap \{h\geqslant 1\}
\bigcup
\{b\geqslant 1\}\cap \{i\geqslant 1\}
\bigcup
\{b\geqslant 1\}\cap \{j\geqslant 1\}
\\
\bigcup
\{c\geqslant 1\}\cap \{f\geqslant 1\}
\bigcup
\{c\geqslant 1\}\cap \{h\geqslant 1\}
\bigcup
\{c\geqslant 1\}\cap \{i\geqslant 1\}
\bigcup
\{c\geqslant 1\}\cap \{j\geqslant 1\}
\\
\bigcup
\{d\geqslant 1\}\cap \{h\geqslant 1\}
\bigcup
\{d\geqslant 1\}\cap \{i\geqslant 1\}
\bigcup
\{d\geqslant 1\}\cap \{j\geqslant 1\}
\\
\bigcup
\{f\geqslant 1\}\cap \{i\geqslant 1\}
\bigcup
\{f\geqslant 1\}\cap \{j\geqslant 1\}
\\
\bigcup
\{g\geqslant 1\}\cap \{j\geqslant 1\}
\endaligned
\]
\end{tiny}

\noindent
Le calcul du compl\'ementaire de cet ensemble est ais\'e, 
et il donne\,:
\[
\footnotesize
\aligned
\Big[
\{a=0\}\cup\{c=e=f=h=i=j=0\}
\Big]
\bigcap
\\
\Big[
\{b=0\}\cup\{e=f=h=i=j=0\}
\Big]
\bigcap
\\
\Big[
\{c=0\}\cup\{f=h=i=j=0\}
\Big]
\bigcap
\\
\Big[
\{d=0\}\cup\{h=i=j=0\}
\Big]
\bigcap
\\
\Big[
\{f=0\}\cup\{i=j=0\}
\Big]
\bigcap
\\
\Big[
\{g=0\}\cup\{j=0\}
\Big],
\ \ \ \
\endaligned
\]
ce qui se simplifie pour donner 7 composantes d\'efinies chacune
par six \'equations\,:
\[
\footnotesize
\aligned
\big\{
0=a=b=c=d=f=g
\big\}
\bigcup
\\
\big\{
0=a=b=c=d=f=j
\big\}
\bigcup
\\
\big\{
0=a=b=c=d=i=j
\big\}
\bigcup
\\
\big\{
0=a=b=c=h=i=j
\big\}
\bigcup
\\
\big\{
0=a=b=f=h=i=j
\big\}
\bigcup
\\
\big\{
0=a=e=f=h=i=j
\big\}
\bigcup
\\
\big\{
0=c=e=f=h=i=j
\big\}.
\ \ \ \
\endaligned
\]
(Incidemment, nous avons \'etabli que l'id\'eal des syzygies restreintes
est une intersection compl\`ete.) Par cons\'equent, l'ensemble de tous
les mon\^omes qu'il nous reste dans l'espace quotient est constitu\'e des
sept listes suivantes\,:
\[
\footnotesize
\aligned
&
\square\ \ \ \ \
\Big(
\Lambda_{1,1,1}^9\big\vert_0
\Big)^e
\Big(
K_{1,1}^{12}\big\vert_0
\Big)^h
\Big(
H_1^{14}\big\vert_0
\Big)^i
\Big(
F_{1,1}^{16}\big\vert_0
\Big)^j,
\\
&
\square\ \ \ \ \
\Big(
\Lambda_{1,1,1}^9\big\vert_0
\Big)^e
\Big(
N^{12}\big\vert_0
\Big)^g
\Big(
K_{1,1}^{12}\big\vert_0
\Big)^h
\Big(
H_1^{14}\big\vert_0
\Big)^i,
\\
&
\square\ \ \ \ \
\Big(
\Lambda_{1,1,1}^9\big\vert_0
\Big)^e
\Big(
M_1^{10}\big\vert_0
\Big)^f
\Big(
N^{12}\big\vert_0
\Big)^g
\Big(
K_{1,1}^{12}\big\vert_0
\Big)^h,
\\
&
\square\ \ \ \ \
\Big(
M^8\big\vert_0
\Big)^d
\Big(
\Lambda_{1,1,1}^9\big\vert_0
\Big)^e
\Big(
M_1^{10}\big\vert_0
\Big)^f
\Big(
N^{12}\big\vert_0
\Big)^g,
\\
&
\square\ \ \ \ \
\Big(
\Lambda_{1,1}^7\big\vert_0
\Big)^c
\Big(
M^8\big\vert_0
\Big)^d
\Big(
\Lambda_{1,1,1}^9\big\vert_0
\Big)^e
\Big(
N^{12}\big\vert_0
\Big)^g,
\\
&
\square\ \ \ \ \
\Big(
\Lambda_1^5\big\vert_0
\Big)^b
\Big(
\Lambda_{1,1}^7\big\vert_0
\Big)^c
\Big(
M^8\big\vert_0
\Big)^d
\Big(
N^{12}\big\vert_0
\Big)^g,
\\
&
\square\ \ \ \ \
\Big(
\Lambda^3\big\vert_0
\Big)^a
\Big(
\Lambda_1^5\big\vert_0
\Big)^b
\Big(
M^8\big\vert_0
\Big)^d
\Big(
N^{12}\big\vert_0
\Big)^g.
\endaligned
\]
Cependant, puisque ces sept listes se recouvrent
partiellement\,\,---\,\,par exemple\,: l'in\-tersection de la premi\`ere et
de la deuxi\`eme ligne est constitut\'ee des mon\^omes de la forme $\big(
\Lambda_{ 1,1,1 }^9 \big \vert_0 \big)^e \big( K_{1,1 }^{12 }\big
\vert_0 \big)^h \big( H_1^{14 }\big \vert_0 \big)^i$\,\,---, nous
devons encore les r\'eorganiser de telle sorte qu'il n'y ait plus aucune
intersection entre elles, et si nous d\'esignons ces listes par les
septs lettres ${\tt A}$, ${\tt B}$, ${\tt C}$, ${\tt D}$, ${\tt F}$ et
${\sf G}$, il nous suffit en fait tout simplement d'\'ecrire\,:
\[
\footnotesize
\aligned
{\tt A}\cup{\tt B}\cup{\tt C}\cup{\tt D}\cup{\tt E}\cup{\tt F}\cup{\tt G}
=
{\tt A}
\bigcup
{\tt B}\big\backslash{\tt A}
\bigcup
{\tt C}\big\backslash
\big({\tt A}\cup{\tt B}\big)
\bigcup
{\tt D}\big\backslash
\big({\tt A}\cup{\tt B}\cup{\tt C}\big)
\bigcup
{\tt E}\big\backslash
\big({\tt A}\cup{\tt B}\cup{\tt C}\cup{\tt D}\big)
\\
\bigcup
{\tt F}\big\backslash
\big({\tt A}\cup{\tt B}\cup{\tt C}\cup{\tt D}\cup{\tt E}\big)
\bigcup
{\tt G}\big\backslash
\big({\tt A}\cup{\tt B}\cup{\tt C}\cup{\tt D}\cup{\tt E}\cup{\tt G}\big),
\endaligned
\]
ce qui nous donne imm\'ediatement la repr\'esentation \'enonc\'ee dans notre
lemme, au moyen des sept polyn\^omes arbitraires $\mathcal{ P}_0$,
$\mathcal{ R}_0$, $\mathcal{ S}_0$, $\mathcal{ T}_0$, $\mathcal{ U}_0$
et $\mathcal{ V}_0$ dont les quatre arguments se diff\'erencient
successivement d'une unit\'e lorsqu'on saute une ligne (de la
ligne 1 \`a la ligne 7), fait
combinatoire aussi remarquable qu'impr\'evu et que nous aimerions voir
se confirmer, se g\'en\'eraliser et se stabiliser lorsque nous \'etudierons
les bi-invariants pour les jets d'ordre 6\,\,---\,\,projet ambitieux s'il en
est.

\smallskip
Afin d'\'etablir la deuxi\`eme assertion du lemme, supposons maintenant
qu'une relation du type\,:
\[
0
\equiv
\mathcal{P}_0
+
\Lambda_{1,1}^7\big\vert_0\,\mathcal{Q}_0
+
\Lambda_{1,1,1}^9\big\vert_0\,\mathcal{R}_0
+
M_1^{10}\big\vert_0\,\mathcal{S}_0
+
K_{1,1}^{12}\big\vert_0\,\mathcal{T}_0
+
H_1^{14}\big\vert_0\,\mathcal{U}_0
+
F_{1,1}^{16}\big\vert_0\,\mathcal{V}_0
\]
est identiquement satisfaite, et rempla\c cons-y alors les six
bi-invariants restreints non fondamentaux par leur expression en
fonction de ceux qui sont fondamentaux, ce qui nous donne\,:
\[
\footnotesize
\aligned
0
&
\equiv
\mathcal{P}_0
\big(
\Lambda^3\big\vert_0,\,
\Lambda_1^5\big\vert_0,\,
M^8\big\vert_0,\,
N^{12}\big\vert_0
\big)
+
{\textstyle{
\frac{\Lambda_1^5\vert_0\,\Lambda_1^5\vert_0}{\Lambda^3\vert_0}
}}\,
\mathcal{Q}_0
\Big(
\Lambda_1^5\big\vert_0,\,
{\textstyle{
\frac{\Lambda_1^5\vert_0\,\Lambda_1^5\vert_0}{\Lambda^3\vert_0}
}},\,
M^8\big\vert_0,\,
N^{12}\big\vert_0
\Big)
+
\\
&
\ \ \ \ \
+
{\textstyle{
\frac{\Lambda_1^5\vert_0\,\Lambda_1^5\vert_0\,\Lambda_1^5\vert_0}
{\Lambda^3\vert_0\,\Lambda^3\vert_0}
}}\,
\mathcal{R}_0
\Big(
{\textstyle{
\frac{\Lambda_1^5\vert_0\,\Lambda_1^5\vert_0}{\Lambda^3\vert_0}
}},\,
M^8\big\vert_0,\,
{\textstyle{
\frac{\Lambda_1^5\vert_0\,\Lambda_1^5\vert_0\,\Lambda_1^5\vert_0}
{\Lambda^3\vert_0\,\Lambda^3\vert_0}
}},\,
N^{12}\big\vert_0
\Big)
+
\\
&
\ \ \ \ \
+
{\textstyle{
\frac{\Lambda_1^5\vert_0\,M^8\vert_0}{\Lambda^3\vert_0}
}}\,
\mathcal{S}_0
\Big(
M^8\big\vert_0,\,
{\textstyle{
\frac{\Lambda_1^5\vert_0\,\Lambda_1^5\vert_0\,\Lambda_1^5\vert_0}
{\Lambda^3\vert_0\,\Lambda^3\vert_0}
}},\,
{\textstyle{
\frac{\Lambda_1^5\vert_0\,M^8\vert_0}{\Lambda^3\vert_0}
}},\,
N^{12}\big\vert_0
\Big)
+
\\
&
\ \ \ \ \
+
{\textstyle{
\frac{\Lambda_1^5\vert_0\,\Lambda_1^5\vert_0\,M^8\vert_0}
{\Lambda^3\vert_0\,\Lambda^3\vert_0}
}}\,
\mathcal{T}_0
\Big(
{\textstyle{
\frac{\Lambda_1^5\vert_0\,\Lambda_1^5\vert_0\,\Lambda_1^5\vert_0}
{\Lambda^3\vert_0\,\Lambda^3\vert_0}
}},\,
{\textstyle{
\frac{\Lambda_1^5\vert_0\,M^8\vert_0}{\Lambda^3\vert_0}
}},\,
N^{12}\big\vert_0,\,
{\textstyle{
\frac{\Lambda_1^5\vert_0\,\Lambda_1^5\vert_0\,M^8\vert_0}
{\Lambda^3\vert_0\,\Lambda^3\vert_0}
}}
\Big)
+
\\
&
\ \ \ \ \
+
{\textstyle{
\frac{\Lambda_1^5\vert_0\,N^{12}\vert_0}{\Lambda^3\vert_0}
}}\,
\mathcal{U}_0
\Big(
{\textstyle{
\frac{\Lambda_1^5\vert_0\,\Lambda_1^5\vert_0\,\Lambda_1^5\vert_0}
{\Lambda^3\vert_0\,\Lambda^3\vert_0}
}},\,
N^{12}\big\vert_0,\,
{\textstyle{
\frac{\Lambda_1^5\vert_0\,\Lambda_1^5\vert_0\,M^8\vert_0}
{\Lambda^3\vert_0\,\Lambda^3\vert_0}
}},\,
{\textstyle{
\frac{\Lambda_1^5\vert_0\,N^{12}\vert_0}{\Lambda^3\vert_0}
}}
\Big)
+
\\
&
\ \ \ \ \
{\textstyle{
\frac{\Lambda_1^5\vert_0\,\Lambda_1^5\vert_0\,N^{12}\vert_0}
{\Lambda^3\vert_0\,\Lambda^3\vert_0}
}}\,
\mathcal{V}_0
\Big(
{\textstyle{
\frac{\Lambda_1^5\vert_0\,\Lambda_1^5\vert_0\,\Lambda_1^5\vert_0}
{\Lambda^3\vert_0\,\Lambda^3\vert_0}
}},\,
{\textstyle{
\frac{\Lambda_1^5\vert_0\,\Lambda_1^5\vert_0\,M^8\vert_0}
{\Lambda^3\vert_0\,\Lambda^3\vert_0}
}},\,
{\textstyle{
\frac{\Lambda_1^5\vert_0\,N^{12}\vert_0}
{\Lambda^3\vert_0}
}},\,
{\textstyle{
\frac{\Lambda_1^5\vert_0\,\Lambda_1^5\vert_0\,N^{12}\vert_0}
{\Lambda^3\vert_0\,\Lambda^3\vert_0}
}}
\Big).
\endaligned
\]
Pour en d\'eduire l'annulation de ces 7 polyn\^omes $\mathcal{ P}_0$,
$\mathcal{ Q}_0$, $\mathcal{ R}_0$, $\mathcal{ S}_0$, $\mathcal{
T}_0$, $\mathcal{ U}_0$ et $\mathcal{ V}_0$, si nous commen\c cons par
multiplier cette identit\'e par le d\'enominateur $\big( \Lambda^3 \vert_0
\big)^\mu$ de son expression r\'eduite, nous obtenons une identit\'e de la
forme\,:
\[
0
\equiv
\big(\Lambda^3\big\vert\big)^\mu\,
\mathcal{P}_0
\big(
\Lambda^3\big\vert_0,\,
\Lambda_1^5\big\vert_0,\,
M^8\big\vert_0,\,
N^{12}\big\vert_0
\big)
+
\big(\Lambda^3\big\vert_0\big)^{\mu-1}\cdot\,
\text{\sf reste polynomial}, 
\]
laquelle montre tout d'abord imm\'ediatement que le premier polyn\^ome
$\mathcal{ P}_0$ s'annule identiquement, puisque $\Lambda^3
\big\vert_0$, $\Lambda_1^5 \big\vert_0$, $M^8 \big\vert_0$ et $N^{
12}\big\vert_0$ sont alg\'ebriquement ind\'ependants. Supprimons donc
$\mathcal{ P}_0$. Ensuite, en multipliant par une puissance
suffisamment \'elev\'ee $\big( \Lambda^3 \vert_0 \big)^\mu$ pour \'eliminer
les puissances n\'egatives de $\Lambda^3 \vert_0$ qui apparaissent dans
$\mathcal{ Q}_0$, $\mathcal{ R}_0$, $\mathcal{ S}_0$, $\mathcal{
T}_0$, $\mathcal{ U}_0$ et $\mathcal{ V}_0$, et en divisant le tout
par le facteur $\Lambda_1^5 \vert_0$ pr\'esent devant chacun des six
polyn\^omes restants, nous obtenons une identit\'e de la forme\,:
\[
\aligned
0
&
\equiv
\big(\Lambda^3\big\vert_0\big)^\mu
\Big[
\Lambda_1^5\big\vert_0\,
\mathcal{Q}_0
+
M^8\big\vert_0\,
\mathcal{S}_0
+
N^{12}\big\vert_0\,
\mathcal{U}_0
\Big]
+
\big(\Lambda^3\big\vert_0\big)^{\mu-1}\cdot\,
\text{\sf reste polynomial},
\endaligned
\]
d'o\`u nous d\'eduisons l'annulation identique du polyn\^ome\,:
\[
\aligned
0
&
\equiv
\big(\Lambda^3\big\vert_0\big)^\mu
\Big[
\Lambda_1^5\big\vert_0\,
\mathcal{Q}_0
\Big(
\Lambda_1^5\big\vert_0,\,
{\textstyle{
\frac{\Lambda_1^5\vert_0\,\Lambda_1^5\vert_0}{\Lambda^3\vert_0}
}},\,
M^8\big\vert_0,\,
N^{12}\big\vert_0
\Big)
+
\\
&
\ \ \ \ \ \ \ \ \ \ \ \ \ \ \ \ \ \ \ \
+
M^8\big\vert_0\,
\mathcal{S}_0
\Big(
M^8\big\vert_0,\,
{\textstyle{
\frac{\Lambda_1^5\vert_0\,\Lambda_1^5\vert_0\,\Lambda_1^5\vert_0}
{\Lambda^3\vert_0\,\Lambda^3\vert_0}
}},\,
{\textstyle{
\frac{\Lambda_1^5\vert_0\,M^8\vert_0}{\Lambda^3\vert_0}
}},\,
N^{12}\big\vert_0
\Big)
+
\\
&
\ \ \ \ \ \ \ \ \ \ \ \ \ \ \ \ \ \ \ \
+
N^{12}\big\vert_0\,
\mathcal{U}_0
\Big(
{\textstyle{
\frac{\Lambda_1^5\vert_0\,\Lambda_1^5\vert_0\,\Lambda_1^5\vert_0}
{\Lambda^3\vert_0\,\Lambda^3\vert_0}
}},\,
N^{12}\big\vert_0,\,
{\textstyle{
\frac{\Lambda_1^5\vert_0\,\Lambda_1^5\vert_0\,M^8\vert_0}
{\Lambda^3\vert_0\,\Lambda^3\vert_0}
}},\,
{\textstyle{
\frac{\Lambda_1^5\vert_0\,N^{12}\vert_0}{\Lambda^3\vert_0}
}}
\Big)
\Big].
\endaligned
\]
Ensuite, si nous introduisons les d\'eveloppement finis de ces trois
polyn\^omes\,:
\[
\small
\aligned
\mathcal{Q}_0
\big(t_1,t_2,t_3,t_4)
=
{\textstyle{\sum}}\,&
{\sf coeff}\cdot
t_1^\alpha\,t_2^\beta\,t_3^\gamma\,t_4^\delta,
\ \ \ \ \ \ \ \ \
\mathcal{S}_0
\big(t_1,t_2,t_3,t_4)
=
{\textstyle{\sum}}\,
{\sf coeff}\cdot
t_1^{\alpha'}\,t_2^{\beta'}\,t_3^{\gamma'}\,t_4^{\delta'},
\\
&
\mathcal{U}_0
\big(t_1,t_2,t_3,t_4)
=
{\textstyle{\sum}}\,
{\sf coeff}\cdot
t_1^{\alpha''}\,t_2^{\beta''}\,t_3^{\gamma''}\,t_4^{\delta''},
\endaligned
\]
nous en d\'eduirons l'annulation de $\mathcal{ Q}_0$, 
de $\mathcal{ S}_0$ et de $\mathcal{ U}_0$ gr\^ace \`a l'observation 
suivante.

\def\theassertion{\!}\begin{assertion} 
Les trois familles de mon\^omes\,:

\begin{footnotesize}
\begin{itemize}

\smallskip\item[{\bf (i)}]\ \ \ \ \
$\big(\Lambda^3\big\vert_0\big)^{\mu-\beta}
\big(\Lambda_1^5\big\vert_0\big)^{1+\alpha+2\beta}
\big(M^8\big\vert_0\big)^{\gamma}
\big(N^{12}\big\vert_0\big)^\delta$,

\smallskip\item[{\bf (ii)}]\ \ \ \ \
$\big(\Lambda^3\big\vert_0\big)^{\mu-2\beta'-\gamma'}
\big(\Lambda_1^5\big\vert_0\big)^{3\beta'+\gamma'}
\big(M^8\big\vert_0\big)^{1+\alpha'+\gamma'}
\big(N^{12}\big\vert_0\big)^{\delta'}$,

\smallskip\item[{\bf (iii)}]\ \ \ \ \
$\big(\Lambda^3\big\vert_0\big)^{\mu-2\alpha''-2\gamma''-\delta''}
\big(\Lambda_1^5\big\vert_0\big)^{3\alpha''+2\gamma''+\delta''}
\big(M^8\big\vert_0\big)^{\gamma''}
\big(N^{12}\big\vert_0\big)^{1+\beta''+\gamma''}$,

\end{itemize}\smallskip
\end{footnotesize}

\noindent
ne contiennent {\sf aucune redondance}, {\rm i.e.} chaque mon\^ome
correspondand \`a un choix de $(\alpha, \beta, \gamma, \delta)$, ou de
$(\alpha', \beta', \gamma', \delta')$, ou encore de $(\alpha'',
\beta'', \gamma'', \delta'')$ appara\^{\i}t {\sf une et une seule fois}.
\end{assertion}

\medskip

En effet, nous v\'erifions tout d'abord pour la premi\`ere famille, que
l'auto-intersec\-tion\,:
\[
\mu-\beta
=
\mu-\underline{\beta},
\ \ \ \ \ \ \
1+\alpha+2\beta
=
1+\underline{\alpha}+2\underline{\beta}
\ \ \ \ \ \ \
\gamma
=
\underline{\gamma}
\ \ \ \ \ \ \
\delta
=
\underline{\delta},
\]
est vide, c'est-\`a-dire implique $\alpha = \underline{ \alpha}$, 
$\beta = \underline{ \beta}$, $\gamma = \underline{ \gamma}$, 
$\delta = \underline{ \delta}$, puis de m\^eme pour la
deuxi\`eme famille\,:
\[
\mu-2\beta'-\gamma'
=
\mu-2\underline{\beta}'-\underline{\gamma}',
\ \ \ \ \ \ \
3\beta'+\gamma'
=
3\underline{\beta}'+\underline{\gamma}',
\ \ \ \ \ \ \
\alpha'+\gamma'
=
\underline{\alpha}'+\underline{\gamma}',
\ \ \ \ \ \ \
\delta'
=
\underline{\delta}',
\]
et aussi pour la troisi\`eme famille\,:
\[
\small
\aligned
\mu-2\alpha''-2\gamma''-\delta''
=
\mu-2\underline{\alpha}''-2\underline{\gamma}''-\underline{\delta}'',
\ \ \ \ \ \ \
&
3\alpha''+2\gamma''+\delta''
=
3\underline{\alpha}''+2\underline{\gamma}''+\underline{\delta}'',
\\
&
\gamma''
=
\underline{\gamma}'',
\ \ \ \ \ \ \ 
\beta''+\gamma''
=
\underline{\beta}''+\underline{\gamma}''.
\endaligned
\]
Ensuite, des formules pour l'intersection entre la premi\`ere et la
deuxi\`eme famille\,:
\[
-\beta
=
-2\beta'-\gamma',
\ \ \ \ \ \ \
1+\alpha+2\beta
=
3\beta'+\gamma',
\ \ \ \ \ \ \
\gamma
=
1+\alpha'+\gamma',
\ \ \ \ \ \ \
\delta
=
\delta',
\]
d\'ecoule l'\'equation $\alpha = - \gamma ' - \beta ' - 1$, impossible
parce que l'exposant entier $\alpha$ doit imp\'erativement 
\^etre $\geqslant 0$. De m\^eme, l'intersection entre 
la premi\`ere et la troisi\`eme famille\,:
\[
-\beta
=
-2\alpha''-2\gamma''-\delta'',
\ \ \ \ \ \ \
1+\alpha+2\beta
=
3\alpha''+2\gamma''+\delta'',
\ \ \ \ \ \ \
\gamma
=
\gamma'',
\ \ \ \ \ \ \
\delta
=
1+\beta''+\gamma'',
\]
implique l'\'equation impossible
$\alpha = -1 - \alpha '' - 2\gamma'' -
\delta''$, et enfin aussi, l'intersection entre la deuxi\`eme et la
troisi\`eme famille\,:
\[
\aligned
-2\beta'-\gamma'
=
-2\alpha''-2\gamma''-\delta'',
\ \ \ \ \ \ \
&
3\beta'+\gamma'
=
3\alpha''+2\gamma''+\delta'',
\\
&
1+\alpha'+\gamma'
=
\gamma'',
\ \ \ \ \ \ \
\delta'
=
1+\beta''+\gamma'',
\endaligned
\]
d'o\`u d\'ecoule $\gamma ' = \delta '' + 2\gamma''$, implique l'\'equation
impossible $\alpha ' = - 1 - \gamma '' - \delta''$, ce qui d\'emontre
l'assertion.
\hfill$\square$

\medskip
Ainsi, nous pouvons supprimer $\mathcal{ Q}_0$, 
$\mathcal{ S}_0$ et $\mathcal{ U}_0$, et nous sommes
ramen\'es \`a \'etudier l'identit\'e restante, qui est du type\,:
\[
\aligned
0
&
\equiv
\big(\Lambda^3\big\vert_0\big)^\mu
\Big[
\Lambda_1^5\big\vert_0\,
\mathcal{R}_0
\Big(
{\textstyle{
\frac{\Lambda_1^5\vert_0\,\Lambda_1^5\vert_0}{\Lambda^3\vert_0}
}},\,
M^8\big\vert_0,\,
{\textstyle{
\frac{\Lambda_1^5\vert_0\,\Lambda_1^5\vert_0\,\Lambda_1^5\vert_0}
{\Lambda^3\vert_0\,\Lambda^3\vert_0}
}},\,
N^{12}\big\vert_0
\Big)
+
\\
&
\ \ \ \ \ \ \ \ \ \ \ \ \ \ \ \ \ \ \ \
+
M^8\big\vert_0\,
\mathcal{T}_0
\Big(
{\textstyle{
\frac{\Lambda_1^5\vert_0\,\Lambda_1^5\vert_0\,\Lambda_1^5\vert_0}
{\Lambda^3\vert_0\,\Lambda^3\vert_0}
}},\,
{\textstyle{
\frac{\Lambda_1^5\vert_0\,M^8\vert_0}{\Lambda^3\vert_0}
}},\,
N^{12}\big\vert_0,\,
{\textstyle{
\frac{\Lambda_1^5\vert_0\,\Lambda_1^5\vert_0\,M^8\vert_0}
{\Lambda^3\vert_0\,\Lambda^3\vert_0}
}}
\Big)
\\
& 
\ \ \ \ \ \ \ \ \ \ \ \ \ \ \ \ \ \ \ \
+
N^{12}\big\vert_0\,
\mathcal{V}_0
\Big(
{\textstyle{
\frac{\Lambda_1^5\vert_0\,\Lambda_1^5\vert_0\,\Lambda_1^5\vert_0}
{\Lambda^3\vert_0\,\Lambda^3\vert_0}
}},\,
{\textstyle{
\frac{\Lambda_1^5\vert_0\,\Lambda_1^5\vert_0\,M^8\vert_0}
{\Lambda^3\vert_0\,\Lambda^3\vert_0}
}},\,
{\textstyle{
\frac{\Lambda_1^5\vert_0\,N^{12}\vert_0}
{\Lambda^3\vert_0}
}},\,
{\textstyle{
\frac{\Lambda_1^5\vert_0\,\Lambda_1^5\vert_0\,N^{12}\vert_0}
{\Lambda^3\vert_0\,\Lambda^3\vert_0}
}}
\Big)
\Big].
\endaligned
\]

\def\theassertion{\!}\begin{assertion}
Les trois familles de mon\^omes\,:

\begin{footnotesize}
\begin{itemize}

\smallskip\item[{\bf (iv)}]\ \ \ \ \
$\big(\Lambda^3\big\vert_0\big)^{\mu-\alpha-2\gamma}
\big(\Lambda_1^5\big\vert_0\big)^{1+2\alpha+3\gamma}
\big(M^8\big\vert_0\big)^{\beta}
\big(N^{12}\big\vert_0\big)^\delta$,

\smallskip\item[{\bf (v)}]\ \ \ \ \
$\big(\Lambda^3\big\vert_0\big)^{\mu-2\alpha'-\beta'-2\delta'}
\big(\Lambda_1^5\big\vert_0\big)^{3\alpha'+\beta'+2\delta'}
\big(M^8\big\vert_0\big)^{1+\beta'+\delta'}
\big(N^{12}\big\vert_0\big)^{\gamma'}$,

\smallskip\item[{\bf (vi)}]\ \ \ \ \
$\big(\Lambda^3\big\vert_0\big)^{\mu-2\alpha''-2\beta''-\gamma''-2\delta''}
\big(\Lambda_1^5\big\vert_0\big)^{3\alpha''+2\beta''+\gamma''+2\delta''}
\big(M^8\big\vert_0\big)^{\beta''}
\big(N^{12}\big\vert_0\big)^{1+\gamma''+\delta''}$,

\end{itemize}\smallskip
\end{footnotesize}

\noindent
ne contiennent {\sf aucune redondance}, {\rm i.e.} chaque mon\^ome
correspondand \`a un choix de $(\alpha, \beta, \gamma, \delta)$, ou de
$(\alpha', \beta', \gamma', \delta')$, ou encore de $(\alpha'',
\beta'', \gamma'', \delta'')$ appara\^{\i}t {\sf une et une seule fois}.
\end{assertion}

\medskip

En effet, il est tout d'abord facile de v\'erifier que chacune des trois
auto-intersec\-tions est triviale. Ensuite, des formules pour
l'intersection entre la premi\`ere et la deuxi\`eme famille\,:
\[
\alpha+2\gamma
=
2\alpha'+\beta'+2\delta',
\ \ \ \ \ \ \
1+2\alpha+3\gamma
=
3\alpha'+\beta'+2\delta',
\ \ \ \ \ \ \
\beta
=
1+\beta'+\delta',
\ \ \ \ \ \ \
\delta
=
\gamma',
\]
d\'ecoule l'\'equation $1 + \alpha + \gamma = \alpha'$ dont nous nous
servons pour remplacer $\alpha'$ dans la seconde \'equation, ce qui
conduit \`a l'impossibilit\'e $0 = 2 + \alpha + \beta ' + 2 \delta'$. De
m\^eme, l'intersection entre la premi\`ere et la troisi\`eme famille\,:
\[
\aligned
\alpha+2\gamma
=
2\alpha''+2\beta''+\gamma''+2\delta'',
&
\ \ \ \ \ \ \
1+2\alpha+3\gamma
=
3\alpha''+2\beta''+\gamma''+2\delta'',
\\
&
\ \ \ \ \ \ \
\beta
=
\beta'',
\ \ \ \ \ \ \
\delta
=
1+\gamma''+\delta'',
\endaligned
\]
conduit \`a l'impossibilit\'e $0 = 2 + \alpha + 2\beta '' + \gamma '' + 2
\delta ''$ en rempla\c cant $\alpha '' = 1 + \alpha + \gamma$ dans la
seconde \'equation. Enfin, l'intersection entre la deuxi\`eme et 
la troisi\`eme famille\,: 
\[
\small
\aligned
2\alpha'+\beta'+2\delta'
=
2\alpha''+2\beta''+\gamma''+2\delta'',
&
\ \ \ \ \ \ \
3\alpha'+\beta'+2\delta'
=
3\alpha''+2\beta''+\gamma''+2\delta'',
\\
&
\ \ \ \ \ \ \
1+\beta'+\delta'
=
\beta'',
\ \ \ \ \ \ \
\gamma'
=
1+\gamma''+\delta'',
\endaligned
\]
conduit \`a $\alpha ' = \alpha''$, d'o\`u $\beta ' + 2\delta ' = 2\beta ''
+ \gamma '' + 2\delta''$, puis en r\'e\'ecrivant la troisi\`eme \'equation et
en y rempla\c cant $-\beta ' - 2\delta '$ par $-2\beta '' - \gamma '' -
2\delta ''$, puis $-\beta ''$ par $-1 - \beta' - \delta'$\,:
\[
\small
\aligned
1
&
=
\beta''-\beta'-2\delta'+\delta'
\\
&
=
\beta''-2\beta''-\gamma''-2\delta''+\delta'
\\
&
=
-\beta''-\gamma''-2\delta''+\delta'
\\
&
=-1-\beta'-\gamma''-2\delta'',
\endaligned
\]
\'equation tout aussi impossible que les deux pr\'ec\'edentes.
Ceci ach\`eve la d\'emonstra\-tion de notre lemme principal.
\hfill$\square$

\smallskip\noindent{\bf Syzygies compl\`etes et substitutions 
alg\'ebriques.}
Nous parvenons enfin \`a la derni\`e\-re \'etape de la d\'emonstration de notre
second th\'eor\`eme. Soit ${\sf P}^{ 2 \times {\rm inv}} \big( j^5 f
\big)$ un bi-invariant quelconque, et reprenons son d\'eveloppement en
puissances positives et n\'egatives de $f_1'$\,:
\[
{\sf P}^{2\times{\rm inv}}
=
\sum_{-\frac{4}{5}m\leqslant a\leqslant m}\,
(f_1')^a\,
\mathcal{P}_a
\big(
\Lambda^3,\Lambda_1^5,\Lambda_{1,1}^7,M^8,\Lambda_{1,1,1}^9,
M_1^{10},N^{12},K_{1,1}^{12},H_1^{14},F_{1,1}^{16}
\big),
\] 
que nous avions obtenu pour le repr\'esenter en injectant
(artificiellement) les six bi-invariants fant\^omes dans une expression
initiale qui ne montrait que $\Lambda^3$, $\Lambda_1^5$, $\Lambda_{ 1,
1}^7$ et $\Lambda_{ 1, 1, 1}^9$. Choisissons l'exposant $a$
maximalement n\'egatif et examinons le polyn\^ome\,:
\[
\mathcal{P}_a
\big(
\Lambda^3,\Lambda_1^5,\Lambda_{1,1}^7,M^8,\Lambda_{1,1,1}^9,
M_1^{10},N^{12},K_{1,1}^{12},H_1^{14},F_{1,1}^{16}
\big).
\]
Notre base de Gr\"obner pour les bi-invariants restreints est obtenue
\`a partir de nos 15 syzygies fondamentales restreintes, et ce au moyen
d'un certain nombre d'op\'era\-tions alg\'ebriques \'el\'ementaires\,:
multiplication et addition d'\'equations, calcul de S-poly\-n\^omes et
divisions euclidiennes subs\'equentes, op\'erations autoris\'ees dans toute
structure d'id\'eal alg\'ebrique. Il en d\'ecoule que les m\^emes op\'erations
qui produisent les 21 \'equations normalis\'ees satisfaites sur $\{ f_1' =
0\}$ peuvent aussi \^etre conduites {\it sans poser $f_1' = 0$}, et
alors elles produisent, \`a partir des 15 syzygies compl\`etes, la m\^eme
liste de 21 \'equations dans laquelle chaque identit\'e ``$0 \equiv$''
doit \^etre remplac\'e par ``${\rm O} ( f_1') \equiv$'', le terme ${\rm O}
( f_1')$ d\'esignant un {\sf reste}, variable selon le contexte, qui d\'epend
{\it a priori}\, de tous les onze bi-invariants\footnote{\, Ici, le
raisonnement fonctionne seulement si l'alg\`ebre compl\`ete des
bi-invariant doit co\"{\i}ncider {\it a priori}\, avec l'alg\`ebre
engendr\'ee par crochets\,: sous cette condition, le {\sf reste} derri\`ere
$f_1'$ doit alors n\'ecessairement \^etre fonction des onze bi-invariants.
Toutefois, nous allons voir dans un instant que de nouveaux
bi-invariants fondamentaux ``fant\^omes'' se cachent derri\`ere
$f_1'$, exactement comme nous avions interpr\'et\'e $\Lambda_{ 1, 1,
1}^9$, $M_1^{ 10}$, $N^{ 12}$, $K_{ 1,1}^{ 12}$, $H_1^{ 14}$ et $F_{
1, 1}^{ 16}$, ces nouveaux bi-invariants n'\'etant pas construits par
crochet. } et qui s'annule lorsqu'on fait $f_1 ' = 0$, c'est-\`a-dire
qui est multiple de $f_1'$. Par cons\'equent, si nous appliquons \`a
$\mathcal{ P}_a$ la m\^eme normalisation modulo les syzygies que dans le
lemme principal (ce qui revient \`a substituer toutes les occurences des
mon\^omes de t\^ete), mais sans prendre la restriction \`a $\{ f_1' = 0\}$,
nous obtenons une expression du m\^eme type, et ce, avec un reste\,:
\[
\aligned
\mathcal{P}_a
&
=
\mathcal{P}_a
\big(
\Lambda^3,\,\Lambda_1^5,\,M^8,\,N^{12}
\big)
+
\Lambda_{1,1}^7\,
\mathcal{Q}_a\big(
\Lambda_1^5,\,\Lambda_{1,1}^7,\,M^8,\,N^{12}
\big)
+
\\
&
\ \ \ \ \
+
\Lambda_{1,1,1}^9\,
\mathcal{R}_a\big(
\Lambda_{1,1}^7,\,M^8,\,\Lambda_{1,1,1}^9,\,N^{12}
\big)
+
M_1^{10}\,
\mathcal{S}_a\big(
M^8,\,\Lambda_{1,1,1}^9,\,M_1^{10},\,N^{12}
\big)
+
\\
&
\ \ \ \ \
+
K_{1,1}^{12}\,
\mathcal{T}_a\big(
\Lambda_{1,1,1}^9,\,M_1^{10},\,N^{12},\,K_{1,1}^{12}
\big)
+
H_1^{14}\,
\mathcal{U}_a\big(
\Lambda_{1,1,1}^9,\,N^{12},\,K_{1,1}^{12},\,H_1^{14}
\big)
+
\\
&
\ \ \ \ \
+
F_{1,1}^{16}\,
\mathcal{V}\big(
\Lambda_{1,1,1}^9,\,K_{1,1}^{12},\,H_1^{14},\,F_{1,1}^{16}
\big)+
\\
&
\ \ \ \ \ 
+
f_1'\,{\sf reste}\big(
f_1',\Lambda^3,\Lambda_1^5,\Lambda_{1,1}^7,M^8,
\Lambda_{1,1,1}^9,M_1^{10},N^{12},K_{1,1}^{12},H_1^{14},F_{1,1}^{16}
\big),
\endaligned
\]
{\it a priori}\, non contr\^ol\'e, mais qui est heureusement repouss\'e dans
les puissances sup\'erieures $(f_1 ')^{ a +1}, \, (f_1')^{ a + 2},
\dots$. Ainsi, nous normalisons l'expression du premier polyn\^ome
$\mathcal{ P}_a$, \`a savoir celui qui appara\^{\i}t dans la puissance
maximalement n\'egative de $f_1'$. Et enuite, nous soumettons le nouveau
coefficient de $(f_1')^{ a + 1}$, qui vient de subir l'interf\'erence du
reste, et que nous noterons encore $\mathcal{ P}_a$, nous soumettons
ce nouveau coefficient au m\^eme processus de normalisation modulo les
21 syzygies non restreintes, et ainsi de suite, jusqu'\`a l'exposant
maximalement positif (toujours born\'e par $m$) de $f_1'$, 
ce qui nous donne une expression finale de la forme\,:
\[
\footnotesize
\aligned
{\sf P}^{2\times{\rm inv}}\big(j^5f\big)
&
=
\sum_{-\frac{4}{5}\,m\leqslant a\leqslant m}\,
(f_1')^a\,
\Big[
\mathcal{P}_a
\big(
\Lambda^3,\,\Lambda_1^5,\,M^8,\,N^{12}
\big)
+
\Lambda_{1,1}^7\,
\mathcal{Q}_a\big(
\Lambda_1^5,\,\Lambda_{1,1}^7,\,M^8,\,N^{12}
\big)
+
\\
&
\ \ \ \ \
+
\Lambda_{1,1,1}^9\,
\mathcal{R}_a\big(
\Lambda_{1,1}^7,\,M^8,\,\Lambda_{1,1,1}^9,\,N^{12}
\big)
+
M_1^{10}\,
\mathcal{S}_a\big(
M^8,\,\Lambda_{1,1,1}^9,\,M_1^{10},\,N^{12}
\big)
+
\\
&
\ \ \ \ \
+
K_{1,1}^{12}\,
\mathcal{T}_a\big(
\Lambda_{1,1,1}^9,\,M_1^{10},\,N^{12},\,K_{1,1}^{12}
\big)
+
H_1^{14}\,
\mathcal{U}_a\big(
\Lambda_{1,1,1}^9,\,N^{12},\,K_{1,1}^{12},\,H_1^{14}
\big)
+
\\
&
\ \ \ \ \
+
F_{1,1}^{16}\,
\mathcal{V}_a\big(
\Lambda_{1,1,1}^9,\,K_{1,1}^{12},\,H_1^{14},\,F_{1,1}^{16}
\big)
\Big]. 
\endaligned
\]
Et maintenant enfin, nous pouvons achever la d\'emonstration\,:
s'il existait des puissances n\'egatives de $f_1'$ dans une telle somme,
on multiplierait alors ${\sf P}^{ 2\times {\rm inv}} \big( j^5 f
\big)$ par la puissance positive minimale $(f_1')^{ - a}$ de $f_1'$
qui \'elimine les d\'enominateurs, on poserait $f_1' = 0$ et le lemme
principal\,\,---\,\,{\it which was specially designed on that
purpose}\,\,---\,\, tuerait alors les sept polyn\^omes $\mathcal{ P}_a$,
$\mathcal{ Q}_a$, $\mathcal{ R}_a$, $\mathcal{ S}_a$, $\mathcal{
T}_a$, $\mathcal{ U}_a$ et $\mathcal{ V}_a$, ce qui contredirait le
choix de $a$. Il n'existe donc que des puissances positives de $f_1'$,
et comme tout polyn\^ome de la forme g\'en\'erale
\[
\aligned
&
\mathcal{P}\big(f_1',\Lambda^3,\Lambda_1^5,M^8,N^{12}\big)
+
\Lambda_{1,1}^7\,
\mathcal{Q}\big(f_1',\Lambda_1^5,\Lambda_{1,1}^7,M^8,N^{12}\big)
+
\\
&
\ \ \ \ \
+
\Lambda_{1,1,1}^9\,
\mathcal{R}\big(f_1',\Lambda_{1,1}^7,M^8,\Lambda_{1,1,1}^9,N^{12}\big)
+
M_1^{10}\,
\mathcal{S}\big(f_1',M^8,\Lambda_{1,1,1}^9,M_1^{10},N^{12}\big)
+
\\
&
\ \ \ \ \
+
K_{1,1}^{12}\,
\mathcal{T}\big(f_1',\Lambda_{1,1,1}^9,M_1^{10},N^{12},K_{1,1}^{12}\big)
+
H_1^{14}\,
\mathcal{U}
\big(f_1',\Lambda_{1,1,1}^9,N^{12},K_{1,1}^{12},H_1^{14}\big)
+
\\
&
\ \ \ \ \
+
F_{1,1}^{16}\,
\mathcal{V}
\big(f_1',\Lambda_{1,1,1}^9,K_{1,1}^{12},H_1^{14},F_{1,1}^{16}\big),
\endaligned
\]
constitue trivialement un bi-invariant, \'ecrit qui plus est sous forme
unique gr\^ace au lemme fondamental, notre second th\'eor\`eme est \`a pr\'esent
compl\`etement d\'emontr\'e.
\hfill$\square$

\smallskip\noindent{\bf Observation.}
Dans cette derni\`ere \'etape du raisonnement, les termes de {\sf
reste} ci-dessus ont beau \^etre divisibles par $f_1'$, il ne sont pas
n\'ecessairement polynomiaux en nos onze bi-invariants fondamentaux\,; si
cela avait \'et\'e le cas, la strat\'egie aurait fonctionn\'e comme pour les
jets d'ordre 4. Plus pr\'ecis\'ement, lorsqu'on compare la liste originale
des 15 syzygies restreintes \`a la liste compl\`ete des 21 syzygies
restreintes, les 6 syzygies ajout\'ees sont d\'eduites des 15 initiales en
autorisant \`a diviser une syzygie par tout bi-invariant non
identiquement nul qui est en facteur, notamment par le wronskien
$\Lambda^3$, et c'est pour cette raison que les termes de {\sf
reste} ci-dessus ne sont pas n\'ecessairement polynomiaux en les onze
bi-invariants fondamentaux. Apr\`es un examen d\'etaill\'e, on d\'ecouvre donc
l'existence de exactement 7 bi-invariants ``fant\^omes''
suppl\'ementaires qui ne sont pas obtenus par crochets et qui se
cachent derri\`ere $f_1'$, \`a savoir\,:
\[
\scriptsize
\aligned
X^{18}
:=
&\
\frac{-5\,\Lambda_{1,1,1}^9\,M_1^{10}
+
56\,\Lambda_{1,1}^7\,K_{1,1}^{12}}{
f_1'}
\\
=
&\
f_1'f_1'f_1'\Big(
-
18816\,\Delta^{1,4}\,\big[\Delta^{2,3}\big]^2 - 25088\,
\big[\Delta^{2,3}\big]^3 - 15\,\big[
\Delta^{1,5}\big]^2\,\Delta^{1,2} - 150\,
\Delta^{1,5}\,\Delta^{2,4}\,\Delta^{1,2} 
\\
&
+ 315\,\Delta^{1,5}\,\Delta^{1,4}\,\Delta^{1,3} + 960\,
\Delta^{1,5}\,\Delta^{2,3}\,\Delta^{1,3} - 375\,\big[\Delta^{2,4}\big]^2
\,\Delta^{1,2} + 1575\,\Delta^{2,4}\,\Delta^{1,4}\,\Delta^{1,3}
\\
&
+ 4800\,\Delta^{2,4}\,\Delta^{2,3}\,\Delta^{1,3} - 392\,
\big[\Delta^{1,4}\big]^3 - 4704\,\big[\Delta^{1,4}\big]^2\,\Delta^{2,3}
\Big)
-
f_1'f_1'f_1''\,
\Big(
- 2475\,\Delta^{2,4}\,\Delta^{1,4}\,\Delta^{1,2}
\\
&
- 9900\,\Delta^{2,4}\,\Delta^{2,3}\,\Delta^{1,2} - 2850\,
\Delta^{1,5}\,\big[\Delta^{1,3}\big]^2 + 51330\,\Delta^{1,4}\,
\Delta^{2,3}\,\Delta^{1,3}
\\
&
+ 92760\,\big[\Delta^{2,3}\big]^2\,\Delta^{1,3} - 14250\,\Delta^{2,4}\,
\big[\Delta^{1,3}\big]^2 + 7035\,\big[\Delta^{1,4}\big]^2\,\Delta^{1,3} - 
495\,\Delta^{1,5}\,\Delta^{1,4}\,\Delta^{1,2}
\\
&
- 1980\,\Delta^{1,5}\,\Delta^{2,3}\,\Delta^{1,2}
\Big) 
- 
f_1'f_1'f_1'''\,
\Big(
- 11100\,\Delta^{2,3}\,\big[\Delta^{1,3}\big]^2 - 3150\,
\Delta^{1,4}\,\big[\Delta^{1,3}\big]^2
\Big)
\\
&
+
f_1'f_1''f_1''\,
\Big(
- 109440\,\big[\Delta^{2,3}\big]^2\,\Delta^{1,2} - 19050\,\Delta^{2,3}\,
\big[\Delta^{1,3}\big]^2 - 32325\,\Delta^{1,4}\,\big[\Delta^{1,3}\big]^2 
\\
&
+ 11025\,\Delta^{1,5}\,\Delta^{1,3}\,\Delta^{1,2} + 55125
\,\Delta^{2,4}\,\Delta^{1,3}\,\Delta^{1,2} - 6840\,\big[\Delta^{1,4}\big]^2
\,\Delta^{1,2}
\\
&
- 54720\,\Delta^{1,4}\,\Delta^{2,3}\,\Delta^{1,2}
\Big) 
-
f_1'f_1''f_1'''\,\Big(
+ 30000\,
\big[\Delta^{1,3}\big]^3
\Big) 
- 
f_1''f_1''f_1''\,
\Big(11025\,\Delta^{1,5}\,
\big[\Delta^{1,2}\big]^2
\\
&
- 55125\,\Delta^{2,4}\,
\big[\Delta^{1,2}\big]^2 + 55125\,\Delta^{1,4}\,
\Delta^{1,3}\,\Delta^{1,2} + 110250\,\Delta^{2,3}\,
\Delta^{1,3}\,\Delta^{1,2} 
\\
&
- 49000\,\big[\Delta^{1,3}\big]^3
\Big).
\endaligned
\]
\[
\scriptsize
\aligned
X^{19}
:=
&\
\frac{-5\,M_1^{10}\,M_1^{10}+64\,M^8\,K_{1,1}^{12}}{f_1'}
\\
=
&\,
f_1'\, 
\Big(
1170\,\Delta^{1,5}\,\Delta^{1,4}\,\Delta^{1,3}\,
\Delta^{1,2} - 45\,\big[\Delta^{1,5}\big]^2\,\big[
\Delta^{1,2}\big]^2 - 450\,
\Delta^{1,5}\,\Delta^{2,4}\,\big[\Delta^{1,2}\big]^2 
\\
&
+ 74220\,\big[\Delta^{2,3}\big]^2\,\big[\Delta^{1,3}\big]^2 + 3780\,
\Delta^{1,5}\,\Delta^{2,3}\,\Delta^{1,3}\,\Delta^{1,2} - 1600\,
\Delta^{1,5}\,\big[\Delta^{1,3}\big]^3 
\\
&
- 1125\,\big[\Delta^{2,4}\big]^2\,\big[\Delta^{1,2}\big]^2 + 5850\,
\Delta^{2,4}\,\Delta^{1,4}\,\Delta^{1,3}\,\Delta^{1,2} + 18900\,
\Delta^{2,4}\,\Delta^{2,3}\,\Delta^{1,3}\,\Delta^{1,2} 
\\
&
- 8000\,\Delta^{2,4}\,\big[\Delta^{1,3}\big]^3 - 1344\,
\big[\Delta^{1,4}\big]^3\,
\Delta^{1,2} - 16128\,\big[\Delta^{1,4}\big]^2\,\Delta^{2,3}\,
\Delta^{1,2} + 1995\,\big[\Delta^{1,4}\big]^2\,\big[\Delta^{1,3}\big]^2 
\\
&
- 64512\,\Delta^{1,4}\,\big[\Delta^{2,3}\big]^2\,\Delta^{1,2} + 
27660\,\Delta^{1,4}\,\Delta^{2,3}\,\big[\Delta^{1,3}\big]^2 - 86016\,
\big[\Delta^{2,3}\big]^3\,\Delta^{1,2}
\Big) 
\\
&
+
f_1''\,
\Big(
- 74400\,\Delta^{2,3}\,\big[\Delta^{1,3}\big]^3 - 10800\,\Delta^{2,4}\,
\Delta^{1,4}\,\big[\Delta^{1,2}\big]^2 - 2160\,\Delta^{1,5}\,\Delta^{1,4}
\,\big[\Delta^{1,2}\big]^2 
\\
&
- 8640\,\Delta^{1,5}\,\Delta^{2,3}\,\big[\Delta^{1,2}\big]^2 + 
3600\,\Delta^{1,5}\,\big[\Delta^{1,3}\big]^2\,\Delta^{1,2} + 64800\,
\Delta^{1,4}\,\Delta^{2,3}\,\Delta^{1,3}\,\Delta^{1,2} 
\\
&
- 43200\,\Delta^{2,4}\,\Delta^{2,3}\,\big[\Delta^{1,2}\big]^2 + 
18000\,\Delta^{2,4}\,\big[\Delta^{1,3}\big]^2\,\Delta^{1,2} + 10800\,
\big[\Delta^{1,4}\big]^2\,\Delta^{1,3}\,\Delta^{1,2} 
\\
&
- 27600\,\Delta^{1,4}\,\big[\Delta^{1,3}\big]^3 + 86400\,
\big[\Delta^{2,3}\big]^2\,\Delta^{1,3}\,\Delta^{1,2}
\Big)
+ 
f_1'''\,\Big(
16000\,\big[\Delta^{1,3}\big]^4\Big).
\endaligned
\]
\[
\scriptsize
\aligned
X^{21}
:=
&\
\frac{-5\,M_1^{10}\,N^{12}+8\,M^8\,H_1^{14}}{f_1'}
\\
=
&\, 
- 135\,\big[\Delta^{1,5}\big]^2\,\big[\Delta^{1,2}\big]^3 - 1350
\,\Delta^{1,5}\,\Delta^{2,4}\,\big[\Delta^{1,2}\big]^3 + 1350\,
\Delta^{1,5}\,\Delta^{1,4}\,\Delta^{1,3}\,\big[\Delta^{1,2}\big]^2 
\\
&
+ 2700\,\Delta^{1,5}\,\Delta^{2,3}\,\Delta^{1,3}\,
\big[\Delta^{1,2}\big]^2 - 1200\,\Delta^{1,5}\,\big[\Delta^{1,3}\big]^3\,
\Delta^{1,2} - 3375\,\big[\Delta^{2,4}\big]^2\,\big[\Delta^{1,2}\big]^3 
\\
&
+ 6750\,\Delta^{2,4}\,\Delta^{1,4}\,\Delta^{1,3}\,
\big[\Delta^{1,2}\big]^2 + 13500\,\Delta^{2,4}\,\Delta^{2,3}\,
\Delta^{1,3}\,\big[\Delta^{1,2}\big]^2 - 6000\,\Delta^{2,4}\,
\big[\Delta^{1,3}\big]^3\,
\Delta^{1,2} 
\\
&
- 576\,\big[\Delta^{1,4}\big]^3\,\big[\Delta^{1,2}\big]^2 - 6912\,
\big[\Delta^{1,4}\big]^2\,\Delta^{2,3}\,\big[\Delta^{1,2}\big]^2 - 495\,
\big[\Delta^{1,4}\big]^2\,\big[\Delta^{1,3}\big]^2\,\Delta^{1,2} 
\\
&
- 27648\,\Delta^{1,4}\,\big[\Delta^{2,3}\big]^2\,\big[\Delta^{1,2}\big]^2
+ 9540\,\Delta^{1,4}\,\Delta^{2,3}\,\big[\Delta^{1,3}\big]^2\,
\Delta^{1,2} + 1200\,\Delta^{1,4}\,\big[\Delta^{1,3}\big]^4 
\\
&
- 36864\,\big[\Delta^{2,3}\big]^3\,\big[\Delta^{1,2}\big]^2 + 32580\,
\big[\Delta^{2,3}\big]^2\,\big[\Delta^{1,3}\big]^2\,\Delta^{1,2} - 7200\,
\Delta^{2,3}\,\big[\Delta^{1,3}\big]^4.
\endaligned
\]
\[
\scriptsize
\aligned
X^{23}
:=
&\
\frac{-7\,N^{12}\,K_{1,1}^{12}+M^8\,F_{1,1}^{16}}{f_1'}
\\
=
&\, 
f_1'\,
\Big(432\,\Delta^{1,5}\,\big[\Delta^{1,4}\big]^2\,\big[\Delta^{1,2}\,
\big]^2 + 3456\,\Delta^{1,5}\,\Delta^{1,4}\,\Delta^{2,3}\,
\big[\Delta^{1,2}\big]^2 + 1710\,\Delta^{1,5}\,\Delta^{1,4}\,
\big[\Delta^{1,3}\big]^2\,
\Delta^{1,2} 
\\
&
- 3150\,\Delta^{1,5}\,\Delta^{2,4}\,\Delta^{1,3}\,
\big[\Delta^{1,2}\big]^2 + 540\,\Delta^{1,5}\,\Delta^{2,3}\,\big[\Delta^{1,3}
\big]^2\,\Delta^{1,2} - 1600\,\Delta^{1,5}\,\big[\Delta^{1,3}\big]^4 
\\
&
- 7875\,\big[\Delta^{2,4}\big]^2\,\Delta^{1,3}\,\big[\Delta^{1,2}\big]^2
+ 6912\,\Delta^{1,5}\,\big[\Delta^{2,3}\big]^2\,\big[
\Delta^{1,2}\big]^2 - 8000
\,\Delta^{2,4}\,\big[\Delta^{1,3}\big]^4 
\\
&
- 2352\,\big[\Delta^{1,4}\big]^3\,\Delta^{1,3}\,\Delta^{1,2} - 
23904\,\big[\Delta^{1,4}\big]^2\,\Delta^{2,3}\,\Delta^{1,3}\,\Delta^{1,2}
+ 2205\,\big[\Delta^{1,4}\big]^2\,\big[\Delta^{1,3}\big]^3 
\\
&
- 78336\,\Delta^{1,4}\,\big[\Delta^{2,3}\big]^2\,\Delta^{1,3}\,
\Delta^{1,2} + 34740\,\Delta^{1,4}\,\Delta^{2,3}\,\big[\Delta^{1,3}\big]^3
- 81408\,\big[\Delta^{2,3}\big]^3\,\Delta^{1,3}\,\Delta^{1,2} 
\\
&
+ 72180\,\big[\Delta^{2,3}\big]^2\,\big[\Delta^{1,3}\big]^3 + 2160\,
\Delta^{2,4}\,\big[\Delta^{1,4}\big]^2\,\big[\Delta^{1,2}\big]^2 + 17280\,
\Delta^{2,4}\,\Delta^{1,4}\,\Delta^{2,3}\,\big[\Delta^{1,2}\big]^2 
\\
&
+ 8550\,\Delta^{2,4}\,\Delta^{1,4}\,\big[\Delta^{1,3}\big]^2\,
\Delta^{1,2} + 34560\,\Delta^{2,4}\,\big[\Delta^{2,3}\big]^2\,
\big[\Delta^{1,2}\big]^2 + 2700\,\Delta^{2,4}\,\Delta^{2,3}\,
\big[\Delta^{1,3}\big]^2\,
\Delta^{1,2} 
\\
&
- 315\,\big[\Delta^{1,5}\big]^2\,\Delta^{1,3}\,\big[\Delta^{1,2}\big]^2
\Big)
+
f_1''\,
\Big(
23625\,\big[\Delta^{2,4}\big]^2\,\big[\Delta^{1,2}\big]^3
 - 47250\,\Delta^{2,4}\,\Delta^{1,4}\,\Delta^{1,3}\,\big[\Delta^{1,2}
\big]^2 
\\
&
- 94500\,\Delta^{2,4}\,\Delta^{2,3}\,\Delta^{1,3}\,
\big[\Delta^{1,2}\big]^2 + 42000\,\Delta^{2,4}\,\big[\Delta^{1,3}\big]^3\,
\Delta^{1,2} + 576\,\big[\Delta^{1,4}\big]^3\,\big[\Delta^{1,2}\big]^2 
\\
&
+ 6912\,\big[\Delta^{1,4}\big]^2\,\Delta^{2,3}\,\big[\Delta^{1,2}\big]^2
+ 20745\,\big[\Delta^{1,4}\big]^2\,\big[\Delta^{1,3}\big]^2\,\Delta^{1,2} + 
27648\,\Delta^{1,4}\,\big[\Delta^{2,3}\big]^2\,\big[\Delta^{1,2}\big]^2 
\\
&
+ 945\,\big[\Delta^{1,5}\big]^2\,\big[\Delta^{1,2}\big]^3 + 9450\,
\Delta^{1,5}\,\Delta^{2,4}\,\big[\Delta^{1,2}\big]^3 - 9450\,\Delta^{1,5}
\,\Delta^{1,4}\,\Delta^{1,3}\,\big[\Delta^{1,2}\big]^2 
\\
&
- 18900\,\Delta^{1,5}\,\Delta^{2,3}\,\Delta^{1,3}\,
\big[\Delta^{1,2}\big]^2 + 8400\,\Delta^{1,5}\,\big[\Delta^{1,3}\big]^3\,
\Delta^{1,2} + 71460\,\Delta^{1,4}\,\Delta^{2,3}\,\big[\Delta^{1,3}\big]^2\,
\Delta^{1,2} 
\\
&
- 37200\,\Delta^{1,4}\,\big[\Delta^{1,3}\big]^4 + 36864\,\big[\Delta^{2,3}\,
\big]^3\,\big[\Delta^{1,2}\big]^2 + 48420\,\big[\Delta^{2,3}\big]^2\,
\big[\Delta^{1,3}\big]^2\,\Delta^{1,2} 
\\
&
- 64800\,\Delta^{2,3}\,\big[\Delta^{1,3}\big]^4
\Big)
+ 
f_1'''\,\Big(
16000\,\big[\Delta^{1,3}\big]^5\Big).
\endaligned
\]
\[
\scriptsize
\aligned
Y^{23}
:=
&\
\frac{-8\,N^{12}\,K_{1,1}^{12}+M_1^{10}\,H_1^{14}}{f_1'}
\\
=
&\,
X^{23}.
\endaligned
\]
\[
\scriptsize
\aligned
X^{25}
:=
&\
\frac{-56\,K_{1,1}^{12}\,H_1^{14}+5\,M_1^{10}\,F_{1,1}^{16}}{f_1'}
\\
=
&\, 
f_1'\,f_1'\,
\Big(
-45\,\big[\Delta^{1,5}\big]^2\,\Delta^{1,4}\,
\big[\Delta^{1,2}\big]^2 - 180\,\big[\Delta^{1,5}\big]^2\,\Delta^{2,3}\,
\big[\Delta^{1,2}\big]^2
- 3600\,\big[\Delta^{1,5}\big]^2\,\big[\Delta^{1,3}\big]^2\,\Delta^{1,2} 
\\
&
 - 2800\,\Delta^{1,5}\,\Delta^{1,4}\,\big[\Delta^{1,3}\big]^3 - 
83200\,\Delta^{1,5}\,\Delta^{2,3}\,\big[\Delta^{1,3}\big]^3 - 1125\,
\big[\Delta^{2,4}\big]^2\,\Delta^{1,4}\,\big[\Delta^{1,2}\big]^2 
\\
&
- 4500\,\big[\Delta^{2,4}\big]^2\,\Delta^{2,3}\,\big[\Delta^{1,2}\big]^2
- 90000\,\big[\Delta^{2,4}\big]^2\,\big[\Delta^{1,3}\big]^2\,\Delta^{1,2} - 
14000\,\Delta^{2,4}\,\Delta^{1,4}\,\big[\Delta^{1,3}\big]^3 
\\
&
- 416000\,\Delta^{2,4}\,\Delta^{2,3}\,\big[\Delta^{1,3}\big]^3 - 
150528\,\big[\Delta^{1,4}\big]^3\,\Delta^{2,3}\,\Delta^{1,2} - 903168\,
\big[\Delta^{1,4}\big]^2\,\big[\Delta^{2,3}\big]^2\,\Delta^{1,2} 
\\
&
+ 163800\,\big[\Delta^{1,4}\big]^2\,\Delta^{2,3}\,\big[\Delta^{1,3}\big]^2 
-2408448\,\Delta^{1,4}\,\big[\Delta^{2,3}\big]^3\,\Delta^{1,2} + 
1129500\,\Delta^{1,4}\,\big[\Delta^{2,3}\big]^2\,\big[\Delta^{1,3}\big]^2 
\\
&
- 9408\,\big[\Delta^{1,4}\big]^4\,\Delta^{1,2} + 3675\,\big[\Delta^{1,4}\,
\big]^3\,\big[\Delta^{1,3}\big]^2 - 2408448\,\big[\Delta^{2,3}\big]^4\,
\Delta^{1,2} + 2132400\,\big[\Delta^{2,3}\big]^3\,\big[\Delta^{1,3}\big]^2 
\\
&
- 450\,\Delta^{1,5}\,\Delta^{2,4}\,\Delta^{1,4}\,
\big[\Delta^{1,2}\big]^2 - 1800\,\Delta^{1,5}\,\Delta^{2,4}\,\Delta^{2,3}\,
\big[\Delta^{1,2}\big]^2 
\\
&
- 36000\,\Delta^{1,5}\,\Delta^{2,4}\,\big[\Delta^{1,3}\big]^2\,
\Delta^{1,2} + 11970\,\Delta^{1,5}\,\big[\Delta^{1,4}\big]^2\,
\Delta^{1,3}\,\Delta^{1,2} 
\\
&
+ 187920\,\Delta^{1,5}\,\big[\Delta^{2,3}\big]^2\,\Delta^{1,3}\,
\Delta^{1,2} + 59850\,\Delta^{2,4}\,\big[\Delta^{1,4}\big]^2\,
\Delta^{1,3}\,\Delta^{1,2} 
\\
&
+ 939600\,\Delta^{2,4}\,\big[\Delta^{2,3}\big]^2\,\Delta^{1,3}\,
\Delta^{1,2} + 474300\,\Delta^{2,4}\,\Delta^{1,4}\,\Delta^{2,3}\,
\Delta^{1,3}\,\Delta^{1,2} 
\\
&
+ 94860\,\Delta^{1,5}\,\Delta^{1,4}\,\Delta^{2,3}\,
\Delta^{1,3}\,\Delta^{1,2}
\Big)
+ 
f_1'f_1''\,
\Big(
- 2556600
\,\Delta^{1,4}\,\Delta^{2,3}\,\big[\Delta^{1,3}\big]^3 
\\
\endaligned
\]
\[
\scriptsize
\aligned 
&
- 5014200\,\big[\Delta^{2,3}\big]^2\,\big[\Delta^{1,3}\big]^3 - 187950\,
\big[\Delta^{1,4}\big]^2\,\big[\Delta^{1,3}\big]^3 + 5621760\,\Delta^{1,4}\,
\big[\Delta^{2,3}\big]^2\,\Delta^{1,3}\,\Delta^{1,2}
\\
&
+ 5652480\,\big[\Delta^{2,3}\big]^3\,\Delta^{1,3}\,\Delta^{1,2}
- 2764800\,\Delta^{2,4}\,\big[\Delta^{2,3}\big]^2\,\big[\Delta^{1,2}\big]^2 
\\
&
+ 99000\,\Delta^{2,4}\,\Delta^{2,3}\,\big[\Delta^{1,3}\big]^2\,
\Delta^{1,2} + 500000\,\Delta^{2,4}\,\big[\Delta^{1,3}\big]^4 + 174720\,
\big[\Delta^{1,4}\big]^3\,\Delta^{1,3}\,\Delta^{1,2} 
\\
&
+ 1751040\,\big[\Delta^{1,4}\big]^2\,\Delta^{2,3}\,\Delta^{1,3}\,
\Delta^{1,2} - 276480\,\Delta^{1,5}\,\Delta^{1,4}\,\Delta^{2,3}\,
\big[\Delta^{1,2}\big]^2 
\\
&
- 105300\,\Delta^{1,5}\,\Delta^{1,4}\,\big[\Delta^{1,3}\big]^2\,
\Delta^{1,2} - 552960\,\Delta^{1,5}\,\big[\Delta^{2,3}\big]^2\,
\big[\Delta^{1,2}\big]^2 
\\
&
+ 19800\,\Delta^{1,5}\,\Delta^{2,3}\,\big[\Delta^{1,3}\big]^2\,
\Delta^{1,2} + 100000\,\Delta^{1,5}\,\big[\Delta^{1,3}\big]^4 + 551250\,
\big[\Delta^{2,4}\big]^2\,\Delta^{1,3}\,\big[\Delta^{1,2}\big]^2 
\\
&
- 172800\,\Delta^{2,4}\,\big[\Delta^{1,4}\big]^2\,\big[\Delta^{1,2}\big]^2
- 1382400\,\Delta^{2,4}\,\Delta^{1,4}\,\Delta^{2,3}\,
\big[\Delta^{1,2}\big]^2 
\\
&
- 526500\,\Delta^{2,4}\,\Delta^{1,4}\,\big[\Delta^{1,3}\big]^2\,
\Delta^{1,2} - 34560\,\Delta^{1,5}\,\big[\Delta^{1,4}\big]^2\,
\big[\Delta^{1,2}\big]^2 + 22050\,\big[\Delta^{1,5}\big]^2\,\Delta^{1,3}\,
\big[\Delta^{1,2}\big]^2
\\
&
+ 220500\,\Delta^{1,5}\,\Delta^{2,4}\,\Delta^{1,3}\,
\big[\Delta^{1,2}\big]^2
\Big)
\\
&
+ 
f_1'f_1'''\,
\Big(
28000\,\Delta^{1,4}\,\big[\Delta^{1,3}\big]^4 + 472000\,
\Delta^{2,3}\,\big[\Delta^{1,3}\big]^4
\Big)
\\
&
+
f_1''f_1''\,
\Big(
330750\,\Delta^{1,5}\,\Delta^{1,4}\,\Delta^{1,3}\,\big[\Delta^{1,2}\big]^2
+ 661500\,\Delta^{1,5}\,\Delta^{2,3}\,\Delta^{1,3}\,
\big[\Delta^{1,2}\big]^2 
\\
&
- 294000\,\Delta^{1,5}\,\big[\Delta^{1,3}\big]^3\,\Delta^{1,2} - 
330750\,\Delta^{1,5}\,\Delta^{2,4}\,\big[\Delta^{1,2}\big]^3 
\\
&
+ 1653750\,\Delta^{2,4}\,\Delta^{1,4}\,\Delta^{1,3}\,
\big[\Delta^{1,2}\big]^2 + 3307500\,\Delta^{2,4}\,\Delta^{2,3}\,
\Delta^{1,3}\,\big[\Delta^{1,2}\big]^2 
\\
&
- 1470000\,\Delta^{2,4}\,\big[\Delta^{1,3}\big]^3\,\Delta^{1,2}
- 2880\,\big[\Delta^{1,4}\big]^3\,\big[
\Delta^{1,2}\big]^2 - 34560\,\big[\Delta^{1,4}
\big]^2\,\Delta^{2,3}\,\big[\Delta^{1,2}\big]^2 
\\
&
- 812475\,\big[\Delta^{1,4}\big]^2\,\big[\Delta^{1,3}\big]^2\,\Delta^{1,2}\,
- 138240\,\Delta^{1,4}\,\big[\Delta^{2,3}\big]^2\,\big[\Delta^{1,2}\big]^2 - 
33075\,\big[\Delta^{1,5}\big]^2\,\big[\Delta^{1,2}\big]^3 
\\
&
+ 1446000\,\Delta^{1,4}\,\big[\Delta^{1,3}\big]^4 - 184320\,
\big[\Delta^{2,3}\big]^3\,\Delta^{1,2}\big]^2 
- 3077100\,\big[\Delta^{2,3}\big]^2\,
\big[\Delta^{1,3}\big]^2\,\Delta^{1,2} 
\\
&
+ 2844000\,\Delta^{2,3}\,\big[\Delta^{1,3}\big]^4 - 826875\,
\big[\Delta^{2,4}\big]^2\,\big[\Delta^{1,2}\big]^3 - 3192300\,\Delta^{1,4}\,
\Delta^{2,3}\,\big[\Delta^{1,3}\big]^2\,\Delta^{1,2}
\Big)
\\
&
+
f_1''f_1'''\,
\Big( - 640000\,\big[\Delta^{1,3}\big]^5
\Big).
\endaligned
\]
\[
\scriptsize
\aligned
X^{27}
:=
&\
\frac{-7\,H_1^{14}\,H_1^{14}+5\,N^{12}\,F_{1,1}^{16}}{f_1'}
\\
=
&\,
f_1'\,
\Big( 
-1032192\,\Delta^{1,4}\,\big[\Delta^{2,3}\big]^3\,
\big[\Delta^{1,2}\big]^2 - 186300\,\Delta^{1,4}\,\Delta^{2,3}\,
\big[\Delta^{1,3}\big]^4 - 3375\,\big[\Delta^{2,4}\big]^2\,\Delta^{1,4}\,
\big[\Delta^{1,2}\big]^3
\\
&
+ 5625\,\big[\Delta^{2,4}\big]^2\,\big[\Delta^{1,3}\big]^2\,
\big[\Delta^{1,2}
\big]^2 - 540\,\big[\Delta^{1,5}\big]^2\,\Delta^{2,3}\,
\big[\Delta^{1,2}\big]^3 + 
12705\,\big[\Delta^{1,4}\big]^3\,\big[\Delta^{1,3}\big]^2\,\Delta^{1,2} 
\\
&
+ 1320720\,\big[\Delta^{2,3}\big]^3\,\big[\Delta^{1,3}\big]^2\,
\Delta^{1,2} - 64512\,\big[\Delta^{1,4}\big]^3\,\Delta^{2,3}\,
\big[\Delta^{1,2}\big]^2 + 
225\,\big[\Delta^{1,5}\big]^2\,\big[\Delta^{1,3}\big]^2\,
\big[\Delta^{1,2}\big]^2 
\\
&
- 135\,\big[\Delta^{1,5}\big]^2\,\Delta^{1,4}\,\big[\Delta^{1,2}\big]^3
- 13500\,\big[\Delta^{2,4}\big]^2\,\Delta^{2,3}\,\big[\Delta^{1,2}\big]^3 - 
387072\,\big[\Delta^{1,4}\big]^2\,\big[\Delta^{2,3}\big]^2\,
\big[\Delta^{1,2}\big]^2 
\\
&
- 1350\,\Delta^{1,5}\,\Delta^{2,4}\,\Delta^{1,4}\,
\big[\Delta^{1,2}\big]^3 - 5400\,\Delta^{1,5}\,\Delta^{2,4}\,\Delta^{2,3}
\,\big[\Delta^{1,2}\big]^3 
\\
&
+ 2250\,\Delta^{1,5}\,\Delta^{2,4}\,\big[\Delta^{1,3}\big]^2\,
\big[\Delta^{1,2}\big]^2 + 3510\,\Delta^{1,5}\,\big[\Delta^{1,4}\big]^2\,
\Delta^{1,3}\,\big[\Delta^{1,2}\big]^2 
\\
&
- 10650\,\Delta^{1,5}\,\Delta^{1,4}\,\big[\Delta^{1,3}\big]^3\,
\Delta^{1,2} + 45360\,\Delta^{1,5}\,\big[\Delta^{2,3}\big]^2\,
\Delta^{1,3}\,\big[\Delta^{1,2}\big]^2 
\\
&
- 38100\,\Delta^{1,5}\,\Delta^{2,3}\,\big[\Delta^{1,3}\big]^3\,
\Delta^{1,2} + 17550\,\Delta^{2,4}\,\big[\Delta^{1,4}\big]^2\,
\Delta^{1,3}\,\big[\Delta^{1,2}\big]^2 
\\
&
- 53250\,\Delta^{2,4}\,\Delta^{1,4}\,\big[\Delta^{1,3}\big]^3\,
\Delta^{1,2} + 226800\,\Delta^{2,4}\,\big[\Delta^{2,3}\big]^2\,
\Delta^{1,3}\,\big[\Delta^{1,2}\big]^2 
\\
&
- 190500\,\Delta^{2,4}\,\Delta^{2,3}\,\big[\Delta^{1,3}\big]^3\,
\Delta^{1,2} + 187560\,\big[\Delta^{1,4}\big]^2\,\Delta^{2,3}\,
\big[\Delta^{1,3}\big]^2\,\Delta^{1,2} 
\\
&
+ 877140\,\Delta^{1,4}\,\big[\Delta^{2,3}\big]^2\,\big[\Delta^{1,3}\big]^2
\,\Delta^{1,2} + 8000\,\Delta^{1,5}\,\big[\Delta^{1,3}\big]^5 + 40000\,
\Delta^{2,4}\,\big[\Delta^{1,3}\big]^5 
\\
&
- 4032\,\big[\Delta^{1,4}\big]^4\,\big[\Delta^{1,2}\big]^2 - 9975\,
\big[\Delta^{1,4}\big]^2\,\big[\Delta^{1,3}\big]^4 
- 1032192\,\big[\Delta^{2,3}\big]^4\,
\big[\Delta^{1,2}\big]^2 
\\
&
- 563100\,\big[\Delta^{2,3}\big]^2\,\big[\Delta^{1,3}\big]^4 + 126900\,
\Delta^{2,4}\,\Delta^{1,4}\,\Delta^{2,3}\,\Delta^{1,3}\,
\big[\Delta^{1,2}\big]^2 
\\
&
+ 25380\,\Delta^{1,5}\,\Delta^{1,4}\,\Delta^{2,3}\,
\Delta^{1,3}\,\big[\Delta^{1,2}\big]^2
\Big) 
+
f_1''\,
\Big( 
-259200\,
\Delta^{2,4}\,\Delta^{1,4}\,\Delta^{2,3}\,\big[\Delta^{1,2}\big]^3 
\\
&
- 518400\,\Delta^{2,4}\,\big[\Delta^{2,3}\big]^2\,\big[\Delta^{1,2}\big]^3
+ 432000\,\Delta^{2,4}\,\Delta^{2,3}\,\big[\Delta^{1,3}\big]^2\,
\Delta^{1,2}\big]^2 
\\
&
- 90000\,\Delta^{2,4}\,\big[\Delta^{1,3}\big]^4\,\Delta^{1,2} + 
32400\,\big[\Delta^{1,4}\big]^3\,\Delta^{1,3}\,\big[\Delta^{1,2}\big]^2 
+ 
324000
\,\big[\Delta^{1,4}\big]^2\,
\Delta^{2,3}\,\Delta^{1,3}\,\big[\Delta^{1,2}\big]^2
\\
&
- 136800\,\big[\Delta^{1,4}\big]^2\,\big[\Delta^{1,3}\big]^3\,
\Delta^{1,2} + 1036800\,\Delta^{1,4}\,\big[\Delta^{2,3}\big]^2\,
\Delta^{1,3}\,
\big[\Delta^{1,2}\big]^2 
\\
&
- 878400\,\Delta^{1,4}\,\Delta^{2,3}\,\big[\Delta^{1,3}\big]^3\,
\Delta^{1,2} + 186000\,\Delta^{1,4}\,\big[\Delta^{1,3}\big]^5 + 1036800\,
\big[\Delta^{2,3}\big]^3\,\Delta^{1,3}\,\big[\Delta^{1,2}\big]^2 
\\
&
- 1324800\,\big[\Delta^{2,3}\big]^2\,\big[\Delta^{1,3}\big]^3\,
\Delta^{1,2} + 564000\,\Delta^{2,3}\,\big[\Delta^{1,3}\big]^5 + 108000\,
\Delta^{2,4}\,\Delta^{1,4}\,
\big[\Delta^{1,3}\big]^2\,\big[\Delta^{1,2}\big]^2 
\\
\endaligned
\]
\[
\scriptsize
\aligned 
&
- 103680\,\Delta^{1,5}\,\big[\Delta^{2,3}\big]^2\,\big[\Delta^{1,2}\big]^3
- 6480\,\Delta^{1,5}\,\big[\Delta^{1,4}\big]^2\,\big[\Delta^{1,2}\big]^3 - 
51840\,\Delta^{1,5}\,\Delta^{1,4}\,\Delta^{2,3}\,\big[\Delta^{1,2}\big]^3
\\
&
+ 21600\,\Delta^{1,5}\,\Delta^{1,4}\,\big[\Delta^{1,3}\big]^2\,
\big[\Delta^{1,2}\big]^2 + 86400\,\Delta^{1,5}\,\Delta^{2,3}\,
\big[\Delta^{1,3}\big]^2\,\big[\Delta^{1,2}\big]^2 
\\
&
- 18000\,\Delta^{1,5}\,\big[\Delta^{1,3}\big]^4\,\Delta^{1,2} - 
32400\,\Delta^{2,4}\,\big[\Delta^{1,4}\big]^2\,\big[\Delta^{1,2}\big]^3
\Big)
+
f_1'''\,
\Big(
-
80000\,\Big[\Delta^{1,3}\big]^6\Big).
\endaligned
\]

Toutes les tentatives de calcul alg\'ebrique pour \'egaler l'un de ces 6
bi-invariants suppl\'ementaires \`a un certain polyn\^ome entre les 11
bi-invariants que nous connaissons d\'ej\`a conduisent \`a un \'echec. 
Toutefois, nous v\'erifierons dans un instant que
\[
X^{27}
=
M^8\,X^{19},
\]
et nous \'etablirons que les $11 + 5 = 16$ bi-invariants\,:
\[
\aligned
f_1',\ \ \ \
\Lambda^3,\ \ \ \ \
\Lambda_1^5,\ \ \ \ \
\Lambda_{1,1}^7,\ \ \ \ \
M^8,\ \ \ \ \
&
\Lambda_{1,1,1}^9,\ \ \ \ \
M_1^{10},\ \ \ \ \
N^{12},\ \ \ \ \ 
K_{1,1}^{12},\ \ \ \ \
H_1^{14},\ \ \ \ \
F_{1,1}^{16}
\\
&
\ \ \ \ \ \ \ \ \ \ \ \ \ \ \ \
X^{18},\ \ \ \ \ 
X^{19},\ \ \ \ \
X^{21},\ \ \ \ \ 
X^{23},\ \ \ \ \
X^{25}
\endaligned
\]
sont mutuellement ind\'ependants. De surcro\^{\i}t, nous allons constater que
de nouveaux bi-invariants ``fant\^omes'' doivent encore
n\'ecessairement appara\^{\i}tre.

\smallskip\noindent{\bf D\'eduction de relations impliquant le wronskien
et $X^{ 18}, \dots, X^{ 27}$.} Sans poser $f_1' = 0$, multiplions
l'\'equation ``$\overset{ 10}{ \equiv}$'' par $M^{10}$\,:
\[
0
\equiv
-7\,\underline{\Lambda_1^5}\,\Lambda_{1,1}^7\,\underline{M_1^{10}}
+
3\,\Lambda^3\,\Lambda_{1,1,1}^9\,M_1^{10}
-
f_1'f_1'\,M_1^{10}\,M_1^{10}.
\]
\'Eliminons le bin\^ome soulign\'e $\underline{ \Lambda_1^5 \, M_1^{10}}$
gr\^ace \`a ``$\overset{ 18}{ \equiv}$'' multipli\'ee par $\Lambda_{ 1,
1}^7$\,:
\[
0
\equiv
-\frac{56}{5}\,\Lambda_{1,1}^7\,\Lambda_{1,1}^7\,M^8
-
\frac{7}{5}\,f_1'\,\Lambda_{1,1}^7\,H_1^{14}
+
3\,\Lambda^3\,\Lambda_{1,1,1}^9\,M_1^{10}
-
f_1'f_1'\,M_1^{10}\,M_1^{10}.
\]
Enfin, si nous \'eliminons le bin\^ome soulign\'e $\underline{ \Lambda_{
1, 1}^7 \, M^8}$ gr\^ace \`a l'\'equation (imm\'ediate\-ment d\'eduite de
``$\overset{ 51}{ \equiv}$'' et de ``$\overset{ 18}{
\equiv}$'') suivante\,:
\[
0
\overset{\widetilde{51}}
{\equiv}
-\Lambda_{1,1}^7\,M^8
+
3\,\Lambda^3\,K_{1,1}^{12}
-
f_1'\,H_1^{14}
\] 
multipli\'ee par $\Lambda_{1, 1}^7$, nous obtenons une identit\'e\,:
\[
0
\equiv
\Lambda^3
\big(
-5\Lambda_{1,1,1}^9\,M_1^{10}
+
56\,\Lambda_{1,1}^7\,K_{1,1}^{12}
\big)
+
f_1'
\big(
{\textstyle{\frac{5}{3}}}\,f_1'\,M_1^{10}\,M_1^{10}
-
{\textstyle{\frac{49}{3}}}\,\Lambda_{1,1}^7\,H_1^{14}
\big)
\]
dans laquelle appara\^{\i}t $f_1' X^{ 18}$. 

En proc\'edant de mani\`ere analogue pour $X^{ 19}$, $X^{ 21}$, $X^{ 23}$,
$X^{ 25}$ et $X^{ 27}$, et en \'eliminant \`a la fin le facteur non
identiquement nul $f_1'$, nous obtenons de nouvelles syzygies
impliquant nos six nouveaux bi-invariants\,:
\[
\aligned
0
&
\overset{a}
{\equiv}
-
49\,\Lambda_{1,1}^7\,H_1^{14}
+
3\,\Lambda^3\,X^{18}
+
5\,f_1'\,M_1^{10}\,M_1^{10},
\\
0
&
\overset{b}
{\equiv}
5\,M_1^{10}\,N^{12}
-
56\,M^8\,H_1^{14}
+
3\,\Lambda^3\,X^{19},
\\
0
&
\overset{c}
{\equiv}
5\,N^{12}\,N^{12}
+
64\,M^8\,M^8\,M^8
+
3\,\Lambda^3\,X^{21},
\\
0
&
\overset{d}
{\equiv}
7\,N^{12}\,H_1^{14}
+
8\,M^8\,M^8\,M_1^{10}
+
3\,\Lambda^3\,X^{23},
\\
0
&
\overset{e}
{\equiv}
35\,N^{12}\,F_{1,1}^{16}
-
448\,M^8\,M^8\,K_{1,1}^{12}
+
40\,M^8\,M_1^{10}\,M_1^{10}
+
3\,\Lambda^3\,X^{25},
\\
0
&
\overset{f}
{\equiv}
5\,M^8\,M_1^{10}\,N^{12}
-
56\,M^8\,M^8\,H_1^{14}
+
3\,\Lambda^3\,X^{27}.
\endaligned
\]
Ce ne sont pas les seules syzygies suppl\'ementaires\,: par
exemple\,:
\[
\aligned
0
&
\equiv
\Lambda_{1,1,1}^9\,H_1^{14}
-
8\,\Lambda_{1,1}^7\,F_{1,1}^{16}
+
\Lambda_1^5\,X^{18},
\\
0
&
\equiv
M_1^{10}\,H_1^{14}
-
8\,M^8\,F_{1,1}^{16}
+
\Lambda_1^5\,X^{19},
\\
0
&
\equiv
N^{12}\,H_1^{14}
+
8\,M^8\,M^8\,M_1^{10}
+
\Lambda_1^5\,X^{21},
\\
0
&
\equiv
N^{12}\,F_{1,1}^{16}
+
M^8\,M_1^{10}\,M_1^{10}
+
\Lambda_1^5\,X^{23}.
\endaligned
\]
et d'autres encore peuvent \^etre form\'ees, que nous ne rechercherons pas ici.

\smallskip\noindent{\bf Restriction \`a $\{ f_1' = 0\}$.}
L'expression en fonction de $j^5 f$ de nos six nouveaux bi-invariants
se simplifie lorsqu'on pose
$f_1' = 0$\,:
\[
\small
\aligned
X^{18}\big\vert_0
&
=
f_1''f_1''f_1''\,
\Big(11025\,\Delta_0^{1,5}\,
\big[\Delta_0^{1,2}\big]^2
+ 55125\,\Delta_0^{2,4}\,
\big[\Delta_0^{1,2}\big]^2 - 55125\,\Delta_0^{1,4}\,
\Delta_0^{1,3}\,\Delta_0^{1,2} 
&
\\
& 
\ \ \ \ \ \ \ \ \ \ \ \ \ \ \ \ \ \ \ \ \ \ \
-
110250\,\Delta_0^{2,3}\,
\Delta_0^{1,3}\,\Delta_0^{1,2} 
+ 49000\,\big[\Delta_0^{1,3}\big]^3
\Big),
\endaligned
\]
et l'on a des expressions similaires pour $X^{ 19} \big \vert_0$, 
$X^{ 21} \big \vert_0$, $X^{ 23} \big \vert_0$,
$X^{ 25} \big \vert_0$ et $X^{ 27} \big \vert_0$.
En fait, gr\^ace \`a nos six nouvelles syzygies ``$\overset{ a}{
\equiv}$'', ``$\overset{ b}{ \equiv}$'', ``$\overset{ c}{
\equiv}$'', ``$\overset{ d}{ \equiv}$'', ``$\overset{ e}{
\equiv}$'' et ``$\overset{ f}{ \equiv}$'' faisant intervenir le
wronskien, nous pouvons exprimer ces restrictions en fonction
seulement des quatre bi-invariants restreints alg\'ebriquement
ind\'ependants que sont $\Lambda^3 \big\vert_0$, $\Lambda_1^5 \big
\vert_0$, $M^8 \big\vert_0$ et $N^{ 12} \big\vert_0$ (tout en
rappelant les expressions les expressions que nous connaissons d\'ej\`a)\,:
\[
\footnotesize
\aligned
&
\underline{\Lambda^3}\big\vert_0,
\ \ \ \ \ \ \ \ \ \ \ \ \ \ \ \ 
\underline{\Lambda_1^5},
\ \ \ \ \ \ \ \ \ \ \ \ \ \ \ \ 
\Lambda_{1,1}^7\big\vert_0
=
{\textstyle{\frac{5}{3}\,\frac{\Lambda_1^5\,\Lambda_1^5}{
\Lambda^3\ \ \ \ }}}\Big\vert_0,
\ \ \ \ \ \ \ \ \ \ \ \ \ \ \ \ 
\underline{M^8}\big\vert_0,
\\
&
\Lambda_{1,1,1}^9\big\vert_0
=
{\textstyle{
\frac{35}{9}\,
\frac{\Lambda_1^5\,\Lambda_1^5\,\Lambda_1^5}
{\Lambda^3\,\Lambda^3\ \ \ \ \, }}}\Big\vert_0,
\ \ \ \ \ \ \ \ \ \ \ \ \ \ \ \ 
M_1^{10}\big\vert_0
=
{\textstyle{
\frac{8}{3}\,\frac{M^8\,\Lambda_1^5}
{\Lambda^3\ \ \ \ \ \, }}}\Big\vert_0,
\ \ \ \ \ \ \ \ \ \ \ \ \ \ \ \ 
\underline{N^{12}}\big\vert_0,
\\
&
K_{1,1}^{12}\big\vert_0
=
{\textstyle{
\frac{5}{9}\,\frac{M^8\,\Lambda_1^5\,\Lambda_1^5}
{\Lambda^3\,\Lambda^3\ \ \ \ \, }}}\Big\vert_0,
\ \ \ \ \ \ \ \ \ \ \ \ \ \ \ \ 
H_1^{14}\big\vert_0
=
{\textstyle{
\frac{5}{3}\,\frac{N^{12}\,\Lambda_1^5}
{\Lambda^3\ \ \ \ \, }}}\Big\vert_0,
\ \ \ \ \ \ \ \ \ \ \ \ \ \ \ \ 
F_{1,1}^{16}\big\vert_0
=
{\textstyle{
\frac{35}{9}\,\frac{\Lambda_1^5\,\Lambda_1^5\,N^{12}}
{\Lambda^3\,\Lambda^3\ \ \ \ \,}}}\Big\vert_0,
\endaligned
\]
\[
\footnotesize
\aligned
X^{18}\big\vert_0
&
=
{\textstyle{
1225\,\frac{\Lambda_1^5\,\Lambda_1^5\,\Lambda_1^5\,N^{12}}
{\Lambda^3\,\Lambda^3\,\Lambda^3\ \ \ \ \ \ }
\Big\vert_0}},
\ \ \ \ \ \ \ \ \ \ \ \ \ \ \ \ 
X^{19}\big\vert_0
=
{\textstyle{
\frac{80}{3}\,
\frac{\Lambda_1^5\,M^8\,N^{12}}{\Lambda^3\,\Lambda^3\ \ \ \ \ \ \ }
}}\Big\vert_0,
\\
X^{21}\big\vert_0
&
=
{\textstyle{
-\frac{5}{3}\,
\frac{N^{12}\,N^{12}}{\Lambda^3\ \ \ \ \ \ \ }
}}\Big\vert_0
-
{\textstyle{
\frac{64}{3}\,
\frac{M^8\,M^8\,M^8}{\Lambda^3\ \ \ \ \ \ \ \ \ \ \ }
}}\Big\vert_0,
\ \ \ \ \ 
X^{23}\big\vert_0
=
-
{\textstyle{
\frac{35}{9}\,
\frac{\Lambda_1^5\,N^{12}\,N^{12}}{\Lambda^3\,\Lambda^3\ \ \ \ \ \ \ }
}}\Big\vert_0
-
{\textstyle{
\frac{64}{9}\,\frac{\Lambda_1^5\,M^8\,M^8\,M^8}
{\Lambda^3 \ \ \ \ \ \ \ \ \ \ \ \ \ \ \ }
}}\Big\vert_0,
\\
X^{25}\big\vert_0
&
=
-
{\textstyle{
1225\,
\frac{\Lambda_1^5\,\Lambda_1^5\,N^{12}\,N^{12}}
{\Lambda^3\,\Lambda^3\,\Lambda^3\ \ \ \ \ \ }
}}\Big\vert_0
-
{\textstyle{
\frac{320}{27}\,
\frac{\Lambda_1^5\,\Lambda_1^5\,M^8\,M^8\,M^8}
{\Lambda^3\,\Lambda^3\,\Lambda^3\ \ \ \ \ \ \ \ \ }
}}\Big\vert_0,
\ \ \ \ \ \ \ \ \ 
X^{27}\big\vert_0
=
{\textstyle{
\frac{80}{3}\,
\frac{\Lambda_1^5\,M^8\,M^8\,N^{12}}
{\Lambda^3\,\Lambda^3\,\ \ \ \ \ \ \ \ \ \ }
}}\Big\vert_0.
\endaligned
\]

\def\theassertion{\!}\begin{assertion} 
Le bi-invariant $X^{ 18}$ ne s'exprime pas comme un
certain polyn\^ome en les onze bi-invariants construits par crochets
$f_1'$, $\Lambda^3$, $\Lambda_1^5$, $\Lambda_{ 1, 1}^7$, $M^8$,
$\Lambda_{ 1, 1, 1}^9$, $M_1^{ 10}$, $N^{ 12}$, $K_{ 1, 1}^{ 12}$,
$H_1^{ 14}$ et $F_{ 1, 1}^{ 16}$.
\end{assertion}

\noindent{\em Preuve.}
Par l'absurde, supposons que
\[
\aligned
X^{18}
&
=
\sum\,\text{\sf coeff}\cdot
\big(f_1'\big)^a\,
\big(\Lambda^3\big)^b\,
\big(\Lambda_1^5\big)^c\,
\big(\Lambda_{1,1}^7\big)^d\,
\big(M^8\big)^e\,
\big(\Lambda_{1,1,1}^9\big)^f\,
\big(M_1^{10}\big)^g\,
\\
&
\ \ \ \ \ \ \ \ \ \ \ \ \ \ \ \ \ \ \ \ \ \ \ \ \ \ \ \ \ \ \ \ \ \ \ \ 
\ \ \ \ \ \ \ \ \ \ \ \ \ \ \ \ 
\big(N^{12}\big)^h\,
\big(K_{1,1}^{12}\big)^i\,
\big(H_1^{14}\big)^j\,
\big(F_{1,1}^{16}\big)^k,
\endaligned
\]
avec des exposants entiers $a$, $b$, $c$, $d$, $e$, $f$, $g$, 
$h$, $i$, $j$ et $k$ tous $\geqslant 0$,
et posons $f_1'=0$ pour en d\'eduire une relation de la forme\,:
\[
\aligned
{\textstyle{
\frac{(\Lambda^5)^3\,N^{12}}{(\Lambda^3)^3}}}
\Big\vert_0
&
=
\sum\,\text{\sf coeff}\cdot
\big(\Lambda^3\big)^b\,
\big(\Lambda_1^5\big)^c\,
\Big(
{\textstyle{
\frac{\Lambda_1^5\,\Lambda_1^5}{\Lambda^3\ \ \ \ }}}
\Big)^d\,
\big(M^8\big)^e\,
\Big(
{\textstyle{
\frac{\Lambda_1^5\,\Lambda_1^5\,\Lambda_1^5}{\Lambda^3\,\Lambda^3\ \ \ \ }}}
\Big)^f\,
\Big(
{\textstyle{
\frac{\Lambda_1^5\,M^8}{\Lambda^3\,\ \ \ \ }}}
\Big)^g\,
\\
&
\ \ \ \ \ \ \ \ \ \ \ \ \ \ \ \ \ \ \ \ \ \ \ \ \ \ \ \ \ \ \ \ \ \ 
\big(N^{12}\big)^h\,
\Big(
{\textstyle{
\frac{\Lambda_1^5\,\Lambda_1^5\,M^8}{\Lambda^3\,\Lambda^3\ \ \ \ }}}
\Big)^i\,
\Big(
{\textstyle{
\frac{\Lambda_1^5\,N^{12}}{\Lambda^3\ \ \ \ }}}
\Big)^j\,
\Big(
{\textstyle{
\frac{\Lambda_1^5\,\Lambda_1^5\,N^{12}}{\Lambda^3\,\Lambda^3\ \ \ \ }}}
\Big)^k
\Big\vert_0
\\
&
=
\sum\,\text{\sf coeff}\cdot
\big(\Lambda^3\big)^{b-d-2f-g-2i-j-2k}\,
\big(\Lambda_1^5\big)^{c+2d+3f+g+2i+j+2k}\,
\\
&
\ \ \ \ \ \ \ \ \ \ \ \ \ \ \ \ \ \ \ \ \ \ \ \ 
\big(M^8\big)^{e+g+i}\,
\big(N^{12}\big)^{h+j+k}.
\endaligned
\]
Si nous identifions alors les exposants des quatre quantit\'es
alg\'ebriquement ind\'epen\-dantes\,:
\[
\left\{
\aligned
3
&
=
-b+d+2f+g+2i+j+2k,
\\
3
&
=
c+2d+3f+g+2i+j+2k,
\\
0
&
=
e+g+i,
\\
1
&
=
h+j+k,
\endaligned\right.
\]
nous d\'eduisons $e = g = i = 0$ de la troisi\`eme ligne, puis $0 = c + b
+ d + f$ en soustrayant la premi\`ere de la seconde, d'o\`u $c = b = d = f
= 0$, ce qui fait que premi\`ere et quatri\`eme ligne se simplifient
comme\,:
\[
3=j+2k
\ \ \ \ \ \ \ \ \ \ 
\text{\rm et}
\ \ \ \ \ \ \ \ \ \
1=h+j+k, 
\]
d'o\`u $h = 0$, puis $k = 2$ et enfin $3 = j + 4$, ce qui
est impossible.
\hfill$\square$

\def\theassertion{\!}\begin{assertion} 
Le bi-invariant $X^{ 19}$ ne s'exprime pas comme
un certain polyn\^ome en les douze bi-invariants $f_1'$, $\Lambda^3$,
$\Lambda_1^5$, $\Lambda_{ 1, 1}^7$, $M^8$, $\Lambda_{ 1, 1, 1}^9$,
$M_1^{ 10}$, $N^{ 12}$, $K_{ 1, 1}^{ 12}$, $H_1^{ 14}$, $F_{ 1, 1}^{
16}$ $X^{ 18}$.
\end{assertion}

\noindent{\em Preuve.}
Le m\^eme raisonnement 
que pour l'assertion pr\'ec\'edente
teste l'existence d'une repr\'esentation restreinte
de la forme\,:
\[
\aligned
{\textstyle{
\frac{\Lambda^5\,M^8\,N^{12}}{\Lambda^3\,\Lambda^3\ \ \ \ \ }}}
\Big\vert_0
&
=
\sum\,\text{\sf coeff}\cdot
\big(\Lambda^3\big)^{b-d-2f-g-2i-j-2k-3l}\,
\big(\Lambda_1^5\big)^{c+2d+3f+g+2i+j+2k+3l}\,
\\
&
\ \ \ \ \ \ \ \ \ \ \ \ \ \ \ \ \ \ \ \ \ \ \ \ 
\big(M^8\big)^{e+g+i}\,
\big(N^{12}\big)^{h+j+k+l}
\Big\vert_0,
\endaligned
\]
laquelle conduit au syst\`eme suivant de quatre \'equations entre entiers
$\geqslant 0$\,:
\[
\left\{
\aligned
2
&
-b+d+2f+g+2i+j+2k+3l,
\\
1
&
=
c+2d+3f+g+2i+j+2k+3l,
\\
1
&
e+g+i,
\\
1
&
h+j+k,
\endaligned\right.
\]
et ici, la contradiction se voit imm\'ediatement en soustrayant la
premi\`ere \'equation de la seconde, ce qui nous donne l'\'equation
impossible $-1 = c + b + d + f$.
\hfill$\square$

\smallskip\noindent{\bf Suite.} Par des raisonnements 
similaires, on \'etablit\,\,---\,\,comme annonc\'e\,\,---\,\,qu'{\it il
est r\'eellement n\'ecessaire d'introduire les 5 bi-invariants
suppl\'ementaires{\rm \,:}
\[
X^{18},\ \ \ \ \
X^{19},\ \ \ \ \
X^{21},\ \ \ \ \ 
X^{23},\ \ \ \ \ 
X^{25},
\]
lesquels ne sont pas engendr\'es par crochets}.

\smallskip\noindent{\bf Poursuite du processus d'engendrement.}
Mais ce n'est pas tout\,: gr\^ace \`a la liste des valeurs que prennent nos
16 bi-invariants en $f_1 ' = 0$, nous constatons qu'il nous
faut introduire encore d'autres bi-invariants, notamment\,:
\[
\aligned
&
{\textstyle{
\frac{-7\,\Lambda_{1,1}^7\,F_{1,1}^{16}+\Lambda_1^5\,X^{18}}{f_1'}}},
\ \ \ \ \ \ \ \ \
{\textstyle{
\frac{-5\,\Lambda_{1,1,1}^9\,F_{1,1}^{16}+\Lambda_{1,1}^7\,X^{18}}{f_1'}}},
\ \ \ \ \ \ \ \ \
{\textstyle{
\frac{-49\,K_{1,1}^{12}\,H_1^{14}+M^8\,X^{18}}{f_1'}}},
\\
&
{\textstyle{
\frac{-56\,K_{1,1}^{12}\,F_{1,1}^{16}+M_1^{10}\,X^{18}}{f_1'}}},
\ \ \ \ \ \ \ \ \
{\textstyle{
\frac{-7\,H_1^{14}\,F_{1,1}^{16}+N^{12}\,X^{18}}{f_1'}}},
\ \ \ \ \ \ \ \ \
{\textstyle{
\frac{-5\,F_{1,1}^{16}\,F_{1,1}^{16}+H_1^{14}\,X^{18}}{f_1'}}},
\endaligned
\]
\[
\aligned
&
{\textstyle{
\frac{-48\,M^8\,N^{12}+\Lambda^3\,X^{19}}{f_1'}}},
\ \ \ \ \ \ \ \ \
{\textstyle{
\frac{-6\,M_1^{10}\,H_1^{14}+\Lambda_1^5\,X^{19}}{f_1'}}},
\ \ \ \ \ \ \ \ \
{\textstyle{
\frac{-48\,K_{1,1}^{12}\,H_1^{14}+\Lambda_{1,1}^7\,X^{19}}{f_1'}}},
\\
&
{\textstyle{
\frac{-48\,K_{1,1}^{12}\,F_{1,1}^{16}+\Lambda_{1,1,1}^9\,X^{19}}{f_1'}}},
\endaligned
\]
et ensuite, il nous faut encore soumettre chacune de ces expressions au test
de savoir si elle ne s'exprime pas polynomialement en fonction d'une
liste croissante de bi-invariants connus \`a l'\'etape pr\'ec\'edente.

\smallskip\noindent{\bf Question.} 
L'alg\`ebre $\mathcal{ DS }_2^5$ poss\`ede-t-elle une infinit\'e de
bi-invariants fondamentaux?

\section*{\S8.~Calculs de caract\'eristique d'Euler}

\noindent{\bf Surfaces alg\'ebriques complexes projectives dans $P_3 (
\C)$.} Soit $X \subset P_3 (\C)$ une surface alg\'ebrique complexe
projective lisse de degr\'e $d\geqslant 1$. D'apr\`es~\cite{ de1997},
lorsque $x$ parcourt $X$, la r\'eunion des fibres $\big( \mathcal{ DS}_{
2, m}^\kappa \big)_x$ que nous avons \'etudi\'ees d'un point de vue
purement alg\'ebrique en un point fix\'e, s'organise de mani\`ere coh\'erente
en un sous-fibr\'e $\mathcal{ DS}_{ 2, m}^\kappa T_X^*$ du fibr\'e
$J^\kappa ( \C, X)$ des jets d'ordre $\kappa$ d'applications
holomorphes de $\C$ \`a valeurs dans $X$. Nous renvoyons le lecteur \`a
\cite{ de1997, ro2007} pour de plus amples informations g\'eom\'etriques.
Pour fixer les id\'ees, choisissons maintenant $\kappa = 4$.

\`A chaque mon\^ome bi-invariant parmi les deux listes fournies par le
premier th\'eor\`eme, \`a savoir\,:
\[
(f_1')^a\,
\big(\Lambda^3\big)^b\,
\big(\Lambda_{1,1}^7\big)^d\,
\big(M^8\big)^e
\ \ \ \ \ \ \ \ \ \ \ \
\text{\rm ou}
\ \ \ \ \ \ \ \ \ \ \ \
\Lambda_1^5\,
(f_1')^a\,
\big(\Lambda^3\big)^b\,
\big(\Lambda_{1,1}^7\big)^d\,
\big(M^8\big)^e
\]
correspond alors le {\sl fibr\'e de Schur}\footnote{\, Pour obtenir les
deux entiers $l_1$ et $l_2$ de $\Gamma^{ ( l_1, \, l_2)} T_X^*$,
rappellons qu'il suffit de compter le nombre de fois qu'apparaissent
les indices ``$(\cdot)_1$ et ``$(\cdot)_2$'' dans chacun de ces
mon\^omes, sachant qu'ils apparaissent chacun exactement une fois dans tout
d\'eterminant $\Delta^{ \alpha, \beta}$. }\,:
\[
\Gamma^{(a+b+3d+2e,\,b+d+2e)}\,T_X^*
\ \ \ \ \ \ \ \ \ \ \ \
\text{\rm ou}
\ \ \ \ \ \ \ \ \ \ \ \
\Gamma^{(2+a+b+3d+2e,\,1+b+d+2e)}\,T_X^*,
\] 
de telle sorte que $\mathcal{ DS}_{ 2, m}^4 T_X^*$ est isomorphe \`a la
somme directe de ces deux familles de fibr\'es de Schur, o\`u le
quadruplet d'entiers positifs ou nuls $(a, b, d, e)$ prend toutes les
valeurs telles que $a + 3b + 7 d + 8 e = m$ pour la premi\`ere famille,
et o\`u il prend toutes les valeurs satisfaisant $ 5 + a + 3b + 7d + 8e
= m$ pour la deuxi\`eme famille.

\smallskip\noindent{\bf Retour sur le choix d'une base de Gr\"obner
pour les jets d'ordre 4.} Nous n'avons pas encore fait remarquer que
notre choix de $\underline{ \Lambda_1^5 \, \Lambda_1^5} \big\vert_{
f_1' = 0}$ comme mon\^ome de t\^ete dans l'unique syzygie restreinte qui
existe entre les cinq bi-invariants fondamentaux restreints $\Lambda^3
\big \vert_0$, $\Lambda_1^5 \big \vert_0$, $\Lambda_{1, 1}^7 \big
\vert_0$ et $M^8 \big \vert_0$ au niveau $\kappa = 4$ n'\'etait pas en
harmonie avec le choix d'ordre lexicographique que nous avions fait au
niveau $\kappa = 5$ et qui nous avait fournit une base de 21 syzygies
faisant appara\^{\i}tre un {\it triangle remarquable de mon\^omes de t\^ete}.
En effet, si nous voulions r\'etablir la coh\'erence entre les deux niveaux
$\kappa = 4$ et $\kappa = 5$, nous devrions, au niveau $\kappa = 4$,
choisir plut\^ot l'ordre purement lexicographique d\'eduit de l'ordre
suivant entre bi-invariants restreints\,:
\[
\Lambda^3\,>\,
\Lambda_1^5\,>\,
\Lambda_{1,1}^7\,>\,
M^8
\]
(nous sous-entendons ici la mention ``$(\cdot)_{ f_1' = 0}$''), ce qui
conduit \`a {\it changer de mon\^ome de t\^ete}\, dans l'unique syzygie
restreinte existante\,:
\[
0
\equiv
-5\Lambda_1^5\,\Lambda_1^5
+
3\,\underline{\Lambda^3\,\Lambda_{1,1}^7}\big\vert_{f_1'=0}.
\]

\def\thelemma{\!\!}\begin{lemma}
Avec ce nouveau choix d'ordre qui anticipe une harmonie avec le niveau
suivant $\kappa = 5$, tout polyn\^ome bi-invariant de poids $m$ dans
$\mathcal{ DS}_{ 2, m}^4$ s'\'ecrit de mani\`ere unique{\rm \,:}
\[
\small
\aligned
{\sf P}^{2\times{\rm inv}}\big(j^4f\big)
&
=
\mathcal{P}\big(f_1',\Lambda^3,\Lambda_1^5,M^8\big)
+
\Lambda_{1,1}^7\,
\mathcal{Q}\big(f_1',\Lambda_1^5,\Lambda_{1,1}^7,M^8\big)
\\
&
=
\sum_{a+3b+5c+8e=m}\,{\sf coeff}\cdot
(f_1')^a\,
\big(\Lambda^3)^b\,
\big(\Lambda_1^5\big)^c\,
\big(M^8\big)^e
+
\\
&
\ \ \ \ \ \ \ \ \ \
+
\sum_{7+a+5c+7d+8e=m}\,{\sf coeff}\cdot
\Lambda_{1,1}^7\,
(f_1')^a\,
\big(\Lambda_1^5\big)^c\,
\big(\Lambda_{1,1}^7\big)^d\,
\big(M^8\big)^e,
\endaligned
\]
et par cons\'equent, le fibr\'e de Demailly-Semple $\mathcal{ DS}_{ 2,
m}^4 T_X^*$ est isomorphe \`a la somme directe suivante de fibr\'es de
Schur{\rm \,:}
\[
\small
\aligned
\mathcal{DS}_{2,m}^4T_X^*
&
=
\bigoplus_{a+3b+5c+8e=m}\,
\Gamma^{(a+b+2c+2e,\,b+c+2e)}\,T_X^*
\\
&
\ \ \ \ \ \ \ \ \ \ \ \
\bigoplus_{7+a+5c+7d+8e=m}\,
\Gamma^{(3+a+2c+3d+2e,\,1+c+d+2e)}\,T_X^*.
\endaligned
\]
\end{lemma}

\noindent{\bf Caract\'eristique d'Euler des fibr\'es de Schur.}
Si ${\sf c}_i = {\sf c}_i ( T_X)$, $i = 1, 2$, d\'esignent les classes
de Chern de d'une surface complexe $X$, on a la formule suivante
(\cite{ ro2007})\,:
\[
\boxed{
\aligned
\chi\big(X,\,\Gamma^{(l_1,l_2)}\,T_X^*\big)
&
=
\frac{1}{6}\,{\sf c}_1^2\,
\big[
l_1^3-l_2^3
\big]
-
\frac{1}{6}\,{\sf c}_2\,
(l_1-l_2)^3
+
{\rm O}\big(\vert \lambda\vert^2\big)
\endaligned
}
\]
pour la caract\'eristique d'Euler du fibr\'e de Schur $\Gamma^{ ( l_1, l_2
)} T_X^*$. Mentionnons au passage que pour $X$ de dimension trois, la
formule devient\,:
\[
\small
\aligned
&
-\chi\big(
X,\,
\Gamma^{(l_1,l_2,l_3)}\,T_X^*
\big)
=
\\
&
=
\frac{-{\sf c}_3}{1!\,\,2!\,\,3!}\
\left\vert
\begin{array}{ccc}
l_1\, & l_2\, & l_3\,
\\
l_1^2\, & l_2^2\, & l_3^2\,
\\
l_1^3\, & l_2^3\, & l_3^3\,
\end{array}
\right\vert
+
\frac{-{\sf c}_1{\sf c}_2+{\sf c}_3}{0!\,\,2!\,\,4!}\
\left\vert
\begin{array}{ccc}
1\, & 1\, & 1\,
\\
l_1^2\, & l_2^2\, & l_3^2\,
\\
l_1^4\, & l_2^4\, & l_3^4\,
\end{array}
\right\vert
+
\\
&
\ \ \ \ \ \ \ \
+
\frac{-{\sf c}_1^3+2\,{\sf c}_1{\sf c}_2-{\sf c}_3}{0!\,\,1!\,\,5!}\
\left\vert
\begin{array}{ccc}
1\, & 1\, & 1\,
\\
l_1\, & l_2\, & l_3\,
\\
l_1^5\, & l_2^5\, & l_3^5\,
\end{array}
\right\vert
+
{\rm O}\big(\vert l\vert^5\big),
\endaligned
\]
puis pour $X$ de dimension quatre\,:
\[
\footnotesize
\aligned
&
\chi\big(
X,\,
\Gamma^{(l_1,l_2,l_3,l_4)}\,T_X^*
\big)
=
\frac{-{\sf c}_1^4+3\,{\sf c}_1^2{\sf c}_2-{\sf c}_2^2
-2\,{\sf c}_1{\sf c}_3+{\sf c}_4}{0!\,\,1!\,\,2!\,\,7!}\
\left\vert
\begin{array}{cccc}
1\, & 1\, & 1\, & 1\,
\\
l_1\, & l_2\, & l_3\, & l_4\,
\\
l_1^2\, & l_2^2\, & l_3^2\, & l_4^2\,
\\
l_1^7\, & l_2^7\, & l_3^7\, & l_4^7\,
\end{array}
\right\vert
+
\\
&
+
\frac{-{\sf c}_4}{1!\,\,2!\,\,3!\,\,4!}\
\left\vert
\begin{array}{cccc}
l_1\, & l_2\, & l_3\, & l_4\,
\\
l_1^2\, & l_2^2\, & l_3^2\, & l_4^2\,
\\
l_1^3\, & l_2^3\, & l_3^3\, & l_4^3\,
\\
l_1^4\, & l_2^4\, & l_3^4\, & l_4^4\,
\end{array}
\right\vert
+
\frac{-{\sf c}_1{\sf c}_3+{\sf c}_4}{0!\,\,2!\,\,3!\,\,5!}\
\left\vert
\begin{array}{cccc}
1\, & 1\, & 1\, & 1\,
\\
l_1^2\, & l_2^2\, & l_3^2\, & l_4^2\,
\\
l_1^3\, & l_2^3\, & l_3^3\, & l_4^3\,
\\
l_1^5\, & l_2^5\, & l_3^5\, & l_4^5\,
\end{array}
\right\vert
+
\\
&
+
\frac{-{\sf c}_1^2{\sf c}_2+{\sf c}_2^2+{\sf c}_1{\sf c}_3-{\sf c}_4}
{0!\,\,1!\,\,3!\,\,6!}\
\left\vert
\begin{array}{cccc}
1\, & 1\, & 1\, & 1\,
\\
l_1\, & l_2\, & l_3\, & l_4\,
\\
l_1^3\, & l_2^3\, & l_3^3\, & l_4^3\,
\\
l_1^6\, & l_2^6\, & l_3^6\, & l_4^6\,
\end{array}
\right\vert
+
\frac{{\sf c}_1{\sf c}_3-{\sf c}_2^2}
{0!\,\,1!\,\,4!\,\,5!}\
\left\vert
\begin{array}{cccc}
1\, & 1\, & 1\, & 1\,
\\
l_1\, & l_2\, & l_3\, & l_4\,
\\
l_1^4\, & l_2^4\, & l_3^4\, & l_4^4\,
\\
l_1^5\, & l_2^5\, & l_3^5\, & l_4^5\,
\end{array}
\right\vert
+
{\rm O}\big(\vert l\vert^9\big),
\endaligned
\]
et enfin, ajoutons qu'il ne serait pas difficile de fournir 
la formule g\'en\'erale, valable en dimension quelconque.

\smallskip\noindent{\bf Sommations de caract\'eristiques.}
En revenant \`a la dimension deux, nous d\'eduisons tout d'abord
trivialement du lemme pr\'ec\'edent la formule sommatoire suivante pour la
caract\'eristique d'Euler de la premi\`ere somme directe de fibr\'es de
Schur\,:
\[
\small
\aligned
\chi\bigg(
X,
\bigoplus_{a+3b+5c+8e=m}\,
&
\Gamma^{(a+b+2c+2e,\,b+c+2e)}\,T_X^*
\bigg)
\\
&
=
\sum_{a+3b+5c+8e=m}\,
\chi
\Big(
X,\,\,
\Gamma^{(a+b+2c+2e,\,b+c+2e)}\,T_X^*
\Big),
\endaligned
\]
et ensuite, gr\^ace \`a une table de calculs lin\'eaires destin\'ee
\`a \'eliminer l'exposant $a$\,:
\[
\small
\aligned
m
&
=
a+3b+5c+8e,
\\
l_1
&
=
a+b+2c+2e
\\
&
=
m-2b-3c-6e,
\\
l_2
&
=
\ \ \ \ \ \ \ \ \ \ 
b+c+2e,
\\
l_1-l_2
&
=
m-3b-4c-8e,
\endaligned
\]
ce qui nous permet de remplacer $\sum_{ a + 3b + 5c + 8e = m}$ par
$\sum_{ 3 b + 5c + 8e \leqslant m}$, nous pouvons calculer les deux
coefficients rationnels ${\tt A}_1 \in m^6\cdot \Q$ et ${\tt A}_2 \in
m^6 \cdot \Q$ qui apparaissent devant ${\sf c}_1^2$ et devant $-{\sf
c}_2$ lorsqu'on effectue la premi\`ere somme de caract\'eristiques, que nous
appellerons ``${\tt A}$'' (la seconde s'appellera ``${\tt B}$'')\,:
\[
\footnotesize
\aligned
{\sf A}_1
&
=
\frac{1}{6}\,
\int_0^{\frac{m}{8}}\,de\,
\int_0^{\frac{m-8e}{5}}\,dc\,
\int_0^{\frac{m-5c-8e}{3}}\,db\,
\Big[
(m-2b-3c-6e)^3
-
(b+c+2e)^3
\Big]
+
{\rm O}(m^5)
\\
&
=
\frac{937\,\,m^6}{28\,800\,000},
\\
{\sf A}_2
&
=
\frac{1}{6}\,
\int_0^{\frac{m}{8}}\,de\,
\int_0^{\frac{m-8e}{5}}\,dc\,
\int_0^{\frac{m-5c-8e}{3}}\,db\,
(m-3b-4c-8e)^3
+
{\rm O}(m^5)
\\
&
=
\frac{13\,\,m^6}{900\,000},
\endaligned
\]
o\`u nous avons utilis\'e le fait que les sommes de Riemann sont suffisamment
bien approxim\'ees par des int\'egrales si l'on s'int\'eresse seulement
au coefficient de $m^6$.

\smallskip

Ensuite, si nous proc\'edons de la m\^eme mani\`ere pour 
la deuxi\`eme somme directe de fibr\'es de Schur\,:
\[
\footnotesize
\aligned
\chi\bigg(
X,
&
\bigoplus_{7+a+5c+7d+8e=m}\,
\Gamma^{(3+a+2c+3d+2e,\,1+c+d+2e)}\,T_X^*
\bigg)
\\
&
=
\sum_{7+a+5c+7d+8e=m}\,
\chi
\Big(
X,\,\,
\Gamma^{(3+a+2c+3d+2e,\,1+c+d+2e)}\,T_X^*
\Big)
\\
&
=
\sum_{a+5c+7d+8e=m}\,
\chi
\Big(
X,\,\,
\Gamma^{(a+2c+3d+2e,\,c+d+2e)}\,T_X^*
\Big)
+
{\rm O}(m^5),
\endaligned
\]
en observant que d\'ecalage de 7 dans le poids $m$, ainsi que les deux
d\'ecalages de 3 et de 1 dans $l_1$ et dans $l_2$ ne contribuent en fait
qu'en ${\rm O} ( m^5)$ dans le r\'esultat final, 
ce qui nous permet de les n\'egliger,
nous pouvons dresser
une table de calculs \'el\'ementaires analogue \`a la pr\'ec\'edente\,:
\[
\small
\aligned
m
&
=
a+5c+7d+8e,
\\
l_1
&
=
a+2c+3d+2e
\\
&
=
m-3c-4d-6e,
\\
l_2
&
=
\ \ \ \ \ \ \ \ \ \
c+d+2e,
\\
l_1-l_2
&
=
m-4c-5d-8e,
\endaligned
\]
remplacer $\sum_{ a + 5c + 7d + 8e = m}$ par $\sum_{ 5c + 7d + 8e
\leqslant m}$ apr\`es avoir \'elimin\'e $a$, de telle sorte que nous sommes
ramen\'es \`a calculer les deux int\'egrales suivantes\,:
\[
\footnotesize
\aligned
{\sf B}_1
&
=
\int_0^{\frac{m}{8}}\,de\,
\int_0^{\frac{m-8e}{7}}\,dd\,
\int_0^{\frac{m-7d-8e}{5}}\,dc\,
\Big[
(m-3c-4d-6e)^3
-
(c+d+2e)^3
\Big]
+
{\rm O}(m^5)
\\
&
=
\frac{559\,819\,\,m^6}{34\,574\,400\,000}
+
{\rm O}(m^5),
\\
{\sf B}_2
&
=
\int_0^{\frac{m}{8}}\,de\,
\int_0^{\frac{m-8e}{7}}\,dd\,
\int_0^{\frac{m-7d-8e}{5}}\,dc\,
(m-4c-5d-8e)^3
+
{\rm O}(m^5)
\\
&
=
\frac{36949\,\,m^6}{4\,321\,800\,000}
+
{\rm O}(m^5).
\endaligned
\]
En d\'efinitive, nous obtenons les deux coefficients rationnels totaux
$\mathcal{ C}_1 \in m^6 \cdot \Q$ et $\mathcal{ C}_2 \in m^6 \Q$ de
${\sf c}_1^2$ et de $- {\sf c}_2$\,:
\[
\aligned
\mathcal{C}_1
&
=
{\sf A}_1+{\sf B}_1
=
\frac{1797\,\,m^6}{36\,879\,360}
+
{\rm O}(m^5),
\\
\mathcal{C}_2
&
=
{\sf A}_2+{\sf B}_2
=
\frac{848\,\,m^6}{36\,879\,360}
+
{\rm O}(m^5).
\endaligned
\]

\noindent{\bf Application.} 
Les classes de Chern ${\sf c}_i = {\sf c}_i ( T_X)$, $i = 1, 2$, de
$X$ sont reli\'ees au degr\'e $d$ de $X$ par ${\sf c}_1^2 = (4-d)^2d$ et
${\sf c}_2 = d\, ( d^2 - 4 d + 6)$.

\def\theproposition{\!}\begin{proposition}
En dimension deux pour les jets d'ordre quatre,
la caract\'eristique d'Euler de $\mathcal{ DS}_{ 2, m}^4 
T_X^*$ vaut{\rm \,:}
\[
\chi
\big(
X,\mathcal{DS}_{2,m}^4T_X^*
\big)
=
\frac{m^6}{36\,879\,360}
\big(
1\,797\,{\sf c}_1^2
-
848\,{\sf c}_2
\big)
+
{\rm O}(m^5),
\]
donc si l'on pose{\rm \,:}
\[
\mathcal{C}
:=
{\textstyle{
\frac{1\,797}{848}}}
=
2,119\cdots,
\]
alors en r\'eexprimant le tout en fonction du degr\'e, 
on a l'\'equivalence{\rm \,:}
\[
\chi
\big(
X,\mathcal{DS}_{2,m}^4T_X^*
\big)
\sim
\frac{1\,797\,m^6}{36\,879\,360}\,d\,
{\sf q}_{\mathcal{C}}(d),
\]
quand $m \to \infty$, avec un polyn\^ome quadratique
\[
\boxed{
{\sf q}_{ \mathcal{ C}} (d)
:=
d^2(\mathcal{C}-1)
-
d(8\,\mathcal{C}-4)
+
16\,\mathcal{C}-6
}
\]
qui est positif pour tout degr\'e $d \geqslant {\bf 9}$.
\end{proposition}

\noindent{\bf Remarque.}
Le quotient $\mathcal{ C}$ des coefficients de ${\sf c}_1^2$ et de
${\sf c}_2$ vaut $\frac{ 47}{ 26} = 1,807 \cdots$ pour les jets
d'ordre 3, d'o\`u ${\sf q}_\mathcal{ C} (d)$ est positif pour tout $d
\geqslant {\bf 11}$ ({\it cf.} \cite{ ro2007}), et il vaut $\frac{ 13}{ 9} =
1,444 \cdots$ pour les jets d'ordre 2 ({\it cf.} \cite{ de1997}), d'o\`u 
${\sf q}_\mathcal{ C} (d)$ est positif pour tout $d
\geqslant {\bf 15}$. Demailly a conjectur\'e que ce quotient tend vers
l'infini avec $\kappa$. Les valeurs num\'eriques $1,44 \cdots$, $1,80
\cdots$ et $2, 12 \cdots$ sugg\`erent une certaine lenteur de la
convergence potentielle.

\begin{corollary}
Si $A$ est un fibr\'e en droites ample sur $X$, pour tout $m$
suffisamment grand, il y a des sections globales de $\mathcal{ DS}_{
2, m}^4 T_X^* \otimes A^{ -1}$ lorsque $d \geqslant {\bf 9}$, et toute courbe
enti\`ere $f = \C \to X$ doit satisfaire l'\'equation diff\'erentielle
globale correspondante.
\end{corollary}

\noindent{\bf Passage aux jets d'ordre 5.}
Maintenant, la correspondance entre bi-invariants fondamentaux et
repr\'esentations de Schur\,:
\[
\aligned
&
f_1'
\longleftrightarrow
\Gamma^{(1,0)},
\ \ \ \ \ \ \ \ \ \
\Lambda^3
\longleftrightarrow
\Gamma^{(1,1)},
\\
\Lambda_1^5
\longleftrightarrow
\Gamma^{(2,1)},&
\ \ \ \ \ \ \ \ \ \
\Lambda_{1,1}^7
\longleftrightarrow
\Gamma^{(3,1)},
\ \ \ \ \ \ \ \ \ \
M^8
\longleftrightarrow
\Gamma^{(2,2)},
\\
\Lambda_{1,1,1}^9
\longleftrightarrow
\Gamma^{(4,1)},&
\ \ \ \ \ \ \ \ \ \
M_1^{10}
\longleftrightarrow
\Gamma^{(3,2)},
\ \ \ \ \ \ \ \ \ \
N^{12}
\longleftrightarrow
\Gamma^{(3,3)},
\\
K_{1,1}^{12}
\longleftrightarrow
\Gamma^{(4,2)},&
\ \ \ \ \ \ \ \ \ \
H_1^{14}
\longleftrightarrow
\Gamma^{(4,3)},
\ \ \ \ \ \ \ \ \ \
F_{1,1}^{16}
\longleftrightarrow
\Gamma^{(5,3)},
\endaligned
\]
\'Etant donn\'e qu'il existe des invariants fondamentaux
suppl\'ementaires, nous pourrions attendre encore avant d'entreprendre
un calcul de Riemann-Roch, mais nous avons quand m\^eme l'opportunit\'e de
nous restreindre \`a la sous-alg\`ebre engendr\'ee par les crochets.

\smallskip\noindent{\bf Base de Gr\"obner.}
En choisissant l'ordre purement lexicographique sur les mon\^omes de
$\C\big[ \Lambda^3, \dots, H_1^{ 14}, F_{ 1, 1}^{ 16}, f_1' \big]$ qui
est d\'eduit de l'ordre suivant sur les mon\^omes \'el\'emen\-taires
restreints (noter que $f_1'$ est plac\'e en derni\`ere
position\footnote{\, Sinon, les bases fournies contiennent plus d'une
soixantaine d'\'equations. })\,:
\[
\Lambda^3
>
\Lambda_1^5
>
\Lambda_{1,1}^7
>
M^8
>
\Lambda_{1,1,1}^9
>
M_1^{10}
>
N^{12}
>
K_{1,1}^{12}
>
H_1^{14}
>
F_{1,1}^{16}
>
f_1'
\]
Maple nous donne la base de Gr\"obner r\'eduite suivante pour l'id\'eal
complet des syzygies entre nos onze bi-invariants restreints \`a $\{ f_1 ' =
0 \}$, laquelle est constitu\'ee de 26 \'equations\,:
\[
\scriptsize
\aligned
0
&
=
-5\,f_1'\big(M_1^{10}\big)^2N^{12}\,K_{1,1}^{12}\,F_{1,1}^{16}
+
5\,f_1'\,N^{12}\big(F_{1,1}^{16}\big)^3
-
64\,f_1'\big(M_1^{10}\big)^2\,K_{1,1}^{12}\big(H_1^{14}\big)^2
+
5\,\underline{f_1'\big(M_1^{10}\big)^3H_1^{14}\,F_{1,1}^{16}}
+
\\
&
\ \ \ \ \ 
+
128\,f_1'\,M_1^{10}\,N^{12}\big(K_{1,1}^{12}\big)^2H_1^{14}
-
7\,f_1'\big(H_1^{14}\big)^2\big(F_{1,1}^{16}\big)^2
-
64\,f_1'\big(N^{12}\big)^2\big(K_{1,1}^{12}\big)^3,
\\
0
&
=
15\,f_1'\,M_1^{10}\,K_{1,1}^{12}\,H_1^{14}
-
7\,\Lambda_{1,1,1}^9\big(H_1^{14}\big)^2
-
5\,f_1'\big(M_1^{10}\big)^2F_{1,1}^{16}
+
5\,\underline{\Lambda_{1,1,1}^9\,N^{12}\,F_{1,1}^{16}}
-
8\,f_1'\,N^{12}\big(K_{1,1}^{12}\big)^2,
\\
0
&
=
7\,f_1'\big(M_1^{10}\big)^2K_{1,1}^{12}
+
\Lambda_{1,1,1}^9\,M_1^{10}\,H_1^{14}
+
f_1'\big(F_{1,1}^{16}\big)^2
-
8\,\Lambda_{1,1,1}^9\,N^{12}\,K_{1,1}^{12},
\\
0
&
=
N^{12}\,K_{1,1}^{12}
-
M_1^{10}\,H_1^{14}
+
\underline{M^8\,F_{1,1}^{16}},
\\
0
&
=
64\,\underline{f_1'\,M^8\,M_1^{10}\,K_{1,1}^{12}\,H_1^{14}}
+
7\,f_1'\big(H_1^{14}\big)^2F_{1,1}^{16}
+
5\,f_1'\big(M_1^{10}\big)^2N^{12}\,K_{1,1}^{12}
-
64\,f_1'\,M^8\,N^{12}\,K_{1,1}^{12}
-
\\
&
\ \ \ \ \
-
5\,f_1'\big(M_1^{10}\big)^3H_1^{14}
-
5\,f_1'\,N^{12}\big(F_{1,1}^{16}\big)^2,
\\
0
&
=
f_1'\,H_1^{14}\,F_{1,1}^{16}
+
8\,\underline{M^8\,\Lambda_{1,1,1}^9\,H_1^{14}}
+
5\,f_1'\big(M_1^{10}\big)^3
-
5\,\Lambda_{1,1,1}^9\,M_1^{10}\,N^{12}
-
8\,f_1'\,M^8\,M_1^{10}\,K_{1,1}^{12},
\\
0
&
=
64\,\underline{M^8\,\Lambda_{1,1,1}^9\,N^{12}\,K_{1,1}^{12}}
-
64\,f_1'\,M^8\big(M_1^{10}\big)^2K_{1,1}^{12}
-
5\,\Lambda_{1,1,1}^9\big(M_1^{10}\big)^2N^{12}
+
5\,f_1'\big(M_1^{10}\big)^4
-
\\
&
\ \ \ \ \
-
7\,f_1'\,M_1^{10}\,H_1^{14}\,F_{1,1}^{16}
+
8\,f_1'\,N^{12}\,K_{1,1}^{12}\,F_{1,1}^{16},
\\
\endaligned
\]
\[
\scriptsize
\aligned
0
&
=
7\,f_1'\big(H_1^{14}\big)^2
-
5\,f_1'\,M^8\big(M_1^{10}\big)^2
+
64\,\underline{f_1'\big(M^8\big)^2K_{1,1}^{12}}
-
5\,f_1'\,N^{12}\,F_{1,1}^{16},
\\
0
&
=
64\underline{\big(M^8\big)^2\Lambda_{1,1,1}^9\,K_{1,1}^{12}}
-
5\,M^8\,\Lambda_{1,1,1}^9\big(M_1^{10}\big)^2
-
5\,f_1'\big(M_1^{10}\big)^2F_{1,1}^{16}
+
15\,f_1'\,M_1^{10}\,K_{1,1}^{12}\,H_1^{14}
-
8\,f_1'\,N^{12}\big(K_{1,1}^{12}\big)^2,
\\
0
&
=
-\Lambda_{1,1,1}^9\,H_1^{14}
+
\underline{\Lambda_{1,1}^7\,F_{1,1}^{16}}
+
f_1'\,M_1^{10}\,K_{1,1}^{12},
\\
0
&
=
5\,f_1'\big(M_1^{10}\big)^2
-
5\,\Lambda_{1,1,1}^9\,N^{12}
-
8\,f_1'\,M^8\,K_{1,1}^{12}
+
7\,\underline{\Lambda_{1,1}^7\,H_1^{14}},
\\
0
&
=
56\,\underline{\Lambda_{1,1}^7\,N^{12}\,K_{1,1}^{12}}
-
64\,f_1'\,M^8\,M_1^{10}\,K_{1,1}^{12}
-
5\,\Lambda_{1,1,1}^9\,M_1^{10}\,N^{12}
+
5\,f_1'\big(M_1^{10}\big)^3
-
7\,f_1'\,H_1^{14}\,F_{1,1}^{16},
\\
0
&
=
-8\,M^8\,\Lambda_{1,1,1}^9
+
7\,\underline{\Lambda_{1,1}^7\,M_1^{10}}
-
f_1'\,F_{1,1}^{16},
\\
0
&
=
56\,\underline{\Lambda_{1,1}^7\,M^8\,K_{1,1}^{12}}
-
5\,M^8\,\Lambda_{1,1,1}^9\,M_1^{10}
-
5\,f_1'\,M_1^{10}\,F_{1,1}^{16}
+
7\,f_1'\,K_{1,1}^{12}\,H_1^{14},
\\
0
&
=
49\underline{\big(\Lambda_{1,1}^7\big)^2K_{1,1}^{12}}
-
5\,M^8\big(\Lambda_{1,1,1}^9\big)^2
-
5\,f_1'\,\Lambda_{1,1,1}^9\,F_{1,1}^{16}
+
7\big(f_1'\big)^2\big(K_{1,1}^{12}\big)^2,
\\
0
&
=
-\Lambda_{1,1,1}^9\,N^{12}
+
\underline{\Lambda_1^5\,F_{1,1}^{16}}
+
f_1'\,M_1^{10}\,M_1^{10},
\\
\endaligned
\]
\[
\scriptsize
\aligned
0
&
=
-\Lambda_{1,1}^7\,N^{12}
+
\underline{\Lambda_1^5\,H_1^{14}}
+
f_1'\,M^8\,M_1^{10},
\\
0
&
=
7\,\underline{\Lambda_1^5\,K_{1,1}^{12}}
-
M^8\,\Lambda_{1,1,1}^9
-
f_1'\,F_{1,1}^{16},
\\
0
&
=
-8\,\Lambda_{1,1}^7\,M^8
+
5\,\underline{\Lambda_1^5\,M_1^{10}}
-
f_1'\,H_1^{14},
\\
0
&
=
-7\,\Lambda_{1,1}^7\,\Lambda_{1,1}^7
+
5\,\underline{\Lambda_1^5\,\Lambda_{1,1,1}^9}
-
f_1'f_1'\,K_{1,1}^{12},
\\
0
&
=
-7\,\Lambda_{1,1}^7\,N^{12}
+
3\,\underline{\Lambda^3\,F_{1,1}^{16}}
+
8\,f_1'\,M^8\,M_1^{10},
\\
0
&
=
-5\,\Lambda_1^5\,N^{12}
+
3\,\underline{\Lambda^3\,H_1^{14}}
+
8\,f_1'\,M^8\,M^8,
\\
0
&
=
3\,\underline{\Lambda^3\,K_{1,1}^{12}}
-
\Lambda_{1,1}^7\,M^8
-
f_1'\,H_1^{14},
\\
0
&
=
-8\,\Lambda_1^5\,M^8
+
3\,\underline{\Lambda^3\,M_1^{10}}
-
f_1'\,N^{12},
\\
0
&
=
-7\,\Lambda_1^5\,\Lambda_{1,1}^7
+
3\,\underline{\Lambda^3\,\Lambda_{1,1,1}^9}
-
f_1'f_1'\,M_1^{10},
\\
0
&
=
-5\big(\Lambda_1^5\big)^2
+
3\,\underline{\Lambda^3\,\Lambda_{1,1}^7}
-
f_1'f_1'\,M^8.
\endaligned
\]
D'apr\`es la th\'eorie des bases de Gr\"obner, une base de l'espace
vectoriel\,:
\[
\C\big[
\Lambda^3,\,\Lambda_1^5,\,\dots,
H_1^{14},\,F_{1,1}^{16},\,f_1'\big]
\Big/
\big(
\text{\rm 26 \'equations pr\'ec\'edentes}
\big)
\]
est constitu\'ee de tous les mon\^omes
\[
\big(\Lambda^3\big)^a
\big(\Lambda_1^5\big)^b
\big(\Lambda_{1,1}^7\big)^c
\big(M^8\big)^d
\big(\Lambda_{1,1,1}^9\big)^e
\big(M_1^{10}\big)^f
\big(N^{12}\big)^g
\big(K_{1,1}^{12}\big)^h
\big(H_1^{14}\big)^i
\big(F_{1,1}^{16}\big)^j
\big(f_1'\big)^k
\]
qui n'appartiennent pas \`a l'id\'eal monomial engendr\'e par les 26
mon\^omes de t\^ete de chacun des 26 g\'en\'erateurs que nous avons
soulign\'es. Or un tel mon\^ome appartient \`a cet id\'eal monomial si
et seulement si il est divisible par l'un des 26 mon\^omes de t\^ete,
ce qui revient \`a dire que le muti-indice
\[
\big(a,\,b,\,c,\,d,\,e,\,f,\,g,\,h,\,i,\,j,\,k\big)
\in
\N^{11}
\]
appartient \`a la {\it r\'eunion}\, des 26 sous-ensembles suivants
de $\N^{11}$\,:
\[
\scriptsize
\aligned
\{f\geqslant 3\}\cap\{i\geqslant 1\}\cap\{j\geqslant 1\}\cap\{k\geqslant 1\},
\ \ \ \ \ \ \ \ \ \ \ \ \
\{e\geqslant 1\}\cap\{g\geqslant 1\}\cap\{j\geqslant 1\},
\\
\{e\geqslant 1\}\cap\{f\geqslant 1\}\cap\{i\geqslant 1\},
\\
\endaligned
\]
\[
\scriptsize
\aligned
&
\{d\geqslant 1\}\cap\{j\geqslant 1\},
\\
&
\{d\geqslant 1\}\cap\{f\geqslant 1\}\cap\{h\geqslant 1\}\cap\{i\geqslant 1\}
\cap\{k\geqslant 1\},
\\
&
\{d\geqslant 1\}\cap\{e\geqslant 1\}\cap\{i\geqslant 1\},
\\
&
\{d\geqslant 1\}\cap\{e\geqslant 1\}\cap\{g\geqslant 1\}\cap\{h\geqslant 1\},
\\
&
\{d\geqslant 2\}\cap\{h\geqslant 1\}\cap\{k\geqslant 1\},
\\
&
\{d\geqslant 2\}\cap\{e\geqslant 1\}\cap\{h\geqslant 1\},
\endaligned
\]
\[
\scriptsize
\aligned
\{c\geqslant 1\}\cap\{j\geqslant 1\},
\ \ \ \ \ \ \ \
\{b\geqslant 1\}\cap\{j\geqslant 1\},
\ \ \ \ \ \ \ \
&
\{a\geqslant 1\}\cap\{j\geqslant 1\},
\\
\{c\geqslant 1\}\cap\{i\geqslant 1\},
\ \ \ \ \ \ \ \
\{b\geqslant 1\}\cap\{i\geqslant 1\},
\ \ \ \ \ \ \ \
&
\{a\geqslant 1\}\cap\{i\geqslant 1\},
\\
\{c\geqslant 1\}\cap\{g\geqslant 1\}\cap\{h\geqslant 1\},
\ \ \ \ \ \ \ \
\{b\geqslant 1\}\cap\{h\geqslant 1\},
\ \ \ \ \ \ \ \
&
\{a\geqslant 1\}\cap\{h\geqslant 1\},
\\
\{c\geqslant 1\}\cap\{f\geqslant 1\},
\ \ \ \ \ \ \ \
\{b\geqslant 1\}\cap\{f\geqslant 1\},
\ \ \ \ \ \ \ \
&
\{a\geqslant 1\}\cap\{f\geqslant 1\},
\\
\{c\geqslant 1\}\cap\{d\geqslant 1\}\cap\{h\geqslant 1\},
\ \ \ \ \ \ \ \
\{b\geqslant 1\}\cap\{e\geqslant 1\},
\ \ \ \ \ \ \ \
&
\{a\geqslant 1\}\cap\{e\geqslant 1\},
\\
\{c\geqslant 2\}\cap\{h\geqslant 1\},
\ \ \ \ \ \ \ \ \ \ \ \ \ \ \ \ \ \ \ \ \ \ \ \
\ \ \ \ \ \ \ \ \ \ \ \ \ \ \ \ \ \ \ \ \ \ \ \
&
\{a\geqslant 1\}\cap\{c\geqslant 1\}.
\endaligned
\]
Pour calculer le compl\'ementaire de cette r\'eunion, on proc\`ede comme
dans la Section~7, en regroupant s\'epar\'ement les 6 intersections
commen\c cant par $\{ a \geqslant 1\}$, puis les 5 commen\c cant par $\{
b \geqslant 1\}$, puis les 6 commen\c cant par $\{ c \geqslant 1\}$,
puis les 6 commen\c cant par $\{ d\geqslant 1\}$, puis les 2 comman\c cant
par $\{ e \geqslant 1\}$, et puis enfin la derni\`ere,
qui commence par $\{ f \geqslant
1\}$. Trouver le compl\'ementaire global reviendra donc \`a calculer
l'{\it intersection} de six sous-ensembles de $\N^{ 11}$.

\noindent
Clairement, le premier et le deuxi\`eme compl\'ementaires sont
donn\'es par\,:
\[
\small
\aligned
\mathcal{N}_1
&
:=
\{a=0\}\cup\{c=e=f=g=h=i=j=0\},
\\
\mathcal{N}_2
&
:=
\{b=0\}\cup\{e=f=h=i=j=0\}.
\endaligned
\]

\noindent
Ensuite, calculons le troisi\`eme compl\'ementaire, en simplifiant
progressivement les intersections, et ce, en partant du dernier terme\,:
\[
\footnotesize
\aligned
\mathcal{N}_3
:=
&\
\{c=1\}\cup\{c=0\}\cup\{h=0\}
\bigcap
\{c=0\}\cup\{d=0\}\cup\{h=0\}
\bigcap
\{c=0\}\cup\{f=0\}
\\
&
\ \ \ \ \ \ \ \
\bigcap
\{c=0\}\cup\{g=0\}\cup\{h=0\}
\bigcap
\{c=0\}\cup\{i=0\}
\bigcap
\{c=0\}\cup\{j=0\}
\\
=
&\
\{c=1\}\cup\{c=0\}\cup\{h=0\}
\bigcap
\{c=0\}\cup\{d=0\}\cup\{h=0\}
\bigcap
\{c=0\}\cup\{f=0\}
\\
&
\ \ \ \ \ \ \ \
\bigcap
\{c=0\}\cup\{c=g=i=j=0\}\cup\{h=i=j=0\}
\\
=
&\
\{c=0\}\cup\{f=h=i=j=0\}\cup\{c=1,\,d=f=g=i=j=0\}.
\endaligned
\]
Les calculs suivants donnent\,:
\[
\footnotesize
\aligned
\mathcal{N}_4
:=
&\
\{d=0\}\cup\{e=h=j=0\}\cup\{e=j=k=0\}\cup\{h=i=j=0\}\cup
\\
&
\ \ \ \ \ \ \ \ \ \ \ \ \
\cup\{d=1,\,e=i=j=0\}\cup\{d=1,\,e=f=j=0\}\cup\{d=1,\,g=i=j=0\}\cup
\\
&
\ \ \ \ \ \ \ \ \ \ \ \ \
\cup\{d=1,\,e=h=j=0\}\cup\{d=1,\,h=i=j=0\}\cup\{d=1,\,e=j=k=0\},
\\
\mathcal{N}_5
:=
&\
\{e=0\}\cup\{f=g=0\}\cup\{f=j=0\}\cup\{g=i=0\}\cup\{i=j=0\},
\\
\mathcal{N}_6
:=
&\
\{f=2\}\cup\{f=1\}\cup\{f=0\}\cup\{i=0\}\cup\{j=0\}\cup\{k=0\}.
\endaligned
\]
En d\'eveloppant l'intersection finale\,:
\[
\mathcal{N}_1\cap\mathcal{N}_2\cap\mathcal{N}_3\cap
\mathcal{N}_4\cap\mathcal{N}_5\cap\mathcal{N}_6,
\]
nous pouvons n\'egliger toutes les composantes qui incorporent un nombre
$\geqslant 7$ d'\'equations, puisque dans la sommation de fibr\'es de
Schur, la contribution ne sera qu'en ${\rm O} ( m^7)$, les termes
principaux \'etant multiples rationnels non nuls de $m^8$. Ainsi, en
n\'egligeant de tels termes, nous obtenons exactement 16 composantes de
dimension 5 d\'efinies par 6 \'equations\,:
\[
\aligned
\small
&
\{a=b=c=d=e=0,\,f=2\}
\cup
\{a=b=c=d=e=0,\,f=1\}
\\
&
\cup
\{a=b=c=d=e=f=0\}
\cup
\{a=b=c=d=e=i=0\}
\\
&
\cup
\{a=b=c=d=e=j=0\}
\cup
\{a=b=c=d=e=k=0\}
\\
&
\cup
\{a=b=c=d=f=g=0\}
\cup
\{a=b=c=d=f=i=0\}
\\
&
\cup
\{a=b=c=d=g=i=0\}
\cup
\{a=b=c=d=i=j=0\}
\\
&
\cup
\{a=b=c=e=h=j=0\}
\cup
\{a=b=c=e=j=k=0\}
\\
&
\cup
\{a=b=c=h=i=j=0\}
\cup
\{a=b=f=h=i=j=0\}
\\
&
\cup
\{a=e=f=h=i=j=0\}
\cup
\{c=e=f=h=i=j=0\}.
\endaligned
\]
Les 16 familles de mon\^omes correspondant \`a ces \'equations 
peuvent \^etre rang\'ees dans un tableau\,:
\[
\tiny
\begin{array}{ccccccccccccc}
{\sf A}\,:\ \ \ \ \ \
&
\bullet&\bullet&\bullet&\bullet&\bullet&\big(M_1^{10}\big)^2&
\big(N^{12}\big)^g&\big(K_{1,1}^{12}\big)^h&\big(H_1^{14}\big)^i&
\big(F_{1,1}^{16}\big)^j&
(f_1')^k
\\
{\sf B}\,:\ \ \ \ \ \
&
\bullet&\bullet&\bullet&\bullet&\bullet&M_1^{10}&
\big(N^{12}\big)^g&\big(K_{1,1}^{12}\big)^h&\big(H_1^{14}\big)^i&
\big(F_{1,1}^{16}\big)^j&
(f_1')^k
\\
{\sf C}\,:\ \ \ \ \ \ 
&
\bullet&\bullet&\bullet&\bullet&\bullet&\bullet&
\big(N^{12}\big)^g&\big(K_{1,1}^{12}\big)^h&\big(H_1^{14}\big)^i&
\big(F_{1,1}^{16}\big)^j&
(f_1')^k
\\
{\sf D}\,:\ \ \ \ \ \
&
\bullet&\bullet&\bullet&\bullet&\bullet&\big(M_1^{10}\big)^f&
\big(N^{12}\big)^g&\big(K_{1,1}^{12}\big)^h&\bullet&
\big(F_{1,1}^{16}\big)^j&
(f_1')^k
\\
{\sf E}\,:\ \ \ \ \ \
&
\bullet&\bullet&\bullet&\bullet&\bullet&\big(M_1^{10}\big)^f&
\big(N^{12}\big)^g&\big(K_{1,1}^{12}\big)^h&\big(H_1^{14}\big)^i&\bullet&
(f_1')^k
\\
{\sf F}\,:\ \ \ \ \ \
&
\bullet&\bullet&\bullet&\bullet&\bullet&\big(M_1^{10}\big)^f&
\big(N^{12}\big)^g&\big(K_{1,1}^{12}\big)^h&\big(H_1^{14}\big)^i&
\big(F_{1,1}^{16}\big)^j&
\bullet
\\
{\sf G}\,:\ \ \ \ \ \
&
\bullet&\bullet&\bullet&\bullet&\big(\Lambda_{1,1,1}^9\big)^e&\bullet&
\bullet&\big(K_{1,1}^{12}\big)^h&\big(H_1^{14}\big)^i&
\big(F_{1,1}^{16}\big)^j&
(f_1')^k
\\
{\sf H}\,:\ \ \ \ \ \ 
&
\bullet&\bullet&\bullet&\bullet&\big(\Lambda_{1,1,1}^9\big)^e&\bullet&
\big(N^{12}\big)^g&\big(K_{1,1}^{12}\big)^h&\bullet&
\big(F_{1,1}^{16}\big)^j&
(f_1')^k
\\
{\sf I}\,:\ \ \ \ \ \
&
\bullet&\bullet&\bullet&\bullet&\big(\Lambda_{1,1,1}^9\big)^e&
\big(M_1^{10}\big)^f&
\bullet&\big(K_{1,1}^{12}\big)^h&\bullet&
\big(F_{1,1}^{16}\big)^j&
(f_1')^k
\\
{\sf J}\,:\ \ \ \ \ \
&
\bullet&\bullet&\bullet&\bullet&\big(\Lambda_{1,1,1}^9\big)^e&
\big(M_1^{10}\big)^f&
\big(N^{12}\big)^g&\big(K_{1,1}^{12}\big)^h&\bullet&\bullet&
(f_1')^k
\\
{\sf K}\,:\ \ \ \ \ \
&
\bullet&\bullet&\bullet&\big(M^8\big)^d&\bullet&
\big(M_1^{10}\big)^f&
\big(N^{12}\big)^g&\bullet&\big(H_1^{14}\big)^i&\bullet&
(f_1')^k
\\
{\sf L}\,:\ \ \ \ \ \
&
\bullet&\bullet&\bullet&\big(M^8\big)^d&\bullet&
\big(M_1^{10}\big)^f&
\big(N^{12}\big)^g&\big(K_{1,1}^{12}\big)^h&\big(H_1^{14}\big)^i&\bullet&
\bullet
\\
{\sf M}\,:\ \ \ \ \ \
&
\bullet&\bullet&\bullet&\big(M^8\big)^d&\big(\Lambda_{1,1,1}^9\big)^e&
\big(M_1^{10}\big)^f&
\big(N^{12}\big)^g&\bullet&\bullet&\bullet&
(f_1')^k
\\
{\sf N}\,:\ \ \ \ \ \
&
\bullet&\bullet&\big(\Lambda_{1,1}^7\big)^c&\big(M^8\big)^d&
\big(\Lambda_{1,1,1}^9\big)^e&
\bullet&
\big(N^{12}\big)^g&\bullet&\bullet&\bullet&
(f_1')^k
\\
{\sf O}\,:\ \ \ \ \ \
&
\bullet&\big(\Lambda_1^5\big)^b&\big(\Lambda_{1,1}^7\big)^c&\big(M^8\big)^d&
\bullet&
\bullet&
\big(N^{12}\big)^g&\bullet&\bullet&\bullet&
(f_1')^k
\\
{\sf P}\,:\ \ \ \ \ \
&
\big(\Lambda^3\big)^a&\big(\Lambda_1^5\big)^b&\bullet&\big(M^8\big)^d&
\bullet&
\bullet&
\big(N^{12}\big)^g&\bullet&\bullet&\bullet&
(f_1')^k
\end{array}
\]
Ensuite, lorsqu'on effectue la somme des caract\'eristiques d'Euler des
fibr\'es de Schur correspondants, il n'est pas n\'ecessaire de r\'eorganiser ces
familles de telle sorte qu'elles soient d'intersection vide, puisque de
toute fa\c con, chaque intersection entre deux familles ne contribuera au
final qu'en ${\rm O} ( m^7)$. Nous pouvons donc additionner les seize
couples d'int\'egrales correspondantes. Voici les deux premi\`eres, 
que l'on confie ais\'ement \`a Maple\,:
\[
\footnotesize
\aligned
{\sf A}_1
&
=
\frac{1}{6}\,
\int_0^m\,dk\,\int_0^{\frac{m-k}{16}}\,dj\,
\int_0^{\frac{m-16j-k}{14}}\,di\,\int_0^{\frac{m-14i-16j-k}{12}}\,dh\,
\int_0^{\frac{m-12h-14i-16j-k}{12}}\,dg\,
\\
&
\ \ \ \ \ \ \ \ \ \ \ \ \ \ \ \ \ \ \ \ \ \ \ \ \ \ \ \ \
\big[(m-9g-8h-10i-11j)^3-(3g+2h+3i+3j)^3\big]
+
{\rm O}(m^7)
\\
&
=
\frac{36562817\,\,m^8}
{4933428814282752}
+
{\rm O}(m^7)\,;
\\
{\sf A}_2
&
=
\frac{1}{6}\,
\int_0^m\,dk\,\int_0^{\frac{m-k}{16}}\,dj\,
\int_0^{\frac{m-16j-k}{14}}\,di\,\int_0^{\frac{m-14i-16j-k}{12}}\,dh\,
\int_0^{\frac{m-12h-14i-16j-k}{12}}\,dg\,
\\
&
\ \ \ \ \ \ \ \ \ \ \ \ \ \ \ \ \ \ \ \ \ \ \ \ \ \ \ \ \
(m-12g-10h-13i-14j)^3
+
{\rm O}(m^7)
\\
&
=
\frac{5015441\,\,m^8}
{1233357203570688}
+
{\rm O}(m^7)\,;
\endaligned
\]
les quinze autres fournissent des expressions similaires. Nous
d\'eduisons donc dans chacun des deux cas par sommation de seize nombres
rationnels\,:
\[
\footnotesize
\aligned
\mathcal{C}_1
=
\frac{159897336810563}
{356792619604377600000}
+
{\rm O}(m^7),
\ \ \ \ \ \ \
\mathcal{C}_2
=
\frac{784698232169}
{3303635366707200000}
+
{\rm O}(m^7),
\endaligned
\]
et finalement, le quotient significatif\,:
\[
\mathcal{C}
=
\frac{\mathcal{C}_1}{\mathcal{C}_2}
=
1,\,887\cdots
\]
est inf\'erieur \`a celui de $\mathcal{ DS}_2^4 T_X^*$\,: confirmation
suppl\'ementaire de l'inad\'equation et de l'insuffisance du sous-fibr\'e
engendr\'ee par les crochets.

\smallskip\noindent{\bf Probl\`eme ouvert.} 
Changer d'optique quant au calcul de Riemann-Roch, en tenant compte
d'une \'etude pr\'ealable, reprise \`a partir de z\'ero, de la
structure sp\'ecifique de $\mathcal{ DS}_2^5 T_X^*$. Le proc\'ed\'e de
division par $f_1'$ des g\'en\'erateurs de l'id\'eal des relations
entre invariants restreints \`a $\{ f_1' = 0 \}$ constitue un
proc\'ed\'e ad\'equat et complet d'engendrement qui doit \^etre
pouss\'e au-del\`a de $X^{ 25}$.

\section*{\S9.~Appendice~1\,: jets d'ordre 3 en dimension 3}

\noindent{\bf Expression initiale.}
Comme annonc\'e dans la Section~4, nous d\'etaillons ici le calcul
(d\'elicat) des trois crochets entre les trois invariants $\Lambda_{
1,1}^7$, $\Lambda_{ 1, 2}^7$, $\Lambda_{ 2, 2}^7$ de poids 7:
\[
\small
\aligned
\frac{
\big[
\Lambda_{i,j}^7,\,\Lambda_{k,l}^7
\big]}{7}
&
=
{\sf D}\Lambda_{i,j}^7\cdot\Lambda_{k,l}^7
-
\Lambda_{i,j}^7\cdot{\sf D}\Lambda_{k,l}^7.
\endaligned
\]

\noindent{\bf D\'eveloppement \'econome.}
N'\'ecrivons que le premier produit, en r\'esumant le second par le
symbole $(i, j) \longleftrightarrow (k,l)$, parce qu'il s'en d\'eduit
par ce simple changement d'indices:
\[
\footnotesize
\aligned
&
=
\Big(
\Delta^{1,5}\,f_if_j'
+
5\,\Delta^{2,4}\,f_i'f_j'
-
4\,\Delta^{1,4}\big(f_i''f_j'+f_i'f_j''\big)
-
16\,\Delta^{2,3}\big(f_i''f_j'+f_i'f_j'')
-
\\
&
\ \ \ \ \
-
5\,\Delta^{1,3}\big(f_i'''f_j'+f_i'f_j'''\big)
+
35\,\Delta^{1,3}\,f_i''f_j''
\Big)
\cdot
\Big(
\Delta^{1,4}\,f_k'f_l'
+
4\,\Delta^{2,3}\,f_k'f_l'
-
\\
&
\ \ \ \ \
-
5\,\Delta^{1,3}\big(f_k''f_l'+f_k'f_l''\big)
+
15\,\Delta^{1,2}\,f_k''f_l''
\Big)
-
\\
&
\ \ \ \ \
-
(i,j)\longleftrightarrow(k,l)
\endaligned
\]
\[
\footnotesize
\aligned
&
=
\underline{
\Delta^{1,5}\,\Delta^{1,4}\,f_i'f_j'f_k'f_l'}_{\!\circ}
+
\underline{
5\,\Delta^{2,4}\,\Delta^{1,4}\,f_i'f_j'f_k'f_l'}_{\!\circ}
-
4\,\Delta^{1,4}\,\Delta^{1,4}\big(f_i''f_j'+f_i'f_j'')\,f_k'f_l'
-
\\
&
\ \ \ \ \
-
16\,\Delta^{2,3}\,\Delta^{1,4}\big(f_i''f_j'+f_i'f_j'')\,f_k'f_l'
-
5\,\Delta^{1,4}\,\Delta^{1,3}\big(f_i'''f_j'+f_i'f_j'''\big)\,f_k'f_l'
+
\\
&
\ \ \ \ \
+
35\,\Delta^{1,4}\,\Delta^{1,3}\,f_i''f_j''f_k'f_l'
+
\underline{
4\,\Delta^{1,5}\,\Delta^{2,3}\,f_i'f_j'f_k'f_l'}_{\!\circ}
+
\underline{
20\,\Delta^{2,4}\,\Delta^{2,3}\,f_i'f_j'f_k'f_l'}_{\!\circ}
-
\\
&
\ \ \ \ \
-
16\,\Delta^{1,4}\,\Delta^{2,3}\big(f_i''f_j'+f_i'f_j''\big)\,f_k'f_l'
-
64\,\Delta^{2,3}\,\Delta^{2,3}\big(f_i''f_j'+f_i'f_j''\big)\,f_k'f_l'
-
\\
&
\ \ \ \ \ 
-
20\,\Delta^{2,3}\,\Delta^{1,3}\big(f_i'''f_j'+f_i'f_j'''\big)\,f_k'f_l'
+
140\,\Delta^{2,3}\,\Delta^{1,3}\,f_i''f_j''f_k'f_l'
-
\\
&
\ \ \ \ \
-
5\,\Delta^{1,5}\,\Delta^{1,3}\,f_i'f_j'\big(f_k''f_l'+f_k'f_l''\big)
-
25\,\Delta^{2,4}\,\Delta^{1,3}\,f_i'f_j'\big(f_k''f_l'+f_k'f_l''\big)
+
\\
&
\ \ \ \ \
+
\underline{
20\,\Delta^{1,4}\,\Delta^{1,3}\big(f_i''f_j'+f_i'f_j''\big)
\big(f_k''f_l'+f_k'f_l''\big)}_{\!\circ}
+
\underline{
80\,\Delta^{2,3}\,\Delta^{1,3}\big(f_i''f_j'+f_i'f_j''\big)
\big(f_k''f_l'+f_k'f_l''\big)}_{\!\circ}
+
\\
&
\ \ \ \ \ 
+
25\,\Delta^{1,3}\,\Delta^{1,3}\big(f_i'''f_j'+f_i'f_j'''\big)
\big(f_k''f_l'+f_k'f_l''\big)
-
175\,\Delta^{1,3}\,\Delta^{1,3}\,f_i''f_j''\big(f_k''f_l'+f_k'f_l''\big)
+
\\
&
\ \ \ \ \ 
+
15\,\Delta^{1,5}\,\Delta^{1,2}\,f_i'f_j'f_k''f_l''
+
75\,\Delta^{2,4}\,\Delta^{1,2}\,f_i'f_j'f_k''f_l''
-
60\,\Delta^{1,4}\,\Delta^{1,2}\big(f_i''f_j'+f_i'f_j''\big)\,f_k''f_l''
-
\\
&
\ \ \ \ \
-
240\,\Delta^{2,3}\,\Delta^{1,2}\big(f_i''f_j'+f_i'f_j''\big)\,f_k''f_l''
-
75\,\Delta^{1,3}\,\Delta^{1,2}\big(f_i'''f_j'+f_i'f_j'''\big)\,f_k''f_l''
+
\\
&
\ \ \ \ \
+
\underline{
525\,\Delta^{1,3}\,\Delta^{1,2}\,f_i''f_j''f_k''f_l''}_{\!\circ}
-
\\
&
\ \ \ \ \
-
(i,j)\longleftrightarrow(k,l).
\endaligned
\]
Nous soulignons (en ajoutant un petit cercle) les termes qui
s'annihilent avec ceux qui leur correspondent dans la permutation
$(i,j)\longleftrightarrow(k,l)$. Nous utilisons les relations
pl\"uckeriennes pour remplacer le dernier terme restant, \`a savoir:
\[
\small
\aligned
-
75\,\Delta^{1,3}\,\Delta^{1,2}\big(f_i'''f_j'+f_i'f_j'''\big)\,f_k''f_l''
+
75\,\Delta^{1,3}\,\Delta^{1,2}\big(f_k'''f_l'+f_k'f_l'''\big)\,f_i''f_j''
\endaligned
\]
par:
\[
\small
\aligned
&
-75\,\Delta^{1,3}\,\Delta^{1,3}\,f_i''f_j'f_k''f_l''
+
75\,\Delta^{1,3}\,\Delta^{2,3}\,f_i'f_j'f_k''f_l''
-
\\
&
-
75\,\Delta^{1,3}\,\Delta^{1,3}\,f_i'f_j''f_k''f_l''
+
75\,\Delta^{2,3}\,\Delta^{1,3}\,f_i'f_j'f_k''f_l'',
\endaligned
\]
et nous additionnons tous ces termes en effectuant des regroupements
qui n'impli\-quent que des sommations de nombres entiers:
\[
\scriptsize
\aligned
&
=
-5\,\Delta^{1,5}\,\Delta^{1,3}\,f_i'f_j'\big(f_k''f_l'+f_k'f_l''\big)
-
25\,\Delta^{2,4}\,\Delta^{1,3}\,f_i'f_j'\big(f_k''f_l'+f_k'f_l''\big)
+
15\,\Delta^{1,5}\,\Delta^{1,2}\,f_i'f_j'f_k''f_l''
+
\\
&
\ \ \ \ \
+
75\,\Delta^{2,4}\,\Delta^{1,2}\,f_i'f_j'f_k''f_l''
-
4\,\Delta^{1,4}\,\Delta^{1,4}\big(f_i''f_j'+f_i'f_j''\big)\,f_k'f_l'
-
32\,\Delta^{1,4}\,\Delta^{2,3}\big(f_i''f_j'+f_i'f_j''\big)\,f_k'f_l'
-
\\
&
\ \ \ \ \
-
64\,\Delta^{2,3}\,\Delta^{2,3}\big(f_i''f_j'+f_i'f_j''\big)\,f_k'f_l'
-
5\,\Delta^{1,4}\,\Delta^{1,3}\big(f_i'''f_j'+f_i'f_j'''\big)\,f_k'f_l'
+
35\,\Delta^{1,4}\,\Delta^{1,3}\,f_i''f_j''f_k'f_l'
-
\\
&
\ \ \ \ \
-
20\,\Delta^{2,3}\,\Delta^{1,3}\big(f_i'''f_j'+f_i'f_j'''\big)\,f_k'f_l'
+
140\,\Delta^{2,3}\,\Delta^{1,3}\,f_i''f_j''f_k'f_l'
+
150\,\Delta^{2,3}\,\Delta^{1,3}\,f_i'f_j'f_k''f_l''
-
\\
&
\ \ \ \ \
-
60\,\Delta^{1,4}\,\Delta^{1,2}\big(f_i''f_j'+f_i'f_j''\big)\,f_k''f_l''
-
240\,\Delta^{2,3}\,\Delta^{1,2}\big(f_i''f_j'+f_i'f_j''\big)\,f_k''f_l''
+
\\
&
\ \ \ \ \
+
25\,\Delta^{1,3}\,\Delta^{1,3}\big(f_i'''f_j'+f_i'f_j'''\big)
\big(f_k''f_l'+f_k'f_l''\big)
-
175\,\Delta^{1,3}\,\Delta^{1,3}\,f_i''f_j''\big(f_k''f_l'+f_k'f_l''\big)
-
\\
&
\ \ \ \ \
-
75\,\Delta^{1,3}\,\Delta^{1,3}\big(f_i''f_j'+f_i'f_j''\big)\,f_k''f_l''
-
\\
&
\ \ \ \ \
-
(i,j)\longleftrightarrow(k,l).
\endaligned
\]

\noindent{\bf Synth\`ese de d\'eterminants.}  Maintenant, la
soustraction suivie de la permutation fait appara\^{\i}tre des
d\'eterminants $2 \times 2$: on a en effet cinq relations
imm\'ediatement v\'erifiables par d\'eveloppement:
\[
\small
\aligned
f_i'f_j'\big(f_k''f_l'+f_k'f_l''\big)
-
f_k'f_l'\big(f_i''f_j'+f_i'f_j''\big)
&
=
f_j'f_l'\,\Delta_{i,k}^{1,2}
+
f_i'f_k'\,\Delta_{j,l}^{1,2},
\\
f_i'f_j'f_k''f_l''
-
f_i''f_j''f_k'f_l'
&
=
f_i'f_l''\,\Delta_{j,k}^{1,2}
+
f_k'f_j''\,\Delta_{i,l}^{1,2},
\\
\big(f_i'''f_j'+f_i'f_j'''\big)\,f_k'f_l'
-
\big(f_k'''f_l'+f_k'f_l'''\big)\,f_i'f_j'
&
=
f_j'f_l'\,\Delta_{k,i}^{1,3}
+
f_i'f_k'\,\Delta_{l,j}^{1,3},
\\
\big(f_i''f_j'+f_i'f_j''\big)\,f_k''f_l''
-
\big(f_k''f_l'+f_k'f_l''\big)f_i''f_j''
&
=
f_i''f_k''\,\Delta_{j,l}^{1,2}
+
f_j''f_l''\,\Delta_{i,k}^{1,2},
\endaligned
\]
\[
\small
\aligned
\\
&
\big(f_i'''f_j'+f_i'f_j'''\big)\big(f_k''f_l'+f_k'f_l''\big)
-
\big(f_k'''f_l'+f_k'f_l'''\big)\big(f_i''f_j'+f_i'f_j''\big)
\\
&
\ \ \ \ \ \ \ \ \ \ \ \ \ \ \ \
\ \ \ \ \ \ \ \ \ \ \ \ \ \ \ \
=
f_j'f_l'\,\Delta_{k,i}^{2,3}
+
f_j'f_k'\,\Delta_{l,i}^{2,3}
+
f_i'f_l'\,\Delta_{k,j}^{2,3}
+
f_i'f_k'\,\Delta_{l,j}^{2,3}.
\endaligned
\]

\noindent
En regroupant les termes, on obtient l'expression finale de 
$\frac{ 1}{ 7}\,
\big[ \Lambda_{ i,j}^7, \, \Lambda_{ k,l}^7 \big]$.
\hfill$\square$

\section*{\S10.~Appendice~2\,: jets d'ordre 3 en dimension 3}

\noindent{ Jets d'ordre $\kappa = 3$ en dimension $\nu = 3$.}
Pour terminer, donnons une description ``\`a la main'' des
g\'en\'erateurs de $\mathcal{ DS}_3^3$ qui ne fasse pas appel \`a des
arguments raffin\'es de th\'eorie des invariants ({\it cf.} \cite{
ro2007}).

Recherchons directement les bi-invariants, {\it i.e.} les polyn\^omes
invariants par reparam\'etrisation qui sont aussi invariants par
l'action du sous-groupe unipotent ${\sf U}_3 ( \C) \subset {\sf GL}_3
( \C)$ constitu\'e des matrices de la forme:
\[
\text{\sc u}
:=
\left(
\begin{array}{ccc}
1 & 0 & 0 
\\
u_a & 1 & 0
\\
u_c & u_b & 1
\end{array}
\right),
\]
qui est d\'efinie par $\text{\sc u}\cdot f_1^{ ( \lambda)} := f_1^{ (
\lambda)}$, puis $\text{\sc u} \cdot f_2^{ ( \lambda)} := f_2^{ (
\lambda)} + u_a f_1^{ ( \lambda)}$ et enfin $\text{\sc u} \cdot f_3^{
(\lambda)} := f_3^{ (\lambda)} + u_b \, f_2^{ ( \lambda)} + u_c \,
f_1^{ ( \lambda)}$, pour $\lambda = 1, 2, 3$. Un premier raisonnement
initial enti\`erement similaire \`a celui que nous avons tir\'e de \cite{
ro2007} pour $\mathcal{ DS}_{ 2, m}^4$ fournit une repr\'esentation de
tout ${\sf P} \in \mathcal{ DS}_{3, m}^3$ sous la forme rationnelle
\[
{\sf P}\big(j^3f\big)
=
\sum_{-\frac{2}{3}m\leqslant a\leqslant m}\,
(f_1')^a\,\mathcal{P}_a
\big(f_2',f_3',\Lambda_{1,2}^3,\Lambda_{1,3}^3,\Lambda_{1,2;1}^5,
\Lambda_{1,3;1}^5\big),
\]
o\`u les invariants $\Lambda_{ i,j}^3$ et
$\Lambda_{ i,j; k}^5$ sont simplement d\'efinis par:
\[
\Lambda_{i,j}^3:=\Delta_{i,j}^{1,2}
\ \ \ \ \ \ \ \ \ \ \ \ \ \ \ \ \ \
\text{\rm et}
\ \ \ \ \ \ \ \ \ \ \ \ \ \ \ \ \ \
\Lambda_{i,j;k}^5:=\Delta_{i,j}^{1,3}\,f_k'-3\,\Delta_{i,j}^{1,2}\,f_k''.
\]
Ensuite, si nous consid\'erons le sous-groupe de ${\sf U}_3 ( \C)$
constitu\'e des matrices de la forme:
\[
\overline{\text{\sc u}}
:=
\left(
\begin{array}{ccc}
1 & 0 & 0 
\\
u_a & 1 & 0
\\
u_c & 0 & 1
\end{array}
\right),
\]
lesquelles stabilisent tous les invariants qui apparaissent
dans notre premi\`ere expression rationnelle:
\[
\aligned
\overline{\text{\sc u}}\cdot
\Lambda_{1,2}^3
=
\Lambda_{1,2}^3,
\ \ \ \ \ \ \ \ \ \ \ \ \ \
\overline{\text{\sc u}}\cdot\Lambda_{1,2;1}^5
=
\Lambda_{1,2;1}^5,
\\
\overline{\text{\sc u}}\cdot
\Lambda_{1,3}^3
=
\Lambda_{1,3}^3,
\ \ \ \ \ \ \ \ \ \ \ \ \ \
\overline{\text{\sc u}}\cdot\Lambda_{1,3;1}^5
=
\Lambda_{1,3;1}^5,
\endaligned
\]
mais agissent en perturbant $f_2'$ et $f_3'$ par $\overline{ \text{\sc
u}}\cdot f_2' = f_2' + u_a \, f_1'$ et $\overline{ \text{\sc u}}\cdot
f_3' = f_3' + u_c \, f_1'$, nous voyons que si ${\sf P}$ est aussi
invariant par l'action unipotente, alors chaque polyn\^ome $\mathcal{
P}_a$ ci-dessus doit en fait \^etre ind\'ependant de $f_2'$ et de $f_3'$,
d'o\`u:
\[
{\sf P}^{2\times{\rm inv}}
\big(j^3f\big)
=
\sum_{-\frac{2}{3}m\leqslant a\leqslant m}\,
(f_1')^a\,\mathcal{P}_a\big(\Lambda_{1,2}^3,\,\Lambda_{1,3}^3,\,
\Lambda_{1,2;1}^5,\,\Lambda_{1,3;1}^5\big).
\]
Mais ce n'est pas termin\'e, car il faut encore s'assurer de l'invariance
par l'action du sous-groupe constitu\'e des matrices
de la forme:
\[
\text{\sc u}_b
:=
\left(
\begin{array}{ccc}
1 & 0 & 0
\\
0 & 1 & 0
\\
0 & u_b & 1
\end{array}
\right),
\]
lesquelles agissent comme suit sur les quatre invariants 
alg\'ebriquement ind\'ependant qui apparaissent:
\[
\aligned
&
\text{\sc u}_b\cdot\Lambda_{1,2}^3
=
\Lambda_{1,2}^3,
\ \ \ \ \ \ \ \ \ \ \ \ \ \ \ \ \ \ \ \ \ \ \ \ \ \
\text{\sc u}_b\cdot\Lambda_{1,2;1}^5
=
\Lambda_{1,2;1}^5,
\\
&
\text{\sc u}_b\cdot\Lambda_{1,3}^3
=
\Lambda_{1,3}^3+u_b\,\Lambda_{1,2}^3,
\ \ \ \ \ \ \ \ \ \ \
\text{\sc u}_b\cdot\Lambda_{1,3;1}^5
=
\Lambda_{1,3;1}^5+u_b\,\Lambda_{1,2;1}^5.
\endaligned
\]
D'apr\`es un calcul direct, le troisi\`eme et dernier invariant
fondamental pour cette action, \`a savoir:
\[
\Lambda_{1,2}^3\,\Lambda_{1,3;1}^5
-
\Lambda_{1,3}^3\,\Lambda_{1,2;1}^5
\equiv
f_1'f_1'\,D_{1,2,3}^6
\]
fait na\^{\i}tre, en tant que nouveau bi-invariant ``fant\^ome''
cach\'e derri\`ere $(f_1 ')^2$, le d\'eter\-minant wronskien en
dimension trois:
\[
D_{1,2,3}^6
:=
\left\vert
\begin{array}{ccc}
f_1' & f_2' & f_3'
\\
f_1'' & f_2'' & f_3 ''
\\
f_1''' & f_2''' & f_3'''
\end{array}
\right\vert,
\]
et en injectant ce nouveau bi-invariant, nous obtenons une nouvelle
expression pour le bi-invariant quelconque dont nous \'etions partis:
\[
{\sf P}^{2\times{\rm inv}}
\big(
j^3f
\big)
=
\sum_{-\frac{2}{3}m\leqslant a\leqslant m}\,
(f_1')^a\,\mathcal{P}_a\big(\Lambda_{1,2}^3,\,
\Lambda_{1,2;1}^5,\,D_{1,2,3}^6\big),
\]
laquelle incorpore toujours des puissances n\'egatives de $f_1'$. Mais
pour terminer, nous affirmons qu'il n'existe en fait aucune puissance
n\'egative de $f_1'$, car sinon, apr\`es multiplication par la puissance
maximalement n\'egative de $f_1'$ ({\it cf.} le raisonnement conduit
dans la Section~6) et apr\`es restriction \`a $\{ f_1' = 0\}$, nous
obtiendrions une identit\'e du type:
\[
\aligned
0
&
\equiv
\mathcal{P}_a\big(\Lambda_{1,2}^3,\,\Lambda_{1,2;1}^5,\,D_{1,2,3}^6\big)
\Big\vert_{f_1'=0}
\\
&
=
\mathcal{P}_a\Big(-f_1''f_2',\,3\,f_1''f_2'f_1'',\,
\left\vert
\begin{array}{ccc}
0 & f_2' & f_3 '
\\
f_1'' & f_2'' & f_3 ''
\\
f_1''' & f_2''' & f_3'''
\end{array}
\right\vert
\bigg),
\endaligned
\]
qui impliquerait imm\'ediatement $\mathcal{ P}_a \equiv 0$, par
ind\'ependance alg\'ebrique de ses trois arguments.

\smallskip

En conclusion, nous avons red\'emontr\'e ({\it cf.} \cite{ ro2007}) qu'en
dimension $\nu = 3$ pour les jets d'ordre $\kappa = 3$, les polyn\^omes
bi-invariants s'\'ecrivent:
\[
\mathcal{P}^{2\times{\rm inv}}\big(j^3f\big)
=
\mathcal{P}\big(f_1',\,\Lambda_{1,2}^3,\,
\Lambda_{1,2;1}^5,\,D_{1,2,3}^6\big),
\]
o\`u $\mathcal{P}$ est un polyn\^ome arbitraire, aucune syzygie
n'existant entre ces quatre bi-invariants fondamentaux, et par
cons\'equent, en polarisant les indices\,\,---\,\,ce qui revient \`a
faire agir le groupe complet ${\sf GL}_3 ( \C)$\,\,---, nous
d\'eduisons que les polyn\^omes g\'en\'eraux invariants par
reparam\'etrisation s'\'ecrivent comme polyn\^omes quelconques en
fonction de 16 invariants fondamentaux:
\[
\aligned
{\sf P}\big(j^3f\big)
=
\mathcal{P}
\big(
f_1',\,f_2',\,f_3',\,
&
\Lambda_{1,2}^3,\,\Lambda_{1,3}^3,\,\Lambda_{2,3}^3,\,
\Lambda_{1,2;1}^5,\,
\Lambda_{1,2;2}^5,\,
\Lambda_{1,2;3}^5,\,
\\
&
\Lambda_{1,3;1}^5,\,
\Lambda_{1,3;2}^5,\,
\Lambda_{1,3;3}^5,\,
\Lambda_{2,3;1}^5,\,
\Lambda_{2,3;2}^5,\,
\Lambda_{2,3;3}^5,\,
D_{1,2,3}^6
\big).
\endaligned
\]
On v\'erifie aussi que parmi les 62 syzygies existant entre ces 16 invariants
qui ont \'et\'e obtenues par un calcul sur Maple ({\it cf.} les r\'ef\'erences
dans \cite{ ro2007}), 30 d'entre elles sont fondamentales et qu'elles
proviennent toutes des trois proc\'edures que nous avons d\'ecrites.
\hfill$\square$

{\bf Received: December 20, 2007}


\begin{thebibliography}{99}

\bibitem{clo2007}
{\sc Cox}, D.; {\sc Little}, J.; {\sc O'Shea}, D.:
{\em Ideals, varie\-ties and algorithms}, {\rm An introduction to
computational algebraic geometry and commutative algebra}, 
Third Edition, Undergraduate Texts in Mathematics,
Springer-Verlag, New York, 2007, xvi+551~pp. 

\bibitem{de1997}
{\sc Demailly}, J.-P.:
{\em Algebraic criteria for Kobayashi hyperbolic projective varieties
and jet differentials}, Proc. Sympos. Pure Math., vol. 62, Amer. Math.
Soc., Providence, RI, 1997, 285--360.

\bibitem{de2007}
{\sc Demailly}, J.-P.:
{\em Oral communication at the Conference}, 
{\tt Effective Aspects of Complex Hyperbolic Varieties}, 
Aber Wrac'h, Breizh, 10--14 september 2007.
{\tt stockage.univ-brest.fr/$\sim$} {\tt huisman/rech/clq/eachv}

\bibitem{ol1999}
{\sc Olver}, P.J.:
{\em Classical invariant theory}, London Mathematical Society Student
Texts, 44, Cambridge University Press, 1999, xxi+280~pp.

\bibitem{ro2007}
{\sc Rousseau}, Erwan:
{\em Hyperbolicit\'e des vari\'et\'es complexes}, 
Cours Peccot au Coll\`ege de France, Mai-Juin 2007, 52~pp.,
t\'el\'echargeable \`a l'adresse:
{\tt arxiv.org/abs/0709.3882/}

\end{thebibliography}
\end{document}